\documentclass[onefignum,onetabnum]{siamonline190516}

\usepackage{amsfonts}
\usepackage{array}
\usepackage{algorithm,algorithmic}
\usepackage{color}
\usepackage{verbatim}
\usepackage{booktabs}
\usepackage{threeparttable}
\usepackage{url}

\usepackage{color}
\usepackage{graphicx,epstopdf}
\usepackage{mathrsfs}
\usepackage{color}
\usepackage{longtable}
\usepackage{multirow}
\usepackage{subfigure}
\usepackage{algorithm,algorithmic}
\usepackage{booktabs}
\usepackage{url}

\newcommand{\normmm}[1]{{\vert\kern-0.25ex\vert\kern-0.25ex\vert #1
    \vert\kern-0.25ex\vert\kern-0.25ex\vert}}

\numberwithin{equation}{section}
\allowdisplaybreaks

\newcommand{\wsize}{.20\textwidth}
\newcommand{\hhsize}{.18\textwidth}

\newsiamremark{assumption}{Assumption}
\newsiamremark{example}{Example}
\newsiamremark{remark}{Remark}
\usepackage{enumitem}
\setlist[enumerate]{leftmargin=.5in}
\setlist[itemize]{leftmargin=.5in}

\usepackage{graphicx,epstopdf}
\graphicspath{{./pics/}}
\usepackage[toc,page]{appendix}
\allowdisplaybreaks
\numberwithin{equation}{section}
\numberwithin{figure}{section}

\headers{Imaging Anisotropic Conductivities}{H. Liu, B. Jin, and X. Lu}

\title{Imaging Anisotropic Conductivities from Current Densities\thanks{The work of H.L. is partially supported by the National Science Foundation of China (No. 12101276), the PhD research startup foundation of Jinling Institute of Technology(No. jit-b-202048 and No. jit-fhxm-202117), that of B.J. by UK EPSRC grant EP/T000864/1, and that of X.L. by the National Key Research and Development Program of China (No. 2020YFA0714200) and the National Science Foundation of China (No. 11871385).}}

\author{Huan Liu\thanks{College of Science, Jinling Institute of Technology, Nanjing 211169, P.R. China (huanliumath@jit.edu.cn)} \and Bangti Jin\thanks{Department of Computer Science, University College London, Gower Street, London WC1E 6BT, UK (bangti.jin@gmail.com, b.jin@ucl.ac.uk)}\and
Xiliang Lu\thanks{School of Mathematics and Statistics, Wuhan University, Wuhan 430072, P.R. China, and Hubei Key Laboratory of Computational Science (Wuhan University), Wuhan 430072, China (xllv.math@whu.edu.cn)}}

\ifpdf
\hypersetup{
  pdftitle={Anisotropic Conductivity},
  pdfauthor={H. Liu, B. Jin, and  X. Lu}
}
\fi

\begin{document}

\maketitle

\begin{abstract}
In this paper, we propose and analyze a reconstruction algorithm for imaging an anisotropic
conductivity tensor in a second-order elliptic PDE with a nonzero Dirichlet boundary condition from internal current densities.
It is based on a regularized output least-squares formulation with the standard $L^2(\Omega)^{d,d}$ penalty, which is then discretized
by the standard Galerkin finite element method. We establish the continuity and differentiability of the forward
map with respect to the conductivity tensor in the $L^p(\Omega)^{d,d}$-norms, the existence of minimizers and optimality
systems of the regularized formulation using the concept of H-convergence. Further, we provide a detailed
analysis of the discretized problem, especially the convergence of the discrete approximations with respect
to the mesh size, using the discrete counterpart of H-convergence. In addition, we develop a projected Newton algorithm for solving the first-order
optimality system. We present extensive two-dimensional numerical examples to show the efficiency of
the proposed method.
\end{abstract}


\section{Introduction}

The conductivity value varies widely with soft tissue types \cite{FosterSchwan:1989,Morimoto:1993} and
its imaging can provide valuable information about the physiological and pathological conditions
of tissue \cite{AmmariGarnier:2016}. This underpins several important medical imaging modalities \cite{Borcea:2002,
Ammari:2008,Uhlmann,SeoWoo:2011,WidlakScherzer:2012,AmmariBoulmier:2015,JensenJinKnudsen:2019}.
For example, electrical impedance tomography (EIT) \cite{Borcea:2002,
Uhlmann} aims at recovering the interior conductivity distribution from boundary voltage
measurement. However, it is severely ill-posed, and the attainable resolution is limited,
and the uniqueness for anisotropic conductivity is only for a given conformal class
\cite{Ferreira09,Kenig2011} (see \cite{Hamilton2014} for some numerics).
To overcome the ill-posed nature in EIT, MREIT (magnetic resonance EIT) employs an MRI scanner
to capture the internal magnetic flux density data $\vec{b}$ induced by an externally injected
current \cite{JoyScott:1989,ScottJoyArmstrong:1991,Ider,GambaBayford:1999} and then obtains
the current density $\vec{h}$ according to Ampere's law $\vec{h}=\mu_0^{-1}\nabla \times \vec{b}$, where
$\mu_0$ is the magnetic permeability of free space. This requires measuring all components
of the magnetic flux $\vec{b}$, which may be challenging in practice, as it requires
a rotation of the domain being imaged or of the MRI scanner.

The use of internal data promises much improved image resolution, and the reconstruction problem has
received much attention. Kwon et al \cite{Kwon} proposed a non-iterative reconstruction method for
recovering an isotropic conductivity using equipotential lines, and proved the unique
recovery from one current density vector $\vec{h}$ in two dimension. Later, an iterative algorithm
known as $J$-substitution was proposed \cite{Kwon2002}. With the knowledge of the magnitude of only one
current density magnitude $|\vec{h}|$, the problem was studied in a series of work
\cite{Nachman2007,Nachman2009,YazdanianKnudsen:2021} in the isotropic
case, and also the anisotropic case in a known conformal class in \cite{Hoell2014}.
Ammari et al. \cite{zhangwenlong17} studied the recovery of an anisotropic
conductivity proportional to a known conductivity tensor.

It is widely accepted that most biological tissues have anisotropic conductivity values.  The ratio of the
anisotropy depends on the type of tissue and human skeletal muscle shows anisotropy of up to 10 between
the longitudinal and transversal directions \cite{RushDriscoll:1968}. The conductivity for
cell membrane, muscle fiber, or nerve fiber structure must be investigated using an anisotropic
tensor model \cite{Nicholson:1965,Roth:2000,AmmariGarnier:2016}. It arises naturally also in the mathematical modeling of
boundary deformation  \cite{AlbertiAmmari:2016}. Therefore, in the past few years, there has been a growing interest
in recovering anisotropic conductivities \cite{SeoPyoPark:2004,bal2014,G.Bal2014,Monard:2012,MonardBal:2013,MonardRim:2018}.
The potential of using current densities to
recover anisotropic conductivity $A$ was studied in \cite{SeoPyoPark:2004}. The uniqueness of recovering
an anisotropic resistivity distribution from current densities has been established
\cite{bal2014,G.Bal2014} (see also \cite{Monard:2012,MonardBal:2013} for power densities)
and \cite{MonardRim:2018} for 3D numerical implementation. Bal et al \cite{bal2014} showed
that a minimum number current densities ensures a unique and explicit reconstruction of the
conductivity tensor $A$ locally, and that the reconstruction of $A$ is about the loss of one derivative compared to
errors in the measurement of $\vec{h}$. Ko and Kim \cite{KoKim:2017} gave a resistive
network-based reconstruction method for three sets of internal electrical current densities,
by directly discretizing Faraday's law. In addition, Hsiao and Sprekels \cite{HsiaoSprekels:1988}
investigated the stability of recovering matrices of the form $A=\nabla p\otimes\nabla u$.

In this work, we consider the inverse problem of recovering an anisotropic conductivity
tensor $A(x)$ from internal current densities $\vec{h}=A\nabla u$, where $u\in H^1(\Omega)$ solves
 \begin{equation}\label{pr1}
    \left\{\begin{aligned}
   -\mathrm{div}(A\nabla u)&=f, &&\mbox{in } \Omega,\\
    u&=g, &&\mbox{on } \Gamma,
    \end{aligned}\right.
\end{equation}
where $\Omega\subset \mathbb{R}^d$ ($d=2,3$) is an open bounded Lipschitz domain with a boundary $\Gamma$,
and $(f,g)\in (H^1(\Omega))^\prime\times H^{\frac12}(\Gamma)$.
 The anisotropic conductivity tensor $A(x)=(A_{ij}(x))_{i,j=1}^d$ belongs to the following admissible set
\begin{align}\label{uec}
  \mathcal{A}=\bigg\{A\in L^\infty(\Omega)^{d,d}: A(x)\in \mathcal{S}_d, \alpha|\xi|^2\leq \sum_{i,j=1}^d A_{ij}(x)\xi_i\xi_j \leq \beta |\xi|^2, \forall \xi \in \mathbb{R}^d \ \mbox{a.e. }x\in \Omega\bigg\}
\end{align}
with constants $0<\alpha < \beta <+\infty$, where $\mathcal{S}_d$ denotes the set of all
$d\times d$ symmetric matrices endowed with the Frobenius inner product.
We use  the notation $u(A)(f,g)$ to explicitly indicate the dependence of the solution $u$ on the
conductivity $A$. The $L^p(\Omega)^{d,d}$-norm $(1\leq p<\infty)$ on $\mathcal{A}$ is given by
\begin{equation*}
\normmm{A}_p=\Big(\int_\Omega \sum_{i,j=1}^d |A_{i,j}(x)|^pdx\Big)^{\frac1p},
\end{equation*}
and for the case $p=\infty$, it is defined in terms of supremum as usual. When $p=2$, it
is naturally induced by an inner product, which is denoted by $(\cdot,\cdot)_{L^2(\Omega)^{d,d}}$ below.
For the recovery of the conductivity tensor $A$, suppose that we are given a set of $L$ measurements $\vec{h}^\ell$, $\ell=1,\ldots,L$,
corresponding to the excitations $(f^\ell,g^\ell)$, $\ell=1,\ldots,L$. To handle the ill-posedness of
the inverse problem, we consider the standard Tikhonov regularization \cite{EnglHankeNeubauer:1996,ItoJin:2015},
i.e., the least squares data fitting regularized with an $L^2(\Omega)^{d, d}$ penalty. More precisely, given the
regularization parameter $\gamma>0$, we minimize the following regularized functional over the admissible
set $\mathcal{A}$:
\begin{equation}\label{tik}
 \min_{A\in\mathcal{A}}\Big\{J_\gamma(A)=\frac12\sum_{\ell=1}^L\|A\nabla u^\ell(A)(f^\ell,g^\ell)-\vec{h}^\ell\|^2_{L^2(\Omega)^d}+\frac{\gamma}{2}\normmm{A}^2_2\Big\}.
\end{equation}
Note that the formulation \eqref{tik} involves only the $L^2(\Omega)^{d,d}$ penalty, and it is
suitable for recovering both smooth and discontinuous conductivity tensors. However,
the choice of the $L^2(\Omega)^{d,d}$ penalty introduces certain challenges in the mathematical
and numerical analysis of the regularized problem, due to a lack of weak sequential closeness
of the parameter-to-state map. One way to address this issue is to use a
stronger norm for the penalty, e.g., $\|\cdot\|_{H^1(\Omega)}$, but it can be numerically
cumbersome to treat the ellipticity constraint (i.e., the bounds on the extremal eigenvalues
of the tensor $A$ in the admissible set $\mathcal{A}$, cf. \eqref{uec}).

In this work, we provide a detailed analysis of the reconstruction approach, e.g.,
the parameter-to-state map, especially $L^p(\Omega)^{d,d}$ differentiability, well-posedness of the
variational formulation \eqref{tik} (existence of a minimizer, optimality system and consistency), and the
finite element discretization and its convergence. One distinct challenge in the analysis arises
from the fact that the regularized formulation \eqref{tik} involves only an $L^2(\Omega)^{d,d}$
penalty, which does not induce strong compactness. In order to resolve the challenge, we
resort to the concept of H-convergence \cite{LTartar,Tartar}, and its discrete analogue, i.e., Hd-convergence
\cite{Eymard} (developed for a finite volume scheme). In order to apply the concept, the presence
of a nonzero Dirichlet boundary condition in the governing model \eqref{pr1} necessitates revisiting
known H-convergence results; see Theorems \ref{compact} and \ref{dcompact} for the precise statements.
Further, we present extensive numerical results
to validate the approach, and the numerical results show clearly the efficiency and accuracy of the approach.

Numerical algorithms for recovering matrix parameters have not been extensively studied
\cite{KohnLowe:1988,Deckelnick2011,Deckelnick2012}. Kohn and Lowe \cite{KohnLowe:1988} introduced a variational method involving a
convex functional for recovering a matrix-valued diffusion coefficient, but did not present
numerical experiments. Deckelnick and Hinze
\cite{Deckelnick2011,Deckelnick2012} studied the identification of matrix parameters in elliptic PDEs from
measurements $(z,f) \in Z \times H^{-1}(\Omega)$ ($Z = L^2(\Omega)$ or $Z = H^1_0(\Omega)$) using
the concept of H-convergence. This work is inspired by the prior works \cite{Deckelnick2011,Deckelnick2012}.
However, there are major differences between these works and the present one. First,
the model in \cite{Deckelnick2011,Deckelnick2012} has a zero Dirichlet boundary condition, but the model \eqref{pr1} involves
a nonzero Dirichlet boundary condition. The presence of a nonzero Dirichlet boundary condition,
inherent to the concerned inverse problem, poses big challenges to the mathematical and numerical analysis
of the regularized formulation \eqref{tik}. Most studies deal only with a zero
Dirichlet boundary condition, and the concept of H-convergence and its discrete analogue (i.e., Hd-convergence)
have to be revisited, which in particular requires establishing relevant fundamental results, cf. Theorems \ref{compact}
and \ref{dcompact} for the continuous and discrete H-convergence results with a nonzero Dirichlet boundary
condition, respectively. This represents the major technical novelty of the study. Second,
in \cite{Deckelnick2012}, the authors carry out the optimization by the projected steepest descent
method with Armijo's rule, whereas we solve the resulting optimization problem by a Newton type algorithm
from measurements $(A\nabla u, g)\in L^2(\Omega)^d\times H^{\frac12}(\Gamma)$, which is numerically observed to
be highly efficient and accurate, when coupled with a path-following strategy. Third, we give a
detailed analysis of the continuity and differentiability of parameter-to-state map, which is
essential for rigorously developing the numerical algorithm (see also the work
\cite{YazdanianKnudsen:2021} for relevant results). Last, the works
\cite{Deckelnick2011,Deckelnick2012} are concerned with variational discretization of the conductivity tensor,
which greatly facilitates the convergence analysis, whereas this work
analyzes the Hd-convergence for general discretization of the conductivity tensor $A$.

The rest of the paper is organized as follows. In Section \ref{sec:reg}, we analyze
the well-posedness of the formulation \eqref{tik}, including existence, optimality system and
consistency, using the concept of H-convergence. Then in Section \ref{sec:disc},
we develop the finite element discretization, prove the convergence of the approximations
as the mesh size $h$ tends to zero using the concept of Hd-convergence, and describe a projected Newton method for solving
the smoothed optimality system. Last, in Section \ref{sec:numer}, we present several
two-dimensional examples to illustrate distinct features of the proposed approach.
Throughout, the notation $C$ denotes a generic constant, which may differ at different occurrences
but does not depend on the matrix $A$ and any functions / parameters (e.g., mesh size $h$) involved in the analysis.
The notation $(\cdot,\cdot)$ with suitable subscripts denotes the $L^2(\Omega)$, $L^2(\Omega)^d$ or $L^2(\Omega)^{d, d}$ inner product,
and $\langle\cdot,\cdot\rangle$ denotes the duality pairing, e.g., between $H_0^1(\Omega)$ and its dual $H^{-1}(\Omega)$.

\section{The regularized formulation}\label{sec:reg}

In this section, we study the well-posedness of the regularized formulation \eqref{tik}
using the concept of H-convergence.

\subsection{Preliminary estimates}\label{ssec:basic}
For any fixed $(f,g)\in (H^1(\Omega))^\prime\times H^{\frac12}(\Gamma)$, and any $A\in\mathcal{A}$,
problem \eqref{pr1} has a unique solution $u \in H^1(\Omega)$. The parameter-to-state map $F: A\mapsto u$
is defined by $u=F(A)(f,g)$. We often write $u = F(A)$, by suppressing the dependence on the problem data $(f,g)$.
By Lax-Milgram theorem, we have the following \textit{a priori} regularity estimate of the solution $u$ to problem \eqref{pr1}.
\begin{lemma}\label{lemma1}
For any $A\in \mathcal{A}$ and $(f,g)\in (H^1(\Omega))^\prime\times H^{\frac12}(\Gamma)$, there exists a unique solution
$u=F(A)\in H^1(\Omega)$ to problem \eqref{pr1}, and it satisfies
\begin{equation*}
  \|u\|_{H^1(\Omega)}\leq C(\|f\|_{(H^1(\Omega))^\prime}+\|g\|_{H^{\frac12}(\Gamma)}).
\end{equation*}
\end{lemma}

Lemma \ref{lemma1} implies that $u\equiv F(A)$ is uniformly bounded in $H^1(\Omega)$
for any fixed $(f,g)\in (H^1(\Omega))^\prime\times H^{\frac12}(\Gamma)$. Actually, the
regularity can be slightly improved using Meyers' gradient estimates
\cite[Theorem 1]{Meyers}. We denote by $Q:=Q(\alpha, \beta)\in (2,\infty)$,
defined in \cite[Theorem 1]{Meyers}, by suppressing its dependence on $d$, which satisfies
$Q\rightarrow 2$ as $\frac\alpha\beta\rightarrow 0$ and $Q\rightarrow r$ ($r$ is given in
Theorem \ref{recmeyer}) as $\frac\alpha\beta\rightarrow 1$.
The next result gives an $L^q(\Omega)$ gradient estimate for a nonzero Dirichlet boundary condition.
Throughout, let
\begin{equation*}
Q_d=\left\{\begin{aligned}
\min(Q,4),&\quad d=2,\\
\min(Q,3), &\quad d=3.
\end{aligned}\right.
\end{equation*}


\begin{theorem}\label{recmeyer}
Let $\Omega\subset \mathbb{R}^d(d=2, 3)$ be a bounded $C^r$ domain, with $r\in (2, \infty)$,
and $A\in \mathcal{A}$. For any $f\in L^2(\Omega)$ and $g\in H^1(\Gamma)$,
let $u\in H^1(\Omega)$ be a weak solution of problem \eqref{pr1}.
Then, for any $2<q<Q_d$, problem \eqref{pr1} has a unique solution $u\in W^{1,q}(\Omega)$, and there exists
 $C=C(\alpha, \beta, d, q, \Omega)>0$ such that
\begin{equation*}
\|u\|_{W^{1,q}(\Omega)}\leq C(\|f\|_{L^2(\Omega)}+\|g\|_{H^1(\Gamma)}).
\end{equation*}
\end{theorem}
\begin{proof}
Let $\bar{u}$ be the solution to the Laplace equation $-\triangle \bar{u}=0$ in $\Omega$ with
$\bar{u}|_\Gamma=g\in H^1(\Gamma)$. Then by the standard elliptic regularity theory, we obtain $\bar{u}\in H^{\frac32}(\Omega)$, and
$\|\nabla \bar{u}\|_{H^{\frac12}(\Omega)} \leq C\|g\|_{H^1(\Gamma)}$. For any $q\in Q_d$, by
Sobolev embedding theorem \cite{Evans},
\begin{equation}\label{soboemb}
  H^{\frac12}(\Omega)\hookrightarrow  \left\{
    \begin{array}{ll}
    L^4(\Omega),    &d=2\\
    L^3(\Omega),&d=3
    \end{array}\right.
    \hookrightarrow L^q(\Omega).
\end{equation}
Thus,
$$\|\nabla \bar{u}\|_{L^q(\Omega)}\leq C\|\nabla \bar{u}\|_{H^{\frac12}(\Omega)}
\leq C\|g\|_{H^1(\Gamma)}.$$
Let $w:=u-\bar{u}\in H^1_0(\Omega)$, which satisfies
\begin{equation*}\left\{
 \begin{aligned}
    -\mathrm{div}(A\nabla w)&=\mathrm{div}(A\nabla \bar{u})+f,   && \mbox{in } \Omega,\\
    w&=0, && \mbox{on } \Gamma.
\end{aligned}\right.
\end{equation*}
By Meyers' gradient estimate for a zero Dirichlet boundary problem \cite[Theorem 1]{Meyers},
we have $$
  \|\nabla w\|_{L^q(\Omega)^d} \leq C(\|A\nabla \bar{u}\|_{L^q(\Omega)^d}+\|f\|_{L^2(\Omega)})
  \leq C(\|A\|_{L^\infty(\Omega)^{d,d}}\|\nabla \bar{u}\|_{L^q(\Omega)^d}+\|f\|_{L^2(\Omega)}).$$
This and the triangle inequality complete the proof of the theorem.
\end{proof}

The following H\"{o}lder inequality for matrix- and vector-valued functions is useful.
\begin{lemma}\label{lem:Holder}
If $\frac1p+\frac1q=\frac12$, $p,q\geq1$, then for $A\in L^p(\Omega)^{d,d}$, $u \in W^{1,p}(\Omega)$, and $v\in H^1(\Omega)$, there holds
\begin{equation}\label{inp}
|( A\nabla u\cdot\nabla v)_{L^2(\Omega)^d}|\leq d^{\frac12+\frac1q}\normmm{A}_p\|\nabla u\|_{L^q(\Omega)^d}\|\nabla v\|_{L^2(\Omega)^d}.
\end{equation}
\end{lemma}
\begin{proof}
By the Cauchy-Schwarz inequality and H\"{o}lder inequality with conjugate exponents
satisfying $\frac1m+\frac1n=1$, we have
 \begin{align*}
   &\quad|(A\nabla u\cdot\nabla v)_{L^2(\Omega)^d}|\\
   &\leq \|A\nabla u\|_{L^2(\Omega)^d}\|\nabla v\|_{L^2(\Omega)^d}
   =\Big(\int_\Omega\sum_{i=1}^d(\sum_{j=1}^dA_{ij}\partial_ju)^2dx\Big)^{\frac12}\|\nabla v\|_{L^2(\Omega)^d}\\
   &\leq\Big(\int_\Omega\sum_{i=1}^d(\sum_{j=1}^d|A_{ij}|^m)^{\frac2m}(\sum_{j=1}^d|\partial_ju|^n)^{\frac2n}dx\Big)^{\frac12}\|\nabla v\|_{L^2(\Omega)^d}\\
   &\leq\Big(\int_\Omega\big(\sum_{i=1}^d\big(\sum_{j=1}^d|A_{ij}|^m\big)^{\frac2m}\big)^mdx\Big)^{\frac{1}{2m}}
   \Big(\int_\Omega\big(\sum_{j=1}^d|\partial_ju|^n\big)^2dx\Big)^{\frac{1}{2n}}\|\nabla v\|_{L^2(\Omega)^d}\\
   &\leq \Big(\int_\Omega d^{m-1}\big(\sum_{i=1}^d\big(\sum_{j=1}^d|A_{ij}|^m\big)^2\big) dx\Big)^{\frac{1}{2m}}
   \Big(\int_\Omega d(\sum_{j=1}^d|\partial_ju|^{2n})dx\Big)^{\frac{1}{2n}}\|\nabla v\|_{L^2(\Omega)^d}\\
   &\leq \Big(\int_\Omega d^m\big(\sum_{i,j=1}^d|A_{ij}|^{2m}\big)dx\Big)^{\frac{1}{2m}}
   d^{\frac1{2n}}\Big(\int_\Omega \sum_{j=1}^d|\partial_ju|^{2n}dx\Big)^{\frac{1}{2n}}\|\nabla v\|_{L^2(\Omega)^d}\\
   &=d^{\frac12+\frac1{2n}}\Big(\int_\Omega\sum_{i,j=1}^d|A_{ij}|^{2m}
   dx\Big)^{\frac{1}{2m}}\Big(\int_\Omega\sum_{j=1}^d|\partial_ju|^{2n}dx\Big)^{\frac{1}{2n}}\|\nabla v\|_{L^2(\Omega)^d}.
 \end{align*}
Upon taking $q=2n$ and $p=2m$, the desired inequality follows.
\end{proof}

Now we can derive the
Lipschitz continuity of the parameter-to-state map $F(A)$ with respect to the
$\normmm{\cdot}_p$ -norm for any $p\in(\frac{2Q_d}{Q_d-2}, \infty]$.
\begin{lemma}\label{lem:LipCont}
For any $f\in L^2(\Omega)$ and $g\in H^1(\Gamma)$,
the mapping $F: A\mapsto u(A)$ is Lipschitz continuous from $(\mathcal{A},L^p(\Omega)^{d\times d})$
to $H^1(\Omega)$, for any $p\in (\frac{2Q_d}{Q_d-2}, \infty]$. That is, for any $A_1, A_2\in\mathcal{A}$,
\begin{equation*}
  \|F(A_2)-F(A_1)\|_{H^1(\Omega)}\leq C(\alpha, \beta, d, p,\Omega){(\|f\|_{L^2(\Omega)}+\|g\|_{H^1(\Gamma)})}\normmm{A_2-A_1}_p.
\end{equation*}
\end{lemma}
\begin{proof}
Let $u_1=F(A_1)$ and $u_2=F(A_2)$ be the weak
solutions to problem \eqref{pr1} with $A_1$ and $A_2$, respectively. By the weak formulations, we deduce
\begin{equation}\label{vf2}
  (A_1(\nabla u_2-\nabla u_1),\nabla v)_{L^2(\Omega)^d}=-((A_2-A_1)\nabla u_2,\nabla v)_{L^2(\Omega)^d},\quad \forall v\in H_0^1(\Omega).
\end{equation}
Take $q\in (2,Q_d)$ such that $\frac1p+\frac1q+\frac12=1$, and let $v=u_2-u_1\in H^1_0(\Omega)$ in \eqref{vf2}.
Then the definition of $\mathcal{A}$ and Lemma \ref{lem:Holder} imply
\begin{equation*}
  \alpha\|\nabla(u_2-u_1)\|^2_{L^2(\Omega)^d}\leq d^{\frac12+\frac1q}\normmm{A_2-A_1}_p\|\nabla u_2\|_{L^q(\Omega)^d}\|\nabla (u_2-u_1)\|_{L^2(\Omega)^d}.
\end{equation*}
By Theorem \ref{recmeyer}, there exists a constant $C=C(\alpha, \beta, d, q,\Omega)$ such that
\begin{equation*}
  \|u_2-u_1\|_{H^1(\Omega)}\leq C\normmm{A_2-A_1}_p(\|f\|_{L^2(\Omega)}+\|g\|_{H^1(\Gamma)}),
\end{equation*}
which directly gives the desired assertion.
\end{proof}

Next we show the directional differentiability of the map $F$. Fix $A\in \mathcal{A}$,
and let $H\in L^\infty(\Omega)^{d,d}$ be a feasible direction such that $A+tH\in \mathcal{A}$ for small $t>0$.
Let $ w\in H_0^1(\Omega)$ be the weak solution to the linearized problem
\begin{equation}\label{pro2}
\left\{
    \begin{aligned}
    \mathrm{div}(A\nabla w)&=-\mathrm{div}(H\nabla u), &&\mbox{in } \Omega,\\
    w&=0, &&\mbox{on } \Gamma.
    \end{aligned}\right.
\end{equation}
This allows defining a linear map $DF(A):L^p(\Omega)^{d,d}\mapsto H^1_0(\Omega)$ by $DF(A)[H]=w$.
Now we show that the map $DF(A)$ is bounded, and represents the derivative of $F$ with
respect to $A$ in a generalized sense.
\begin{proposition}\label{prop:deriv}
Let $f\in L^2(\Omega)$, and $g\in H^1(\Omega)$. Then for any $p\in (\frac{2Q_d}{Q_d-2}, \infty]$, the
operator $DF(A):L^p(\Omega)^{d,d}\mapsto H^1_0(\Omega)$ is bounded:
\begin{equation*}
  \|w\|_{H^1(\Omega)}\leq C(\alpha, \beta, d, p, \Omega)\normmm{H}_p(\|f\|_{L^2(\Omega)}+\|g\|_{H^1(\Gamma)}).
\end{equation*}
 Further, for any $p\in(\frac{4Q_d}{Q_d-2},\infty]$, the map $F(A)$ is
differentiable in the sense that for any $A, A+H\in \mathcal{A}$, there holds
\begin{equation}\label{frechet}
  \lim_{ \normmm{H}_p\rightarrow 0}\frac{\|F(A+H)-F(A)-DF(A)[H]\|_{H^1(\Omega)}}{\normmm{H}_p}=0.
\end{equation}
\end{proposition}
\begin{proof}
Fix $H\in L^\infty(\Omega)^{d,d}$, $f\in L^2(\Omega)$ and
$g\in H^1(\Gamma)$. The weak formulation of problem \eqref{pro2} is to find $w\in H_0^1(\Omega)$ such that
\begin{equation*}
  (A\nabla w,\nabla v)_{L^2(\Omega)^d}=-(H\nabla u,\nabla v)_{L^2(\Omega)^d},\quad \forall v\in H_0^1(\Omega).
\end{equation*}
Choose $q\in(2,Q_d)$ such that $\frac1p+\frac1q+\frac12=1$. Then
Lemma \ref{lem:Holder} and letting $v=w\in H^1_0(\Omega)$ give
\begin{equation*}
  \alpha\|\nabla w\|^2_{L^2(\Omega)^d}\leq d^{\frac12+\frac1q}\normmm{H}_p\|\nabla u\|_{L^q(\Omega)^d}\|\nabla w\|_{L^2(\Omega)^d}.
\end{equation*}
By Theorem \ref{recmeyer}, there exists a constant $C=C(\alpha, \beta, d, q,\Omega)$ such that
\begin{equation*}
  \|DF(A)[H]\|_{H^1(\Omega)}\leq C\normmm{H}_p(\|f\|_{L^2(\Omega)}+\|g\|_{H^1(\Gamma)}).
\end{equation*}
This shows the first assertion. Let $R=F(A+H)-F(A)-DF(A)[H].$
Then $R$ satisfies
\begin{equation*}\left\{
    \begin{aligned}
    \mathrm{div}((A+H)\nabla R)&=-\mathrm{div}(H\nabla w),   &&\mbox{in } \Omega,\\
    R&=0, &&\mbox{on } \Gamma.
    \end{aligned}\right.
\end{equation*}
With $\frac1p+\frac1q=\frac12$, the preceding argument leads to
$\|\nabla R\|_{L^2(\Omega)^d}\leq C\normmm{H}_p\|w\|_{L^q(\Omega)^d}$.
Since $w$ satisfies \eqref{pro2}, by \cite[Theorem 1]{Meyers}, we deduce
\begin{equation*}
  \|\nabla w\|_{L^q(\Omega)^d}\leq C\|H\nabla u\|_{L^q(\Omega)^d}.
\end{equation*}
The choice of $p\in(\frac{4Q_d}{Q_d-2},\infty]$ implies $q\in(2, \frac{4Q_d}{2+Q_d})\subset (2, Q_d)$.
By choosing $\varepsilon>0$ such that $\bar{q}=q+\varepsilon=\frac{2q}{4-q}<Q_d$, we have, again by Theorem \ref{recmeyer},
\begin{align*}
  \|H\nabla u\|^q_{L^q(\Omega)^d}&\leq \normmm{H}^q_{\frac{q\bar{q}}\varepsilon}\|\nabla u\|^q_{L^{\bar{q}}(\Omega)^d}
  \leq C\normmm{H}^{\frac{p\varepsilon}{\bar{q}}}_p\|\nabla u\|^q_{L^{\bar{q}}(\Omega)^d}
  \leq C\normmm{H}^q_p(\|f\|_{L^2(\Omega)}+\|g\|_{H^1(\Gamma)})^q.
\end{align*}
The preceding estimates together give
$$\|R\|_{H^1(\Omega)} \leq C\normmm{H}^2_p(\|f\|_{L^2(\Omega)}+\|g\|_{H^1(\Gamma)})
,$$
completing the proof of the proposition.
\end{proof}

\subsection{Well-posedness of problem \eqref{tik}}

Note that the regularized functional $J_\gamma$ involves only the $L^2(\Omega)^{d,d}$ penalty
$\normmm{A}^2_2$, which is weaker than the commonly used Sobolev smooth penalty
or total variation penalty, and does not induce very strong compactness to ensure
the weak sequential closeness of the parameter-to-state map, which is commonly used for
analyzing the well-posedness of the optimization problem \cite{ItoJin:2015}. To justify the approach,
we employ the concept of H-convergence. First, we recall the concept of
H-convergence \cite{Tartar} for elliptic problems with a zero Dirichlet
boundary condition. Note that the H-limit is unique.
\begin{definition}\label{def1}
A sequence $\{A^n\}_{n\in\mathbb{N}} \subset \mathcal{A}$ is said to be
H-convergent to $A\in \mathcal{A}$, denoted by $A^n\xrightarrow{H} A$,
if for all $f\in (H_0^1(\Omega))'$, we have $u_n\rightharpoonup u$ in $H^1(\Omega)$, and
$A^n\nabla u_n\rightharpoonup A\nabla u$ in $L^2(\Omega)^d$, where $u_n$ is the solution to
 \begin{equation*}\left\{
    \begin{aligned}
   -\mathrm{div}(A^n\nabla u_n)&=f, &&\mbox{in } \Omega,\\
    u_n&= 0, &&\mbox{on } \Gamma,
    \end{aligned}\right.
\end{equation*}
and $u$ solves the problem with $A^n$ replaced by $A$.
\end{definition}

The following verion of div-curl lemma plays an important role in the analysis.
\begin{lemma}\label{curl-div}
Let $\{U_n = \nabla u_n\}_{n\in\mathbb{N}}$ and $\{V_n\}_{n\in\mathbb{N}}$ be two bounded sets
in $L^2(\Omega)^d$ such that
\begin{equation*}
 (u_n,V_n) \rightharpoonup (u_\infty, V_\infty) \mbox{ in } H^1(\Omega)\times L^2(\Omega)^d\quad\mbox{and}\quad
 \mathrm{div} V_n\rightarrow \mathrm{div} V_\infty\mbox{ in }H^{-1}(\Omega).
\end{equation*}
Then for any $\phi\in C_0^\infty(\Omega)$, we have
\begin{equation*}
\lim_{n\rightarrow \infty} \int_\Omega \phi \nabla u_n \cdot V_n dx= \int_\Omega \phi \nabla u_\infty \cdot V_\infty dx.
\end{equation*}
\end{lemma}
\begin{proof}
The proof can be found in \cite[Page 91]{LTartar}. We sketch the proof for completeness. Since $\phi u_n  \rightharpoonup \phi u_\infty $ in $H^1_0(\Omega)$, the strongly-weakly convergence in $\langle \cdot, \cdot\rangle_{H^{-1}(\Omega),H^1_0(\Omega)}$ implies
\begin{equation*}
\lim_{n\rightarrow \infty}\langle \mathrm{div}V_n, \phi u_n\rangle_{H^{-1}(\Omega),H^1_0(\Omega)} = \langle \mathrm{div}V_\infty, \phi u_\infty\rangle_{H^{-1}(\Omega),H^1_0(\Omega)}.
\end{equation*}
This directly gives
\begin{equation*}
\lim_{n\rightarrow \infty}( V_n, \phi \nabla u_n + u_n\nabla \phi)_{L^2(\Omega)^d} =( V_\infty, \phi \nabla u_\infty + u_\infty\nabla \phi)_{L^2(\Omega)^d} .
\end{equation*}
Meanwhile, by compact Sobolev embedding $H^1(\Omega)\hookrightarrow L^2(\Omega)$,
$( V_n, u_n\nabla \phi)_{L^2(\Omega)^d}$ converges to $( V_\infty, u_\infty\nabla \phi)_{L^2(\Omega)^d}$, and the proof is completed.
\end{proof}

To prove the existence of a minimizer of problem \eqref{tik}, we need  the following
compactness result, which is the counterpart of H-convergence for a nonzero Dirichlet
boundary condition.
\begin{theorem}\label{compact}
For any sequence $\{A^n\}_{n\in \mathbb{N}}\subset \mathcal{A}$,
there exists a subsequence, still denoted by $\{A^{n}\}$,
and an element $A\in \mathcal{A}$ such that
for every $(f,g)\in (H^1(\Omega))^\prime\times H^{\frac12}(\Gamma)$, there holds
\begin{equation}\label{Hconvergent}
  u_{n}\rightharpoonup u \ in \ H^1(\Omega) \ and \ A^{n}\nabla u_{n}\rightharpoonup A\nabla u \ in \ L^2(\Omega)^d,
\end{equation}
where $u_n$ and $u$ be solutions of problem \eqref{pr1} with conductive matrices $A^n$ and $A$, respectively.
\end{theorem}
\begin{proof}
Since $\{A^n\}_{n\in\mathbb{N}}\subset \mathcal{A}$, there exists a subsequence, still denote by $A^n$, and $A\in\mathcal{A}$
such that $A^n\xrightarrow{H} A$ {\cite[Theorem 6.5]{LTartar}}.
First, we show that for any fixed $g\in H^{\frac12}(\Gamma)$, there exists a subsequence of $\{A^{n}\}_{n\in \mathbb{N}}$, still denoted by $\{A^{n}\}_{n\in \mathbb{N}}$,
such that \eqref{Hconvergent} holds, following the argument of \cite[Lemma 10.4]{LTartar}.
We define the bilinear forms $a_n(\cdot,\cdot), a(\cdot, \cdot): H^1(\Omega)\times H^1(\Omega)\to \mathbb{R}$ by
\begin{align*}
a_n(u,v)=( A^n\nabla u, \nabla v )_{L^2(\Omega)^d} \quad \mbox{and}\quad a(u,v)=(A\nabla u, \nabla v)_{L^2(\Omega)^d}, \quad \forall  u,v\in H^1(\Omega).
\end{align*}
The weak formulation of problem \eqref{pr1} (with $A^n$) is to find $u_n\in K:=\{u\in H^1(\Omega):
u=g \mbox{ on } \Gamma\}$ such that
\begin{equation*}
a_n(u_n,v)= \langle f,v\rangle_{(H^1(\Omega))',H^1(\Omega)}, \quad \forall v\in H^1_0(\Omega).
\end{equation*}
By Lemma \ref{lemma1}, the sequence $\{u_n\}_{n\in\mathbb{N}}$ is uniformly bounded in $H^1(\Omega)$.
Thus, we can extract a subsequence, again denoted by $u_n$, such that $u_n \rightharpoonup u$ in $H^1(\Omega)$.
Since the set $K$ is closed and convex, it is weakly closed, and we deduce $u\in K$. Since the space $C^\infty_c(\Omega)$
is dense in $L^2(\Omega)$, for any $\phi\in L^2(\Omega)^d$, there exists a sequence $\{\phi_\varepsilon\}_{\varepsilon>0}
\subset C^\infty_c(\Omega)^d$ such that $\phi_\varepsilon\rightarrow \phi$ in $L^2(\Omega)^d$ as $\varepsilon\to0^+$.
For all $\phi \in L^2(\Omega)^d$,
\begin{align*}
|(\phi,A^n\nabla u_n-A\nabla u)_{L^2(\Omega)^d}|
&\leq \|\phi-\phi_\varepsilon\|_{L^2(\Omega)^d}(\|A^n\nabla u_n\|_{L^2(\Omega)^d}+
\|A\nabla u\|_{L^2(\Omega)^d}) \\
&\quad +|(\phi_\varepsilon, A^n\nabla u_n-
A\nabla u)_{L^2(\Omega)^d}|.
\end{align*}
By \cite[Lemma 10.3]{LTartar}, we obtain
\begin{equation}\label{wc}
A^n\nabla u_n\rightharpoonup A \nabla u\quad \mbox{in } L^2(\Omega)^d.
\end{equation}
Next we show
\begin{equation}\label{wc2}
a_n(u_n,u_n-v) \rightarrow a(u,u-v),\quad  \forall  v\in H^1(\Omega).
\end{equation}
Since $A^n\nabla u_n\rightharpoonup A \nabla u,
\nabla u_n\rightharpoonup \nabla u$ in $L^2(\Omega)^d$ and $\mathrm{div}\;(A^n \nabla u_n) = f$, and $\mathrm{curl}\;(\nabla u_n)=0$,
 the div-curl lemma (cf. Lemma \ref{curl-div}) implies
\begin{equation}
( \phi A^n\nabla u_n,\nabla u_n)_{L^2(\Omega)^d}\rightarrow (\phi A \nabla u, \nabla u)_{L^2(\Omega)^d}\quad \forall \phi\in C_c(\Omega).
\end{equation}
Repeating the argument of \cite[Lemma 10.4]{LTartar} and noting
$$a_n(u_n,v) \rightarrow a(u,v),$$
we obtain the assertion desired \eqref{wc2}. Thus we have
\begin{equation*}
  a(u,v)=\langle f, v \rangle_{(H^1(\Omega))',H^1(\Omega)},\quad \forall v \in H^1_0(\Omega).
\end{equation*}
Next we apply a density / diagonal argument to show that for all
$(f,g)\in (H^1(\Omega))^\prime\times H^{\frac12}(\Gamma)$, the assertion
\eqref{Hconvergent} holds. Since the space $H^{\frac12}(\Gamma)$ is separable,
there exists a countable dense subset $G$ of $H^{\frac12}(\Gamma)$. Let $\{0_n\} = \{n\}$.
For any fixed $\{g^k\}\subset {G}$ ($k\geq 1$), the preceding argument ensures the existence of a subsequence
 $\{(k)_{n}\}\subset \{(k-1)_n\}$ such that
 $$u_{k_n}^{k}\rightharpoonup u^k  \;\mbox{ in }\; H^1(\Omega) \quad  \textrm{and} \quad
A^{k_n}\nabla u_{k_n}^{k}\rightharpoonup A\nabla u^k \;\mbox{ in }\; L^2(\Omega)^d,$$
where $u_{k_n}^{k}$ and $u^k$ satisfy
\begin{eqnarray*}
\left\{\begin{aligned}
  -\mathrm{div}(A^{k_n}\nabla u_{k_n}^{k})&=f && \mbox{in } \Omega,\\
   u_{k_n}^{k}&=g^k,&& \mbox{on }\Gamma,
\end{aligned}\right. \quad \mbox{and}\quad
\left\{\begin{aligned}
  -\mathrm{div}(A\nabla u^k)&=f, &&\mbox{in } \Omega,\\
   u^k&=g^k,&& \mbox{on }\Gamma,
\end{aligned}
\right.
\end{eqnarray*}
respectively. To show the convergence for any $(f,g)\in (H^1(\Omega))^\prime\times H^{\frac12}(\Gamma)$, we apply a diagonal argument, i.e., choosing the subsequence $\{k_k\}$, still denoted by $A^n$, i.e., $A^n \equiv A^{n_n}$. Since $G$ is dense in $H^{\frac12}(\Gamma)$, for any
$g\in H^{\frac12}(\Gamma)$ there exists $\{g^k\}\subset G$ such that ${g^k}\rightarrow g$ in $H^{\frac12}(\Gamma)$. For any $\xi\in (H^1(\Omega))^\prime$, we have
\begin{align*}
&\quad\big|\langle \xi, u_n\rangle_{(H^1(\Omega))',H^1(\Omega)} - \langle \xi, u\rangle_{(H^1(\Omega))',H^1(\Omega)} \big|\\
&= \big|\langle \xi, u_n-u_{n}^{k}\rangle_{(H^1(\Omega))',H^1(\Omega)} +\langle \xi, u_{n}^{k}-u^k\rangle_{(H^1(\Omega))',H^1(\Omega)} +\langle \xi, u^k-u\rangle_{(H^1(\Omega))',H^1(\Omega)} \big| \\
&\leq \|\xi\|_{(H^1(\Omega))'}(\|u_n-u_{n}^{k}\|_{H^1(\Omega)}+\|u^k-u\|_{H^1(\Omega)})+
{|\langle \xi, u_{n}^{k}-u^k\rangle_{(H^1(\Omega))',H^1(\Omega)}|} \\
&\leq C\|g-g^k\|_{H^{\frac12}(\Gamma)}+{|\langle \xi, u_{n}^{k}-u^k\rangle_{(H^1(\Omega))',H^1(\Omega)}|},
\end{align*}
which yields $u_n\rightharpoonup u$ in $H^1(\Omega)$.
Since $\|A^n\nabla u_n-A^n\nabla {u_{n}^k}\|_{L^2(\Omega)^d}\leq C\|u_n-{u_{n}^k}\|_{H^1(\Omega)}$ and
repeating the argument gives
$A^n\nabla u_n\rightharpoonup A \nabla u$ in $L^2(\Omega)^d$. This completes the proof of the theorem.
\end{proof}

Theorem \ref{compact} gives the the H-convergence for a nonzero Dirichlet boundary data, and further,
the H-limit for the nonzero Dirichlet boundary condition case is identical with the zero case.
The next lemma gives the norm inequality for the H-limit \cite{Tartar}.
\begin{lemma}\label{compare-cont}
Let the sequence $\{A^n\}_{n\in \mathbb{N}}\subset \mathcal{A}$ be
$A^n\xrightarrow{H} A$ and $A^n  \overset{\ast}{\rightharpoonup} A^0$ in $L^\infty(\Omega)^{d,d}$.
Then
\begin{equation}\label{liminequ}
 A\leq A^0 \mbox{ a.e. in } \Omega\quad \mbox{and}\quad \normmm{A}^2_2\leq \normmm{A^0}^2_2\leq \liminf_{n\rightarrow \infty}\normmm{A^n}^2_2.
\end{equation}
\end{lemma}

With the compactness result in Theorem \ref{compact}, we can state the existence of a minimizer.
\begin{theorem}\label{uniqsol}
There exists at least one minimizer to problem \eqref{tik}.
\end{theorem}
\begin{proof}
The functional $J_\gamma$ is bounded from below by zero, and thus we can find a minimizing sequence
$\{A^n\}_{n\in\mathbb{N}}\subset \mathcal{A}$ such that
$$\lim_{n\rightarrow\infty}J_\gamma(A^n)=\inf_{A\in \mathcal{A}}J_\gamma(A).$$
By Theorem \ref{compact}, since the sequence $\{A^n\}_{n\in \mathbb{N}}$ is bounded
in $L^\infty(\Omega)^{d,d}$, there exists a subsequence $\{A^{n_k}\}_{k\in\mathbb{N}}$ and
some $A\in \mathcal{A}$ such that $A^{n_k}\xrightarrow{H} A$ in $L^\infty(\Omega)^{d,d}$.
Letting $u^{\ell}_{n_k}=u^{\ell}(A^{n_k}), \ell=1,...,L$, we have
$u^{\ell}_{n_k}\rightharpoonup u^{\ell}(A), \ell=1,...,L$ in $H^1(\Omega)$.
Hence, the weak lower semi-continuity of the $L^2(\Omega)$-norm and Lemma \ref{compare-cont} give
\begin{align*}
J_\gamma(A)&=\frac12\sum_{\ell=1}^L\|A\nabla u^{\ell}-\vec{h}^{\ell}\|^2_{L^2(\Omega)^d}
+\frac{\gamma}{2}\normmm{A}^2_2\\
&\leq \frac12\liminf_{n^\prime\rightarrow\infty}\sum_{\ell=1}^L\|A^{n^\prime} \nabla u^{\ell}_{n^\prime}-\vec{h}^{\ell}\|^2_{L^2(\Omega)^d}+\frac\gamma 2\liminf_{n^\prime\rightarrow\infty}\normmm{A^{n^\prime}}^2_2\\
&\leq \liminf_{n^\prime\rightarrow\infty}J_\gamma(A^{n^\prime})=\inf_{A\in \mathcal{A}}J_\gamma(A).
\end{align*}
Thus, $A$ is a global minimizer. This completes the proof of the theorem.
\end{proof}

Next we derive the expression of the gradient $J_\gamma'$ of the functional $J_\gamma$ and
the first-order necessary optimality system using the adjoint technique. The former allows
applying popular gradient-descent type algorithms, whereas the latter is useful for designing
Newton type methods. Let the adjoint variable $\bar p^\ell$, $\ell=1,...,L$, solve
\begin{equation}\label{adjoint}
\left\{\begin{aligned}
  -\mathrm{div}(A\nabla \bar p^\ell)&= \mathrm{div}(A(A\nabla u^\ell-\vec{h}^\ell)), &&\mbox{in } \Omega,\\
   \bar{p}^\ell &= 0, &&\mbox{on } \Gamma.
  \end{aligned}\right.
\end{equation}
We show the differentiability of $J_\gamma(A)$ in the $L^p(\Omega)^{d,d}$
topology. For two vectors $a, b\in \mathbb{R}^d$, we denote the symmetrized tensor product by
$(a\otimes b)_{ij}=\frac12(a_ib_j+a_jb_i)$, $i,j=1,...,d$, and the Frobenius inner product between
two matrices $A,B\in\mathbb{R}^d$ by $A\cdot B$. The lengthy but routine proof is deferred to
the appendix.

\begin{theorem}\label{thm:diff-tikh}
Let $\bar{p}^\ell$ be defined in \eqref{adjoint},
for any  $H\in L^\infty(\Omega)^{d,d}$ such that $A+tH\in \mathcal{A}$ for sufficiently small $t>0$,
the directional derivative $J_\gamma^\prime(A)[H]$ of $J_\gamma$ is given by
\begin{equation}\label{dj}
 J_\gamma^\prime(A)[H] =\int_\Omega\big(\sum_{\ell=1}^L(\nabla u^\ell\otimes \nabla \bar{p}^\ell+\nabla u^\ell\otimes (A\nabla u^\ell-\vec{h}^\ell)+\gamma A\big)\cdot Hdx.
\end{equation}
Moreover, if $A$ is an interior point of $\mathcal{A}$, then $J_\gamma(A)$ is Fr\'{e}chet
differentiable in the $L^p(\Omega)^{d,d}$ topology, for any $p\in(\frac{4Q_d}{Q_d-2},\infty]$.
\end{theorem}

Next we derive the first-order necessary optimality system of problem \eqref{tik}.
Recall the subset $K:=\{A\in \mathcal{S}_d|\alpha I\leq A\leq \beta I\}
\subset\mathcal{S}_d$ given in \eqref{uec}, with $A\leq B$ indicating that $B-A$ is symmetric
positive semidefinite. Since $K$ is a convex and closed subset of
$\mathcal{S}_d$, we can define an orthogonal projection
$P_K: \mathcal{S}_d\rightarrow K$, characterized by the following variational inequality
\begin{equation*}
  (A-P_K(A))\cdot (B-P_K(A))\leq 0 ,\quad \forall  B\in K.
\end{equation*}
\begin{theorem}\label{kkt}
Let $A^*\in \mathcal{A}$ solve problem \eqref{tik}, then the tuple $(A^*,u^{\ell}, \bar{p}^{\ell})$ satisfies the following
optimality system, for every $\lambda>0$,
\begin{align*}
  & \left\{\begin{aligned}
    -\mathrm{div}(A^*\nabla u^{\ell})&=f^{\ell}   &&\mbox{in }\Omega,\\
    u^{\ell}&=g^{\ell}   &&\mbox{on } \Gamma,
    \end{aligned}\right. \tag*{primal}\\
   & \left\{\begin{aligned}
    -\mathrm{div}(A^*\nabla {p}^{\ell})&=\mathrm{div}(A^*(A^*\nabla u^{\ell}-\vec{h}^{\ell})),  &&\mbox{in } \Omega,\\
    {p}^{\ell}&=0, &&\mbox{on } \Gamma,
    \end{aligned}\right. \tag*{dual}\\
    &A^*(x)=P_K\Big(A^*-\lambda\sum_{\ell=1}^L(\nabla u^{\ell}\otimes (A^*\nabla u^{\ell}-\vec{h}^{\ell})+\nabla u^{\ell}\otimes \nabla {p}^{\ell}+\gamma A^*)\Big) \ a.e. \ in \ \Omega.\tag*{complementarity}
\end{align*}
\end{theorem}
\begin{proof}
The optimality of $A^*$ implies
$J_\gamma^\prime(A^*)(B-A^*)\geq 0$ for all $ B\in \mathcal{A}.$
In view of \eqref{dj}, it can be rewritten as
\begin{equation*}
  \int_\Omega\Big(A^*-\lambda\sum_{\ell=1}^L\big(\nabla u^{\ell}\otimes (A^*\nabla u^{\ell}-\vec{h}^{\ell})+\nabla u^{\ell}\otimes \nabla {p}^{\ell}\big)+\gamma A^*\big)-A^*\Big)\cdot(B-A^*)dx\leq 0,\quad \forall  B\in \mathcal{A}.
\end{equation*}
The inequality is equivalent to
\begin{equation}\label{funspace}
 A^*=P_\mathcal{A}\Big(A^*-\lambda\sum_{\ell=1}^L\big(\nabla u^{\ell}\otimes (A^*\nabla u^{\ell}-\vec{h}^{\ell})+\nabla u^{\ell}\otimes \nabla {p}^{\ell}\big)+\gamma A^*\big)\Big),
\end{equation}
since $P_\mathcal{A}$ is a projection onto a closed convex subset
in a Hilbert space. Now for every $x\in \Omega$,
\begin{align*}
&P_K\Big(A^*(x)-\lambda\sum_{\ell=1}^L\big(\nabla u^{\ell}(x)\otimes (A^*(x)\nabla u^{\ell}(x)-\vec{h}^{\ell}(x))+\nabla u^{\ell}(x)\otimes \nabla {p}^{\ell}(x)\big)+\gamma A^*(x)\big)\Big) \\
&=\arg\min_{B\in K} \Big \vert\kern-0.25ex\Big\vert\kern-0.25ex\Big\vert A^*(x)-\lambda\sum_{\ell=1}^L\big(\nabla u^{\ell}(x)\otimes (A^*(x)\nabla u^{\ell}(x)-\vec{h}^{\ell}(x))+\nabla u^{\ell}(x)\otimes \nabla {p}^{\ell}(x)\big)\\
  &\quad +\gamma A^*(x)\big)-B\Big\vert\kern-0.25ex\Big\vert\kern-0.25ex\Big\vert_p.
\end{align*}
This and the uniqueness of a projection onto convex sets imply
that the relation \eqref{funspace} is indeed equivalent to the pointwise projection
$
 A^*(x)=P_K\big(A^*(x)-\lambda\sum_{\ell=1}^L\big(\nabla u^{\ell}(x)\otimes (A^*(x)\nabla u^{\ell}(x)-\vec{h}^{\ell}(x))+\nabla u^{\ell}(x)\otimes \nabla {p}^{\ell}(x)\big)+\gamma A^*(x)\big)\big)
$
a.e. in $\Omega$,  which directly implies the desired complementary condition.
\end{proof}

\begin{remark}
Choosing $\lambda=\frac1 \gamma$ in the complementarity condition gives
\begin{equation}\label{A}
 A^*(x)=P_K\Big(-\gamma^{-1} \sum_{\ell=1}^L\nabla u^{\ell}\otimes (A^*\nabla u^{\ell}-\vec{h}^{\ell})+\nabla u^{\ell}\otimes \nabla {p}^{\ell}\Big) \quad\mbox{ a.e. \ in }\ \Omega.
\end{equation}
This choice is employed in the numerical implementation in Section 5.
\end{remark}

Last, we state a convergence result with respect to the noise level. Suppose that we are given a set of noisy
measurements $\vec{h}^{\ell,\delta} \in L^2(\Omega)^d, \ell=1,...,L$, with a noise level $\delta$, i.e.,
\begin{equation}\label{delta}
  \|\vec{h}^{\ell}-\vec{h}^{\ell,\delta}\|_{L^2(\Omega)^d}\leq \delta,\quad \ell=1,\ldots,L.
\end{equation}
Accordingly, consider the following optimization problem:
\begin{equation}\label{probdel}
\min_{A\in \mathcal{A}} \bigg\{J_\gamma^\delta(A):=\frac12\sum_{\ell=1}^L\|A\nabla u^{\ell}(A)-\vec{h}^{\ell,\delta}\|^2_{L^2(\Omega)^d}+\frac{\gamma}{2}\normmm{A}^2_2\bigg\}.
\end{equation}

The next result shows the convergence of regularized solutions to a minimum norm solution as
the noise level $\delta$ tends to zero, i.e., regularizing property of problem \eqref{tik}.
The minimum norm solution $A^\dag$ is defined by
$A^\dag\in\textrm{arg}\min_{A\in \mathcal{I}_{\mathcal{A}}(\vec{h}^L)}\normmm{A}_2^2,$
where the set $\mathcal{I}_{\mathcal{A}}(\vec{h},L)=\{A\in\mathcal{A}|\vec{h}^\ell=A\nabla u^\ell(A), \ell = 1,...,L\}$
is assumed to be nonempty. {This result ensures that the solutions to the regularized problem
\eqref{tik} don't differ too much from the reference solution, provided that the noise level $\delta$
is sufficiently small and the regularization parameter $\gamma$ is chosen properly.}

\begin{theorem}\label{stability}
 Let $\vec{h}^{\ell,\delta}$ satisfy \eqref{delta}, $A^\delta\in \mathcal{A}$ be a solution of problem \eqref{probdel} and suppose that
 $\gamma\to 0^+$ and $\frac{\delta}{\sqrt{\gamma}}\to 0$, as $\delta\to 0^+$.
 Then up to a subsequence, there holds $A^{\delta}\rightarrow A$ in $L^2(\Omega)^{d,d}$ as $\delta\rightarrow 0^+$, for some minimum norm solution $A$.
\end{theorem}
\begin{proof}
Given a sequence $\{\delta_k\}_{k\in \mathbb{N}}\subset\mathbb{R}_+$ with $\lim_{k\rightarrow \infty}\delta_k=0^+$, we choose $\gamma_k>0$ such that $\lim_{k\to\infty}\gamma_k=0$ and $\lim_{k\to\infty}\gamma_k^{-\frac12}\delta_k=0$.
Further, let $A_k:=A^{\delta_k}_{\gamma_k}\in \mathcal{A}$ be a minimizer of $J_{\gamma_k}^{\delta_k}(A)$ over $\mathcal{A}$.
For any minimum norm solution $A^\dag\in \mathcal{A}$, there  holds
$J_{\gamma_k}^{\delta_k}(A_k)\le J_{\gamma_k}^{\delta_k}(A^\dag).$ That is,
\begin{align}\label{ineq}
&\quad\frac12\sum_{\ell=1}^L\|A_k\nabla u^{\ell}(A_k)-\vec{h}^{\ell,\delta_k}\|^2_{L^2(\Omega)^d}
+\frac{\gamma_k}{2}\normmm{A_k}_2^2\nonumber\\
&\le \frac12\sum_{\ell=1}^L\|A^\dag\nabla u^{\ell}(A^\dag)-\vec{h}^{\ell,\delta_k}\|^2_{L^2(\Omega)^d}
+\frac{\gamma_k}{2}\|A^\dag\|^2_{L^2(\Omega)^{d,d}}\\
&\le \frac L2\delta_k^2 + \frac{\gamma_k}{2}\normmm{A^\dag}_2^2.\nonumber
\end{align}
This and the condition on $\gamma_k$ imply
\begin{equation}\label{form1}
\limsup_{k\to\infty}\normmm{A_k}_2^2
\le \limsup_{k\to\infty}(L\gamma_k^{-1}\delta_k^2 + \normmm{A^\dag}_2^2)
=\normmm{A^\dag}_2^2.
\end{equation}
By Theorem \ref{compact}, there exists a subsequence, again denoted by
 $\{A_k\}_{k\in \mathbb{N}}$, $A,A^0\in \mathcal{A}$ such that
$A_k\xrightarrow{H} A$ and $ A_k\to A^0\mbox{ weak $\ast$ in } L^\infty(\Omega)^{d,d}.$
 Lemma \ref{compare-cont} implies
\begin{equation}\label{inequ}
  A\leq A^0 \mbox{ a.e. in } \Omega \quad \mbox{and} \quad \normmm{A}_2^2\le \normmm{A_0}_2^2\le \liminf_{k\rightarrow \infty} \normmm{A_k}_2^2.
\end{equation}
By the definition of H-convergence, we have $u^{\ell}(A_k)\rightharpoonup u^{\ell}(A)$ in $H^1(\Omega)$, $A_k\nabla u^{\ell}(A_k)\rightharpoonup A\nabla u^{\ell}(A)$ in $L^2(\Omega)^d$.
Since $\|\vec{h}^{\ell}-\vec{h}^{\ell, \delta_k}\|_{L^2(\Omega)^d}\leq \delta_k$, from the estimate \eqref{ineq}, we deduce
 \begin{equation}\label{wcon}
   \sum_{\ell=1}^L\|A\nabla u^{\ell}(A)-\vec{h}^{\ell}\|_{L^2(\Omega)^d}\le
   \liminf_{k\rightarrow \infty}\sum_{\ell=1}^L\|A_k\nabla u^{\ell}(A_k)-\vec{h}^{\ell,\delta_k}\|_{L^2(\Omega)^d}=0.
 \end{equation}
Therefore, $A\in \mathcal{I}_{\mathcal{A}}(\vec{h},L)$, and it is a minimum norm solution.
From \eqref{form1} it follows that $A_k\to A^0\mbox{ weak $\ast$ in } L^\infty(\Omega)^{d,d}$ and $A \le A^0$ a.e. in $\Omega$. Then we have
\begin{align*}
\limsup_{k\to\infty}\normmm{A_k-A}_2^2
&=\limsup_{k\to\infty}(\normmm{A_k}^2_2+\normmm{A}^2_2)
-2(A_k,A)_{L^2(\Omega)^{d,d}}\\
&\le 2\normmm{A}^2_2-2(A_0,A)_{L^2(\Omega)^{d,d}}
\le 0
\end{align*}
In conclusion, $A_k\rightarrow A$ in $L^2(\Omega)^{d,d}$ as $k\rightarrow \infty$, and $A$ is a minimum norm solution.
\end{proof}
\begin{remark}
Note that the minimum norm solution $A^\dag$ is generally nonunique, and it does not necessarily
coincide with the exact conductivity tensor. Nonetheless, it coincides with the exact conductivity
tensor if the set $\mathcal{I}_{\mathcal{A}}(\vec{h},L)$ consists of one singleton. Under suitable
regularity conditions on the problem data and given a sufficient number of current densities, the
anisotropic conductivity tensor is uniquely determined \cite{bal2014,G.Bal2014}. Then the
set $\mathcal{I}_{\mathcal{A}}(\vec{h},L)$ consists of one single element, and accordingly,
the minimum norm solution is unique, and also by the standard subsequence argument,
the whole sequence converges to the unique minimum norm solution as $\delta \rightarrow 0^+$.
\end{remark}

\section{Numerical discretization and convergence}\label{sec:disc}
In practice, numerically solving problem \eqref{tik} requires suitable discretization.
Since the problem involves variable coefficients, it is most conveniently carried out using
the Galerkin finite element method \cite{Ciarlet,Brenner}. Throughout, we assume that the domain $\Omega$
is a polygon in 2D or polyhedron in 3D. Let $\mathcal{T}$ be a quasi-uniform
triangulation of the domain $\Omega$ with a maximum mesh size $h$.
We employ the piecewise linear finite element spaces $V_h\subset H^1(\Omega)$
and $X_h\subset L^2(\Omega)$ defined respectively by
\begin{align*}
  V_h&=\{v_h\in C(\overline{\Omega}): v_h|_T\in P_1(T), \forall T\in \mathcal{T}\}, \\
  X_h&=\{A_h: (A_h)_{i,j}\in L^2(\Omega), (A_h)_{i,j}|_T\in P_1(T), i,j=1,...,d, \forall T\in \mathcal{T}\},
\end{align*}
where $P_1(T)$ consists of all linear functions over $T$. The discrete admissible set $\mathcal{A}_h$ is given by
$\mathcal{A}_h=\mathcal{A}\cap X_h.$ Let $V_{h,0}=V_h\cap H_0^1(\Omega)\subset V_h$,
whose elements vanish on the boundary $\Gamma$, and denote by $V_h(\Gamma)$ the trace space of $V_h$
on $\Gamma$. Let $P_h: L^2(\Omega)\rightarrow V_{h,0}$ be the standard $L^2(\Omega)$
orthogonal projection. Then it is $H^1(\Omega)$-stable \cite{BX91}, i.e.,
$$\|P_h(u-\tilde{u})\|_{H^1_0(\Omega)}\leq C\|u-\tilde{u}\|_{H^1_0(\Omega)}.$$

To formulate the finite element approximation of problem \eqref{pr1}, we fix a function
$g_h\in V_h(\Gamma)$ which approximates the given Dirichlet boundary condition $g$.
Then the finite element problem for problem \eqref{pr1} reads: find
$u_h\in V_h$ such that $u_h=g_h$ on $\Gamma$ and
 \begin{equation}\label{ph}
  (A_h\nabla u_h,\nabla v_h)_{L^2(\Omega)^d}= (f,v_h)_{L^2(\Omega)}, \quad \forall v_h\in V_{h,0}.
\end{equation}
This defines a discrete forward operator $F_h: A_h\rightarrow u_h$ (again suppressing the dependence on the problem data $(f,g)$).

\begin{remark}
Since the space $V_h(\Gamma)$ is dense in $H^{\frac12}(\Gamma)$, there exists $g_h\in V_h(\Gamma)$
such that $\|g_h-g\|_{H^{\frac12}(\Gamma)}\rightarrow 0$ as $h\rightarrow 0^+$.
If $g\in H^{\frac32}(\Gamma)$, a standard choice for the Dirichlet data $g_h\in V_h(\Gamma)$
is the Lagrange interpolation $I_h g$ of $g$ \cite{Ciarlet,Brenner}.
 If $g\in H^{\frac12}(\Gamma)$, we may let $g_h$ be the $L^2(\Gamma)$-projection of
  $g$ onto $V_h(\Gamma)$ \cite{Fix1983,Bartels04}.
\end{remark}

Throughout, we take $g_h=\tilde{u}_h|_\Gamma$, with $\tilde{u}_h\in V_h$ and
$\|\tilde{u}_h-\tilde{u}\|_{H^1(\Omega)}\rightarrow 0$, where $\tilde{u}$ satisfies $-\Delta
\tilde{u}=0$ in $\Omega$ and $\tilde{u}|_\Gamma=g$ (i.e., $\tilde u$ is the harmonic
extension of $g$ from the boundary $\Gamma$ to the domain $\Omega$). Then by \cite[p. 143]{Ciarlet} and \cite[p. 200]{fix73}
and the density of $V_h$ in $H^1(\Omega)$, $V_{h,0}$ in $H^1_0(\Omega)$, we have
\begin{align*}
&\|u-u_h\|_{H^1(\Omega)}\leq \beta^\frac12 \alpha^{-\frac12}\inf_{v_h\in V_{h,0}}\|u-(\tilde{u}_h+v_h)\|_{H^1(\Omega)}\\
\leq&\beta^\frac12 \alpha^{-\frac12}\big(\inf_{v_h\in V_{h,0}}\|u-\tilde{u}-v_h\|_{H_0^1(\Omega)}
+\|\tilde{u}-\tilde{u}_h\|_{H^1(\Omega)}\big).
\end{align*}
Upon choosing $v_h=P_h(u-\tilde{u})\in V_{h,0}$, then we have
\begin{align*}
&\|u_h\|_{H^1(\Omega)}\leq \|u-u_h\|_{H^1(\Omega)}+\|u\|_{H^1(\Omega)}\\
\leq& \beta^\frac12 \alpha^{-\frac12}(2\|u-\tilde{u}\|_{H_0^1(\Omega)}
+\|\tilde{u}-\tilde{u}_h\|_{H^1(\Omega)})+\|u\|_{H^1(\Omega)}.
\end{align*}
Thus, $\|u_h\|_{H^1(\Omega)}$ is uniformly bounded.
This fact will be frequently used below.

Now we can formulate the discrete counterpart of problem \eqref{tik}:
\begin{equation}\label{discret}
\min_{A_h\in \mathcal{A}_h}\bigg\{ J_{\gamma,h}(A_h):=\frac12\sum_{\ell=1}^L\|A_h\nabla u^{\ell}_h(A_h)
  -\vec{h}^{\ell}\|^2_{L^2(\Omega)^d}+\frac{\gamma}{2}\normmm{A_h}^2_2\bigg\}.
\end{equation}

\subsection{Convergence analysis}

To analyze the convergence of problem \eqref{discret}, we use the concept of Hd-convergence defined below.
It was introduced in \cite{Eymard} and \cite{Deckelnick2011} in the context of finite volume and finite element
discretizations of boundary value problems, respectively, to adapt H-convergence results
to a discrete setting.
\begin{definition}\cite{Deckelnick2011}
A sequence $\{A_h\} \subset \mathcal{A}$ is said to be
Hd-convergent to $A\in \mathcal{A}$, denoted by $A_h\xrightarrow{Hd} A$,
if for all $f\in H^{-1}(\Omega)$, we have $u_h\rightharpoonup u$ in $H^1(\Omega)$, and
$A_h\nabla u_h\rightharpoonup A\nabla u$ in $L^2(\Omega)^d$, where $u_h\in V_{h,0}$
is the solution to the finite element problem
 \begin{equation*}
   (A_h\nabla u_h, \nabla v_h)_{L^2(\Omega)^d}=\langle f,v_h\rangle_{H^{-1}(\Omega),H^1_0(\Omega)},  \quad \forall v_h\in V_{h,0},
 \end{equation*}
and $u$ solves $-\mathrm{div}(A\nabla u)=f$ in $\Omega$ and $u=0$ on the boundary $\Gamma$.
\end{definition}

Next we establish the Hd-convergence for a nonzero Dirichlet boundary condition, which
will be useful for the convergence analysis of the discrete problem \eqref{discret}.
\begin{theorem}\label{dcompact}
Assume that $\{A_h\in \mathcal{A}_h\}_{h>0}\subset\mathcal{A}$ with $A_h\xrightarrow{Hd} A$, for any given
$f\in H^{-1}(\Omega)$, $g\in H^\frac{1}{2}(\Gamma)$, $g_h\in V_h(\Gamma)$ with $\|g_h - g\|_{H^\frac{1}{2}(\Gamma)} \rightarrow 0$.
Let $u_h\in V_h$ and $u\in H^1(\Omega)$ be the solutions of the following finite element problem and the elliptic problem \eqref{pr1}, respectively:
$u_h\in V_h$ with $u_h|_{\Gamma} = g_h$ with
 \begin{equation*}
   (A_h\nabla u_h, \nabla v_h)_{L^2(\Omega)^d}=\langle f,v_h\rangle_{H^{-1}(\Omega),H^1_0(\Omega)},  \quad \forall v_h\in V_{h,0},
 \end{equation*}
and $u$ solves $-\mathrm{div}(A\nabla u)=f$ in $\Omega$ and $u=g$ on the boundary $\Gamma$.
Then there hold
\begin{equation*}
  u_h \rightharpoonup u\ in\ H^1(\Omega)\quad\mbox{and}\quad A_{h}\nabla u_{h}\rightharpoonup A\nabla u \ in
  \ L^2(\Omega)^d.
\end{equation*}
\end{theorem}

\begin{proof}
Let $\bar{u} \in H^1(\Omega)$ with $\bar{u}|_{\Gamma} = g$. Then
\begin{equation*}
(A \nabla (u-\bar{u}), \nabla v)_{L^2(\Omega)^d} = \langle f + \mathrm{div} A\nabla \bar{u}, v\rangle_{H^{-1}(\Omega),H^1_0(\Omega)}, \quad \forall v\in H_0^1(\Omega).
\end{equation*}
Let $z_h\in V_{h,0}$ be the solution to the finite element problem
\begin{equation*}
(A_h \nabla z_h, \nabla v_h)_{L^2{\Omega}^d} = \langle f + \mathrm{div} A\nabla \bar{u}, v_h\rangle_{H^{-1}(\Omega),H^1_0(\Omega)}, \quad \forall v_h\in V_{h,0}.
\end{equation*}
By the Hd-convergence \cite[Theorem 3.1]{Deckelnick2011}, we have (up to a subsequence)
\begin{equation}\label{eq:tm1}
z_h \rightharpoonup u - \bar{u} \mbox{ in } H_0^1(\Omega) \quad \mbox{and}\quad A_h\nabla z_h \rightharpoonup A \nabla (u-\bar{u}) \mbox{ in } L^2(\Omega)^d.
\end{equation}
We denote by $\bar{u}_h \in V_h$ with $\bar{u}_h|_{\Gamma} = g_h$ and
\begin{equation*}
(A_h \nabla \bar{u}_h, \nabla v_h)_{L^2{\Omega}^d} = \langle - \mathrm{div} A\nabla \bar{u}, v_h\rangle_{H^{-1}(\Omega),H^1_0(\Omega)} = (A\nabla \bar{u}, \nabla v_h)_{L^2{\Omega}^d}, \quad \forall v_h\in V_{h,0}.
\end{equation*}
By definition, $A_h \nabla \bar{u}_h$ is the $L^2(\Omega)$-projection of $A\nabla \bar{u}$ onto the space $\nabla V_{h,0}:=\{\nabla v_h:v_h\in V_{h,0}\}$.
Since $\nabla V_{h,0}$ is the standard piecewise constant finite element space and it is dense in $L^2(\Omega)^d$ space, we deduce
\begin{equation}\label{eq:tm2}
A_h\nabla \bar{u}_h \rightharpoonup A\nabla \bar{u} \quad \mbox{in } L^2(\Omega)^d.
\end{equation}
Meanwhile, by the definition of $u_h$, we can split $u_h$ into $u_h = z_h + \bar{u}_h$, and furthermore, in view of \eqref{eq:tm1}-\eqref{eq:tm2}, we have
\begin{equation*}
 A_h\nabla u_h \rightharpoonup A\nabla u\mbox{ in }L^2(\Omega)^d.
\end{equation*}
It remains to show that $u_n \rightharpoonup u$ in $H^1(\Omega)$.
Since the sequence $\{u_h\}_{h>0}$ is bounded in $H^1(\Omega)$, there exists a weak accumulation point $\hat{u}\in H^1(\Omega)$, and it suffices to prove $\hat{u} = u$. Then up to a subsequence, let $u_h \rightharpoonup \hat{u}$ in $H^1_0(\Omega)$.
For any $\hat{f}\in H^{-1}(\Omega)$, let $y_h\in V_{h,0}$ and $y\in H^1_0(\Omega)$ are solutions to the finite element discretization and the continuous PDE, respectively, i.e.,
\begin{align*}
  (A_h\nabla y_h, \nabla v_h)_{L^2(\Omega)^d}&=\langle \hat{f},v_h\rangle_{H^{-1}(\Omega),H^1_0(\Omega)},  \quad \forall v_h\in V_{h,0},\\
  (A\nabla y, \nabla v)_{L^2(\Omega)^d}&=\langle \hat{f},v\rangle_{H^{-1}(\Omega),H^1_0(\Omega)},  \quad \forall v\in H^1_0(\Omega).
\end{align*}
Then the Hd-convergence \cite[Theorem 3.1]{Deckelnick2011} implies
\begin{equation*}
 y_h\rightharpoonup y\mbox{ in }H^1_0(\Omega)\quad\mbox{and}\quad A_h\nabla y_h \rightharpoonup A \nabla y\mbox{ in }L^2(\Omega)^d.
\end{equation*}
For any given $\phi\in C_0^\infty(\Omega)$, we define
\begin{equation*}
E_h = (\phi A_h \nabla y_h, \nabla u_h)_{L^2(\Omega)^d}.
\end{equation*}
Then direct computation gives
\begin{eqnarray*}
E_h & = &  (A_h\nabla y_h, \phi \nabla u_h)_{L^2(\Omega)^d} \\
 &= &(A_h\nabla y_h, \nabla(\phi  u_h))_{L^2(\Omega)^d} - (A_h\nabla y_h,  u_h\nabla \phi)_{L^2(\Omega)^d} \\
& = &  (A_h\nabla y_h, \nabla I_h(\phi  u_h))_{L^2(\Omega)^d} - (A_h\nabla y_h,  u_h\nabla \phi)_{L^2(\Omega)^d} \\
  &&+ (A_h\nabla y_h, \nabla ( \phi  u_h - I_h(\phi  u_h)))_{L^2(\Omega)^d}.
\end{eqnarray*}
Then using the weak formulation of $y_h$, we obtain
\begin{align}
E_h =& \langle \hat{f},I_h(\phi u_h)\rangle_{H^{-1}(\Omega),H^1_0(\Omega)} - (A_h\nabla y_h,  u_h\nabla \phi)_{L^2(\Omega)^d} \nonumber\\
     & + (A_h\nabla y_h, \nabla ( \phi  u_h - I_h(\phi  u_h)))_{L^2(\Omega)^d},\label{eqn:Eh-splitting}
\end{align}
where $I_h$ denotes the standard Lagrange interpolation operator. Then there hold
\begin{equation*}
I_h(\phi u_h) \rightharpoonup \phi \hat{u}
\mbox{ in }H^1_0(\Omega)\quad\mbox{and}\quad \|\phi  u_h - I_h(\phi  u_h)\|_{H^1(\Omega)} \rightarrow 0.
\end{equation*}
Indeed, by the standard interpolation estimate \cite{Brenner}, we have
\begin{equation*}
  \|\phi u_h-I_h(\phi u_h)\|_{H^1(T)}\leq Ch\|D^2(\phi u_h)\|_{L^2(T)}\leq ch\|\phi\|_{W^{2,\infty}(T)}\|u_h\|_{H^1(T)}, \quad \forall T\in \mathcal{T},
\end{equation*}
so that $\phi u_h-I_h(\phi u_h)\to0 $ in $H_0^1(\Omega)$ as $h\to0^+$ since $\|u_h\|_{H^1(\Omega)}\leq C$.
This and the weak convergence $\phi u_h\rightharpoonup \phi \hat u$ imply $I_h(\phi u_h) \rightharpoonup \phi \hat{u}$.
Upon passing to the limit $h\to0^+$ in the identity \eqref{eqn:Eh-splitting}, we obtain
\begin{eqnarray*}
\lim_{h\rightarrow 0^+}E_h  & = &  \langle \hat{f},\phi \hat{u}\rangle_{H^{-1}(\Omega),H^1_0(\Omega)}- (A\nabla y,  \hat{u} \nabla \phi)_{L^2(\Omega)^d} \\
 & = & \langle -\mathrm{div}(A\nabla y),\phi \hat{u}\rangle_{H^{-1}(\Omega),H^1_0(\Omega)}- (A\nabla y,  \hat{u} \nabla \phi)_{L^2(\Omega)^d}\\
 & = & (A\nabla y, \phi\nabla \hat{u})_{L^2(\Omega)^d}.
\end{eqnarray*}
Likewise, repeating the preceding argument yields
\begin{eqnarray*}
E_h & = &  (A_h\nabla u_h, \phi \nabla y_h)_{L^2(\Omega)^d} \\
 &=& (A_h\nabla u_h, \nabla(\phi  y_h))_{L^2(\Omega)^d} - (A_h\nabla u_h,  y_h\nabla \phi)_{L^2(\Omega)^d} \\
& = &  (A_h\nabla u_h, \nabla I_h(\phi  y_h))_{L^2(\Omega)^d} - (A_h\nabla u_h,  y_h\nabla \phi)_{L^2(\Omega)^d} \\
  &&+ (A_h\nabla u_h, \nabla ( \phi  y_h - I_h(\phi  y_h)))_{L^2(\Omega)^d} \\
& = &  \langle f,I_h(\phi y_h)\rangle_{H^{-1}(\Omega),H^1_0(\Omega)}  - (A_h\nabla u_h,  y_h\nabla \phi)_{L^2(\Omega)^d}\\
 &&+ (A_h\nabla u_h, \nabla ( \phi  y_h - I_h(\phi  y_h)))_{L^2(\Omega)^d}.
\end{eqnarray*}
Passing to the limit $h\to0^+$ in the last identity gives
\begin{eqnarray*}
\lim_{h\rightarrow 0^+}E_h  & = & \langle f,\phi y\rangle_{H^{-1}(\Omega),H^1_0(\Omega)}- (A\nabla u,  y \nabla \phi)_{L^2(\Omega)^d} \\
 & = & \langle -\mathrm{div}(A\nabla u),\phi y\rangle_{H^{-1}(\Omega),H^1_0(\Omega)}- (A\nabla u,  y \nabla \phi)_{L^2(\Omega)^d} \\
 & = & (A\nabla u, \phi\nabla y)_{L^2(\Omega)^d}.
\end{eqnarray*}
Since the two limits are identical, we conclude
\begin{equation*}
(A\nabla y, \phi\nabla \hat{u})_{L^2(\Omega)^d} = (A\nabla u, \phi\nabla y)_{L^2(\Omega)^d}.
\end{equation*}
Since the choice of the functions $\hat{f}$ and $\phi$ is arbitrary, we have $\nabla u = \nabla \hat{u}$.
Since $u_h\rightharpoonup\hat{u}$ weakly in $H^1(\Omega)$, by Sobolev embedding theorem and
trace theorem, $u_h$ converges to $u$ strongly in $L^2(\Gamma)$, i.e., $\|\hat{u}|_{\Gamma}
- g_h\|_{L^2(\Gamma)} \rightarrow 0 $. This directly gives $\hat{u}|_{\Gamma} = g = u|_{\Gamma}$.
Together with the identity $\nabla u = \nabla \hat{u}$, we obtain $ u = \hat{u}$, which completes
the proof of the theorem.
\end{proof}

The next result is a discrete analogue of Lemma \ref{compare-cont}  \cite[Corollary 3.1]{Deckelnick2011}.
\begin{lemma}\label{dineq}
Let the sequence $\{A_h\}_{h>0}$ be
$A_h\xrightarrow{Hd} A$ and $A_h \overset{\ast}{\rightharpoonup}  A^0$ in $L^\infty(\Omega)^{d,d}$.
Then
\begin{equation*}
A\leq A^0\quad \mbox{a.e. in }\Omega\quad\mbox{and}\quad \normmm{A}^2_2\leq \normmm{A^0}^2_2\leq \liminf_{h\rightarrow 0}\normmm{A_h}^2_2.
\end{equation*}
\end{lemma}

Now we can state a convergence result for a sequence $\{A_h\}_{h>0}$
of solutions of the discrete problem \eqref{discret}. It is a finite element version of the result
obtained in \cite{Eymard}.

\begin{theorem}
There exists at least one solution $A_h\in \mathcal{A}_h$ to problem \eqref{discret}.
Further, there exists a subsequence $\{A_{h^\prime}\}_{h^\prime>0}$ and $A\in \mathcal{A}$
such that $A_{h^\prime}\rightarrow A$ in $L^2(\Omega)^{d,d}$ and $A$ is a solution of
problem \eqref{tik}.
\end{theorem}
\begin{proof}
Problem \eqref{discret} is a finite-dimensional optimization problem, and the existence of a minimizer $A_h$
follows directly from the coercivity and continuity of the functional $J_{\gamma,h}$. It remains to prove the convergence.
By Theorem \ref{dcompact} and Lemma \ref{dineq}, there exists a subsequence
$\{A_{h^\prime}\in \mathcal{A}_{h'}\}_{h^\prime>0}$ and $A\in \mathcal{A}$ such that
 $A_{h^\prime}\xrightarrow{Hd} A$
and $A_{h^\prime} \overset{\ast}{\rightharpoonup} A^0$ in $L^\infty(\Omega)^{d,d}$, and
\begin{equation}\label{inequ2}
  A\leq A^0\mbox{ a.e. in } \Omega, \quad \mbox{and}\quad \normmm{A}^2_2\leq \liminf_{k\rightarrow \infty}\normmm{A_{h^\prime}}^2_2.
\end{equation}
Let $u^{\ell}_{h^\prime}(A_{h^\prime})=F^{\ell}_{h^\prime}(A_{h^\prime})$,
$u^{\ell}(A)=F^{\ell}(A)$, $i=1,\ldots,M$. Then we have
\begin{equation*}
u^{\ell}_{h^\prime}(A_{h^\prime})\rightharpoonup u^{\ell}(A) \mbox{ in }H^1(\Omega)\quad \mbox{and}\quad A_{h^\prime}\nabla u^{\ell}_{h^\prime}\rightharpoonup A\nabla u^{\ell}\mbox{ in }L^2(\Omega)^d.
\end{equation*}
The argument of Theorem \ref{uniqsol} implies
\begin{equation*}
J_{\gamma}(A)\leq \liminf_{h^\prime\rightarrow 0} J_{\gamma,h^\prime}(A_{h^\prime}).
\end{equation*}
 Next, Theorem \ref{uniqsol}
 implies that problem \eqref{tik} has a solution $\tilde{A}\in \mathcal{A}$.
 Then we have
 \begin{align*}
   J_\gamma(\tilde{A})&\leq J_\gamma(A)\leq \liminf_{h^\prime\rightarrow 0}J_{\gamma,h^\prime}(A_{h^\prime})\leq \limsup_{h^\prime\rightarrow 0}J_{\gamma,h^\prime}(A_{h^\prime})\\
   &\leq \limsup_{h^\prime\rightarrow 0}J_{\gamma,h^\prime}(\tilde{A})=J_\gamma(\tilde{A}).
 \end{align*}
 Consequently,
\begin{equation}\label{JAh}
  \lim_{h^\prime\rightarrow 0}J_{\gamma,h^\prime}(A_{h^\prime})=J_\gamma(A)=J_\gamma(\tilde{A}).
\end{equation}
In particular, $A$ is a minimizer of the functional $J_\gamma$.
Furthermore, direct computation gives
\begin{align*}
&\quad \frac12\sum_{\ell=1}^L\|A\nabla u^{\ell}(A)-A_{h^\prime}\nabla u^{\ell}_{h^\prime}(A_{h^\prime})\|^2_{L^2(\Omega)^d}
+\frac \gamma 2\normmm{A-A_{h^\prime}}^2_2\\
&=J_\gamma(A)+J_{\gamma,h^\prime}(A_{h^\prime})+\sum_{\ell=1}^L(A\nabla u^{\ell}(A)-\vec{h}^{\ell},\vec{h}^{\ell}-A_{h^\prime}\nabla u^{\ell}_{h^\prime}(A_{h^\prime}))_{L^2(\Omega)^d}-\gamma(A, A_{h^\prime})_{L^2(\Omega)^{d,d}}.
\end{align*}
Now \eqref{inequ2} and the convergence $A_{h^\prime} \overset{\ast}{\rightharpoonup} A^0$ in $L^\infty(\Omega)^{d,d}$ imply
\begin{align*}
  \lim_{h'\to0^+}  (A\nabla u^{\ell}(A)-\vec{h}^{\ell},\vec{h}^{\ell}-A_{h^\prime}\nabla u^{\ell}_{h^\prime}(A_{h^\prime}))_{L^2(\Omega)^d} &= \|A\nabla u^{\ell}(A)-\vec{h}^{\ell}\|^2_{L^2(\Omega)^d},\\
  \lim_{h'\to0^+}  (A, A_{h^\prime})_{L^2(\Omega)^{d,d}} &= (A,A^0)_{L^2(\Omega)^{d,d}} \geq \normmm{A}^2_2.
\end{align*}
These identities together with \eqref{JAh} imply
\begin{align*}
  &\lim_{h'\to0^+} \frac12\sum_{\ell=1}^L\|A\nabla u^{\ell}(A)-A_{h^\prime}\nabla u^{\ell}_{h^\prime}(A_{h^\prime})\|^2_{L^2(\Omega)^d}
+\frac \gamma 2\normmm{A-A_{h^\prime}}^2_2\\
 \leq& 2J_\gamma(A) -\sum_{\ell=1}^L\|A\nabla u^{\ell}(A)-\vec{h}^{\ell}\|^2_{L^2(\Omega)^d} - \gamma\normmm{A}^2_2  =0,
\end{align*}
which completes the proof of the theorem.
\end{proof}

\subsection{Projected Newton algorithm}
Based on the necessary optimality system in Theorem \ref{kkt}, there are several
different ways to develop algorithms for solving the regularized formulation \eqref{tik}. One
direct choice is gradient descent, which was explored for a related inverse conductivity
problem in \cite{Deckelnick2012}. Generally, gradient type methods are known to converge
steadily but often slowly, especially when the sought-for conductivity tensor is nonsmooth.
Thus, it is still of much interest to develop efficient algorithms. In this part we develop
a projected Newton algorithm. First we derive the necessary optimal condition of the finite
element problem \eqref{discret}. The argument in Section \ref{sec:reg} shows that the
directional derivative $J'_{\gamma,h}(A_h)[H]$ of $J_{\gamma,h}(A_h)$ is given by
\begin{equation}\label{DFre}
  J^\prime_{\gamma,h}(A_h)[H]=\int_\Omega(\sum_{\ell=1}^L\nabla u^{\ell}_h(A_h)\otimes (A_h\nabla u^{\ell}_h(A_h)-\vec{h}^{\ell})+\sum_{\ell=1}^L\nabla u^{\ell}_h(A_h)\otimes \nabla \bar{p}^{\ell}_h+\gamma A_h)\cdot Hdx,
\end{equation}
for any feasible direction $H\in X_h$ such that $A_h+H\in \mathcal{A}_h$,
where $u^{\ell}_h(A_h)=F^{\ell}_h(A_h)\in V_h$ and $\bar{p}^{\ell}_{h}\in V_{h,0}$, $\ell=1,\ldots, L$, solves
\begin{equation}\label{adj}
  (A_h\nabla\bar{p}^{\ell}_h,\nabla v_h)_{L^2(\Omega)^d}+(A_h(A_h\nabla u^{\ell}_h(A_h)-\vec{h}^{\ell}),
  \nabla v_h)_{L^2(\Omega)^d}=0,\quad \forall v_h\in V_h(\mathcal{T}).
\end{equation}
The optimality condition can be interpreted pointwise as
\begin{equation}\label{dsol}
  A_h(x)=P_K\Big(-\frac1\gamma\big(\sum_{\ell=1}^L\nabla u^{\ell}_h(A_h)\otimes (A_h\nabla u^{\ell}_h-\vec{h})+\sum_{\ell=1}^L\nabla u^{\ell}_h(A_h)\otimes\nabla\bar{p}^{\ell}_h\big)\Big) \quad \mbox{a.e. in } \Omega.
\end{equation}
In sum, the first-order necessary optimality system in the variational form is given by
\begin{equation*}
\left\{\begin{aligned}
  &(A_h\nabla u^{\ell}_h,\nabla v_h)_{L^2(\Omega)^d} -(f^{\ell},v_h)_{L^2(\Omega)^d}=0,\quad  \forall v_h\in V_{h,0},\quad\ell=1,\ldots,L, \\
  &(A_h\nabla \bar{p}^{\ell}_h,\nabla v_h)_{L^2(\Omega)^d}+ (A_h(A_h\nabla u^{\ell}_h-\vec{h}^{\ell})
   ,\nabla v_h)_{L^2(\Omega)^d} =0, \quad \forall v_h\in V_{h,0}, \quad\ell=1,\ldots,L, \\
  &\Big(\sum_{\ell=1}^L\nabla u^{\ell}_h\otimes(A_h\nabla u^{\ell}_h-\vec{h}^{\ell})+ \nabla u^{\ell}_h\otimes \nabla \bar{p}^{\ell}_h+\gamma A_h,D_h-A_h\Big)_{L^2(\Omega)^{d,d}}\geq 0,\quad \forall  D_h\in \mathcal{A}_h.
\end{aligned}\right.
\end{equation*}

To apply the Newton method, one crucial step is to derive the Newton update $(\bar u^\ell,\bar p^\ell,\bar A)$
for the unknowns $(u^\ell, p^\ell, A)$. This can be achieved by solving
\begin{equation}\label{eqn:Newton}
        \left\{
            \begin{aligned}
   (\bar A \nabla u_h^{\ell,n},\nabla v_h)_{L^2(\Omega)^d}+ (A_h^n\nabla \bar u^{\ell},\nabla v_h)_{L^2(\Omega)^d} = (A_h^n\nabla u_h^{\ell,n},\nabla v_h)_{L^2(\Omega)^d},&\\
        \quad  \forall v_h\in V_{h,0}, \ell=1,\ldots,L, &\\
    (\bar A (\nabla p_h^{\ell,n}+A_h^n\nabla u_h^{\ell,n}-\vec{h}^{\ell})+A_h^n(\bar A\nabla u_h^{\ell,n}),\nabla v_h)_{L^2(\Omega)^d}
      +
   (A_h^n(A_h^n\nabla \bar u^\ell,\nabla v_h)_{L^2(\Omega)^d}&\\
   +(A_h^n\nabla \bar p^{\ell},\nabla v_h)_{L^2(\Omega)^d}
   =  (A_h^n\nabla {p}_h^{\ell,n} + A_h^n(A_h^n\nabla u_h^{\ell,n}-\vec{h}^{\ell}),\nabla v_h)_{L^2(\Omega)^d},&\\
   \quad  \forall v_h\in V_{h,0},\ell=1,\ldots,L, &\\
   \sum_{\ell=1}^L (\nabla u_h^{\ell,n}\otimes(\bar A \nabla u_h^{\ell,n})+\gamma \bar A+\nabla \bar u^{\ell}\otimes(A_h^n\nabla u_h^{\ell,n}-\vec{h}^{\ell}),D_h)_{L^2(\Omega)^{d,d}} &\\
+  \sum_{\ell=1}^L (\nabla u_h^{\ell,n}\otimes A_h^n\nabla \bar u^{\ell}+\nabla \bar u^{\ell,n}\otimes
  \nabla p_h^{\ell,n} + \nabla \bar p^{\ell}\otimes \nabla u_h^{\ell,n}, D_h)_{L^2(\Omega)^{d,d}}& \\
    = \sum_{\ell=1}^L \big(\nabla u_h^{\ell,n}\otimes(A_h^n\nabla u_h^{\ell,n}-\vec{h}^{\ell})+ \nabla u_h^{\ell,n}\otimes \nabla p_h^{\ell,n}+\gamma
    A_h^n, D_h\big)_{L^2(\Omega)^{d,d}}, &\\
     \quad \forall D_h\in \mathcal{A}_h.&
\end{aligned}\right.
\end{equation}
The solution of the coupled system \eqref{eqn:Newton} is denoted by $\bar A_h^n, \bar u_h^{\ell,n}, \bar p_h^{\ell,n}$,
$\ell=1,\ldots,L$. When formulating the Newton method for the KKT system, we do not treat directly
the variational inequality (i.e., the box-constraint on the extremal eigenvalues of $A_h$, or equivalently the projection operator in
\eqref{dsol}). That is, we use the Newton iteration to update the solution to \eqref{ph}, \eqref{adj}
and \eqref{dsol} without the pointwise projection $P_K$. After computing the updates, we project the updated solution
onto the convex set $\mathcal{A}_h$, which can be performed pointwise. This leads to a projected Newton algorithm, whose detail is listed in
Algorithm 1. Note that each iteration of the method requires solving one coupled linear system in
the state variable $\bar u^{\ell}_h$, adjoint variable $\bar{p}^{\ell}_h$ and conductivity tensor $\bar A_h$.
The stopping criterion in line 6 is taken so that the relative error
\begin{equation*}
 \max\bigg\{\frac{\normmm{A_h^{n+1}-A_h^n}_2}{\normmm{A_h^{n+1}}_2},
 \frac{\|u_h^{\ell,n+1}-u_h^{\ell,n}\|_{L^2(\Omega)}}{\|u_h^{\ell,n+1}\|_{L^2(\Omega)}},
 \frac{\|{p}_h^{\ell,n+1}-{p}_h^{\ell,n}\|_{L^2(\Omega)}}{\|\bar{p}_h^{\ell,n+1}\|_{L^2(\Omega)}}, \ell=1,..., L\bigg\}
\end{equation*}
falls below a given tolerance $\tau$ or the max iteration number exceeds a pre-determined number.

\begin{algorithm}
\caption{Projected Newton algorithm}
\label{alg:A}
{\small{
\begin{algorithmic}[1]
\STATE Given initial guess $A_h^0, u_h^{\ell,0},{p}_h^{\ell,0}$, $\ell=1,\ldots,L$.
\FOR { $n=0,1,2,...$}
\STATE Obtain the increments $\bar A^n, \bar u^{\ell,n}, \bar p^{\ell,n}$ of $A_h^n, u_h^{\ell,n}, {p}_h^{\ell,n}$ by solving the Newton system \eqref{eqn:Newton}.
\STATE Update $A_h^n, u_h^{\ell,n}, \bar{p}_h^{\ell,n}$ by
 \begin{equation*}
   A_h^{n+1} = A_h^n-\bar A_h^n,\quad u_h^{\ell,n+1}=u_h^{\ell,n}-\bar u_h^{\ell,n}, \quad p_h^{\ell,n+1}={p}_h^{\ell,n}-\bar p_h^{\ell,n},\quad \ell=1,\ldots,L.
 \end{equation*}
\STATE Project $A_h^{n+1}$ onto $\mathcal{A}_h$, and reset $A_h^{n+1}$.
\STATE Check the stopping criterion.
\ENDFOR
\end{algorithmic}}}
\end{algorithm}

It is well known that Newton type algorithms converge very fast, if a good initial guess
is provided, which however is generally nontrivial. Meanwhile, the choice of the regularization parameter $\gamma$ is very important in order
to obtain satisfactory reconstructions. To choose a suitable regularization parameter
and to provide a good initial guess for the Newton algorithm simultaneously, we adopt an easy-to-implement yet very powerful path-following strategy,
which has been successfully applied in many applications \cite{JiaoJinLu:2015}.
Specifically, fix a decreasing factor $\rho\in (0,1)$, and we apply Algorithm 1 with $\gamma_n=\gamma_0\rho^n$.
with the initial guess given by the solution of the $\gamma_{n-1}-$ problem, i.e.,
\begin{equation*}
     \min_{A\in \mathcal{A}}\bigg\{ J_{\gamma_{n-1}}(A)=\frac12\sum_{\ell=1}^L\|A\nabla u^\ell(A)(f^\ell,g^\ell)-\vec{h}^\ell\|^2_{L^2(\Omega)^d}+\frac{\gamma_{n-1}}{2}\normmm{A}^2_2\bigg\}.
\end{equation*}
This step warm starts the Newton update, and is to fully exploit the fast local convergence of
Newton type methods. The final regularization parameter is determined by the classical discrepancy
principle \cite{EnglHankeNeubauer:1996,ItoJin:2015}, i.e., determining the smallest $k^*\in\mathbb{N}$ such that
\begin{equation*}
  \sum_{\ell=1}^L\|A_{\gamma_{k^*},h}\nabla u_h^\ell(A_{\gamma_{k^*},h})(f^\ell,g^\ell)-\vec{h}^\ell\|^2_{L^2(\Omega)^d} \leq L\delta^2,
\end{equation*}
where $A_{\gamma,h}$ denotes the minimizer for the discrete functional $J_{\gamma,h}$. In summary, there are
two loops in the algorithm: one is the inner iteration as shown in Algorithm \ref{alg:A}, and the
other is the outer iteration, which performs the path-following strategy over the penalty parameter $\gamma$. Due to
the fast local convergence of inner Newton iterations and proper initial guess from the path-following
strategy, one often needs only one or two inner iterations to ensure convergence, and hence, Algorithm \ref{alg:A} is
expected to be highly efficient, which is also confirmed by the extensive numerical experiments in Section \ref{sec:numer}.

\section{Numerical experiments and discussions}\label{sec:numer}
Now we present several two-dimensional numerical experiments to demonstrate the accuracy and efficiency
of the algorithm. All the experiments are carried out using FreeFem++ \cite{freefem} on a personal laptop.
The domain $\Omega$ is taken to be the square $(-1,1)^2$. We use an $(N+1)\times(N+1)$ uniform
square grid with $N=20$ and mesh size $h=2/N$. The current densities $h^\ell = A\nabla u^\ell$ for reconstructing
the conductivity tensor $A$ are simulated by solving problem \eqref{pr1} using the Galerkin finite
element method with a finer mesh. The noisy data $\vec{h}^{\ell,\delta}$ are generated by perturbing
the exact data $\vec{h}^\ell$ pointwise as ${\vec{h}}^{\ell,\delta}(x)= \vec{h}^\ell(x)(1+\delta\xi(x))$,
where $\xi$ is uniformly distributed on $[-1,1]$, and $\delta>0$ denotes the relative noise level.
In the numerical experiments, $V_h(\mathcal{T})$ is chosen as the $P_3$ conforming finite element space
and $X_h(\mathcal{T})$ is the $P_1$ discontinuous finite element space. Due to the use of cubic finite
elements, the choice $N=20$ is sufficient for ensuring reasonably accurate reconstructions for noisy data.
Numerically this choice is observed to significantly outperform low-order finite elements. To investigate
distinct features of the approach under different problem settings, we consider examples with either
diagonal or non-diagonal anisotropic conductivities, and unless otherwise stated,
the data correspond to five (i.e. $L=5$) different Dirichlet data $(g_1,g_2,g_3,g_4,g_5)$, given by
$(x_1+x_2,x_2+0.5x_1^2,x_1-0.1x_2^2,0.1(\cos(10x_2)-\cos(10x_1)),x_1x_2)$ and a vanishing source
$f\equiv0$, This choice is motivated by the fact that in the two-dimensional case, four judiciously chosen
excitations are sufficient to ensure the unique recovery \cite{bal2014,G.Bal2014}. The reconstruction procedure
is initialized with a constant matrix $A^0=\left[
  \begin{array}{cc}
       2 & -1 \\
      -1 & 2 \\
  \end{array}
\right],$
and in order to warm start the algorithm, we take $\gamma_0=10$ and reduce its value by a factor $\rho = 0.7$ after each inner loop.
Below we use $A^\dag$ (with entries $A_{ij}$, $i,j=1,2$) to denote the exact conductivity tensor.
The test cases and their distinct features are listed in Table \ref{tab:exam}.

\begin{table}[hbt!]
  \centering\caption{The test cases for the conductivity tensor. The functions
  $\xi:=2+\sin(\pi x_1)\sin(\pi x_2)$ and $\zeta:=0.5\sin(2\pi x_1)$, and $l(x)=
  1.5\chi_{S_1}+0.5\chi_{S_2}+2\chi_{S_3}+\chi_{S_4}$ with $S_1=[-0.4, 0.6] \times
  [-0.4, 0.6]$, $S_2=[-1, -0.7] \times [-0.4, 0.1]$, $S_3=[0.7, 1] \times [-0.8, -0.3]$
  and $S_4=\overline{\Omega}\setminus\cup_{i=1}^3 S_i$.\label{tab:exam}}
  \begin{tabular}{ccccc}
  \toprule
   Example  & $A_{11}$ & $A_{12}$ & $A_{22}$ & feature\\
  \midrule
    1       & 1        & 0 & 1 & constant\\
    2       & $1+x_1^2+x_2^2$ & 0 & $1+x_1^2+x_2^2$ & smooth\\
    3       &   $\xi$   & 0 & $1+\zeta^2$ & smooth, oscillatory\\
    4       &   $\xi$  & $\zeta$ & $1+x_1^2+x_2^2$ & nondiagonal\\
    5       & $1+|x_1|+|x_2|$ & $0$ & $2+|\sin(\pi x_1)\sin(\pi x_2)|$ & nonsmooth\\
    6       & $1+|x_1|+|x_2|$ & $|x_1x_2|$ & $2+|\sin(\pi x_1)\sin(\pi x_2)|$& nonsmooth, nondiagonal\\
    7       & $l(x)\xi$ & 0 & $l(x)(1+\zeta^2)$ & discontinuous\\
  \bottomrule
  \end{tabular}
\end{table}

\begin{table}[!htbp]
\tabcolsep 0.9mm
\begin{center}
\caption{The numerical results (accuracy error $e$) for Examples 1--7.}\label{table2}
\begin{tabular}{c|cccc}
\toprule
Example $\backslash$$\delta$&$0\%$&$1\%$&$5\%$&$10\%$\\
\midrule
1&3.490e-4&1.538e-3&4.180e-3&8.456e-3\\
2&1.492e-3&3.571e-3&8.452e-3&1.237e-2\\
3&2.731e-3&4.498e-3&1.256e-2&1.963e-2 \\
4&5.496e-3&1.110e-2&2.065e-2&3.249e-2\\
5&1.771e-3&3.566e-3&9.463e-3& 1.438e-2\\
6&2.499e-3&3.975e-3&9.999e-3&1.509e-2\\
7&4.970e-2&4.965e-2&5.029e-2&5.283e-2 \\
\bottomrule
\end{tabular}
\end{center}
\end{table}

\begin{figure}[!htbp]
\centering
\setlength{\tabcolsep}{0pt}
\begin{tabular}{cccc}
 \includegraphics[width=\wsize,height=\hhsize,trim={1.5cm 1cm 1cm 0.5cm},clip]{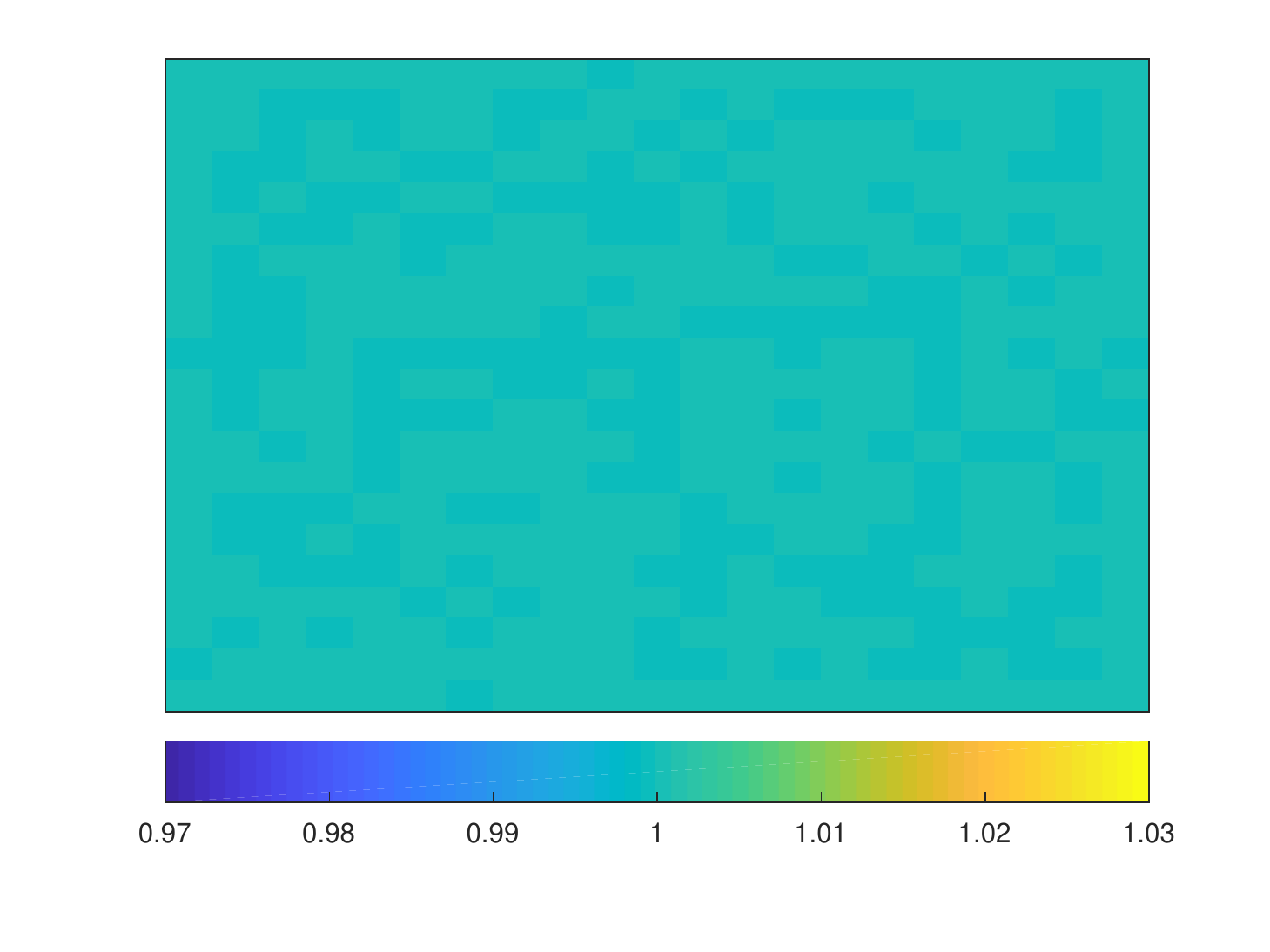}
&\includegraphics[width=\wsize,height=\hhsize,trim={1.5cm 1cm 1cm 0.5cm},clip]{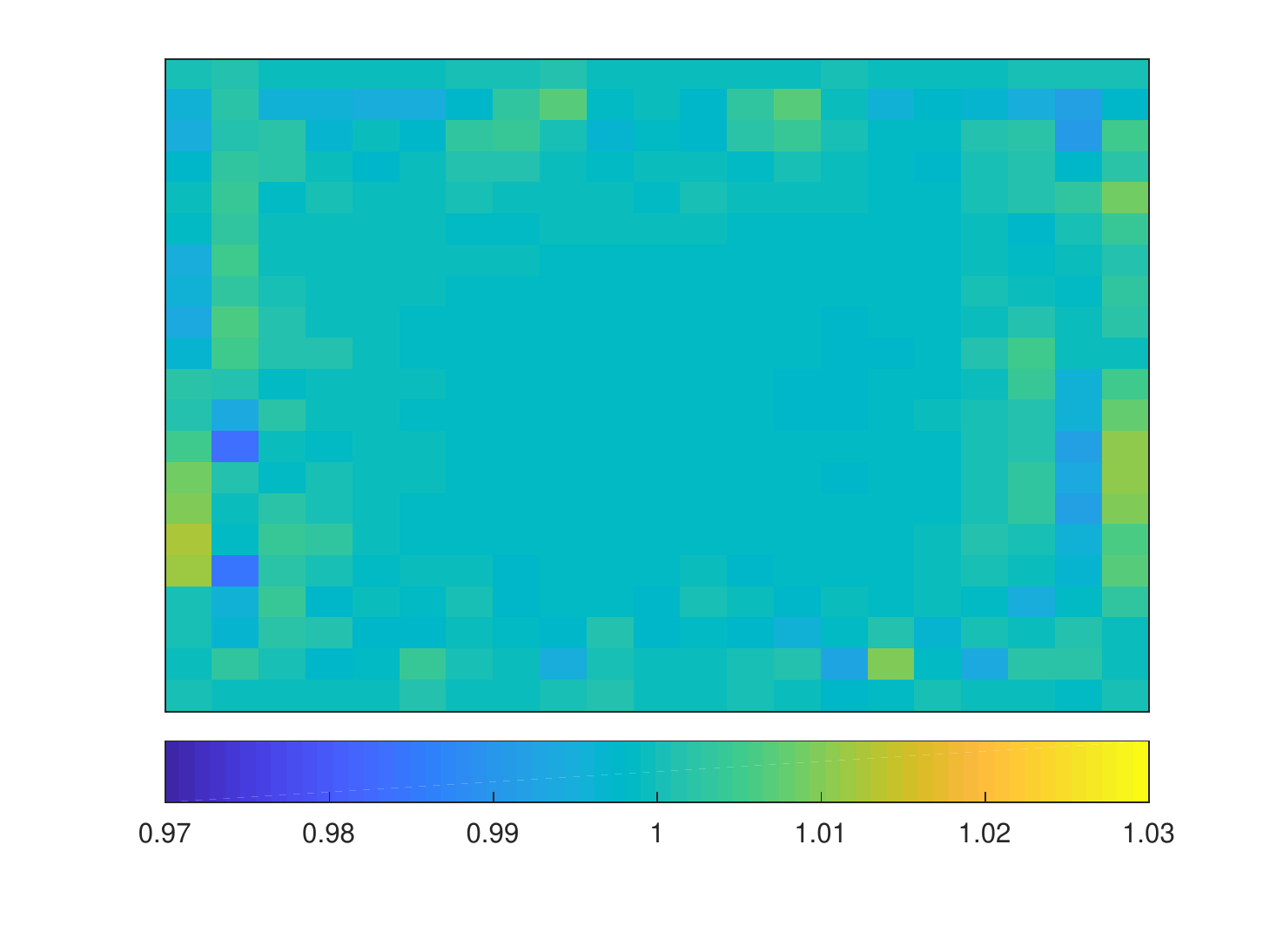}
&\includegraphics[width=\wsize,height=\hhsize,trim={1.5cm 1cm 1cm 0.5cm},clip]{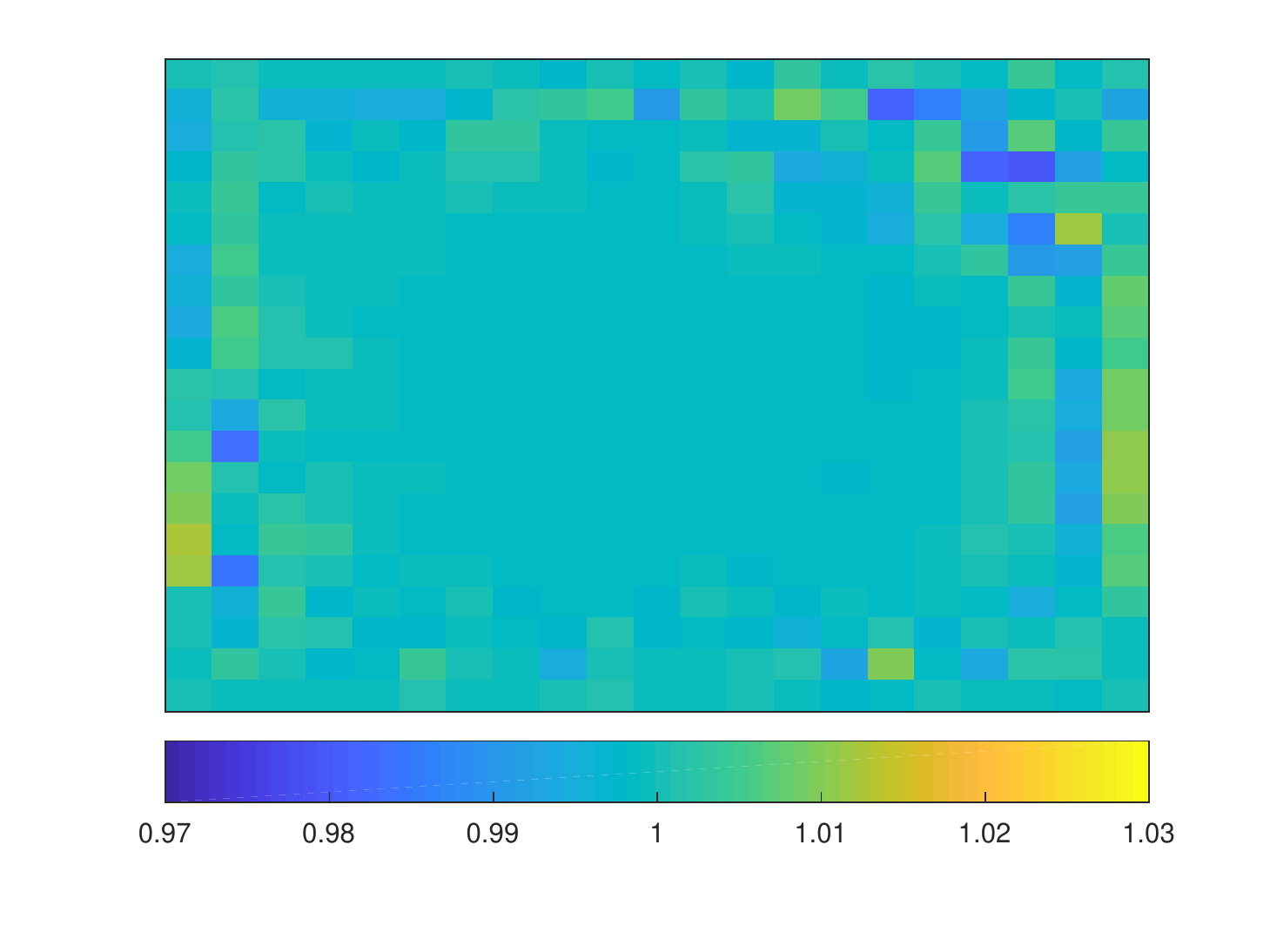}
&\includegraphics[width=\wsize,height=\hhsize,trim={1cm .5cm .5cm .5cm },clip]{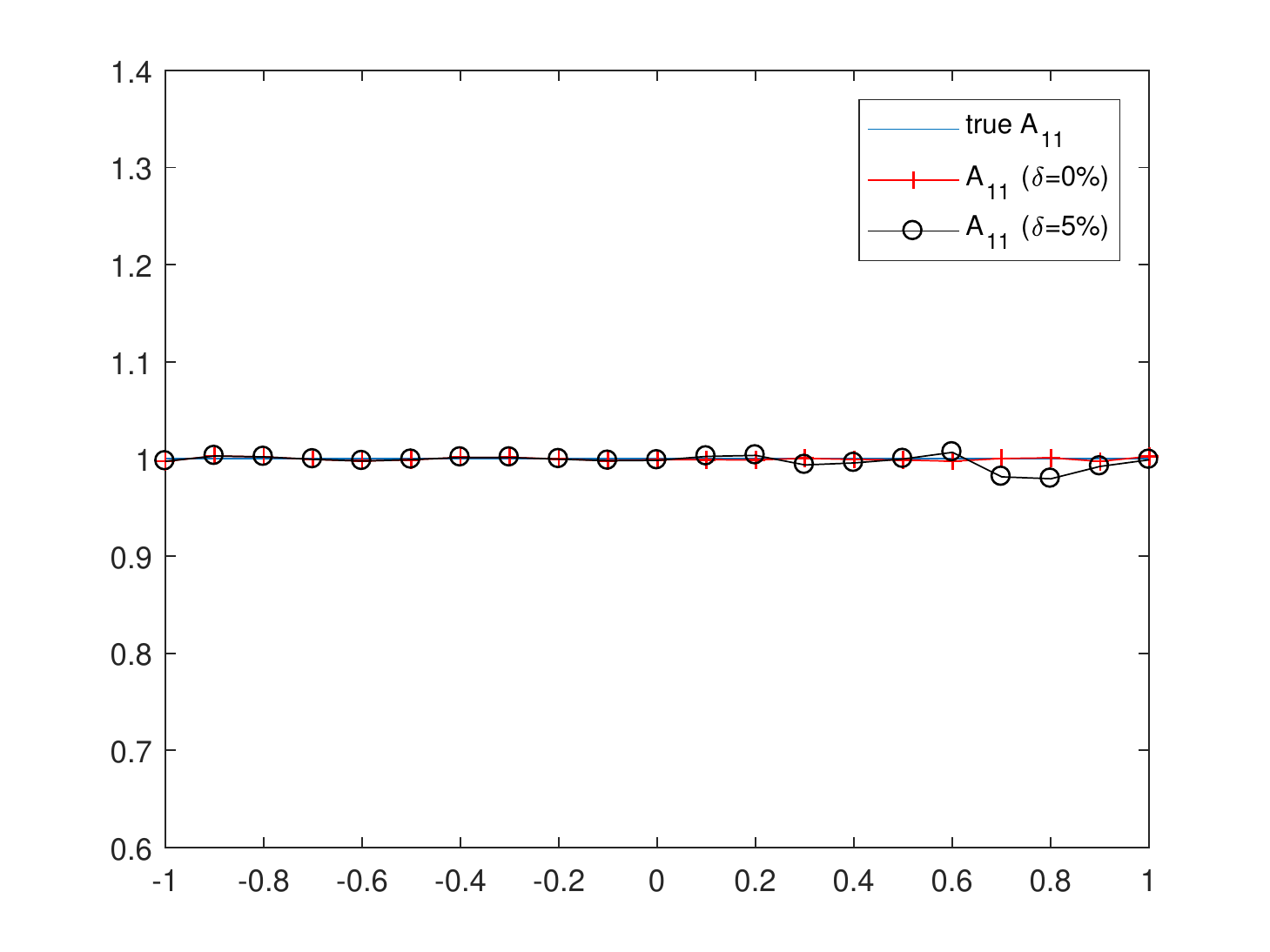}\\
 \includegraphics[width=\wsize,height=\hhsize,trim={1.5cm 1cm 1cm 0.5cm},clip]{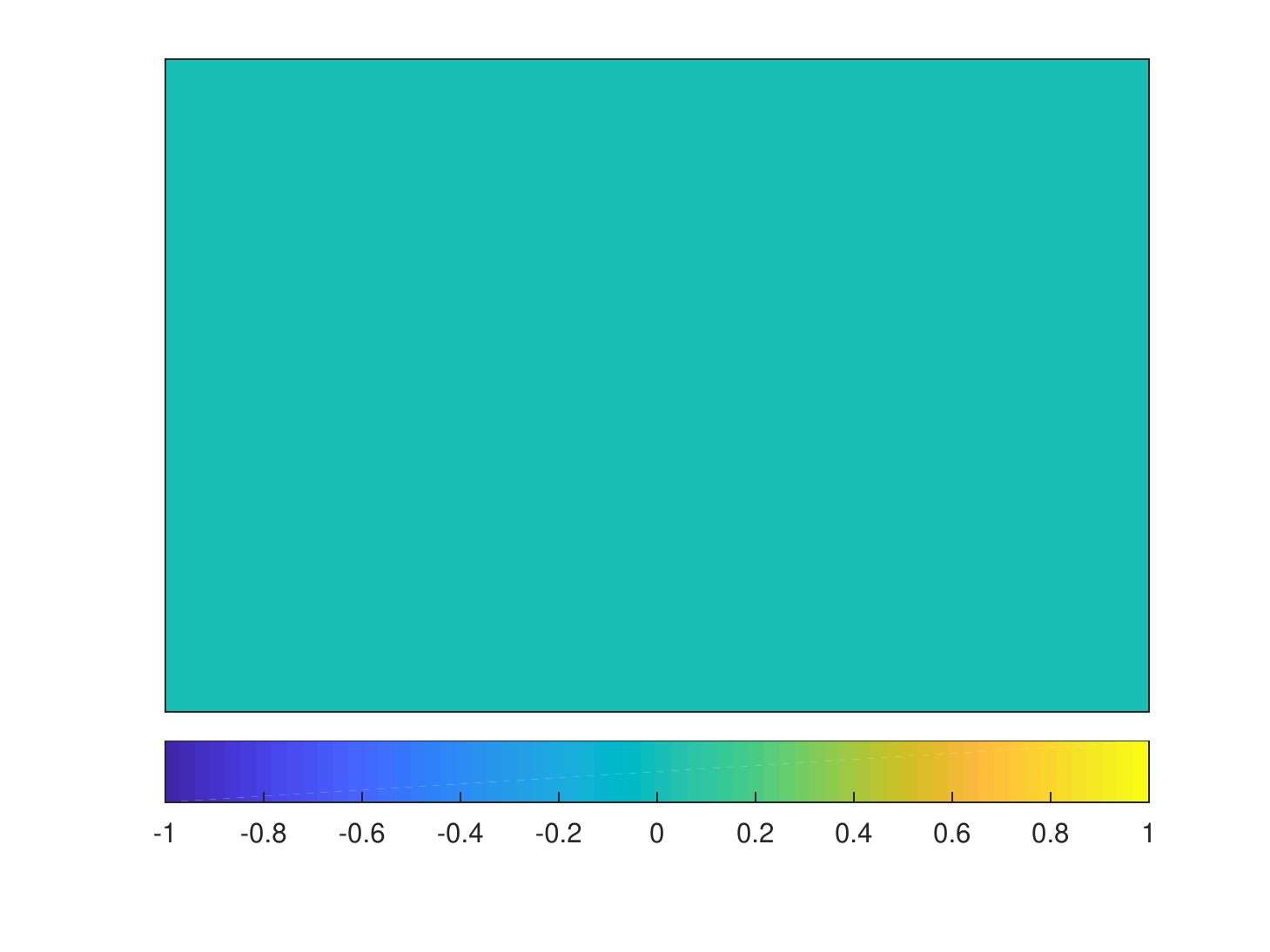}
&\includegraphics[width=\wsize,height=\hhsize,trim={1.5cm 1cm 1cm 0.5cm},clip]{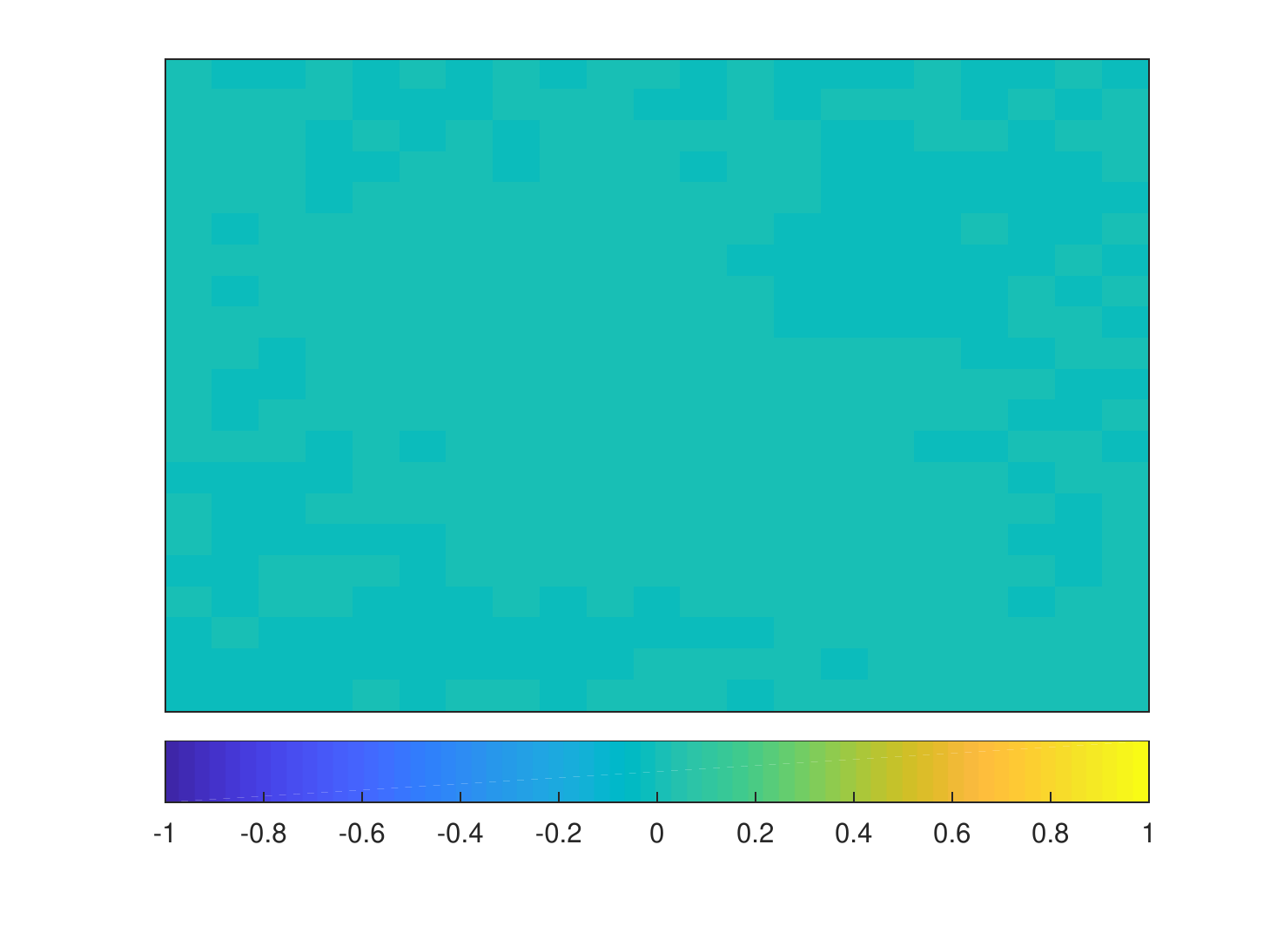}
&\includegraphics[width=\wsize,height=\hhsize,trim={1.5cm 1cm 1cm 0.5cm},clip]{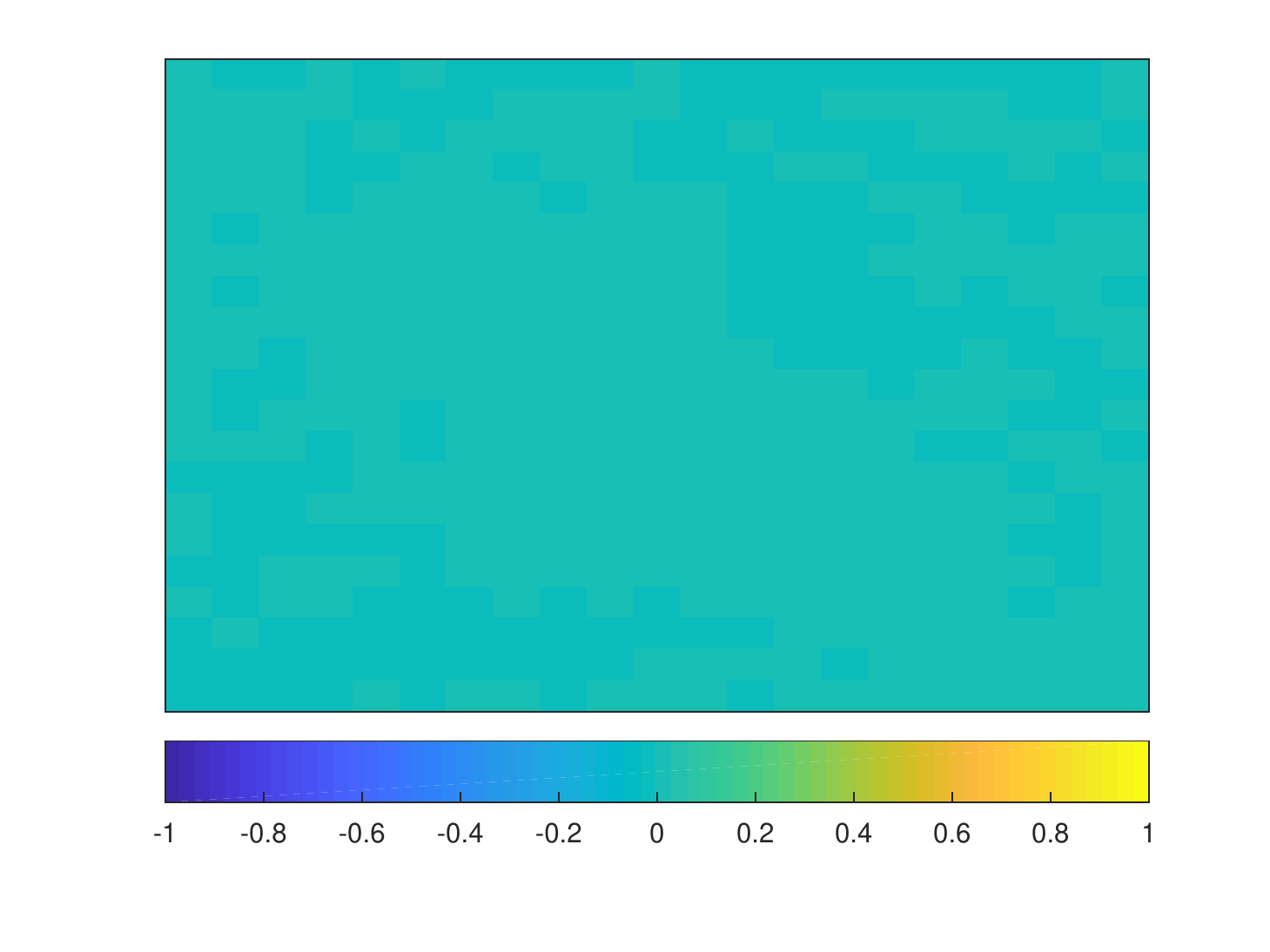}
&\includegraphics[width=\wsize,height=\hhsize,trim={1cm 0.5cm .5cm .5cm},clip]{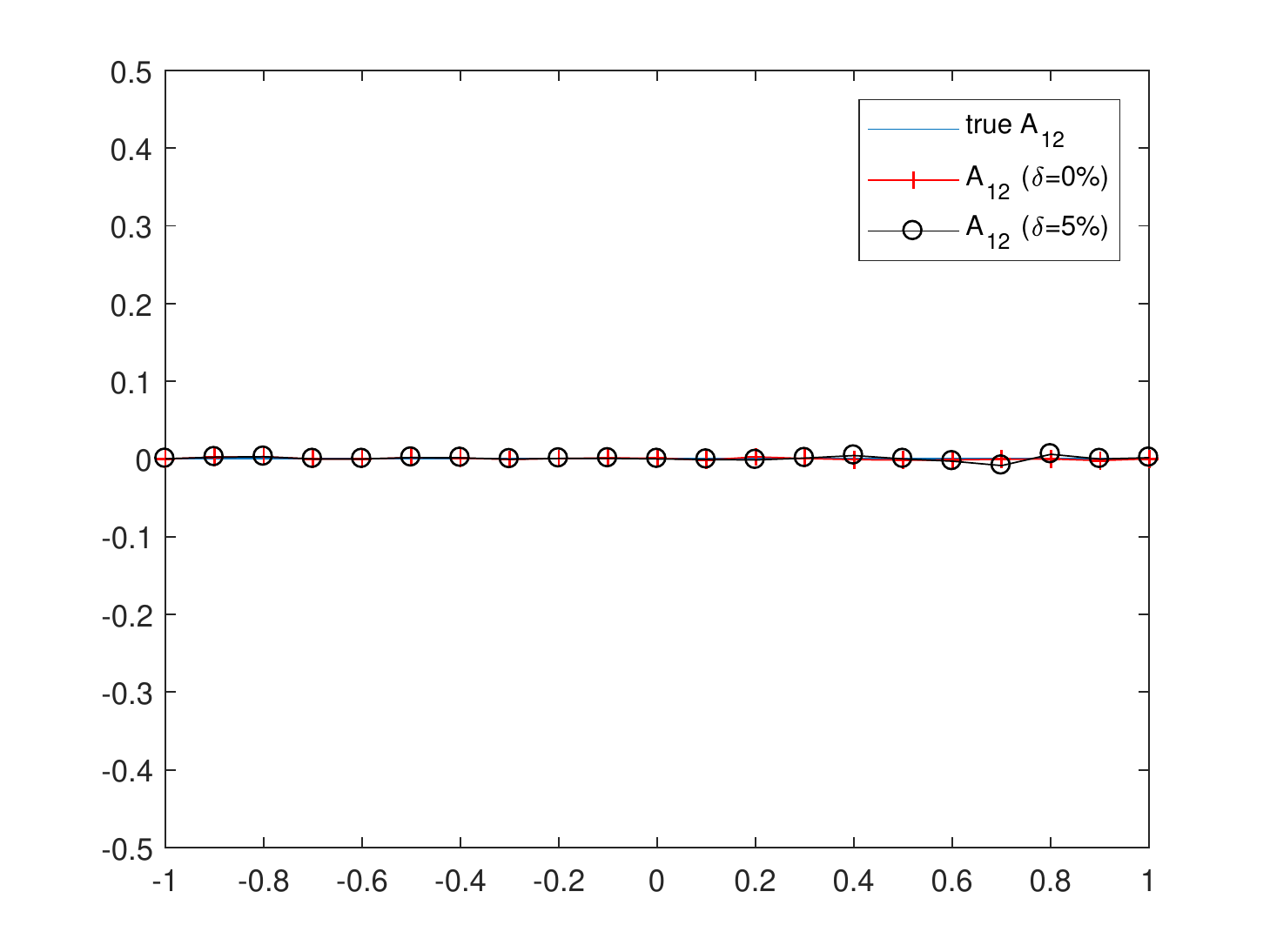}\\
 \includegraphics[width=\wsize,height=\hhsize,trim={1.5cm 1cm 1cm 0.5cm},clip]{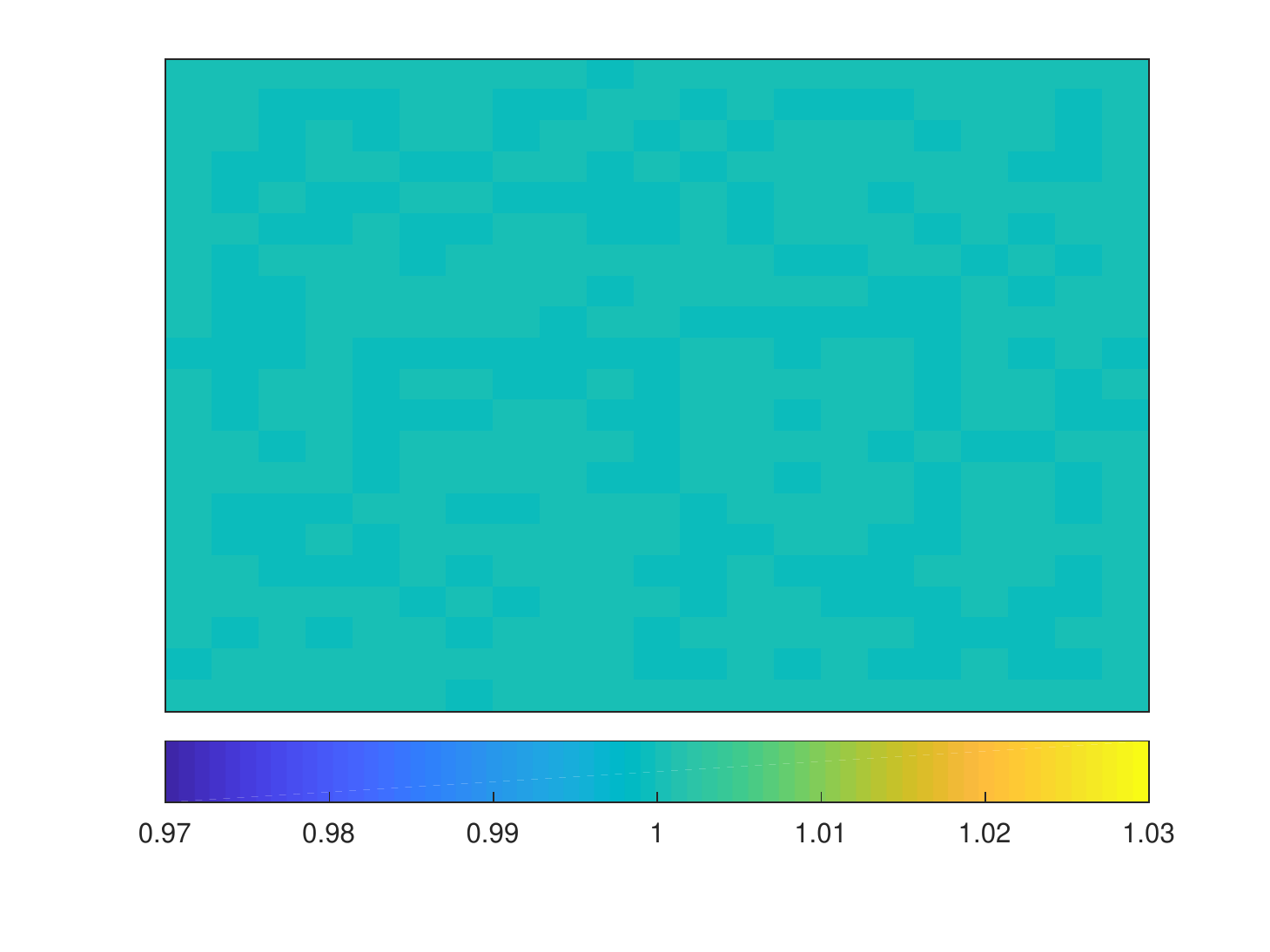}
&\includegraphics[width=\wsize,height=\hhsize,trim={1.5cm 1cm 1cm 0.5cm},clip]{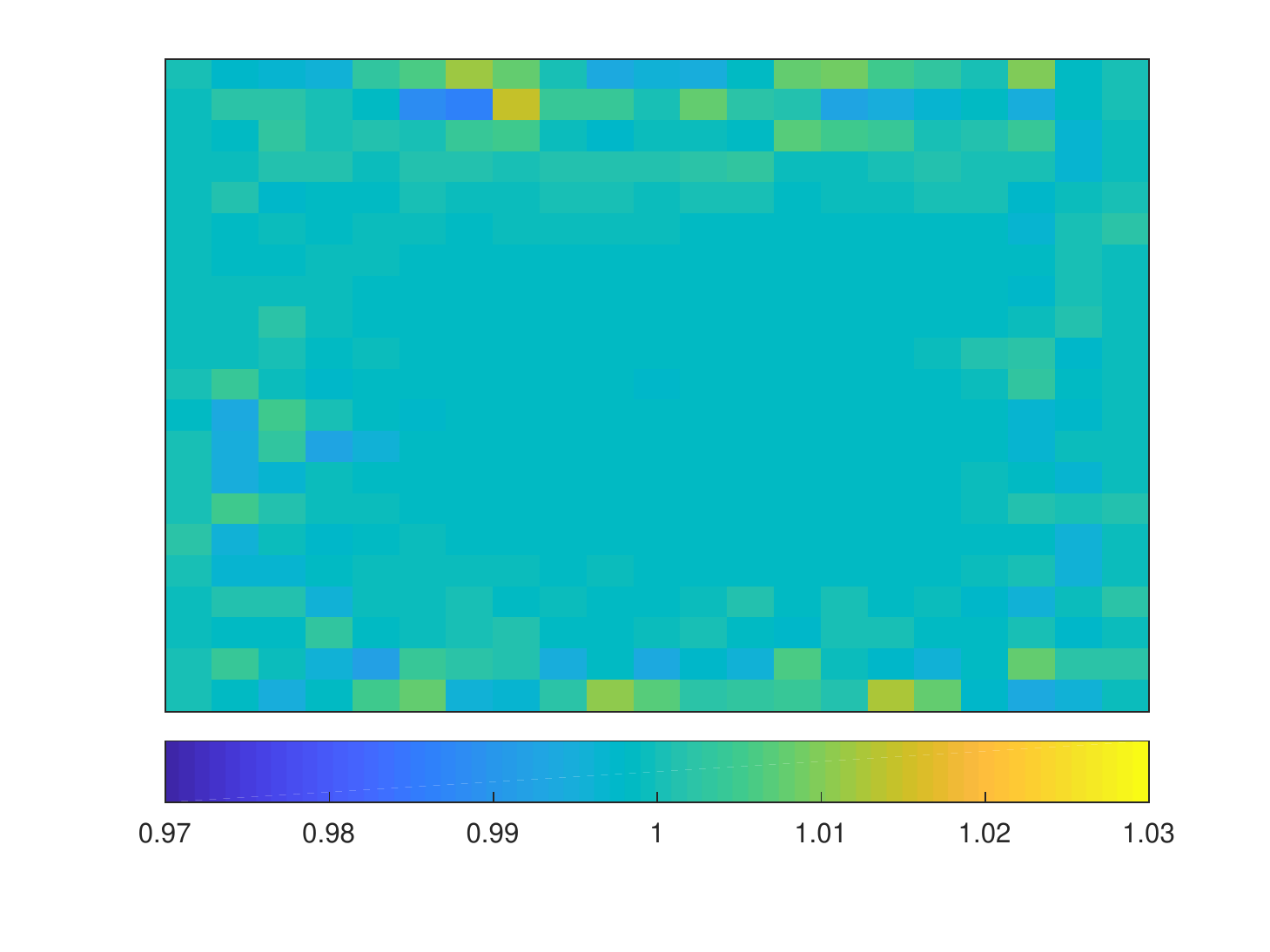}
&\includegraphics[width=\wsize,height=\hhsize,trim={1.5cm 1cm 1cm 0.5cm},clip]{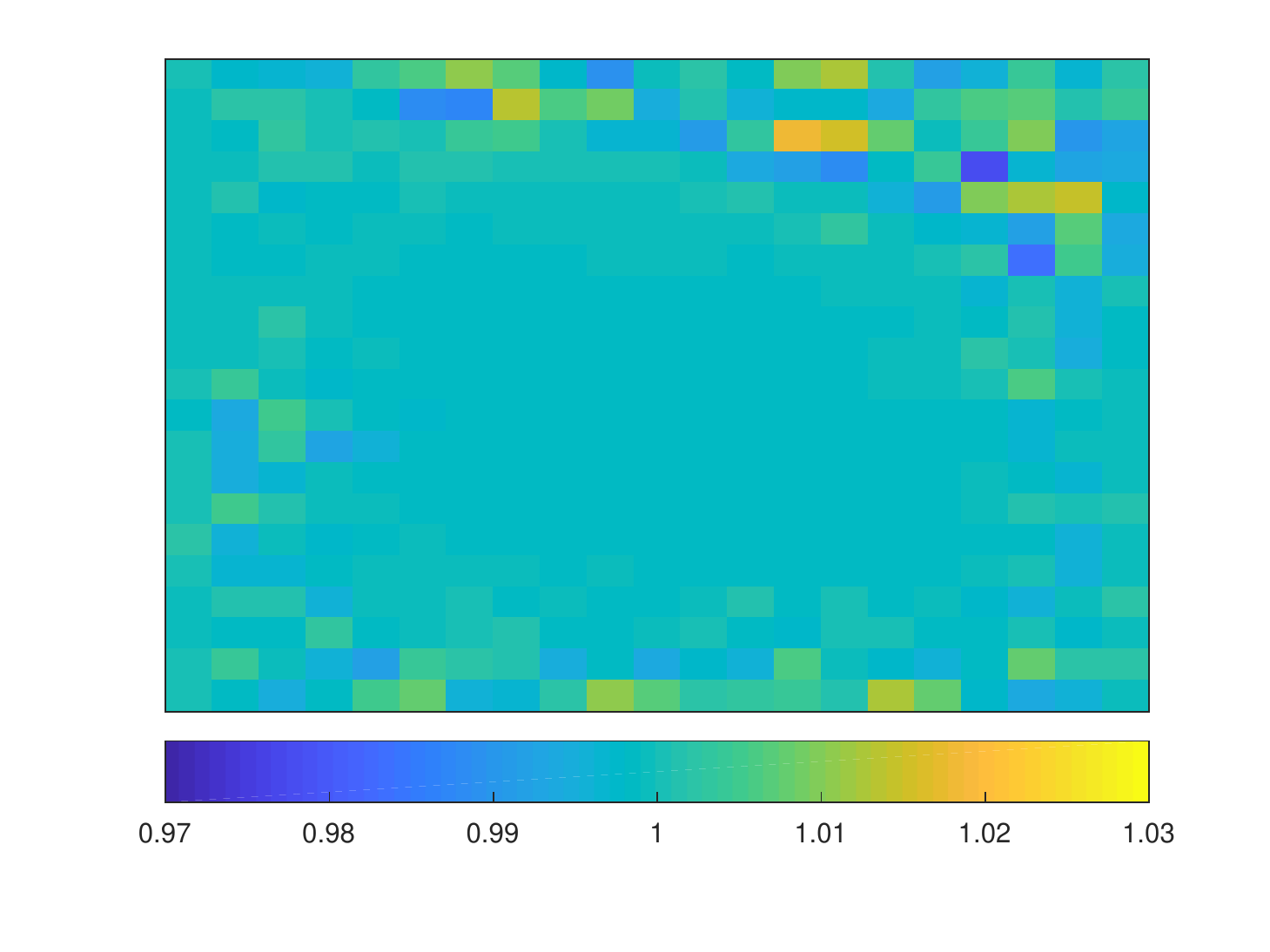}
&\includegraphics[width=\wsize,height=\hhsize,trim={1cm .5cm .5cm .5cm },clip]{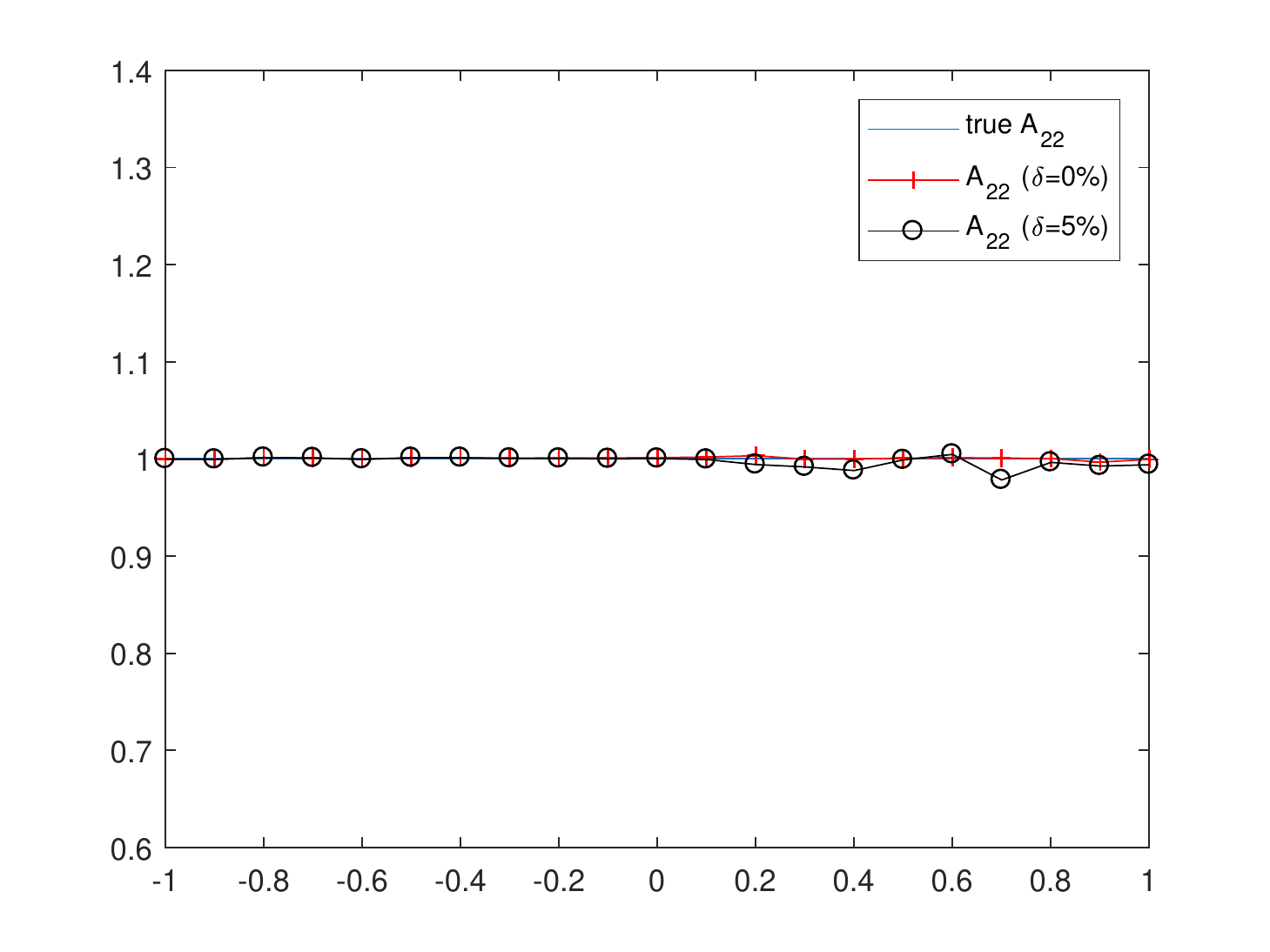}\\
\end{tabular}
\caption{The reconstructions for the conductivity tensor for Example 1: the top, middle and bottom rows refer to $A_{11}$, $A_{12}$ and $A_{22}$,
respectively. Column 1 is for exact conductivity, Column 2 for reconstruction with exact data, Column 3 for noisy data with $\delta=5\%$, and column 4 for the cross section along \{$x_2=-0.7$\}.}
\label{exam1}
\end{figure}

\begin{figure}[!htbp]
\centering
\setlength{\tabcolsep}{0pt}
\begin{tabular}{cccc}
 \includegraphics[width=\wsize,height=\hhsize,trim={1.5cm 1cm 1cm 0.5cm},clip]{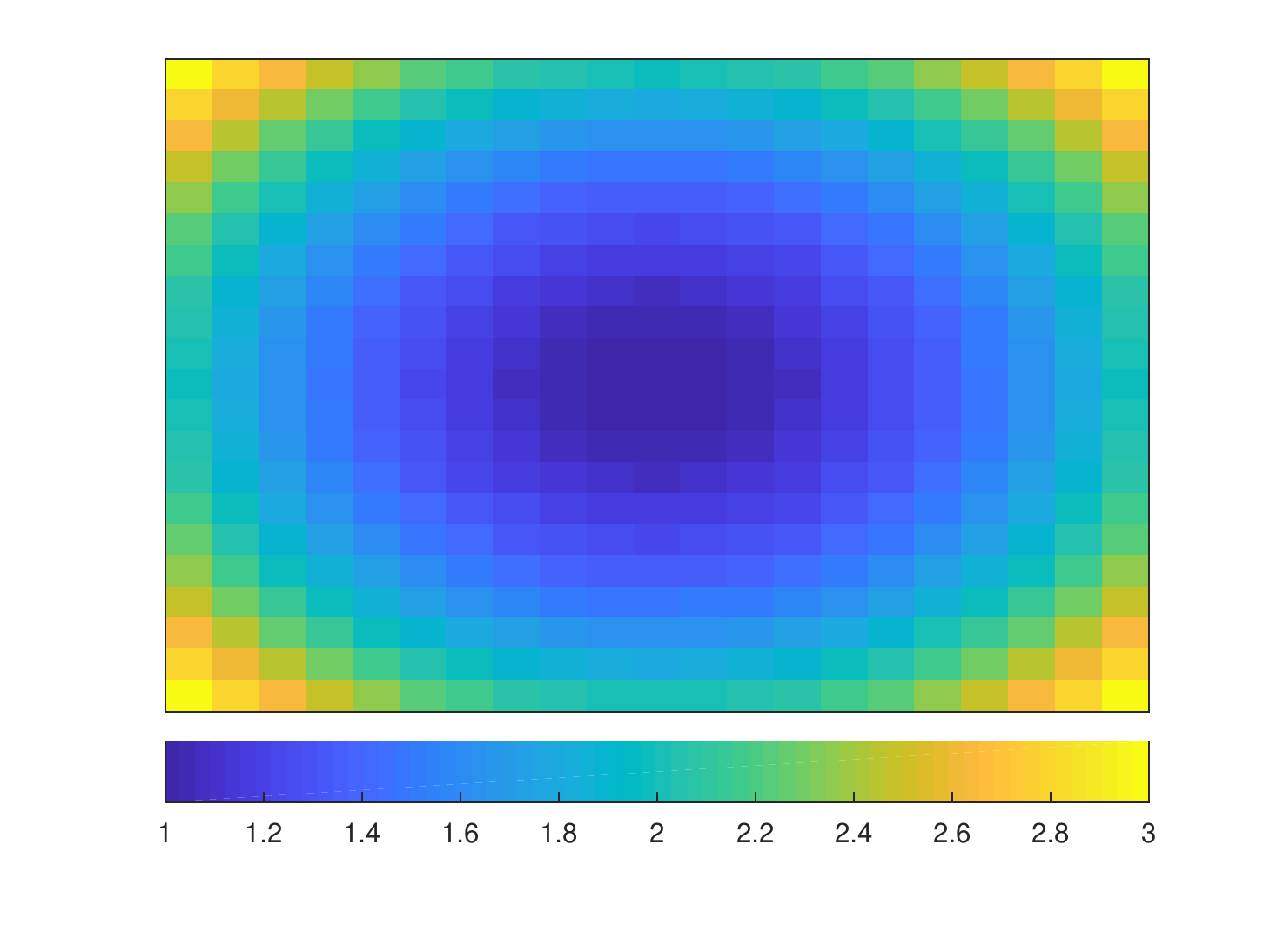}
&\includegraphics[width=\wsize,height=\hhsize,trim={1.5cm 1cm 1cm 0.5cm},clip]{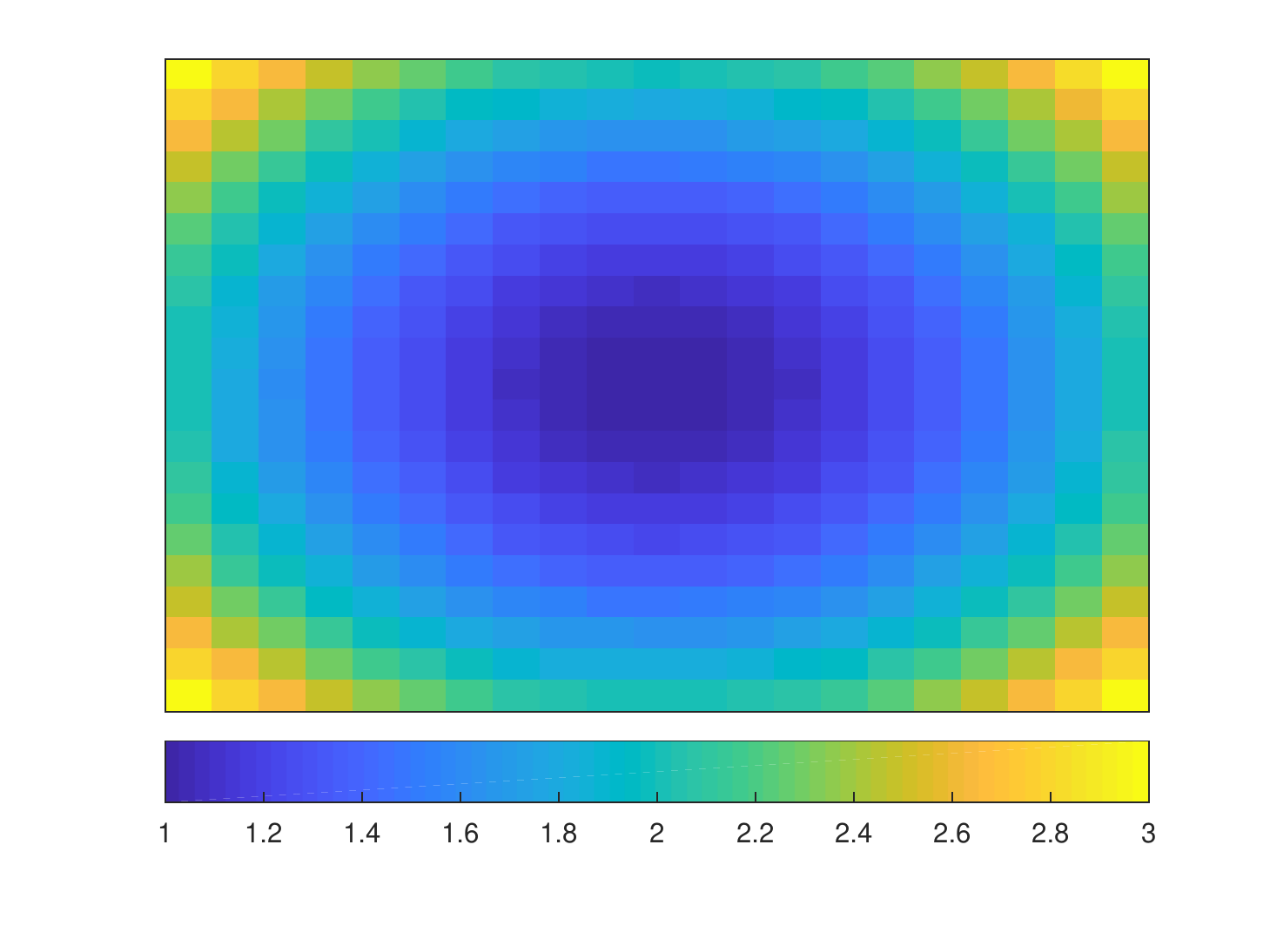}
&\includegraphics[width=\wsize,height=\hhsize,trim={1.5cm 1cm 1cm 0.5cm},clip]{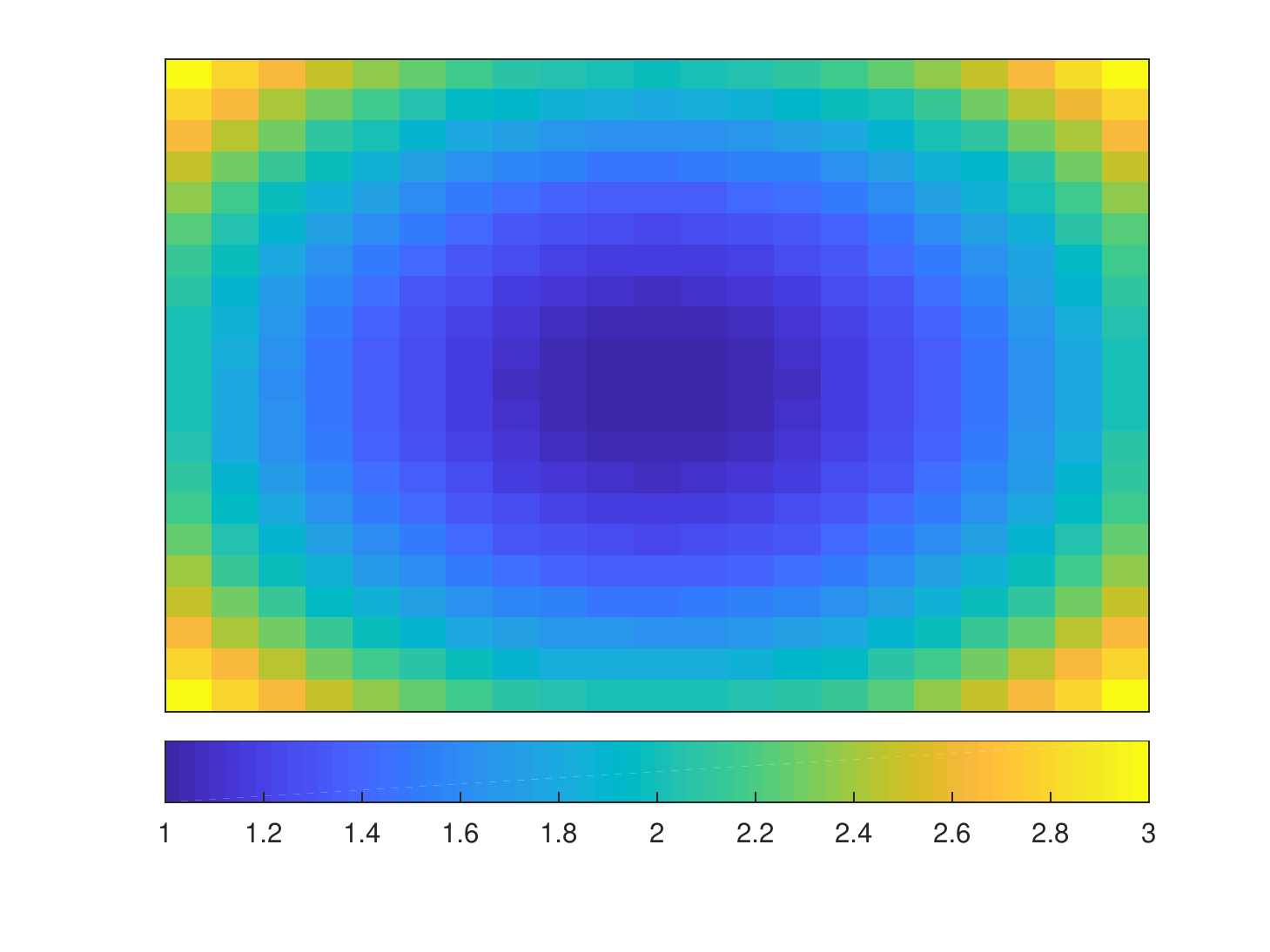}
&\includegraphics[width=\wsize,height=\hhsize,trim={1cm .5cm .5cm .5cm },clip]{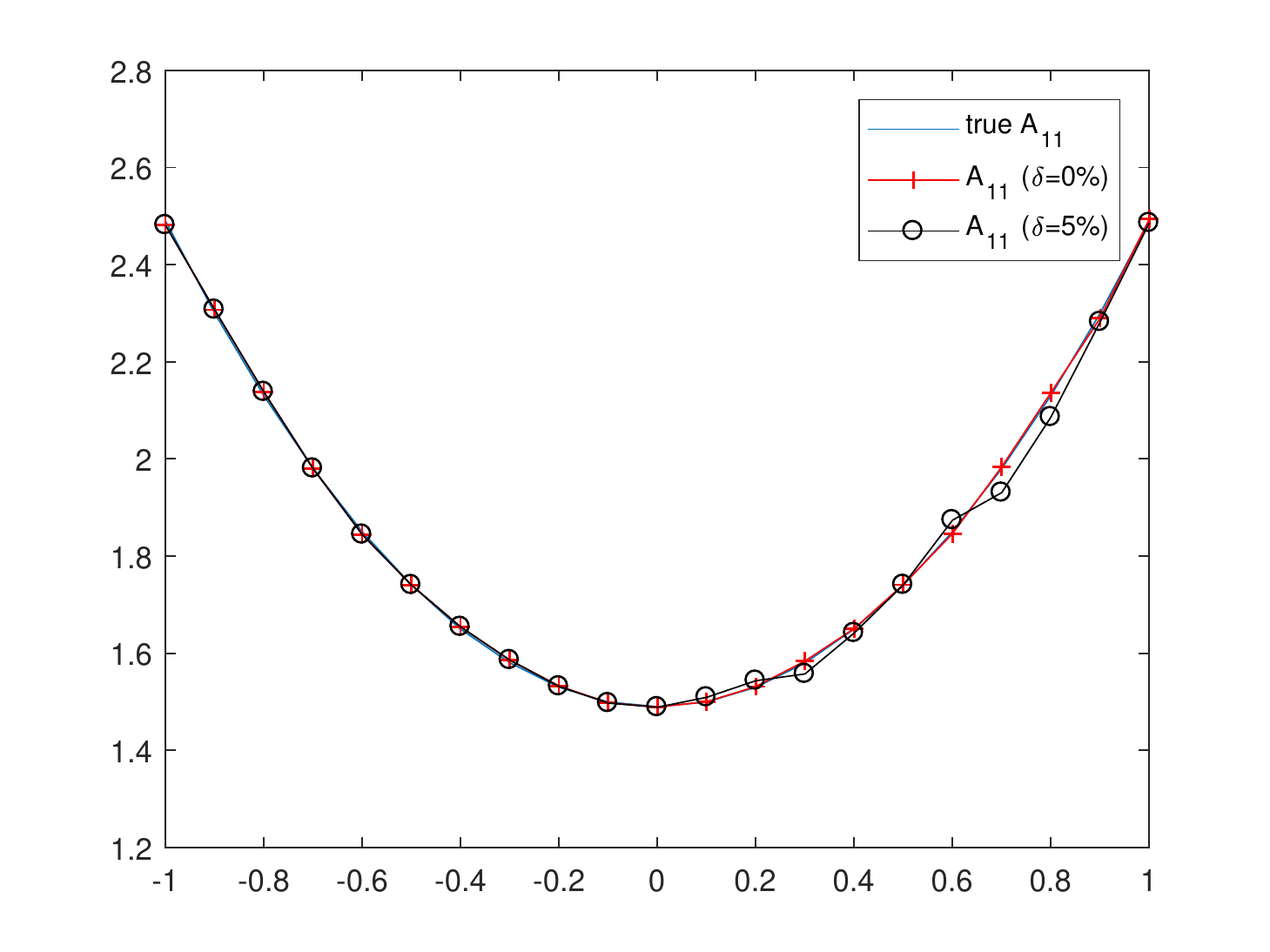}\\
 \includegraphics[width=\wsize,height=\hhsize,trim={1.5cm 1cm 1cm 0.5cm},clip]{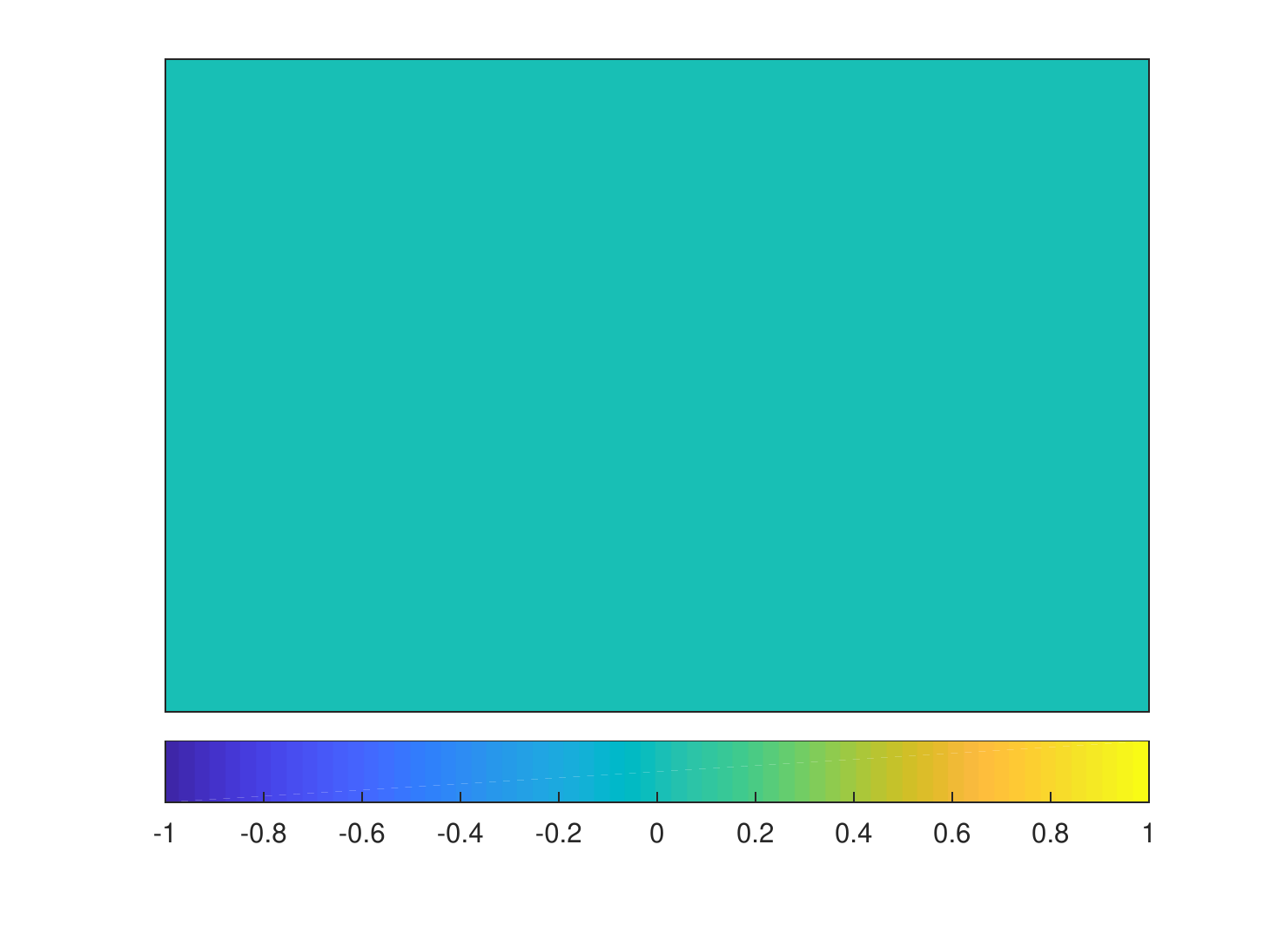}
&\includegraphics[width=\wsize,height=\hhsize,trim={1.5cm 1cm 1cm 0.5cm},clip]{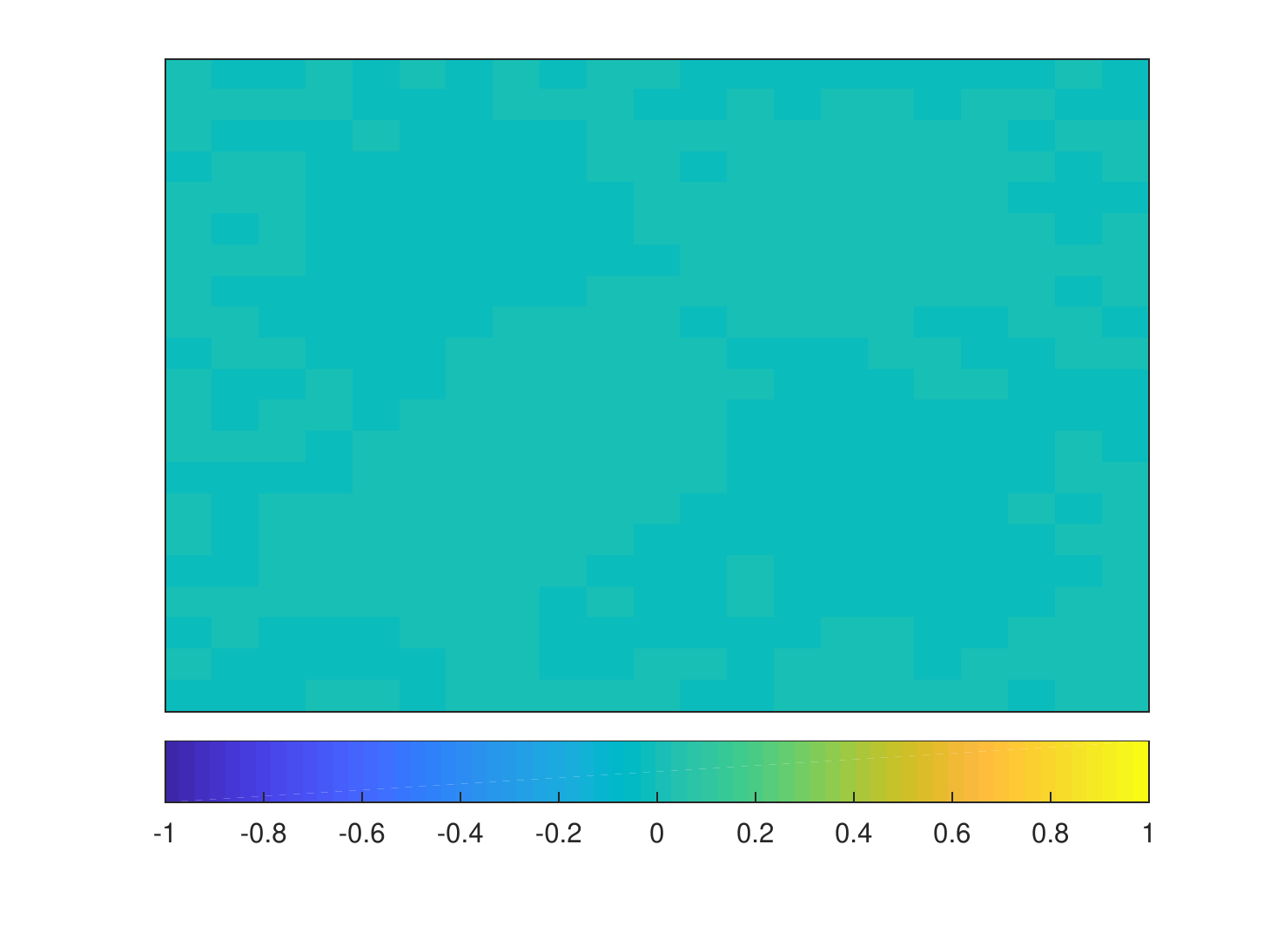}
&\includegraphics[width=\wsize,height=\hhsize,trim={1.5cm 1cm 1cm 0.5cm},clip]{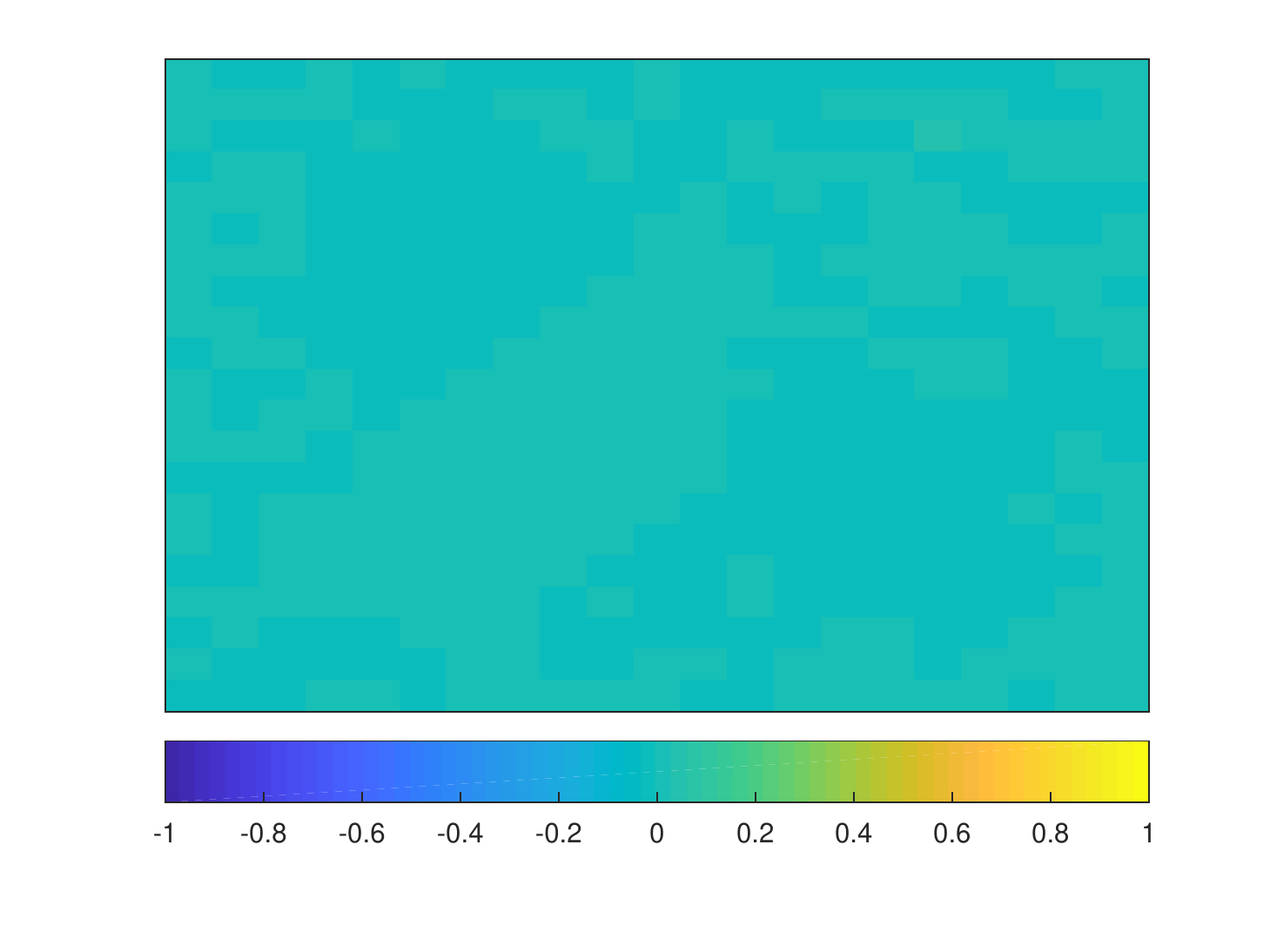}
&\includegraphics[width=\wsize,height=\hhsize,trim={1cm 0.5cm .5cm .5cm},clip]{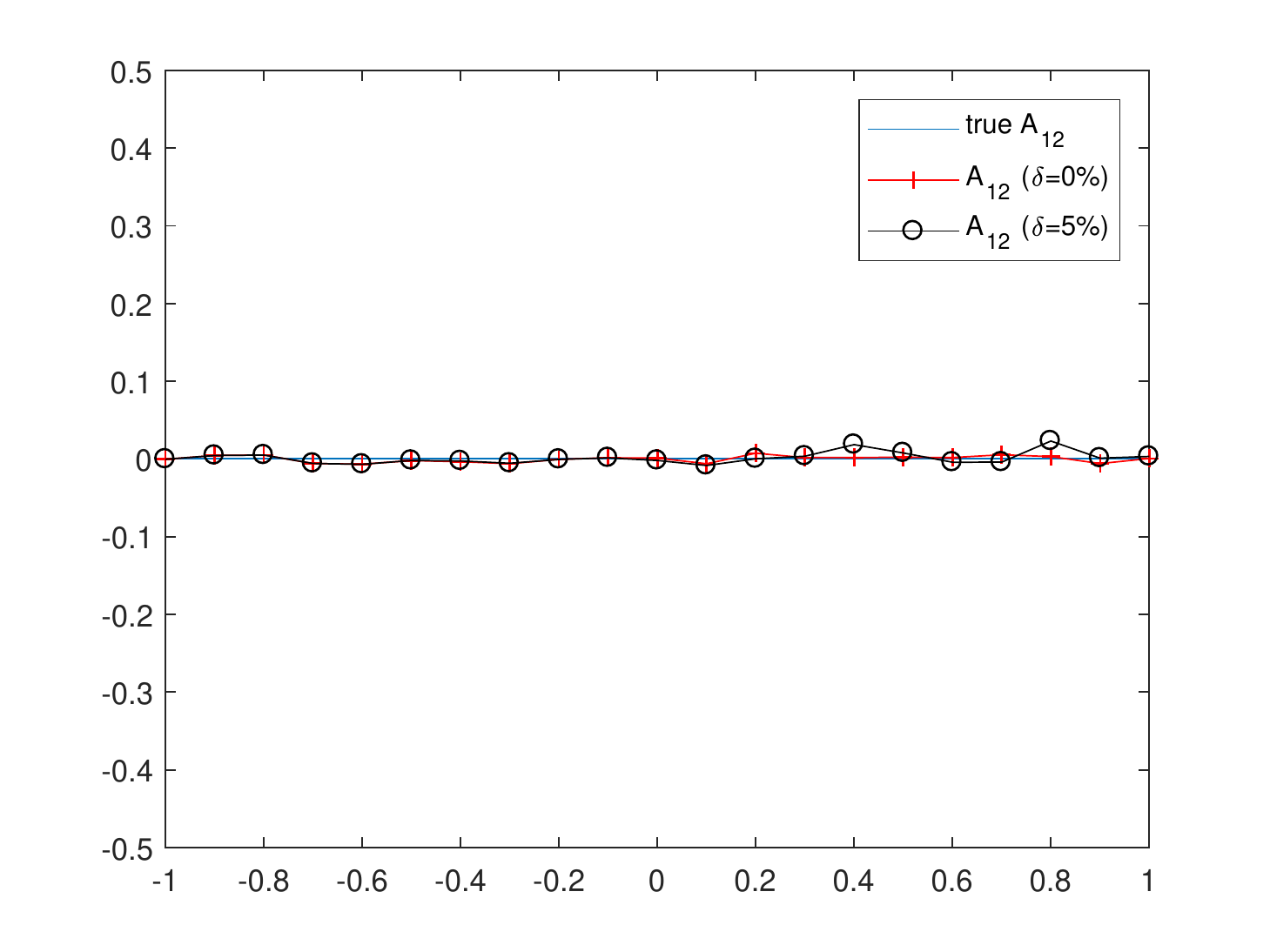}\\
 \includegraphics[width=\wsize,height=\hhsize,trim={1.5cm 1cm 1cm 0.5cm},clip]{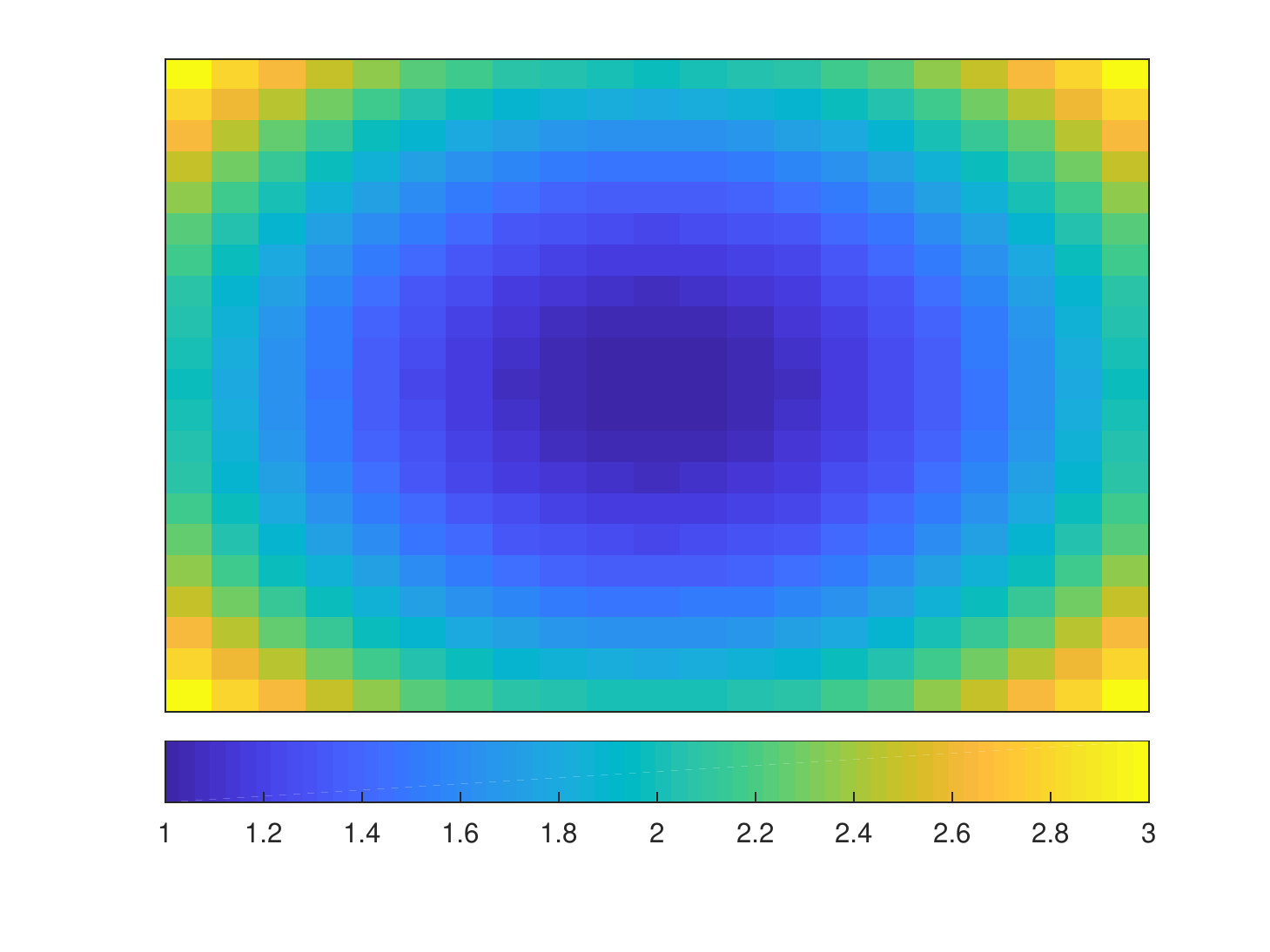}
&\includegraphics[width=\wsize,height=\hhsize,trim={1.5cm 1cm 1cm 0.5cm},clip]{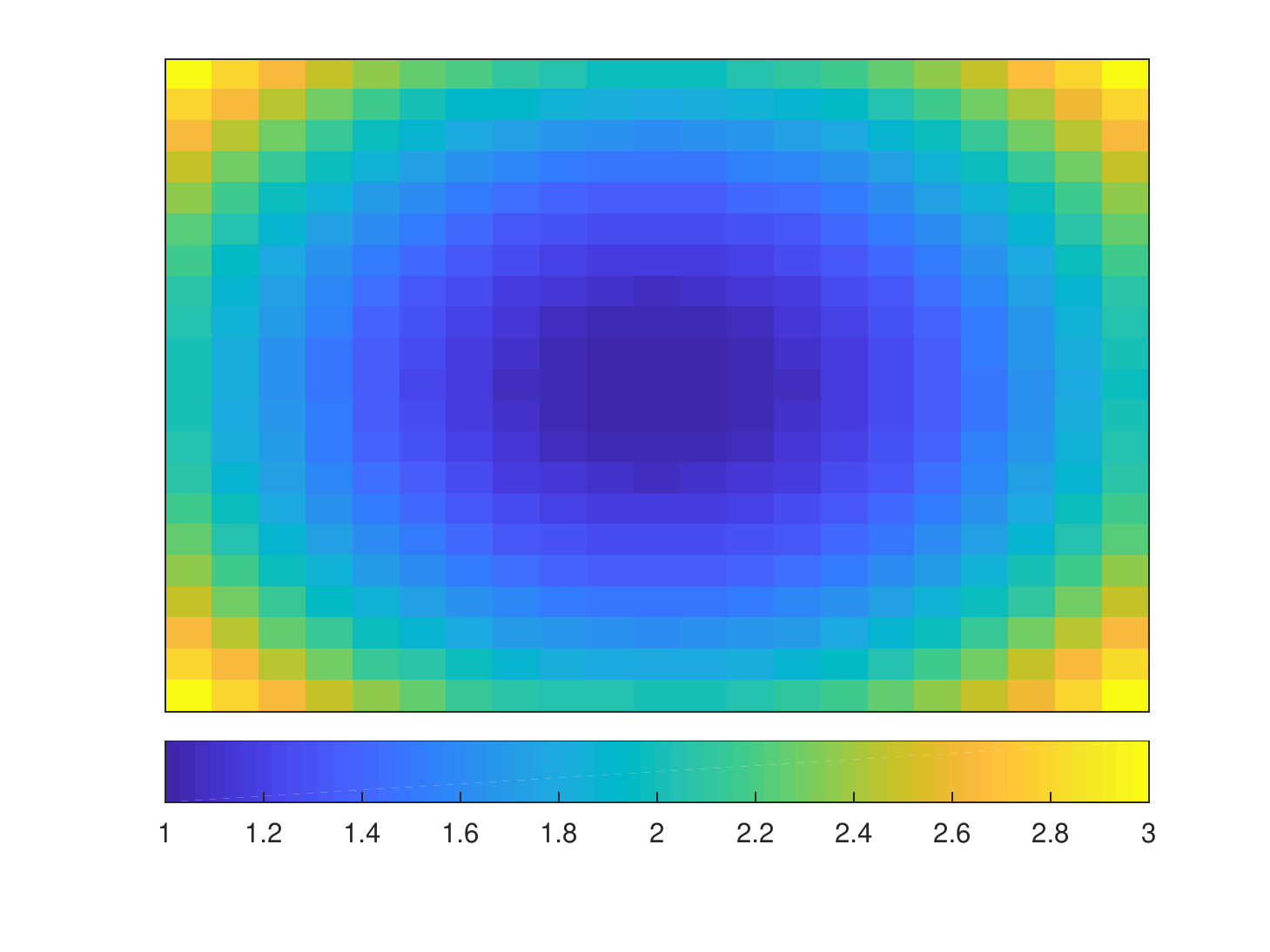}
&\includegraphics[width=\wsize,height=\hhsize,trim={1.5cm 1cm 1cm 0.5cm},clip]{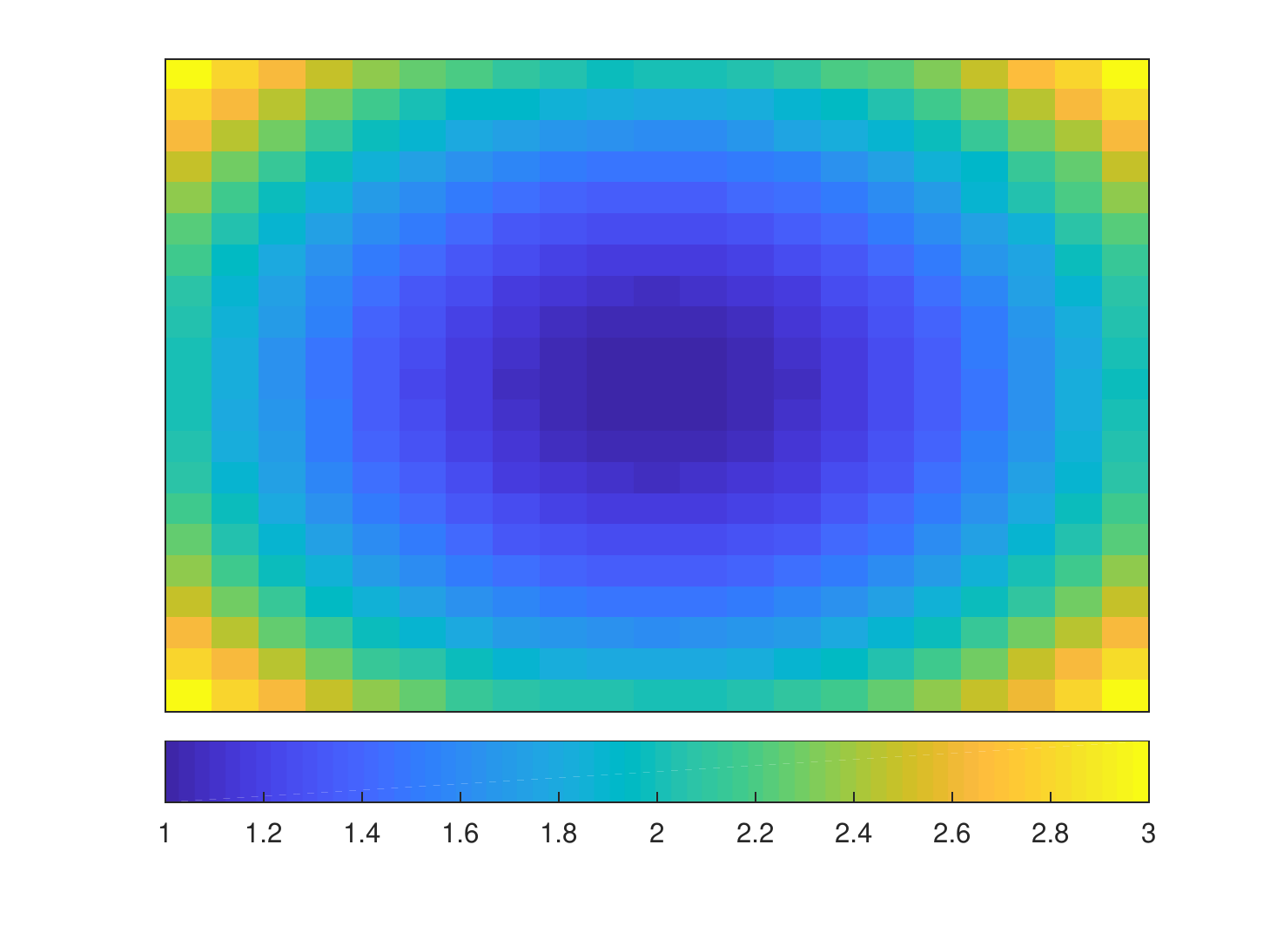}
&\includegraphics[width=\wsize,height=\hhsize,trim={1cm .5cm .5cm .5cm },clip]{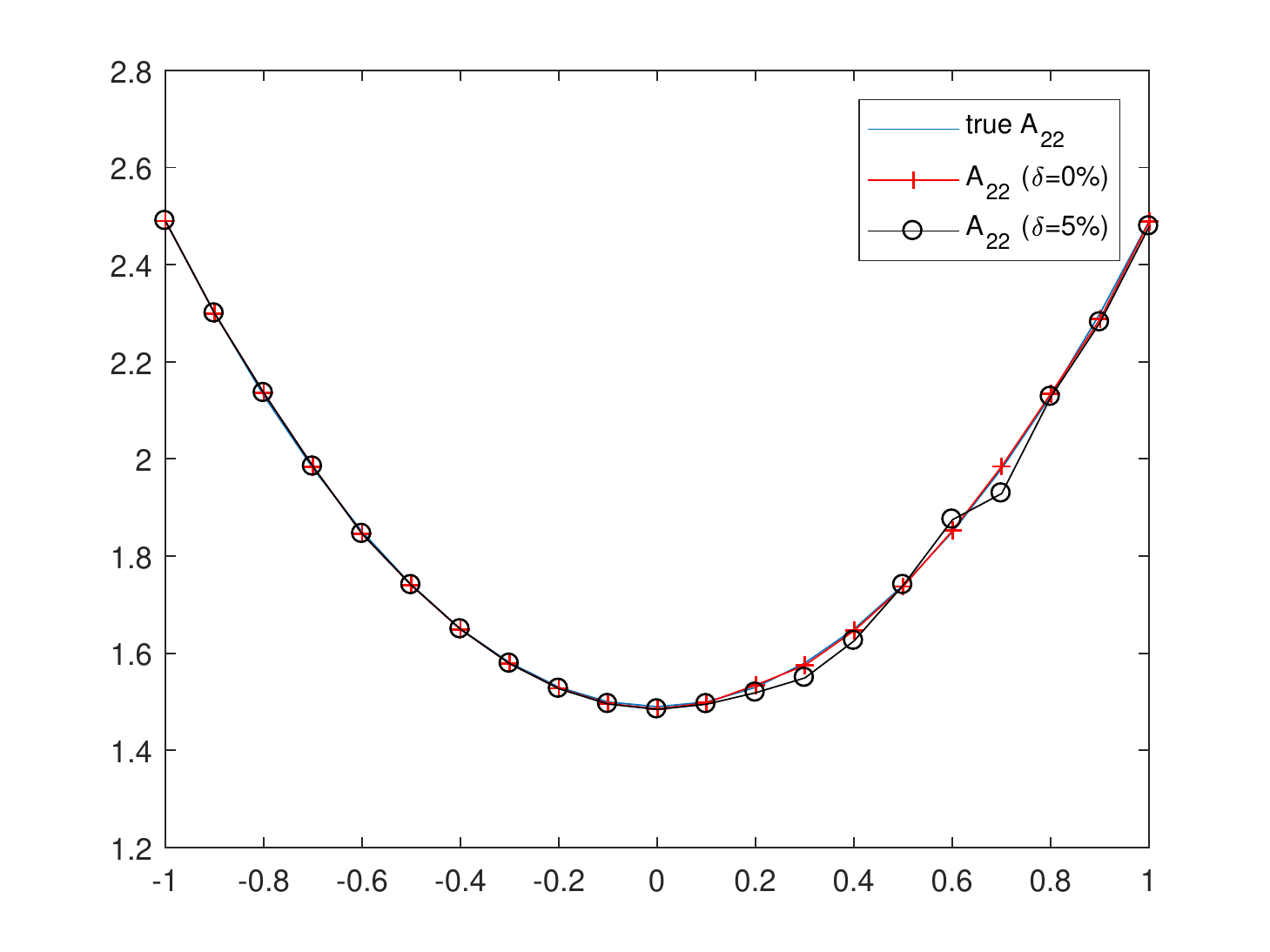}
\end{tabular}
\caption{The reconstructions for the conductivity tensor for Example 2: the top, middle and bottom rows refer to $A_{11}$, $A_{12}$ and $A_{22}$,
respectively. Column 1 is for exact conductivity, Column 2 for reconstruction with exact data, Column 3 for noisy data with $\delta=5\%$, and column 4 for the cross section along \{$x_2=-0.7$\}.}
\label{exam2}
\end{figure}

\begin{figure}[!htbp]
\centering
\setlength{\tabcolsep}{0pt}
\begin{tabular}{cccc}
 \includegraphics[width=\wsize,height=\hhsize,trim={1.5cm 1cm 1cm 0.5cm},clip]{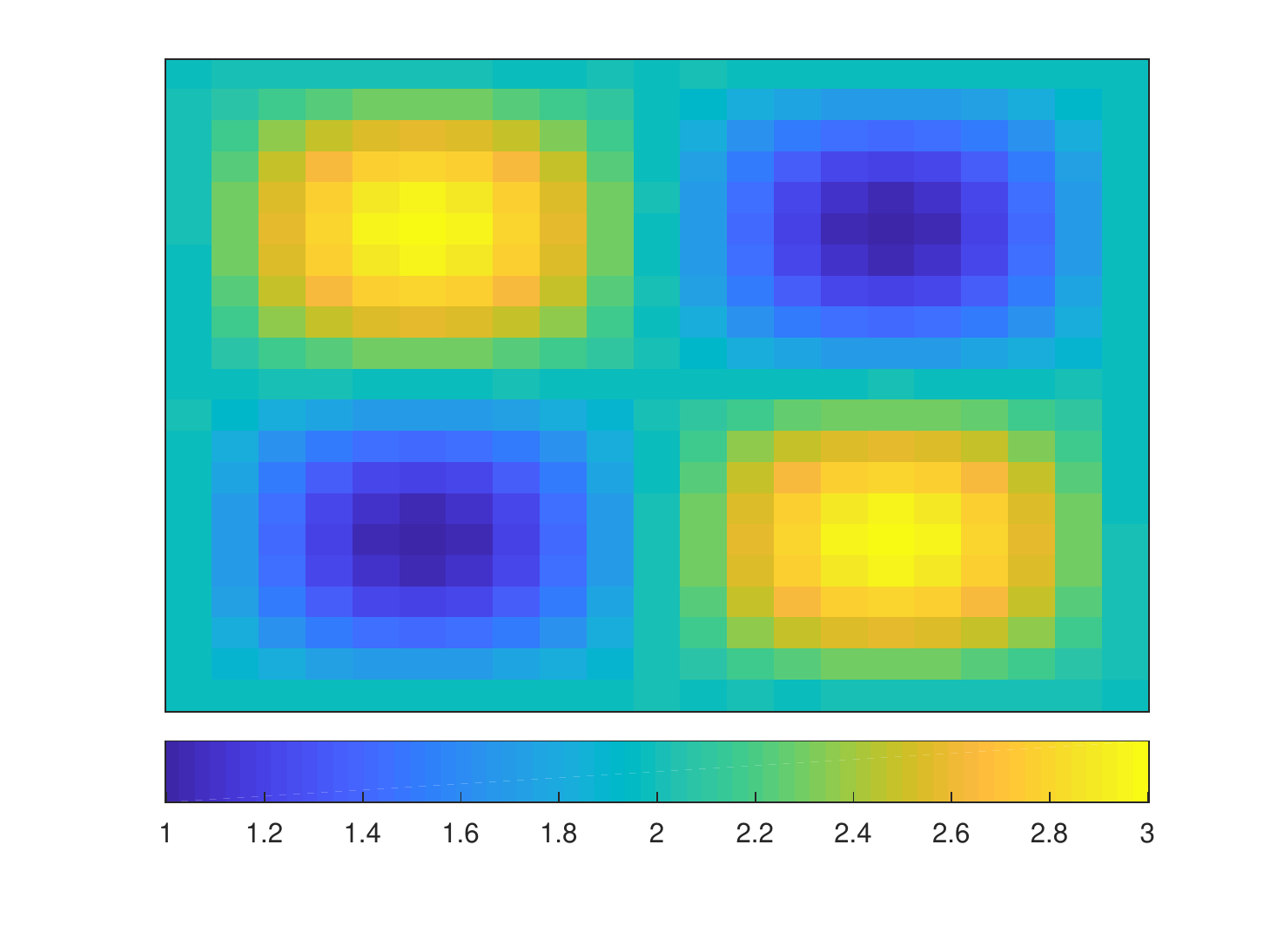}
&\includegraphics[width=\wsize,height=\hhsize,trim={1.5cm 1cm 1cm 0.5cm},clip]{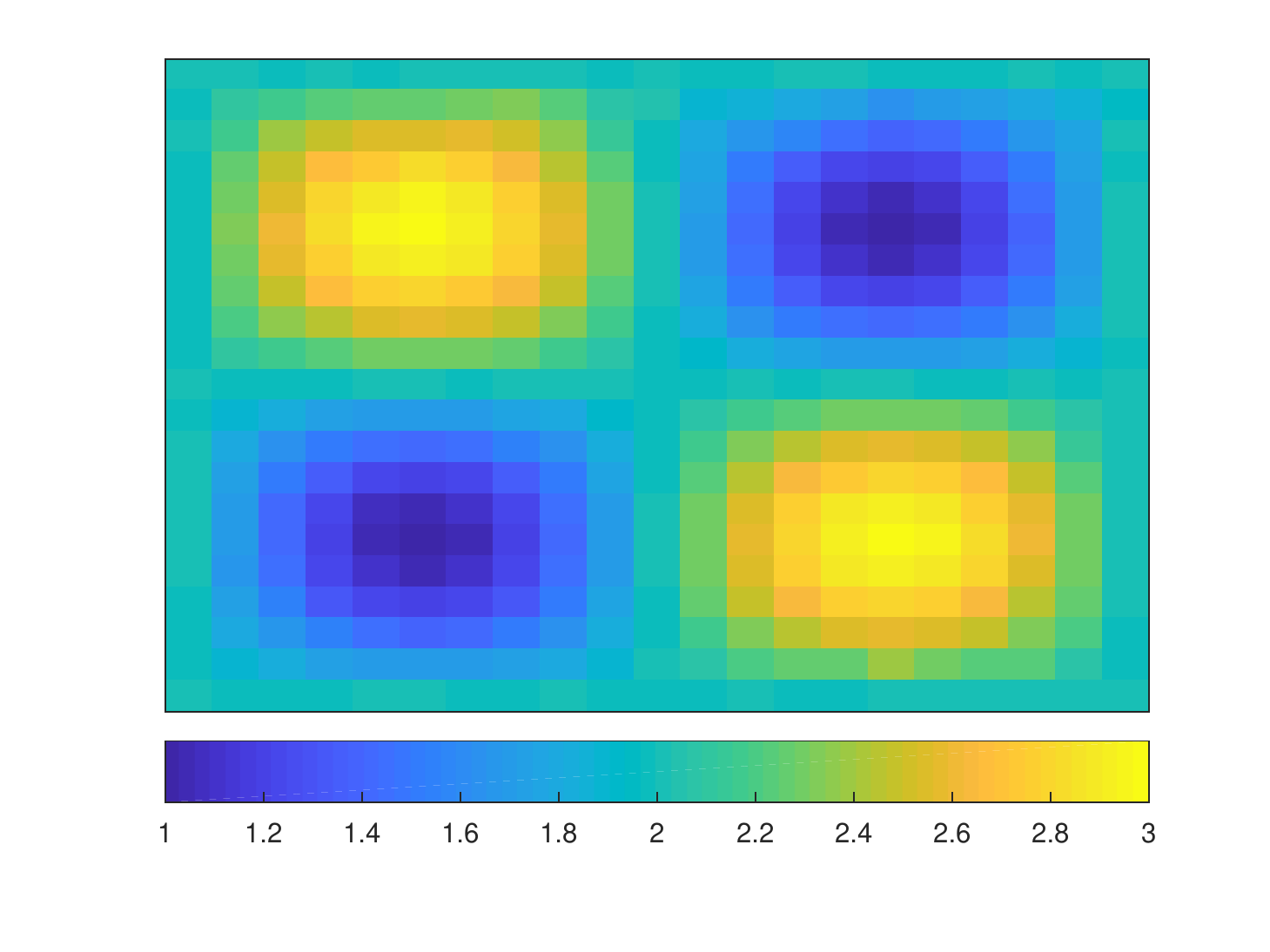}
&\includegraphics[width=\wsize,height=\hhsize,trim={1.5cm 1cm 1cm 0.5cm},clip]{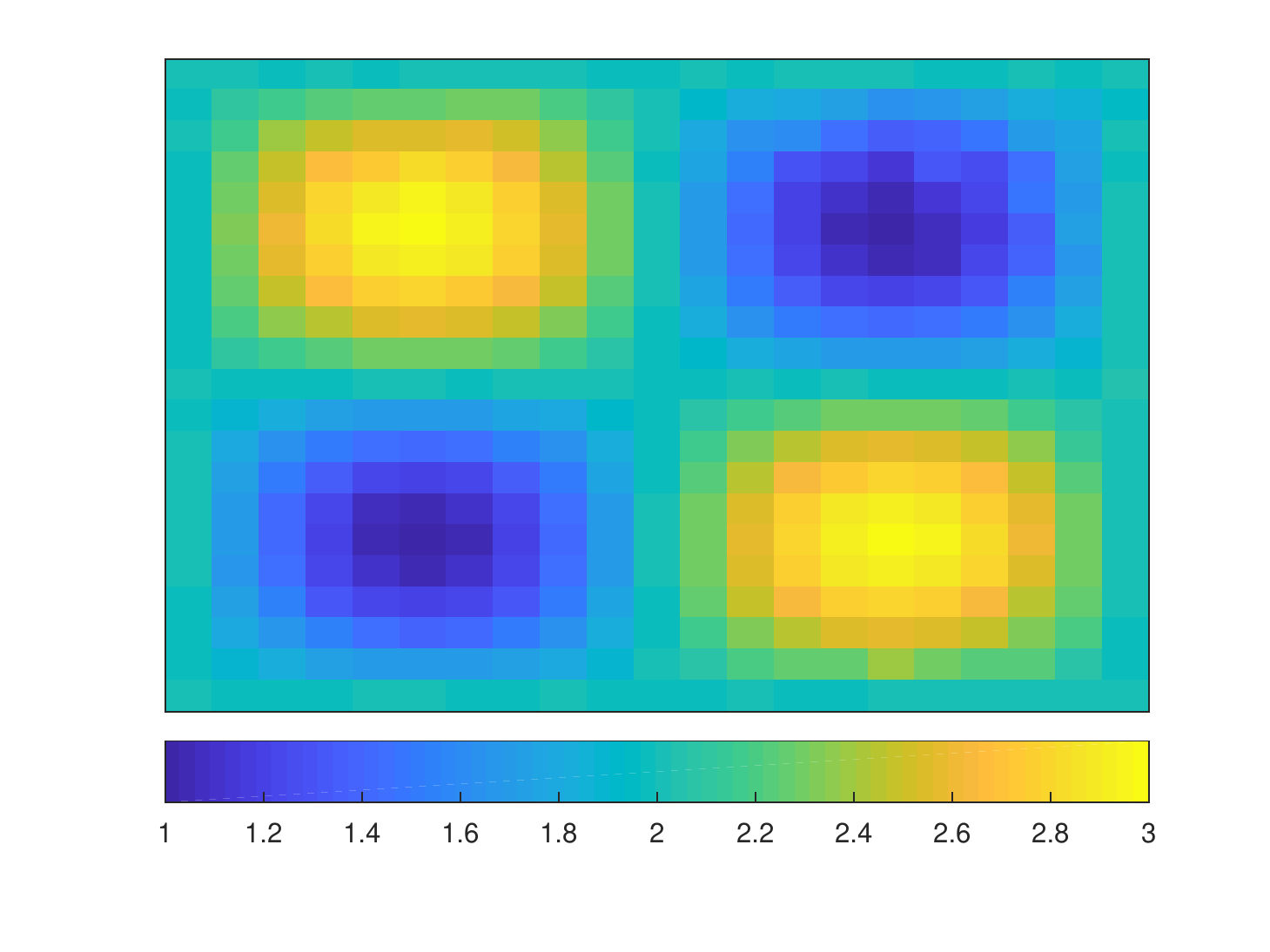}
&\includegraphics[width=\wsize,height=\hhsize,trim={1cm .5cm .5cm .5cm },clip]{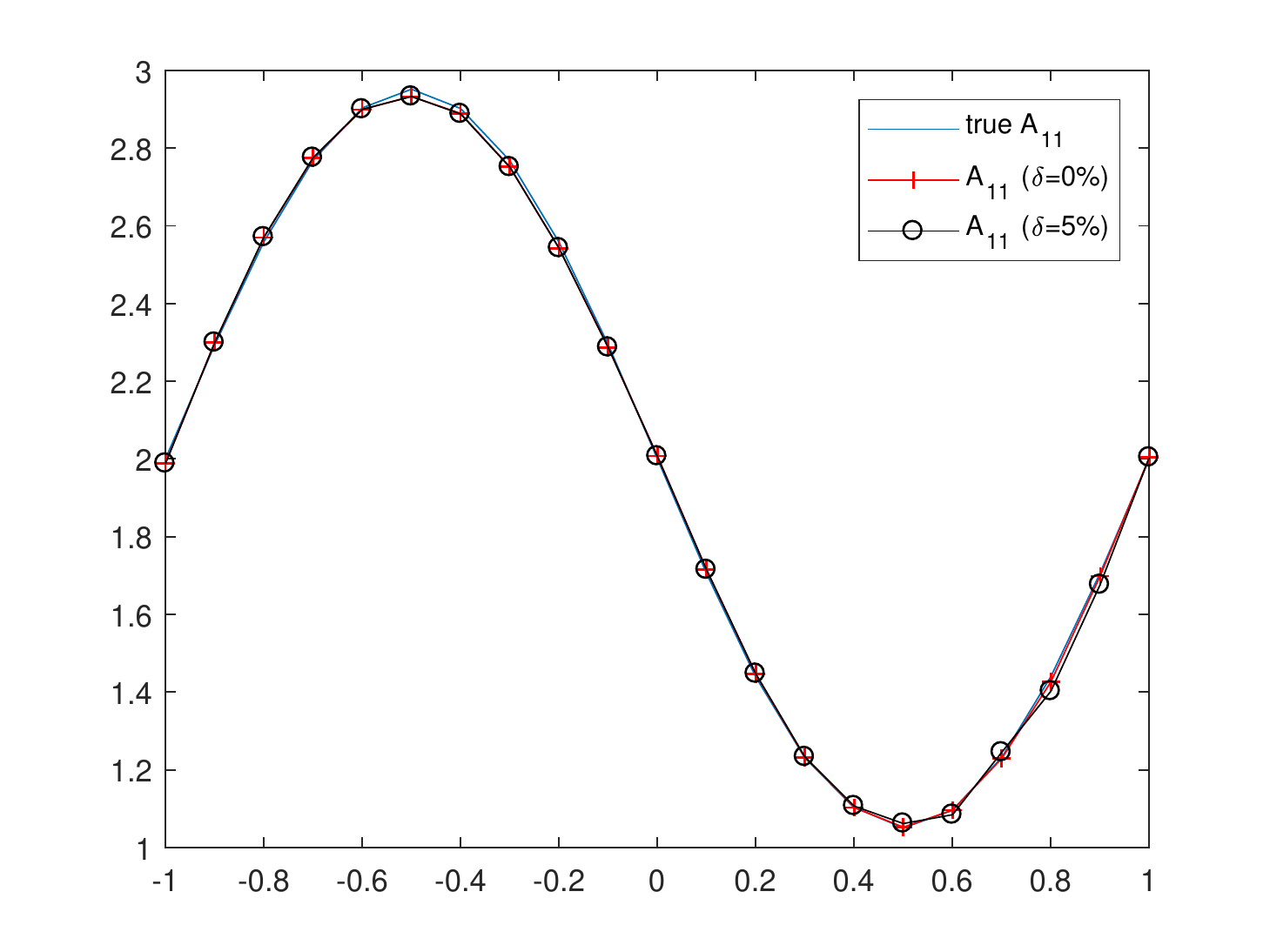}\\
 \includegraphics[width=\wsize,height=\hhsize,trim={1.5cm 1cm 1cm 0.5cm},clip]{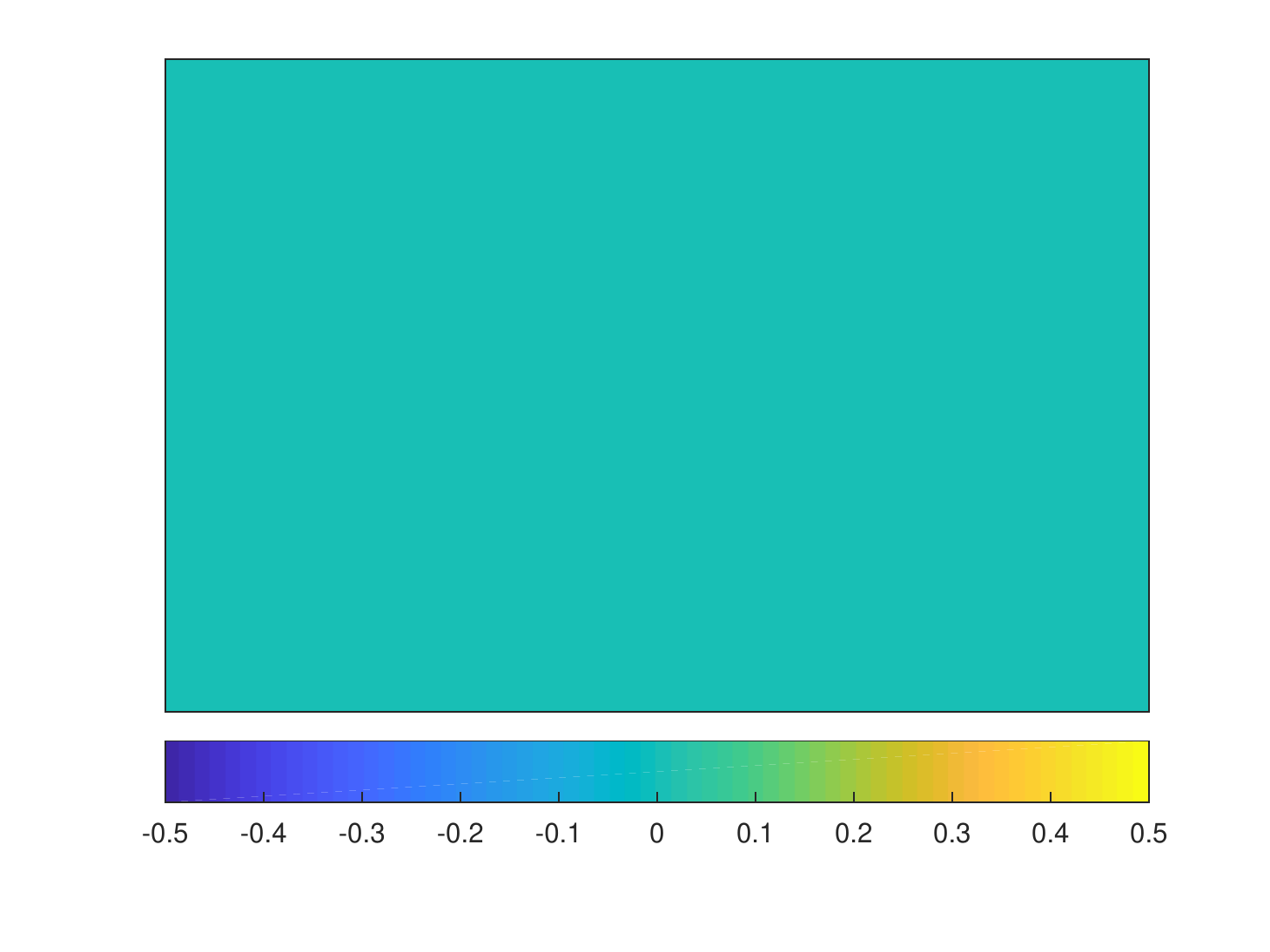}
&\includegraphics[width=\wsize,height=\hhsize,trim={1.5cm 1cm 1cm 0.5cm},clip]{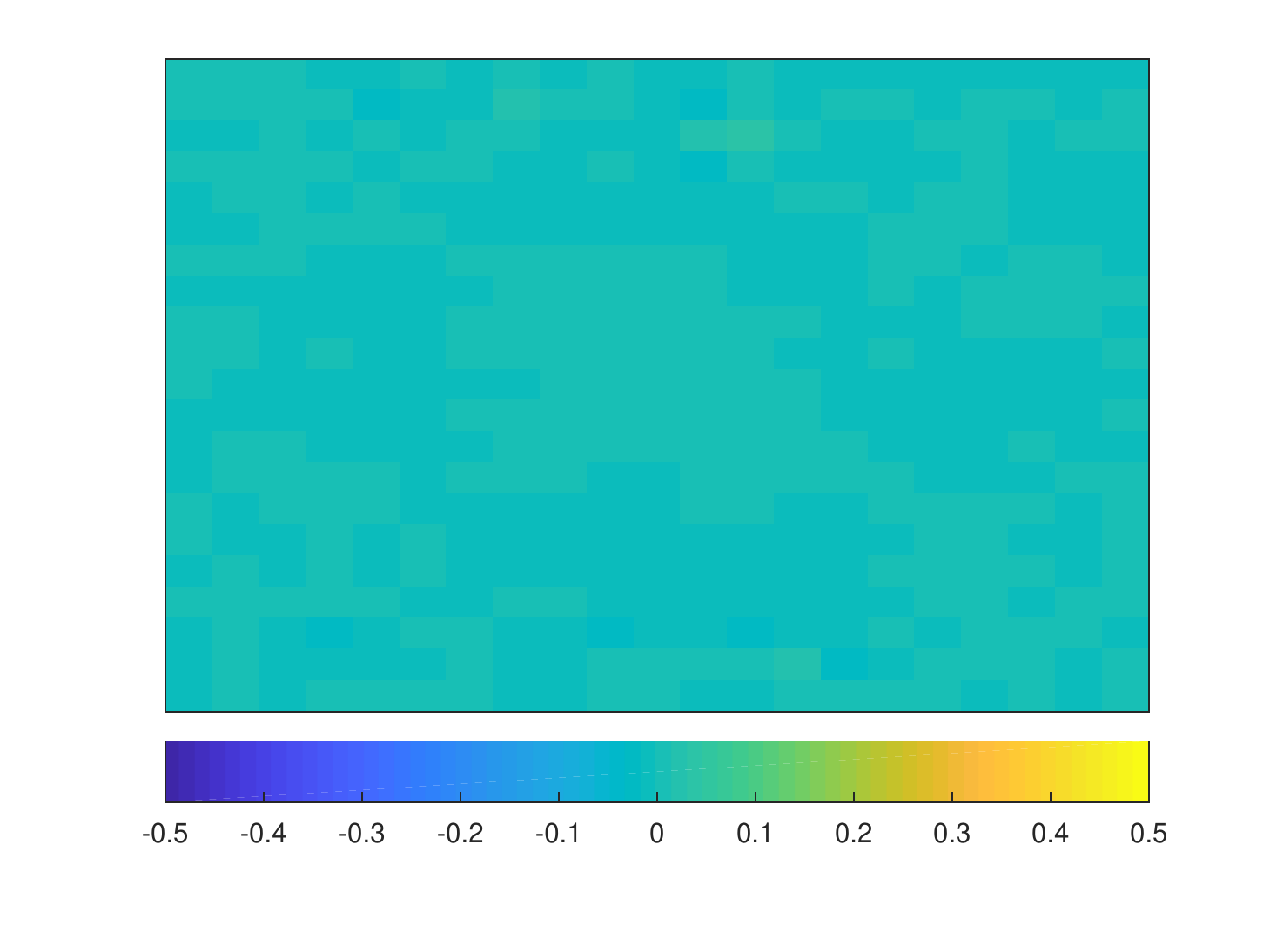}
&\includegraphics[width=\wsize,height=\hhsize,trim={1.5cm 1cm 1cm 0.5cm},clip]{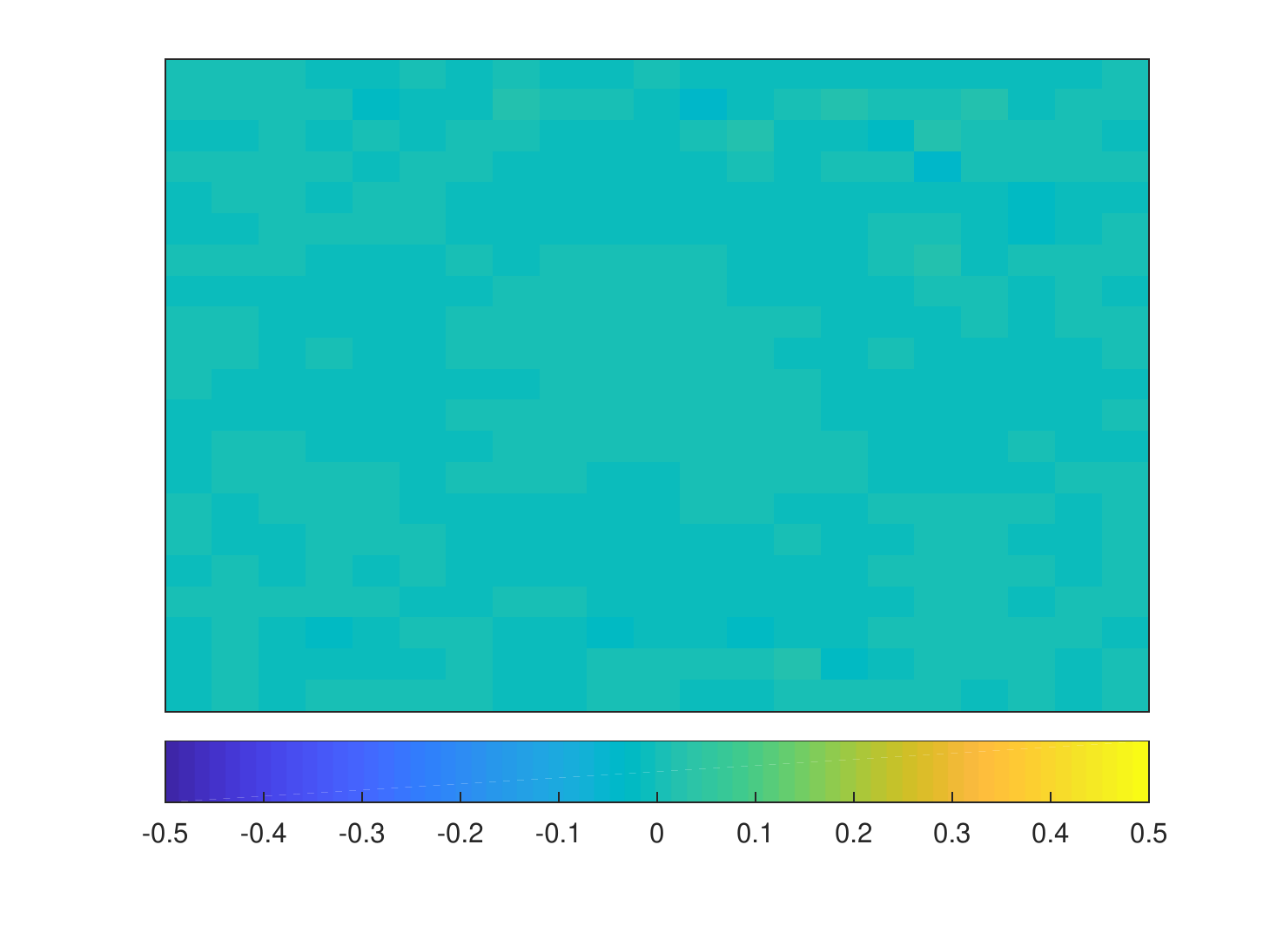}
&\includegraphics[width=\wsize,height=\hhsize,trim={1cm 0.5cm .5cm .5cm},clip]{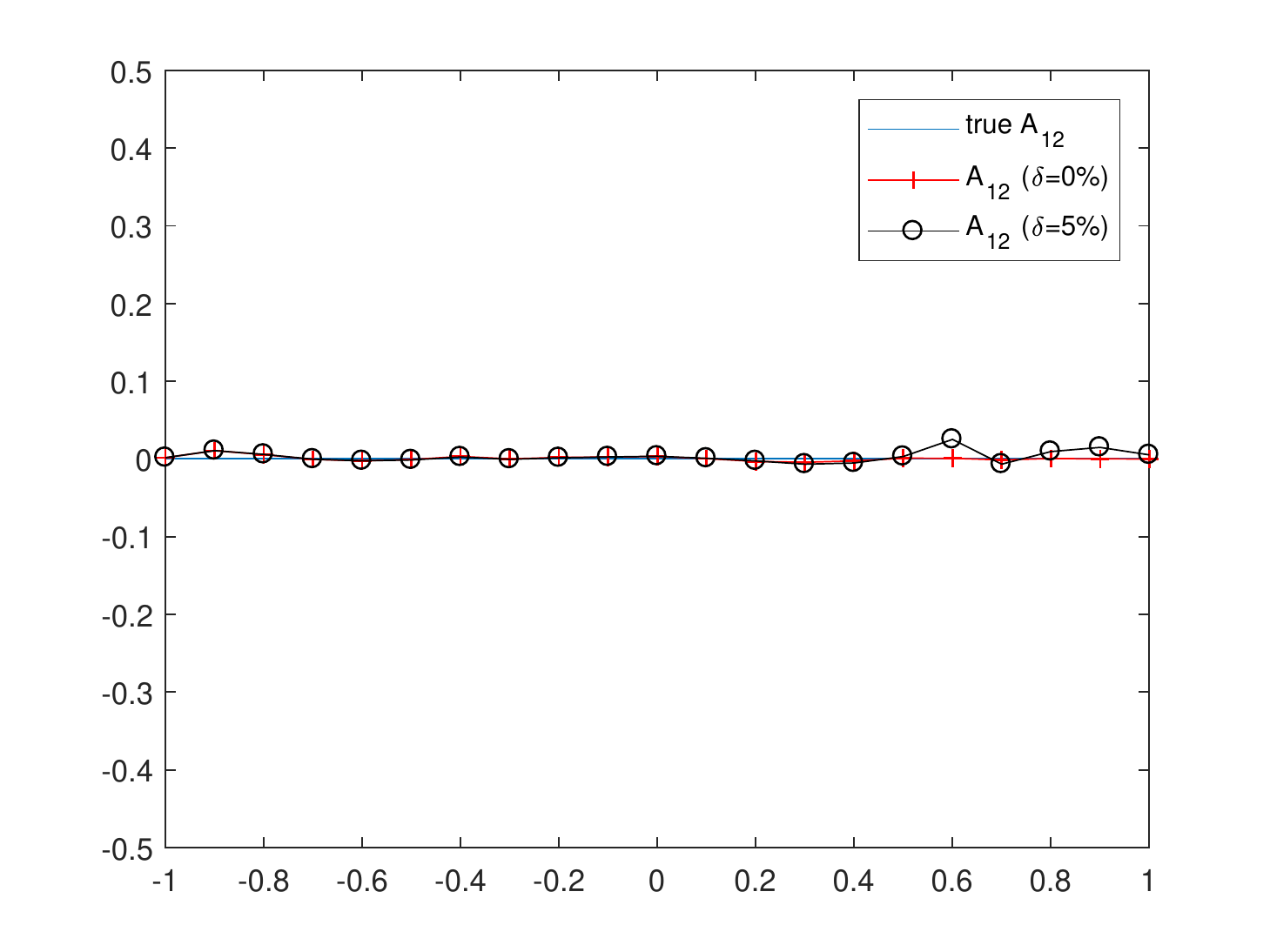}\\
 \includegraphics[width=\wsize,height=\hhsize,trim={1.5cm 1cm 1cm 0.5cm},clip]{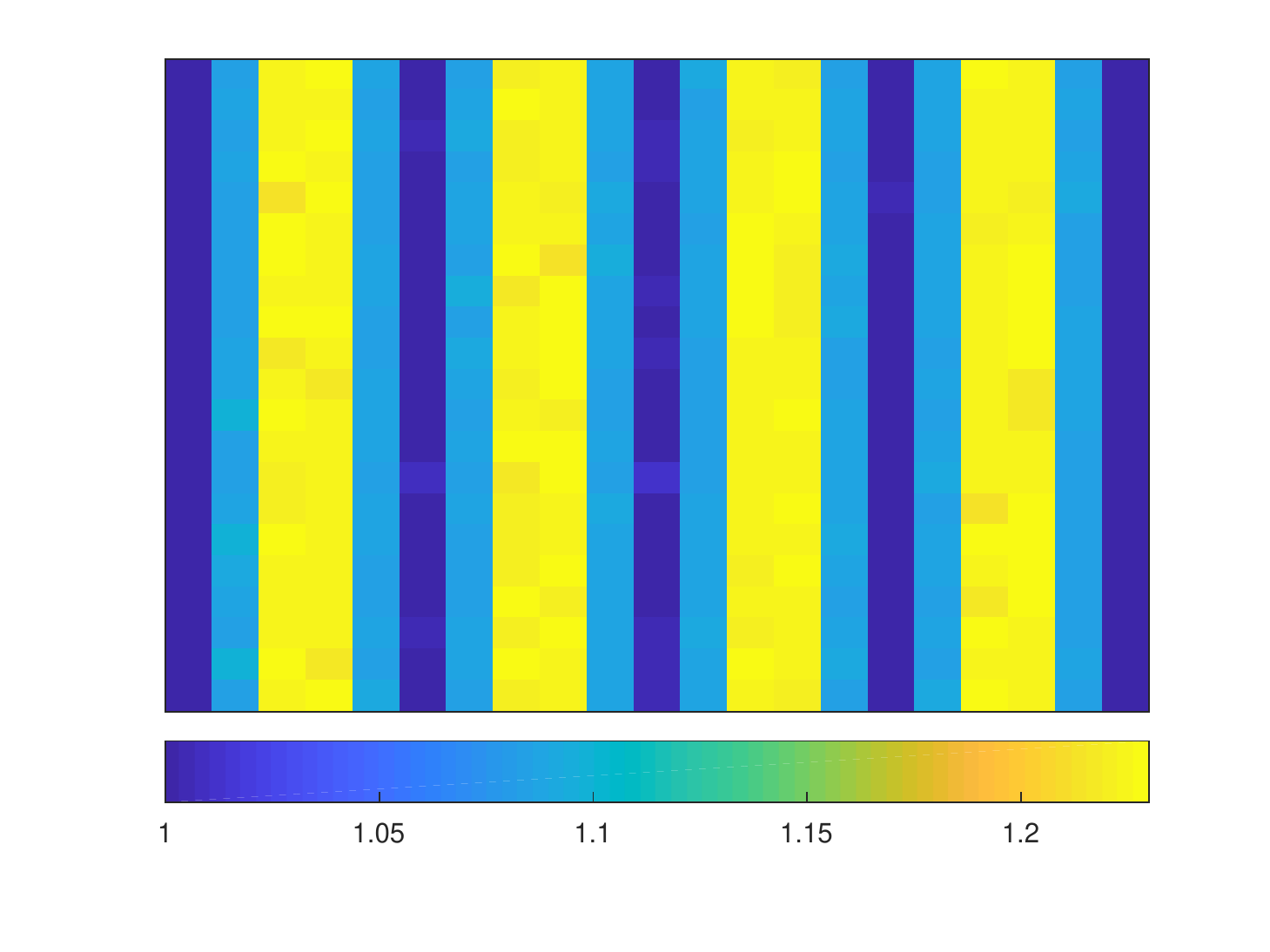}
&\includegraphics[width=\wsize,height=\hhsize,trim={1.5cm 1cm 1cm 0.5cm},clip]{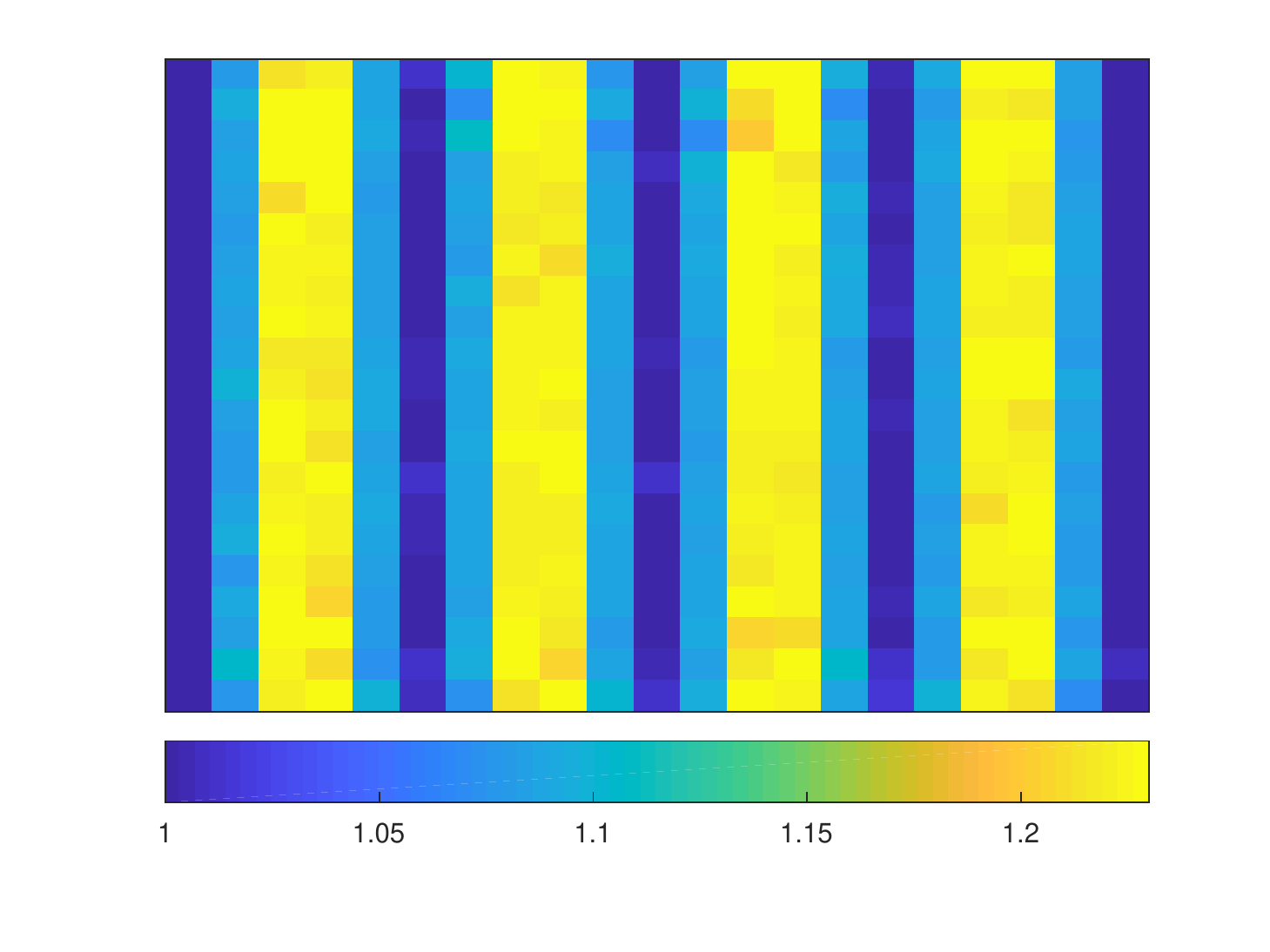}
&\includegraphics[width=\wsize,height=\hhsize,trim={1.5cm 1cm 1cm 0.5cm},clip]{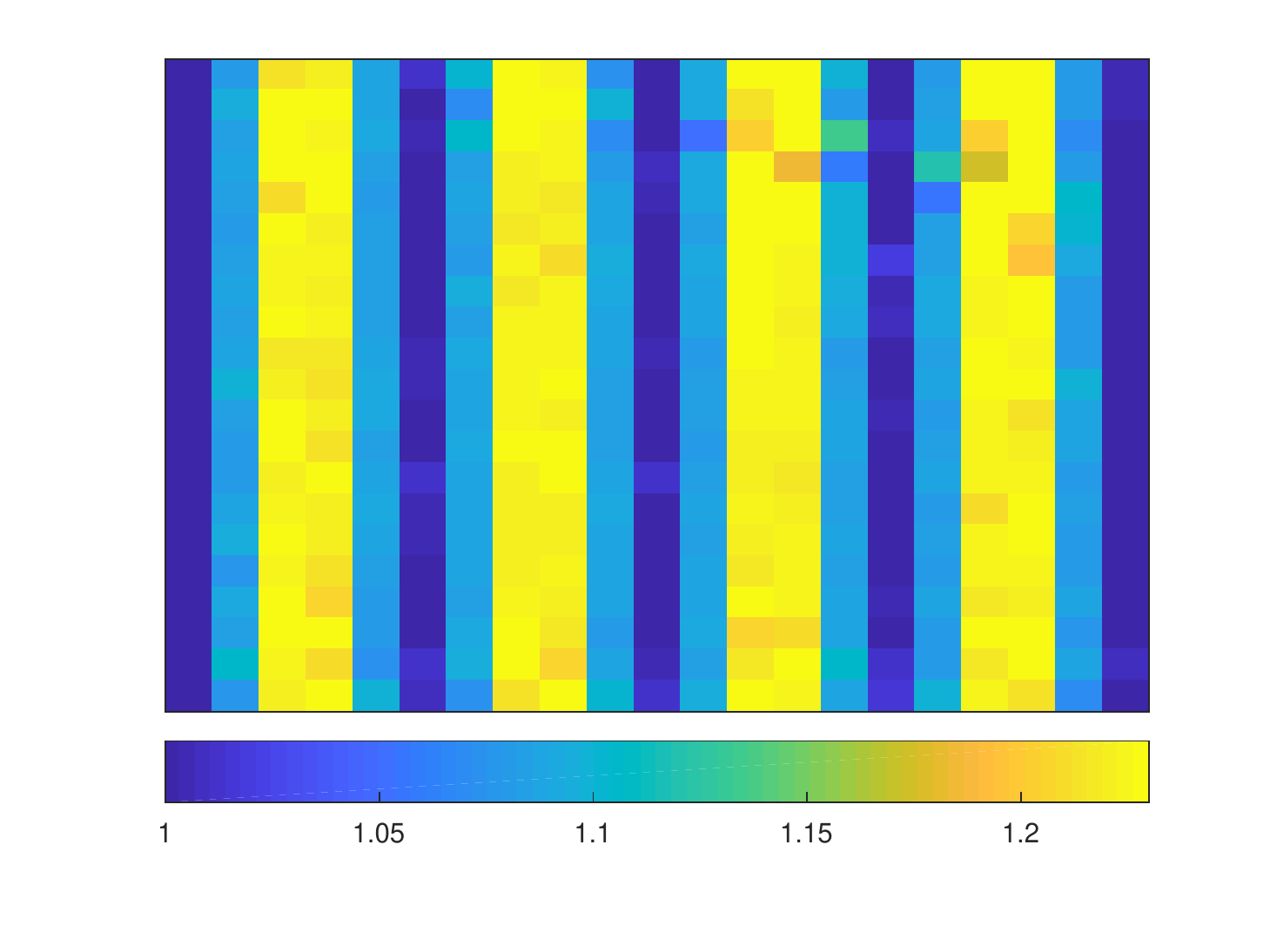}
&\includegraphics[width=\wsize,height=\hhsize,trim={1cm .5cm .5cm .5cm },clip]{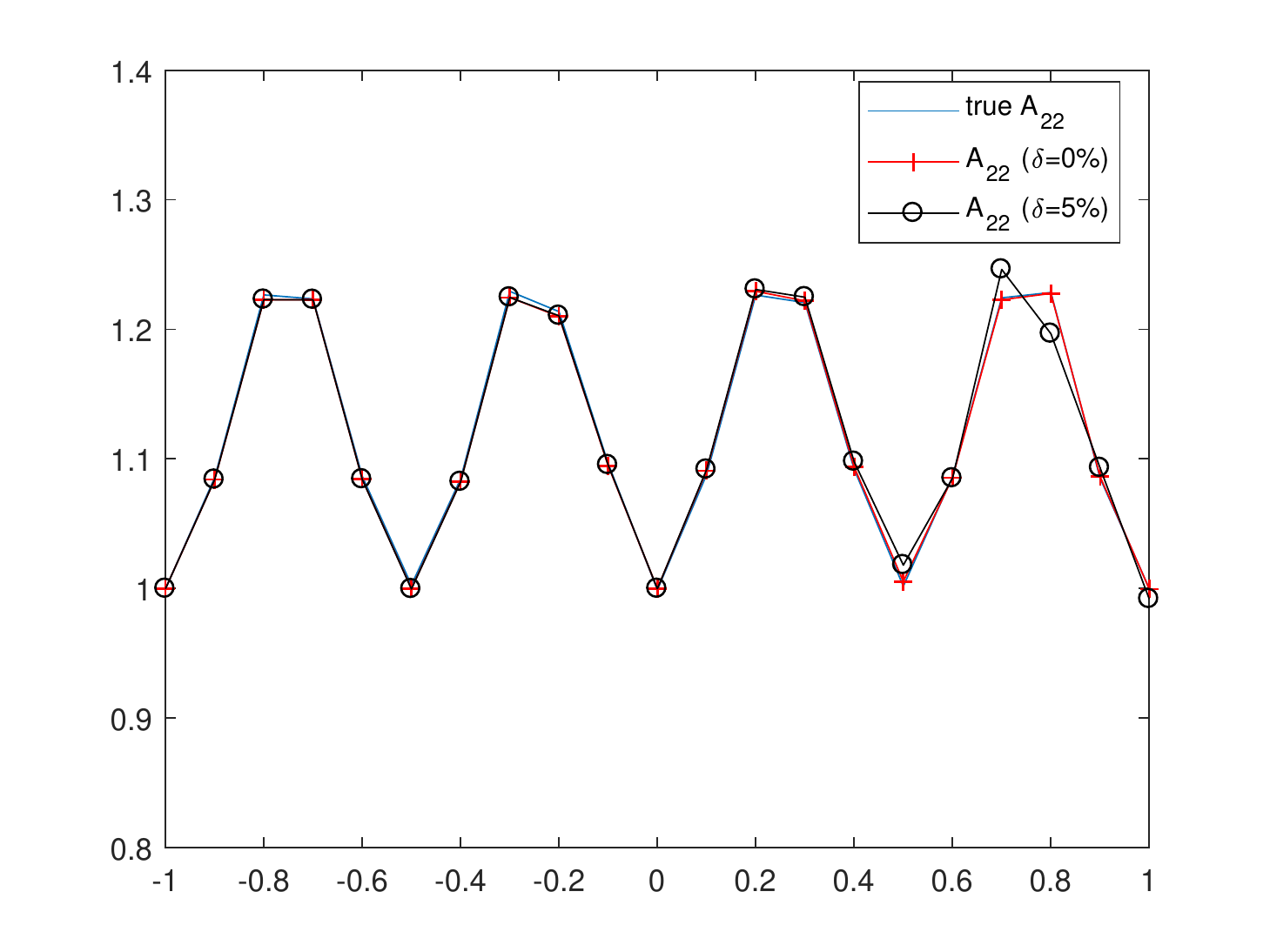}\\
\end{tabular}
\caption{The reconstructions for the conductivity tensor for Example 3: the top, middle and bottom rows refer to $A_{11}$, $A_{12}$ and $A_{22}$,
respectively. Column 1 is for exact conductivity, Column 2 for reconstruction with exact data, Column 3 for noisy data with $\delta=5\%$, and column 4 for the cross section along \{$x_2=-0.4$\}.}
\label{exam3}
\end{figure}

\begin{figure}[!htbp]
\centering
\setlength{\tabcolsep}{0pt}
\begin{tabular}{cccc}
 \includegraphics[width=\wsize,height=\hhsize,trim={1.5cm 1cm 1cm 0.5cm},clip]{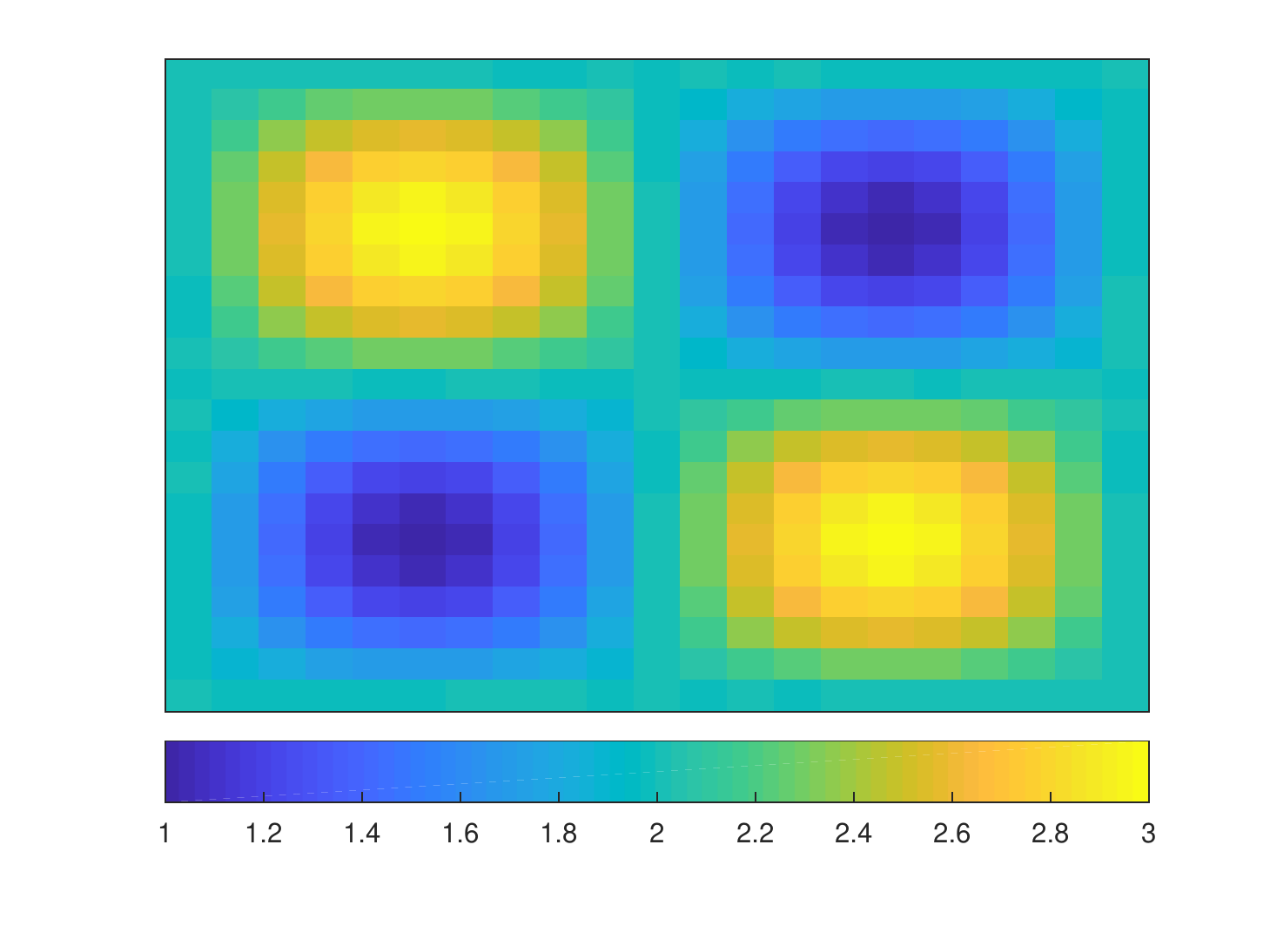}
&\includegraphics[width=\wsize,height=\hhsize,trim={1.5cm 1cm 1cm 0.5cm},clip]{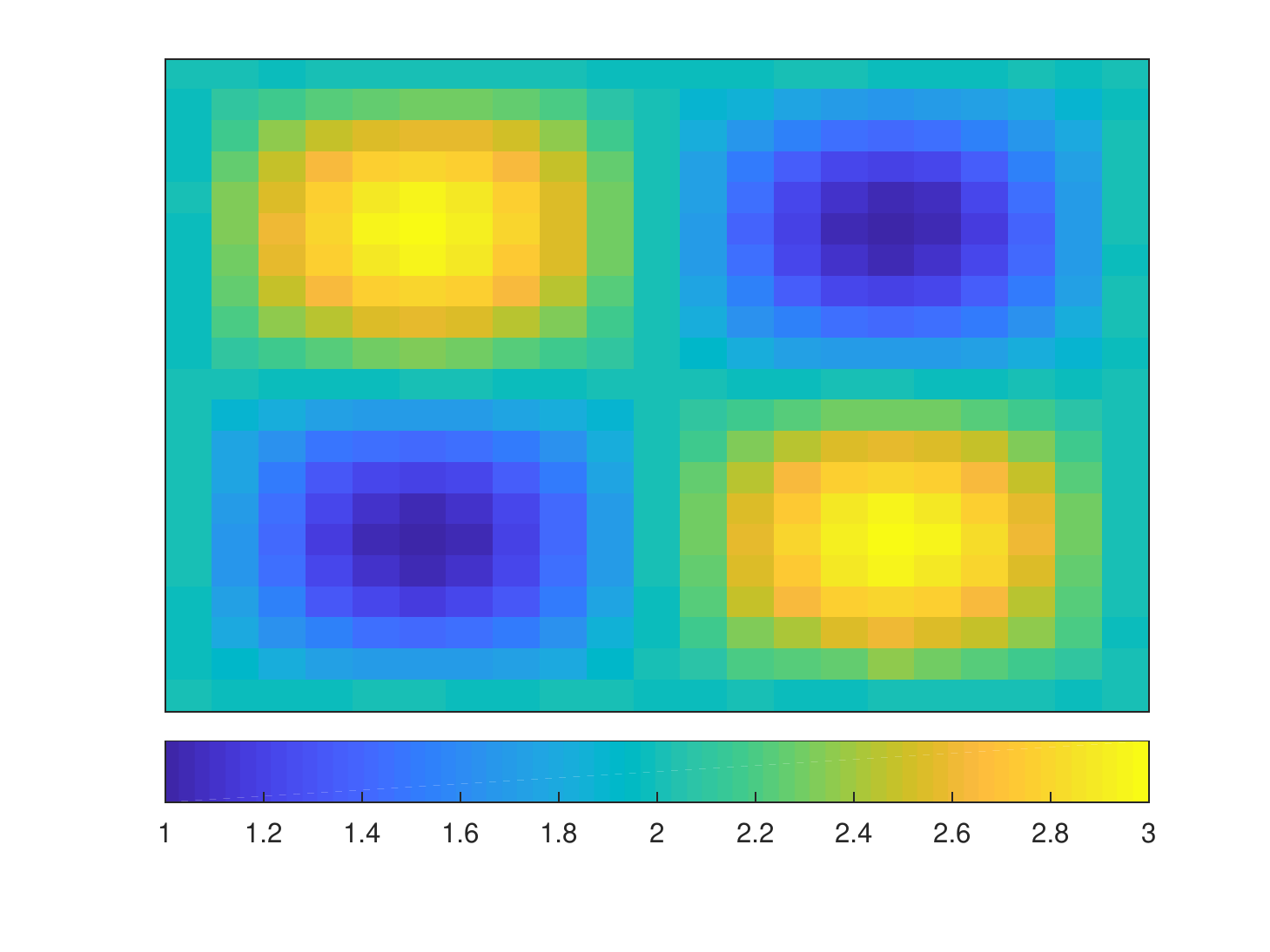}
&\includegraphics[width=\wsize,height=\hhsize,trim={1.5cm 1cm 1cm 0.5cm},clip]{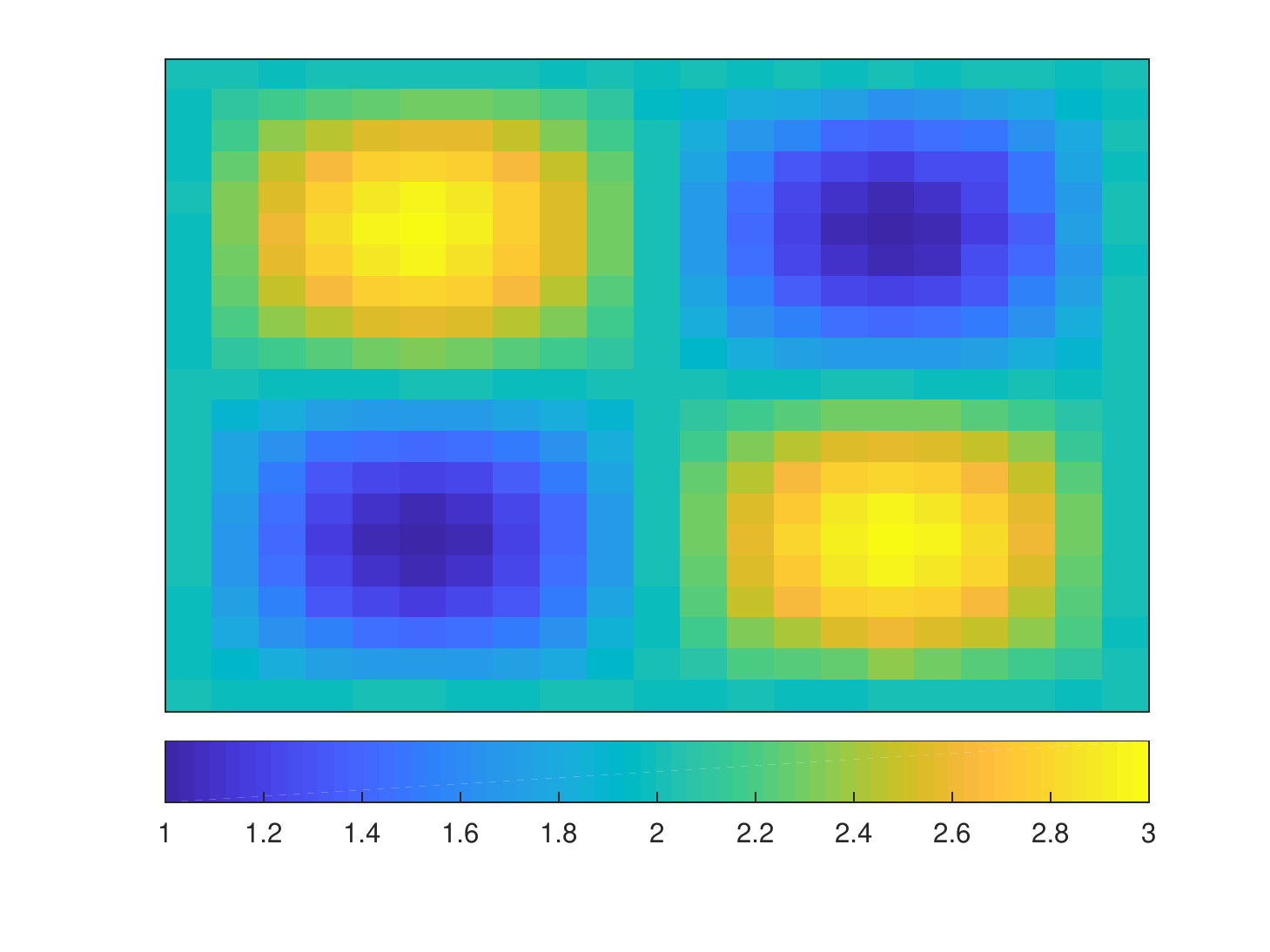}
&\includegraphics[width=\wsize,height=\hhsize,trim={1cm .5cm .5cm .5cm },clip]{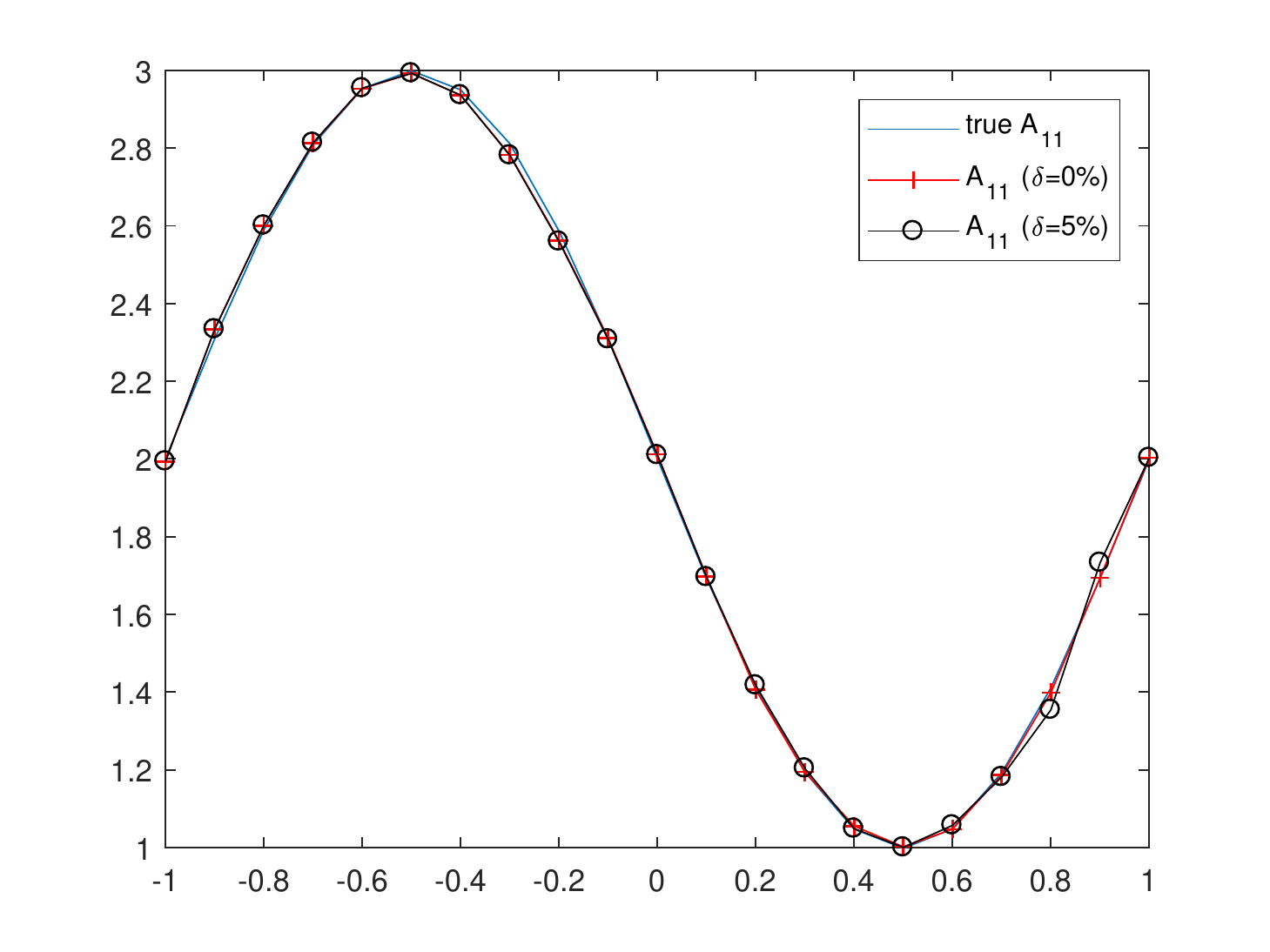}\\
 \includegraphics[width=\wsize,height=\hhsize,trim={1.5cm 1cm 1cm 0.5cm},clip]{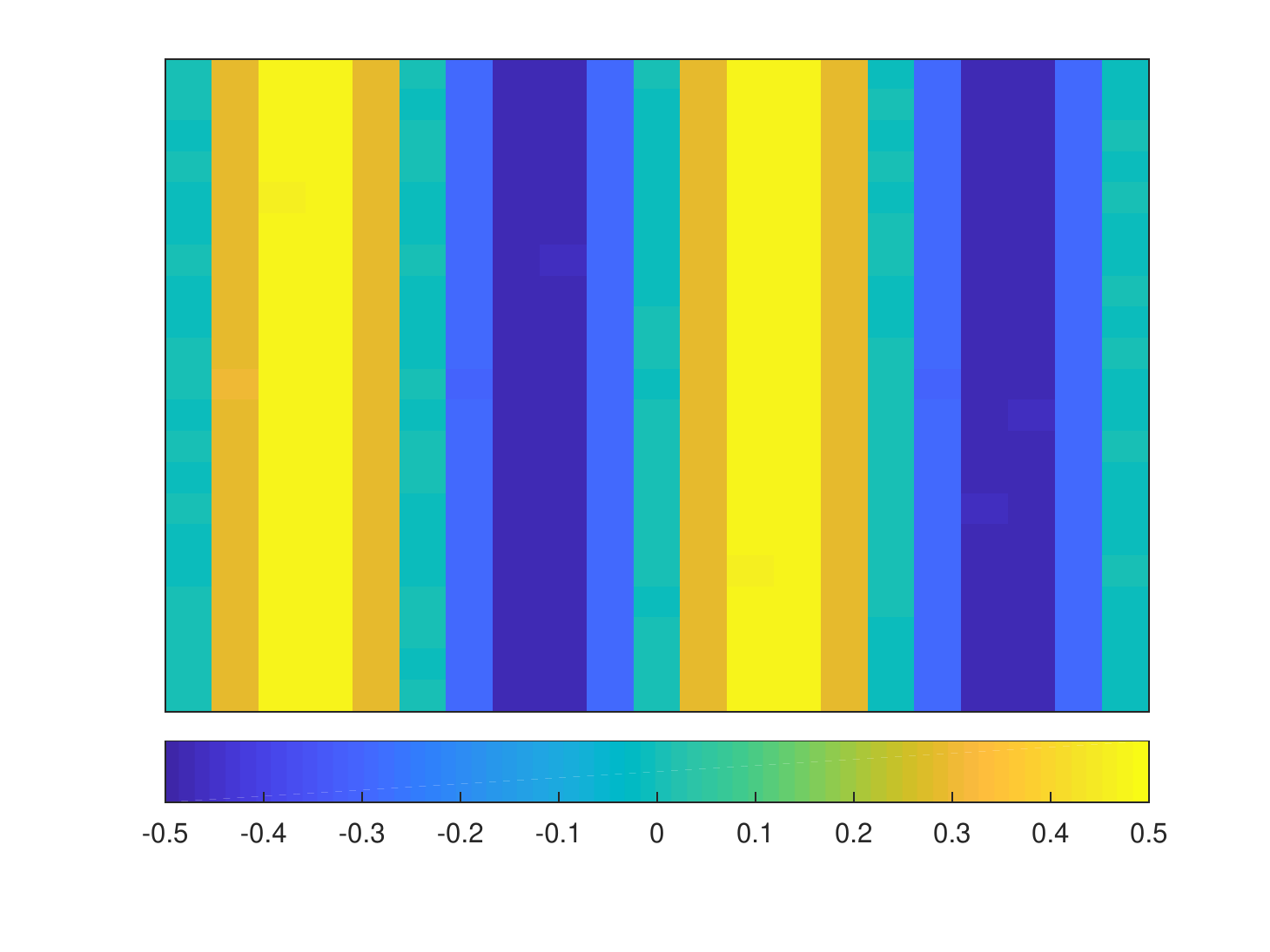}
&\includegraphics[width=\wsize,height=\hhsize,trim={1.5cm 1cm 1cm 0.5cm},clip]{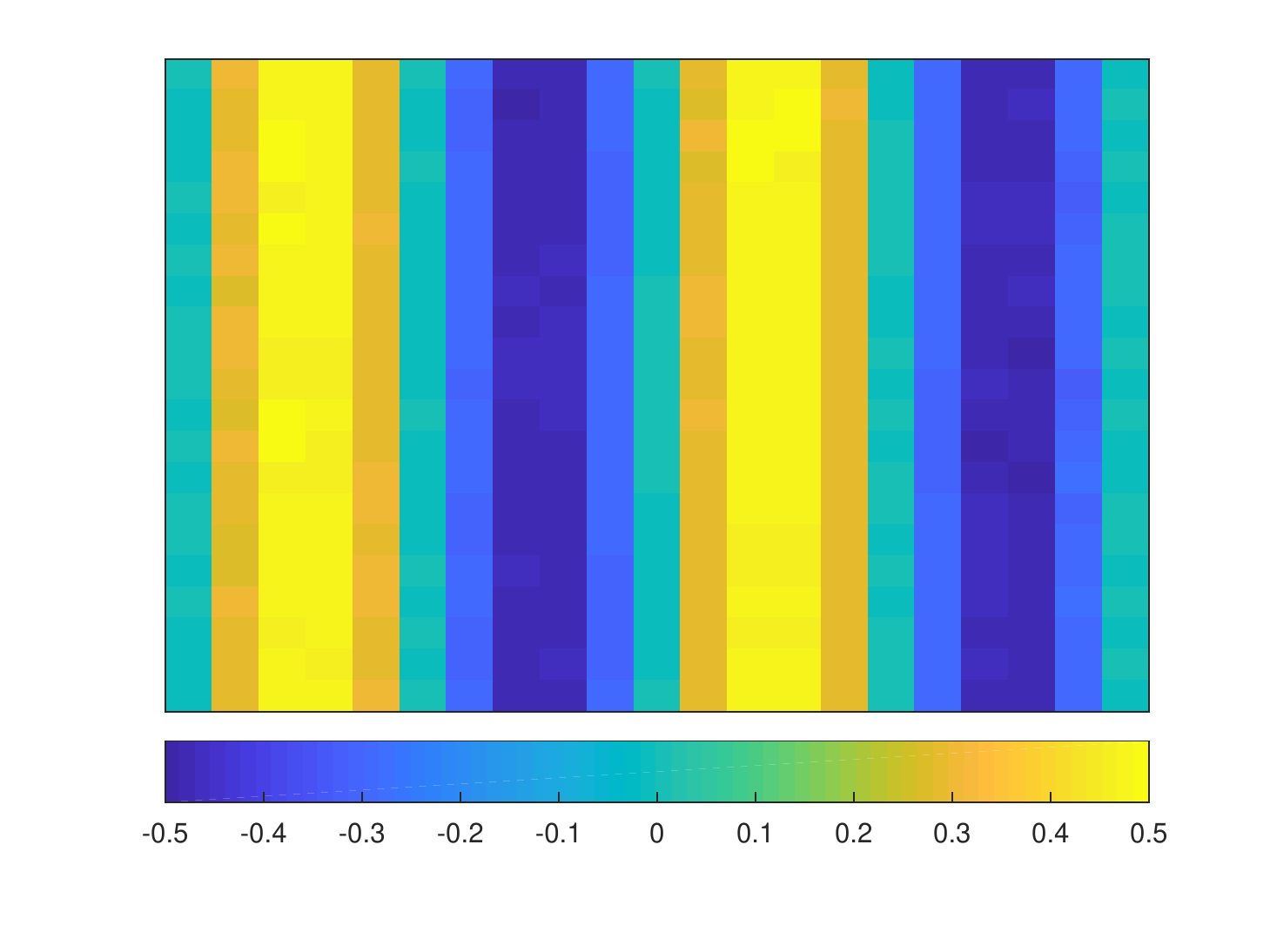}
&\includegraphics[width=\wsize,height=\hhsize,trim={1.5cm 1cm 1cm 0.5cm},clip]{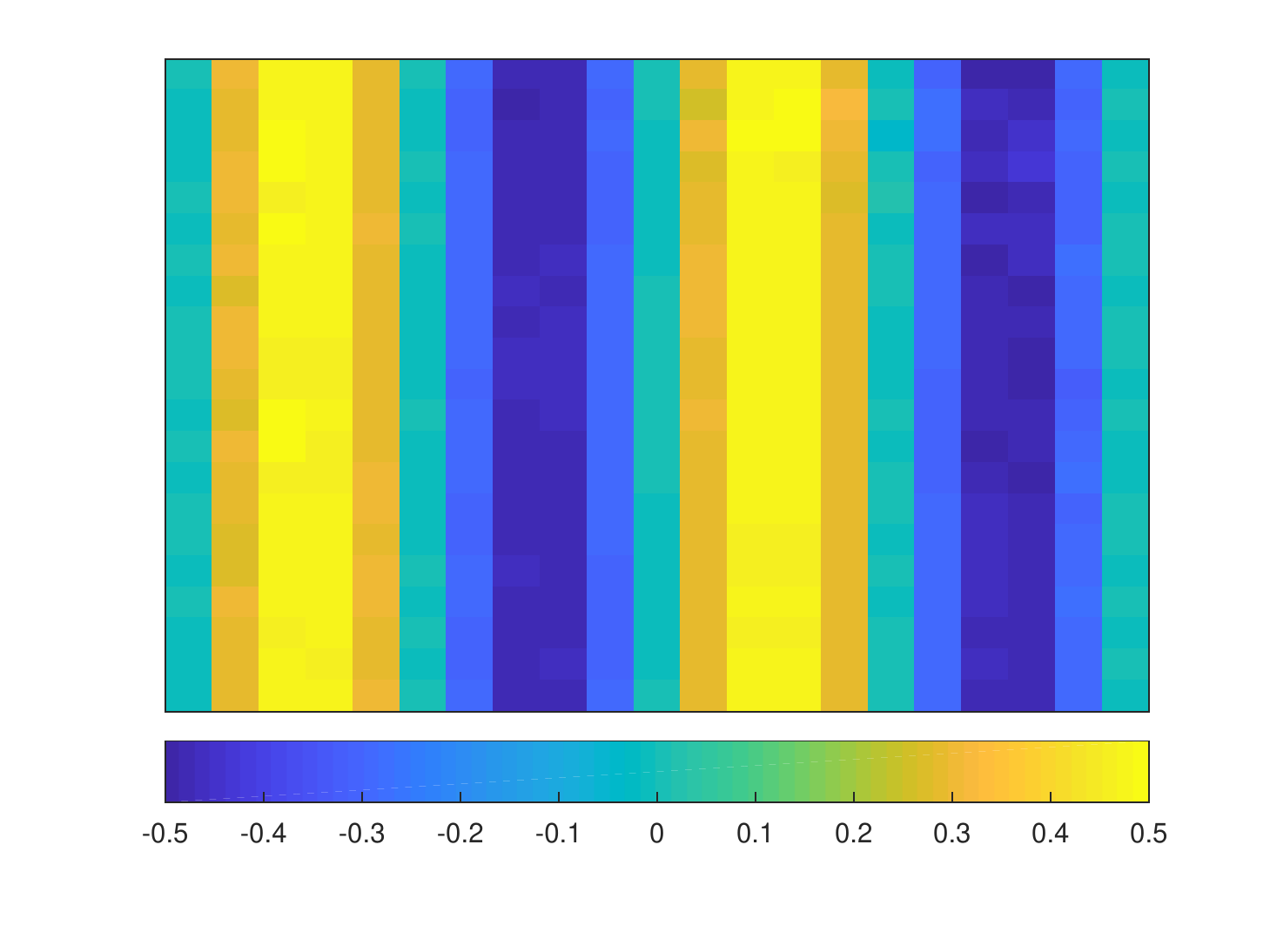}
&\includegraphics[width=\wsize,height=\hhsize,trim={1cm 0.5cm .5cm .5cm},clip]{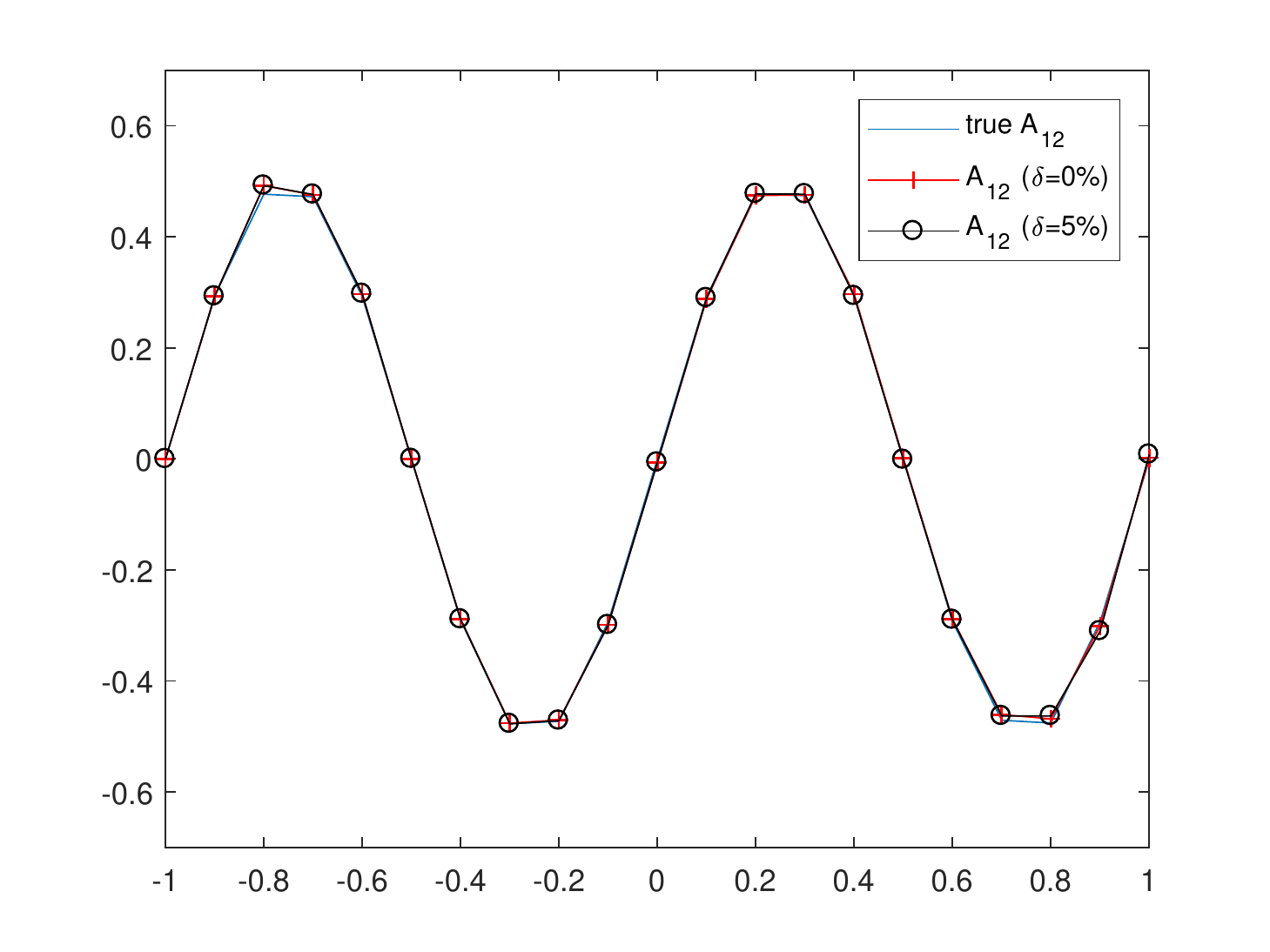}\\
 \includegraphics[width=\wsize,height=\hhsize,trim={1.5cm 1cm 1cm 0.5cm},clip]{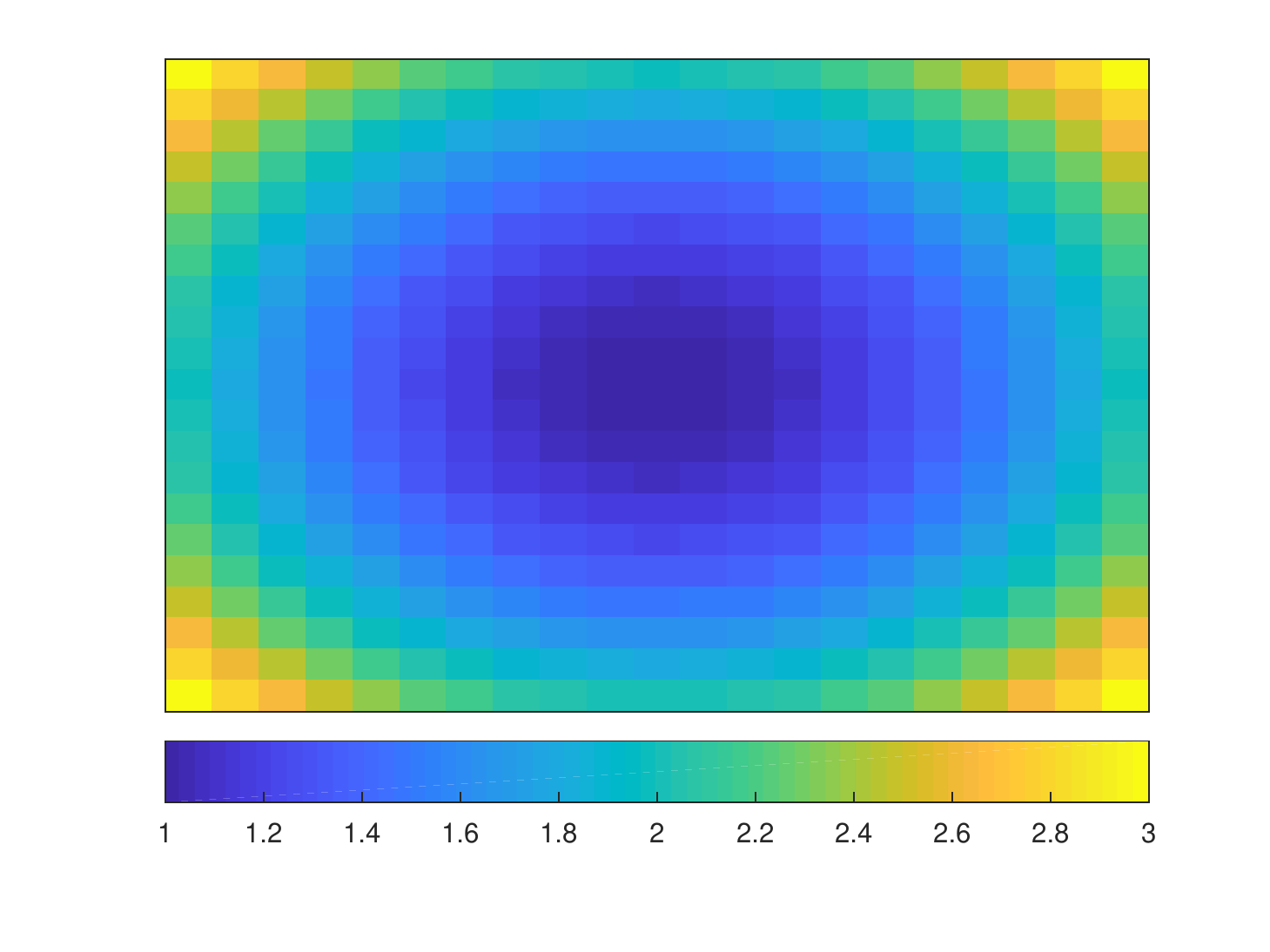}
&\includegraphics[width=\wsize,height=\hhsize,trim={1.5cm 1cm 1cm 0.5cm},clip]{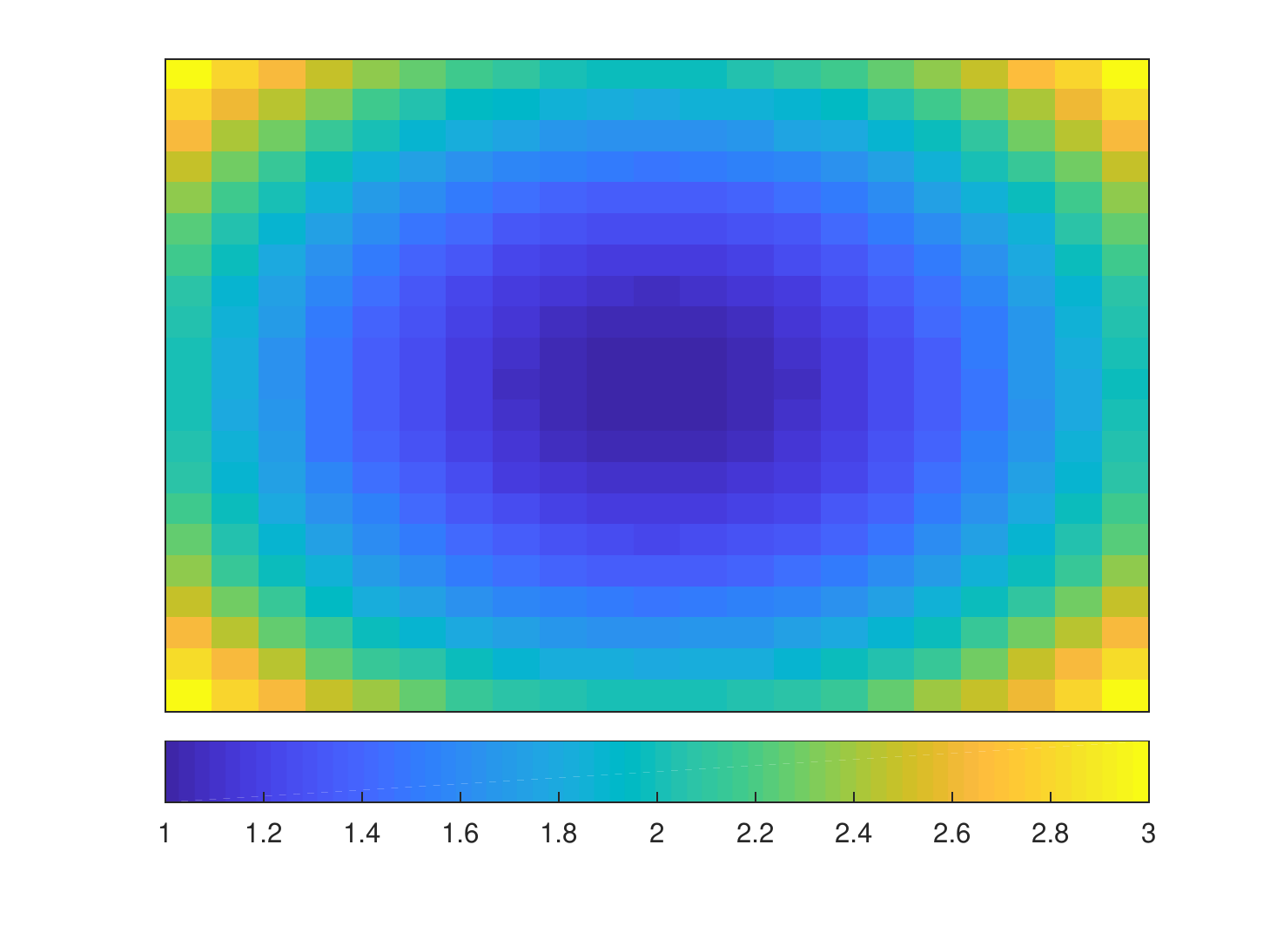}
&\includegraphics[width=\wsize,height=\hhsize,trim={1.5cm 1cm 1cm 0.5cm},clip]{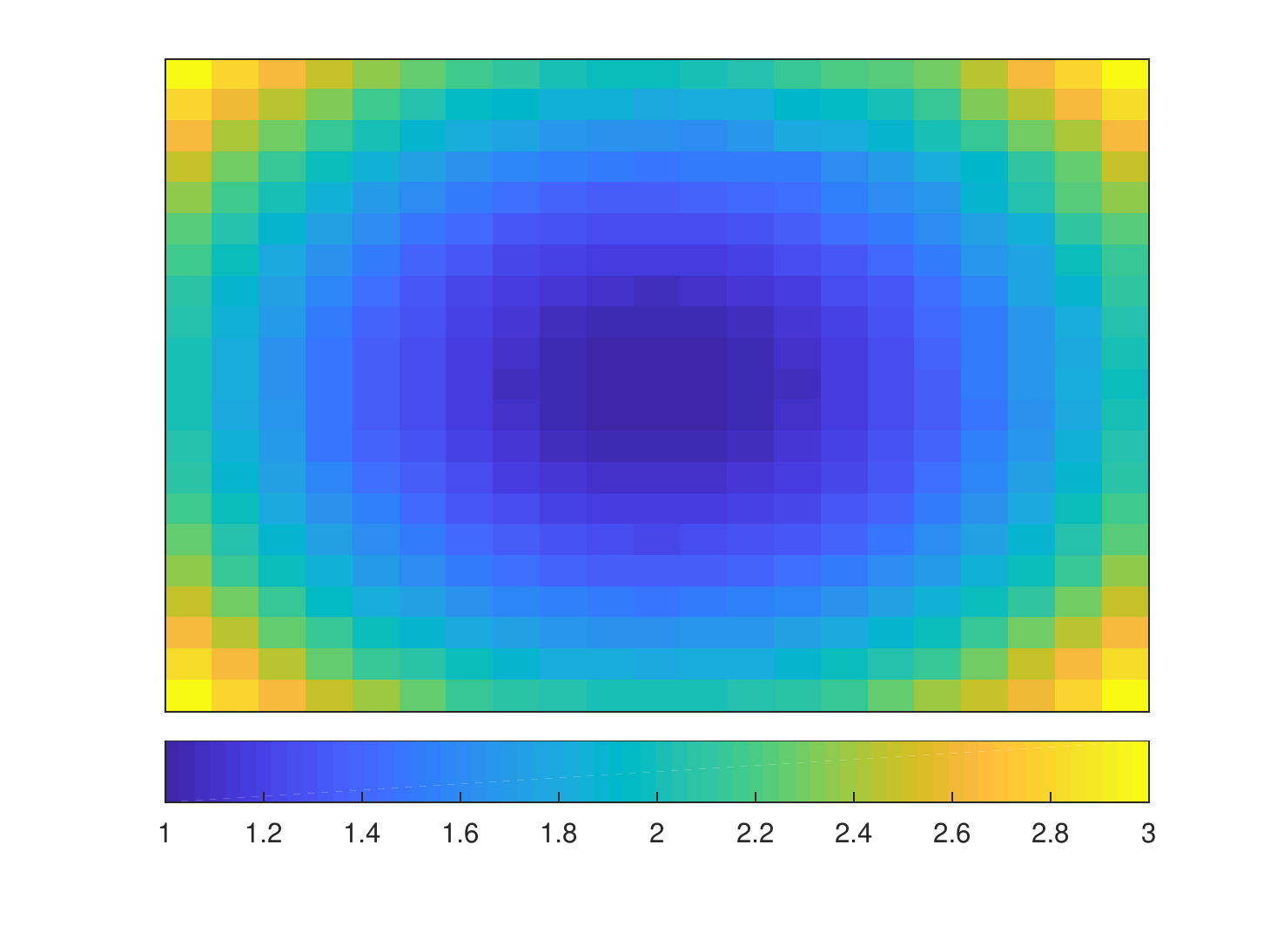}
&\includegraphics[width=\wsize,height=\hhsize,trim={1cm .5cm .5cm .5cm },clip]{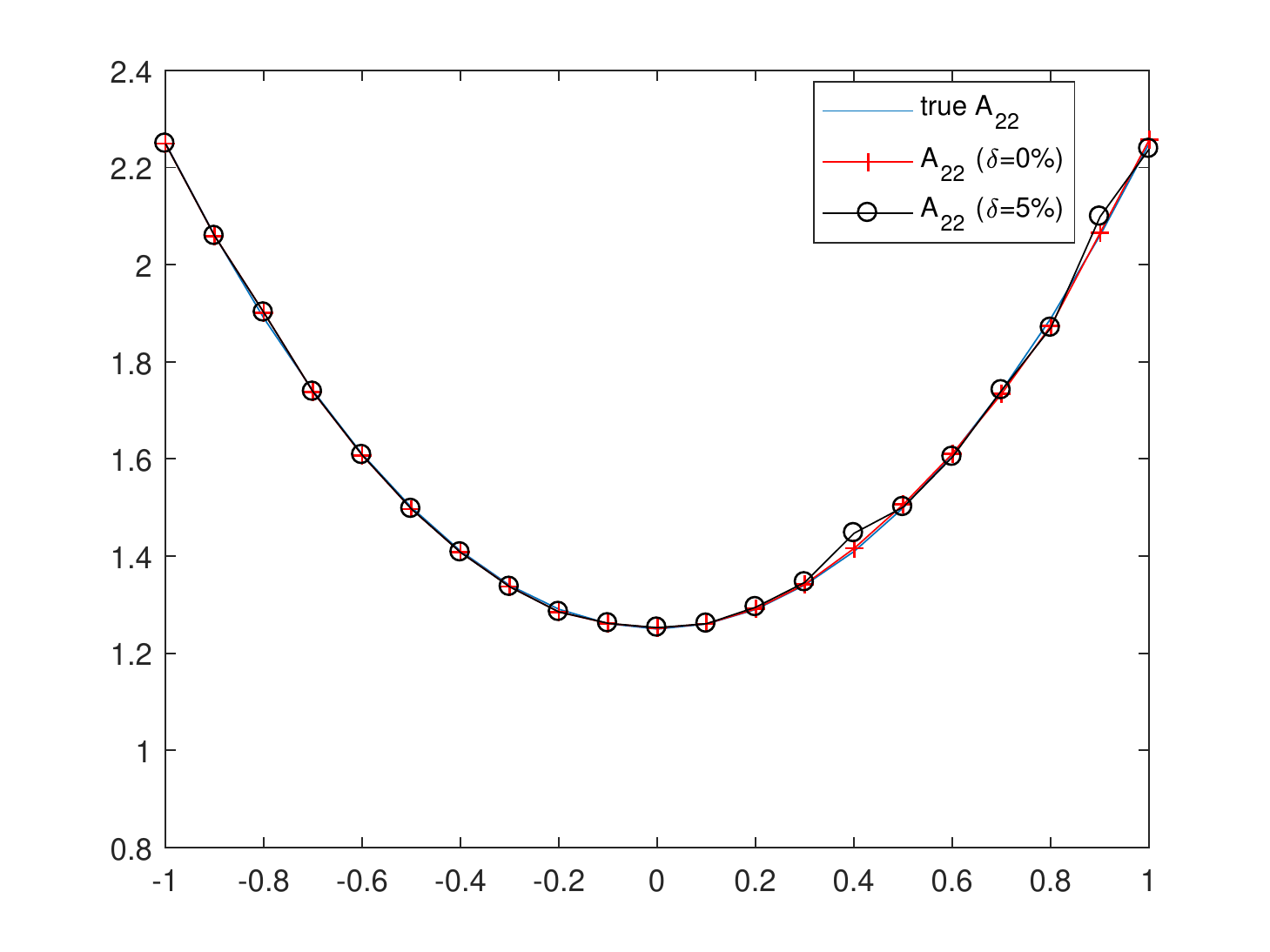}\\
\end{tabular}
\caption{The reconstructions for the conductivity tensor for Example 4: the top, middle and bottom rows refer to $A_{11}$, $A_{12}$ and $A_{22}$,
respectively. Column 1 is for exact conductivity, Column 2 for reconstruction with exact data, Column 3 for noisy data with $\delta=5\%$, and column 4 for the cross section along \{$x_2=-0.5$\}.}
\label{exam4}
\end{figure}

\begin{figure}[!htbp]
\centering
\setlength{\tabcolsep}{0pt}
\begin{tabular}{cccc}
 \includegraphics[width=\wsize,height=\hhsize,trim={1.5cm 1cm 1cm 0.5cm},clip]{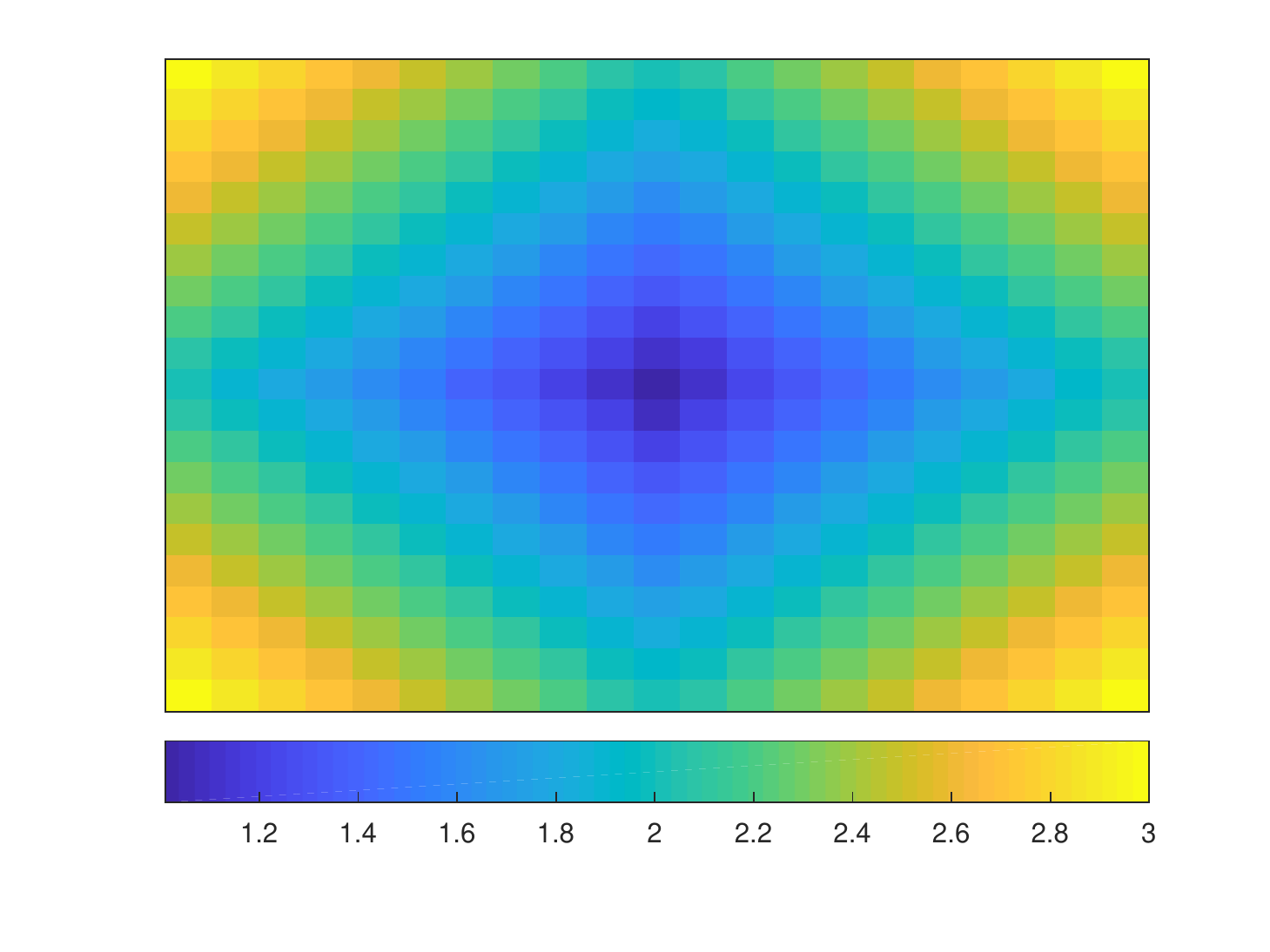}
&\includegraphics[width=\wsize,height=\hhsize,trim={1.5cm 1cm 1cm 0.5cm},clip]{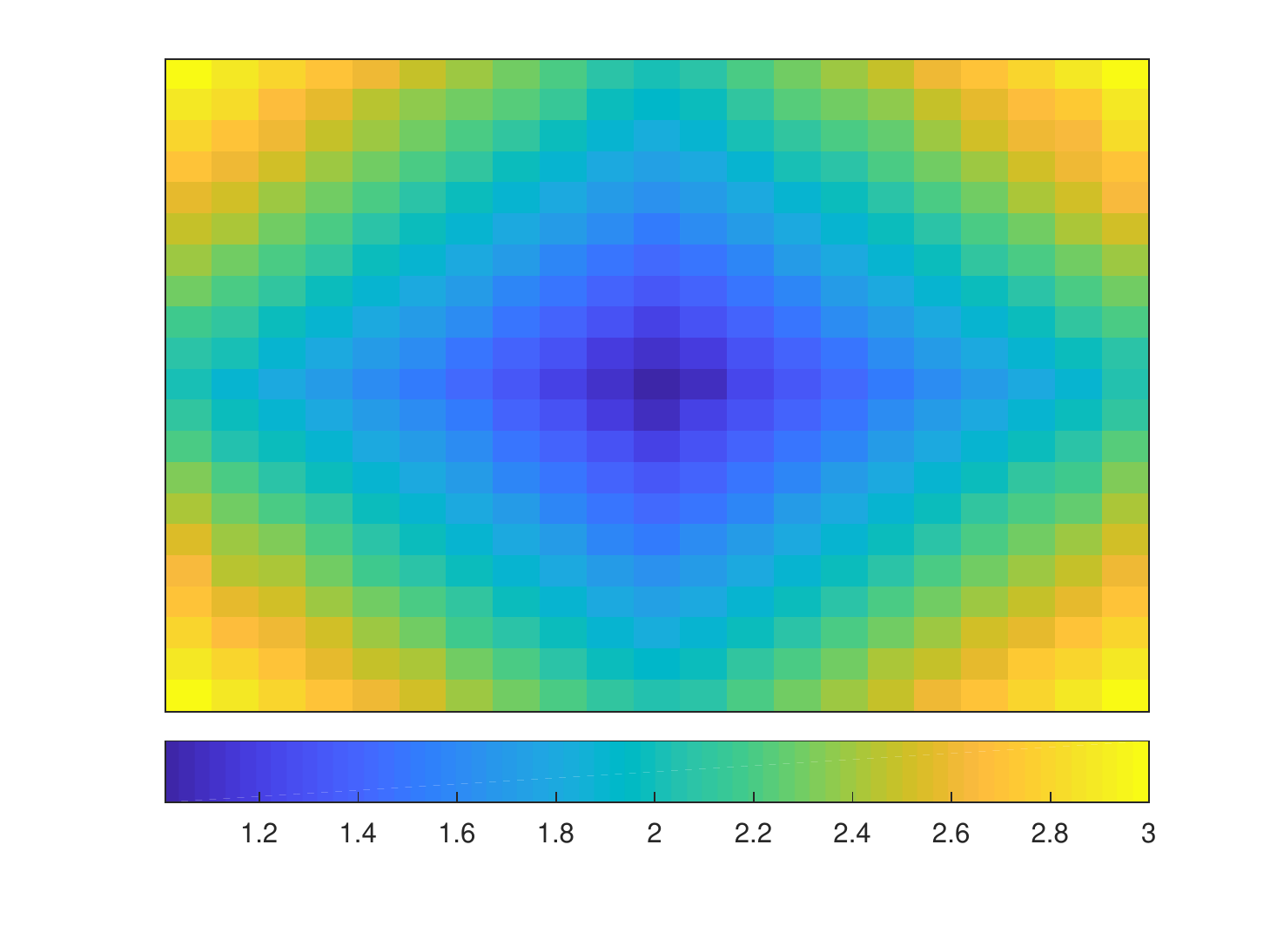}
&\includegraphics[width=\wsize,height=\hhsize,trim={1.5cm 1cm 1cm 0.5cm},clip]{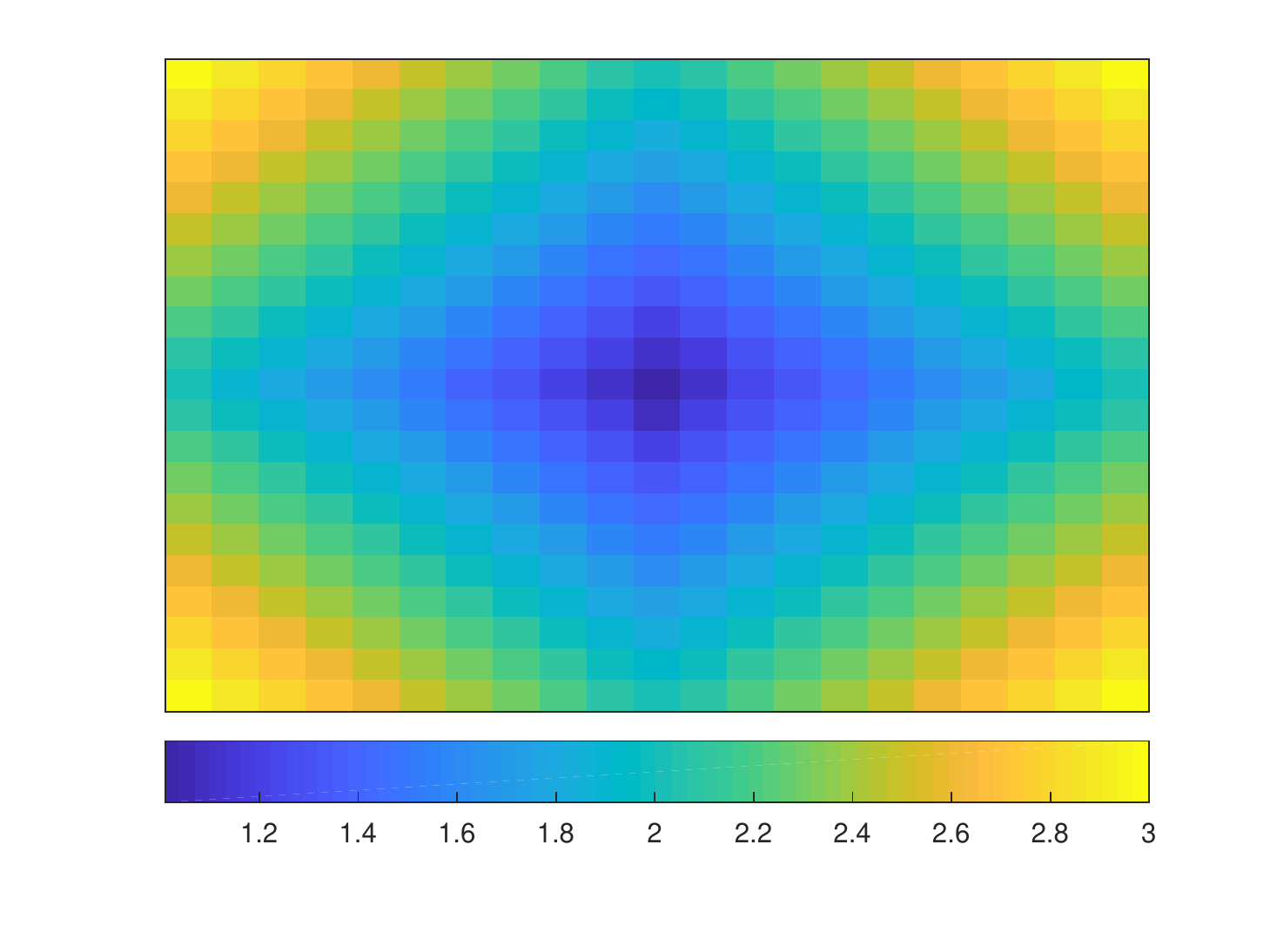}
&\includegraphics[width=\wsize,height=\hhsize,trim={1cm .5cm .5cm .5cm },clip]{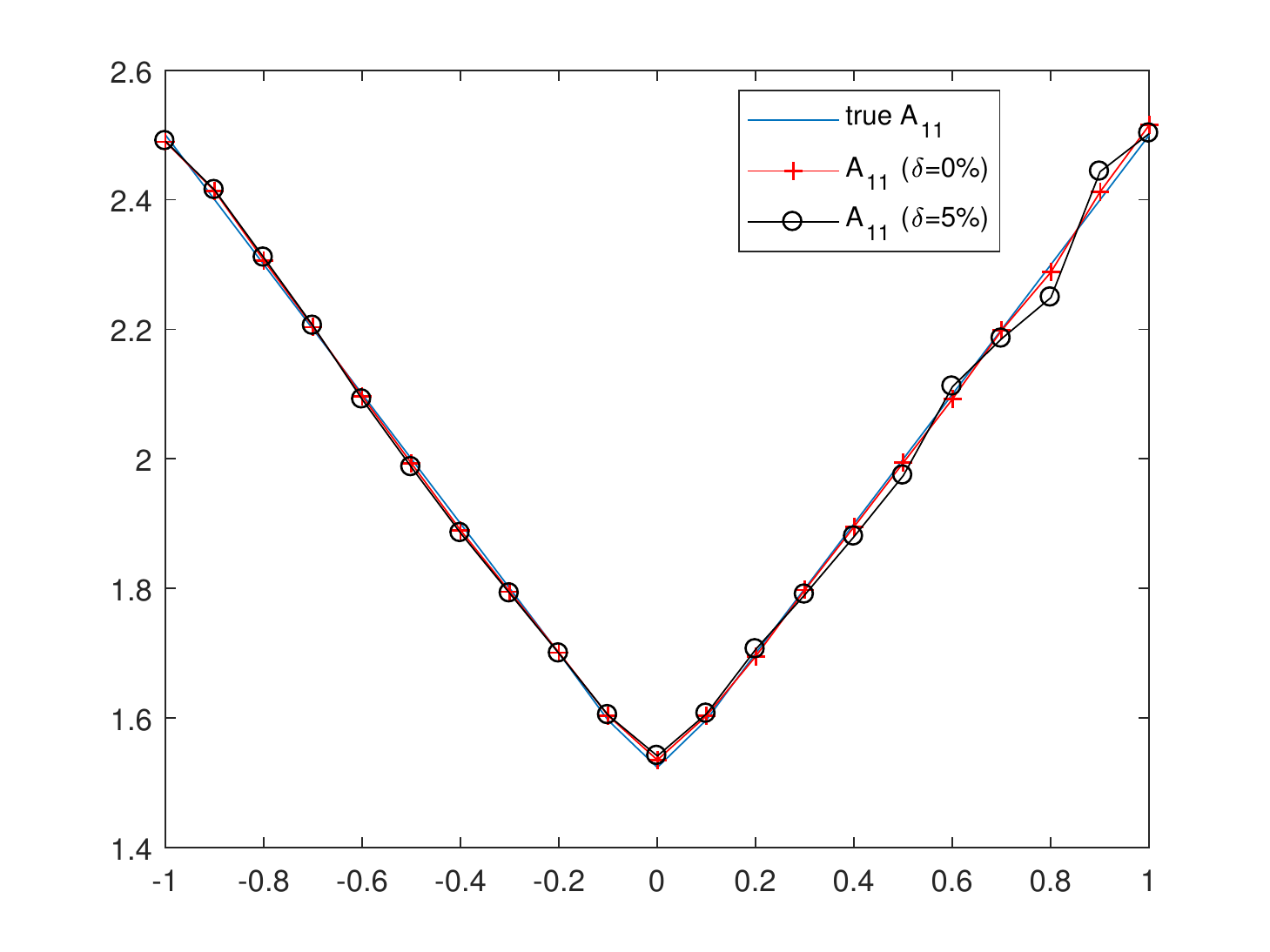}\\
 \includegraphics[width=\wsize,height=\hhsize,trim={1.5cm 1cm 1cm 0.5cm},clip]{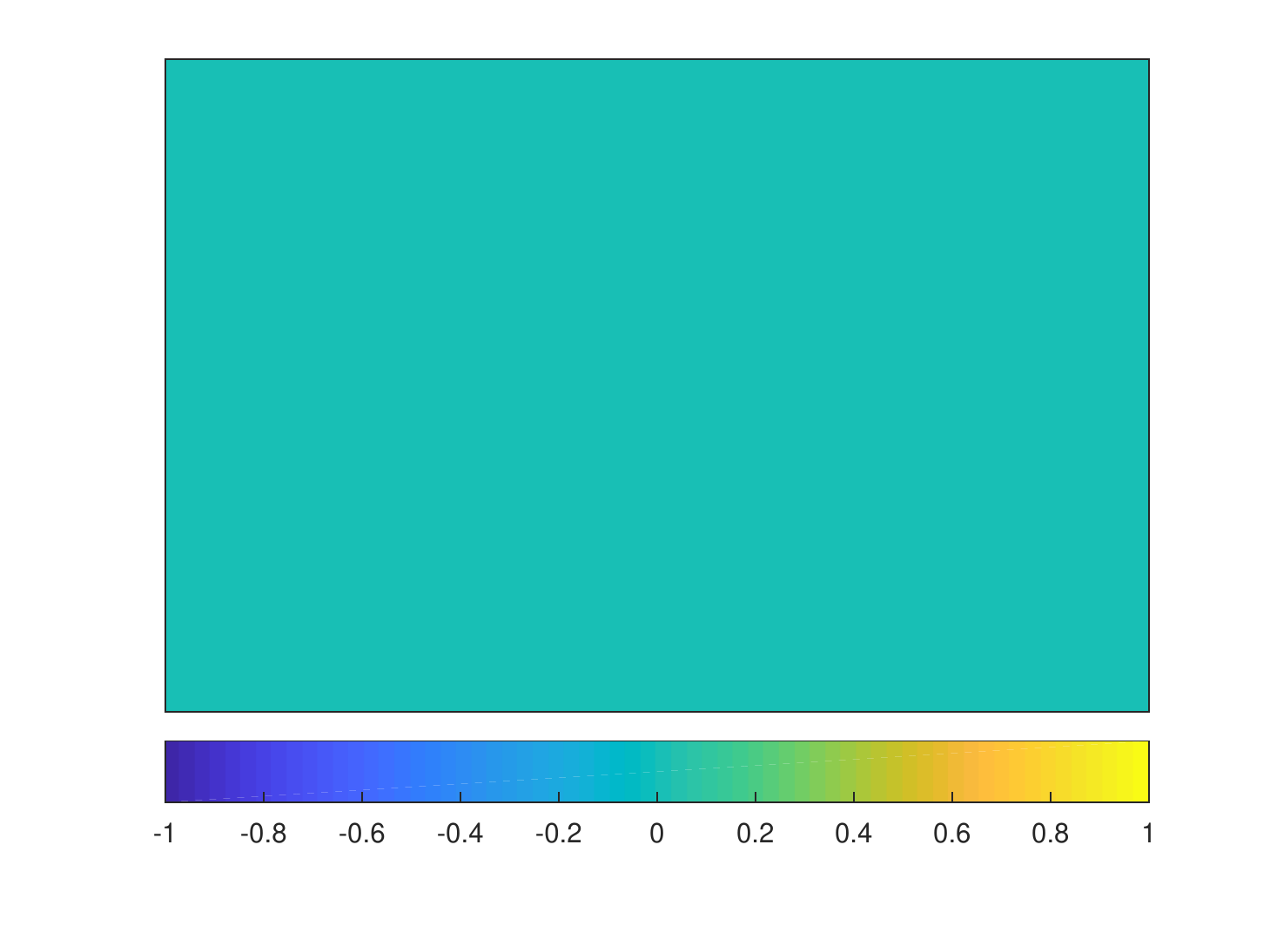}
&\includegraphics[width=\wsize,height=\hhsize,trim={1.5cm 1cm 1cm 0.5cm},clip]{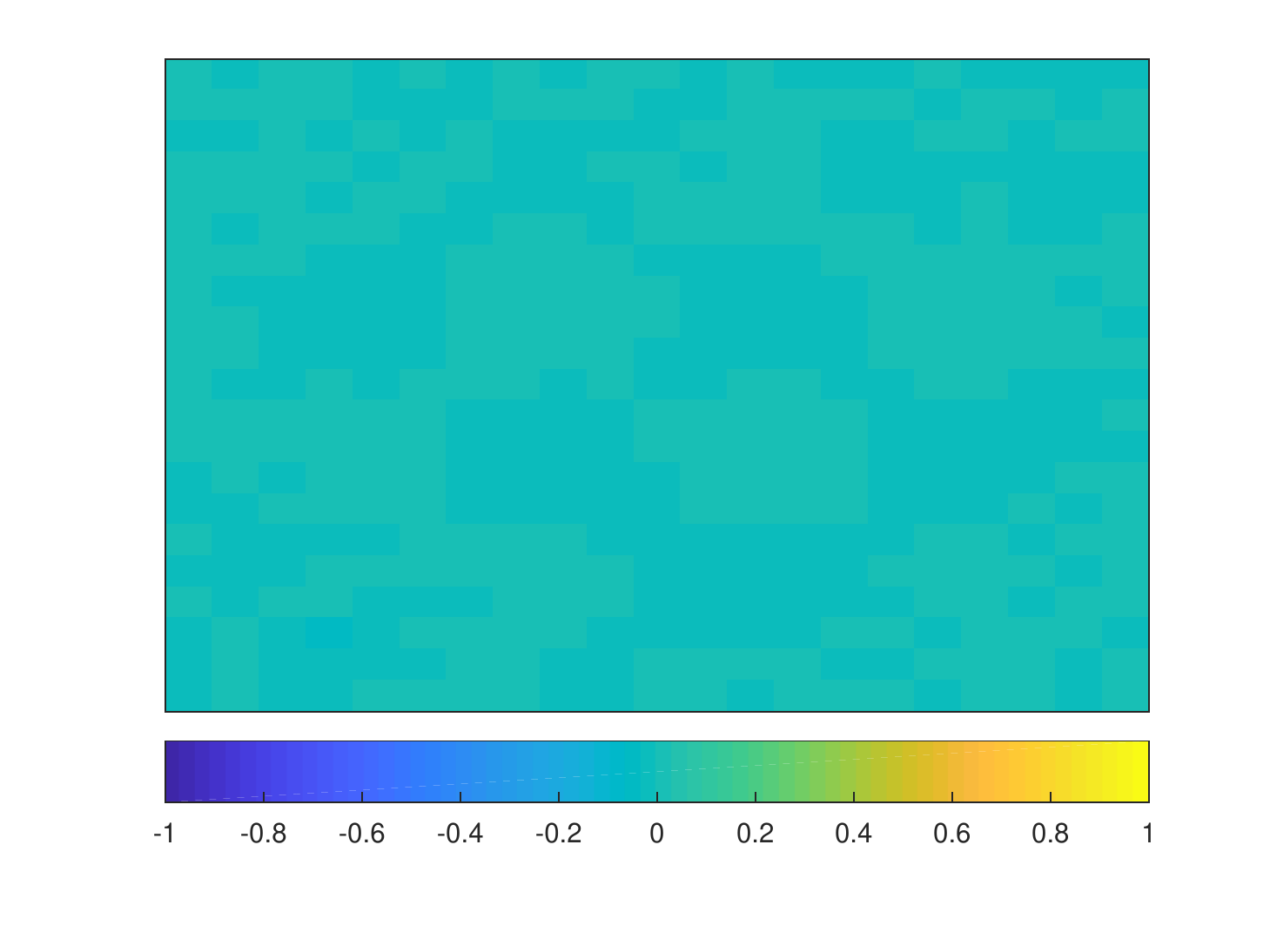}
&\includegraphics[width=\wsize,height=\hhsize,trim={1.5cm 1cm 1cm 0.5cm},clip]{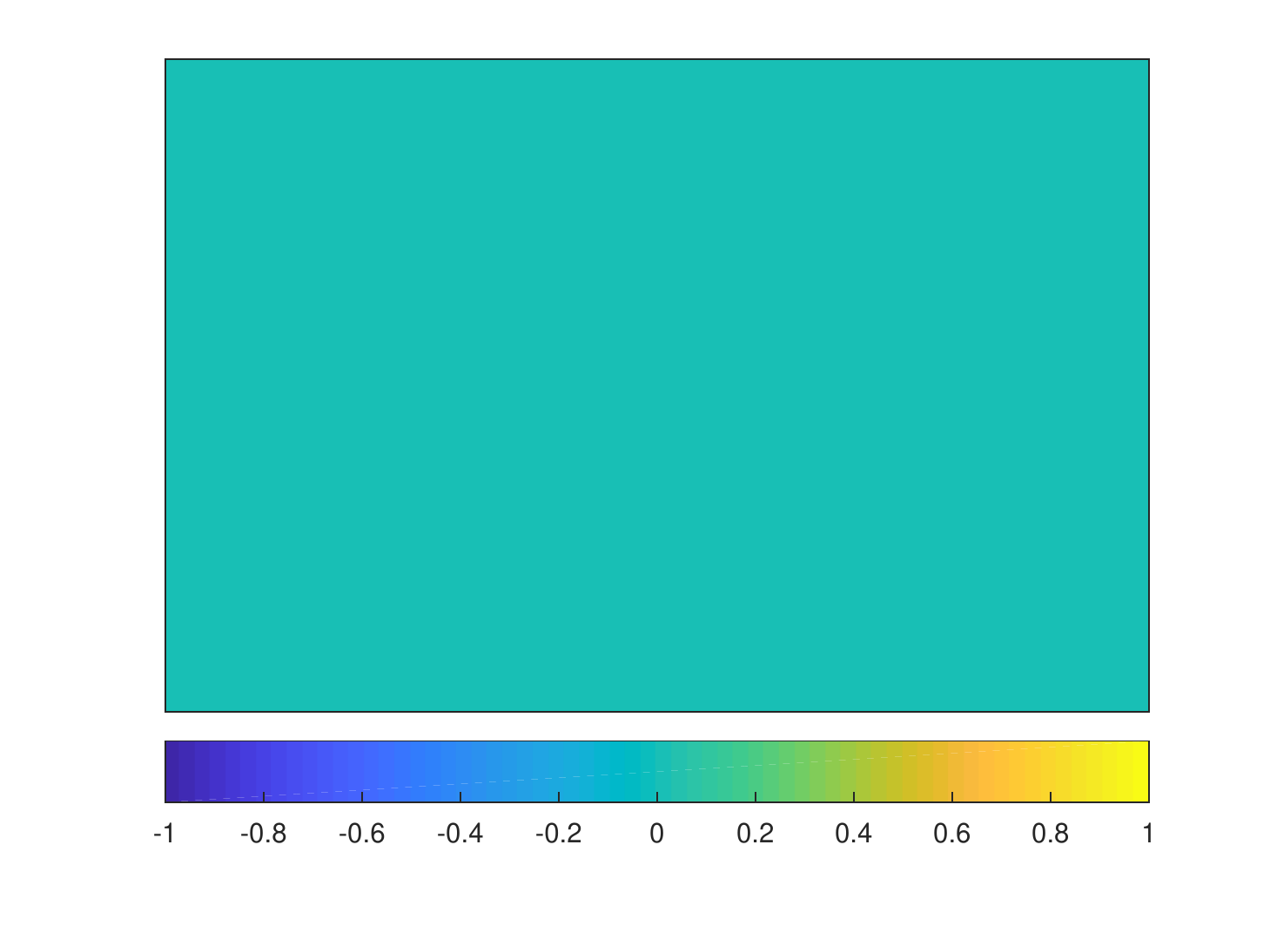}
&\includegraphics[width=\wsize,height=\hhsize,trim={1cm 0.5cm .5cm .5cm},clip]{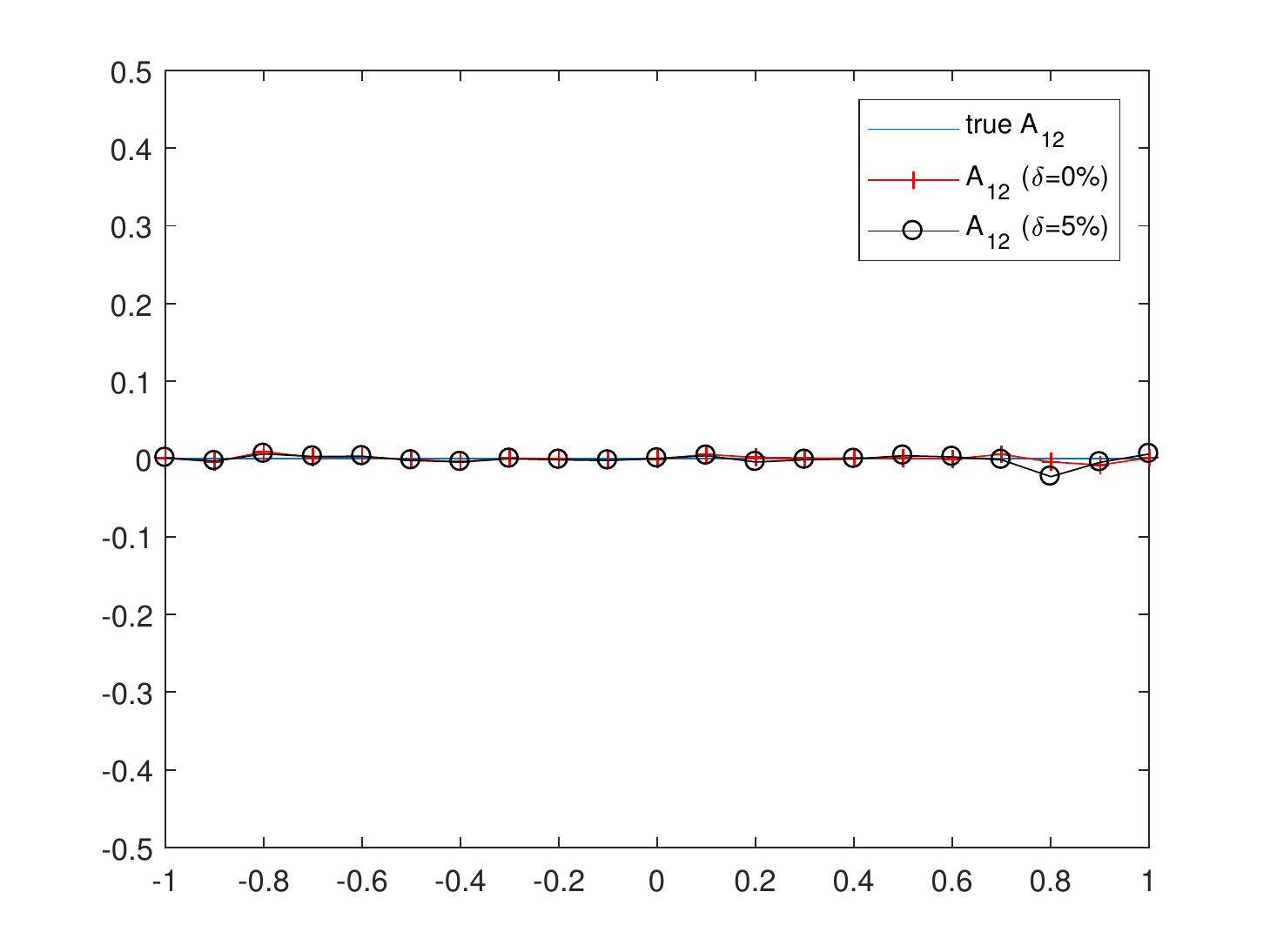}\\
 \includegraphics[width=\wsize,height=\hhsize,trim={1.5cm 1cm 1cm 0.5cm},clip]{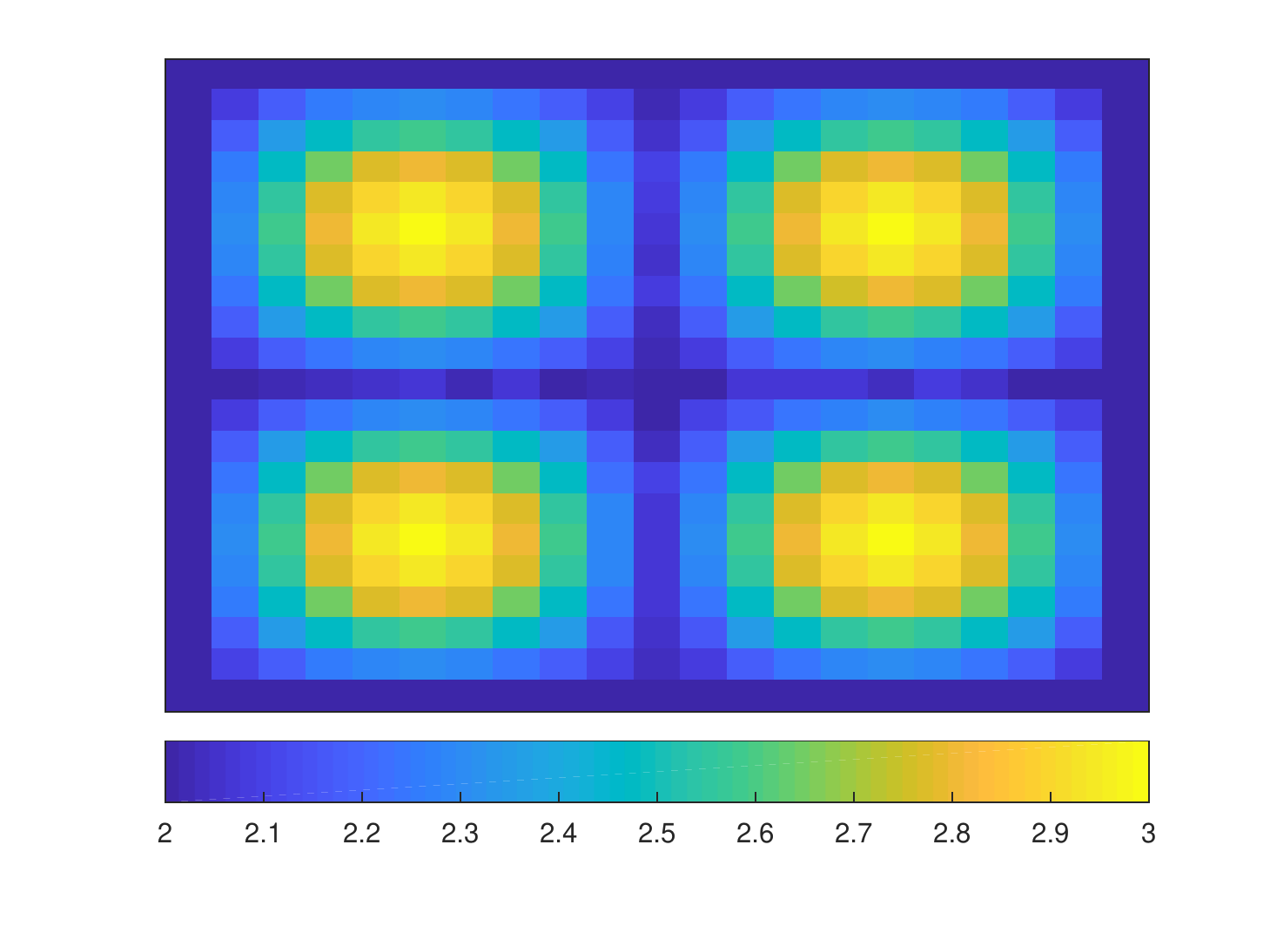}
&\includegraphics[width=\wsize,height=\hhsize,trim={1.5cm 1cm 1cm 0.5cm},clip]{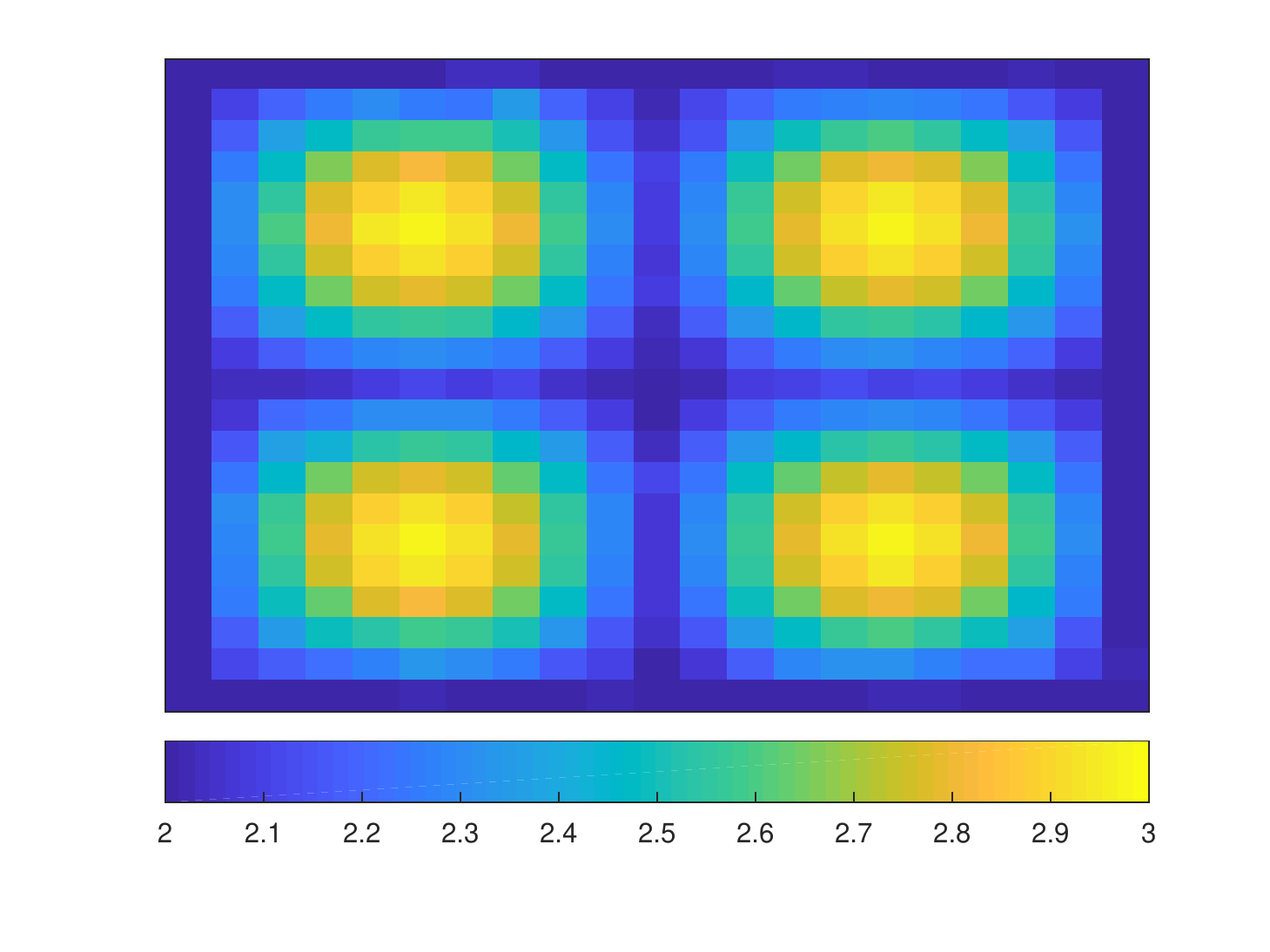}
&\includegraphics[width=\wsize,height=\hhsize,trim={1.5cm 1cm 1cm 0.5cm},clip]{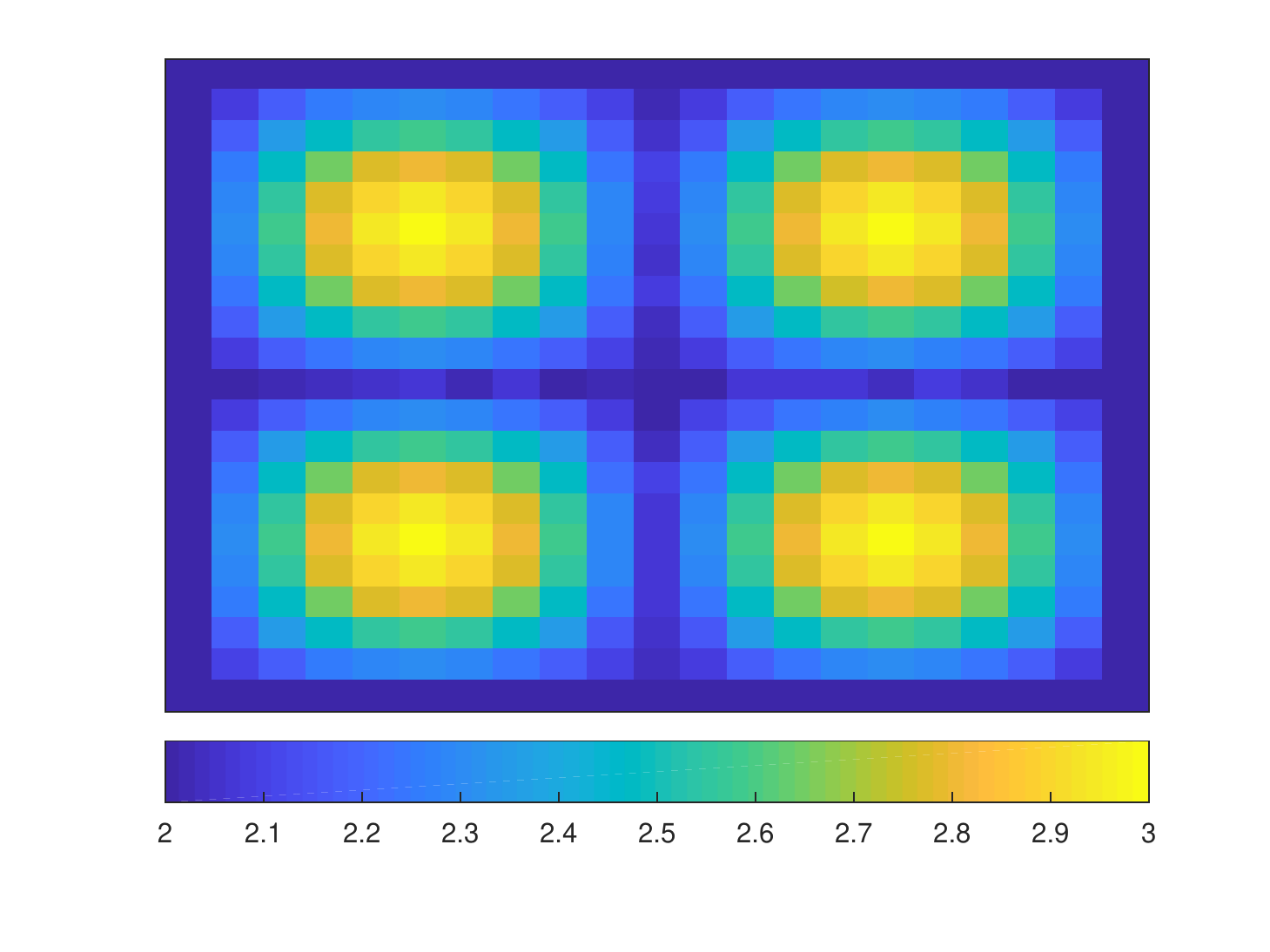}
&\includegraphics[width=\wsize,height=\hhsize,trim={1cm .5cm .5cm .5cm },clip]{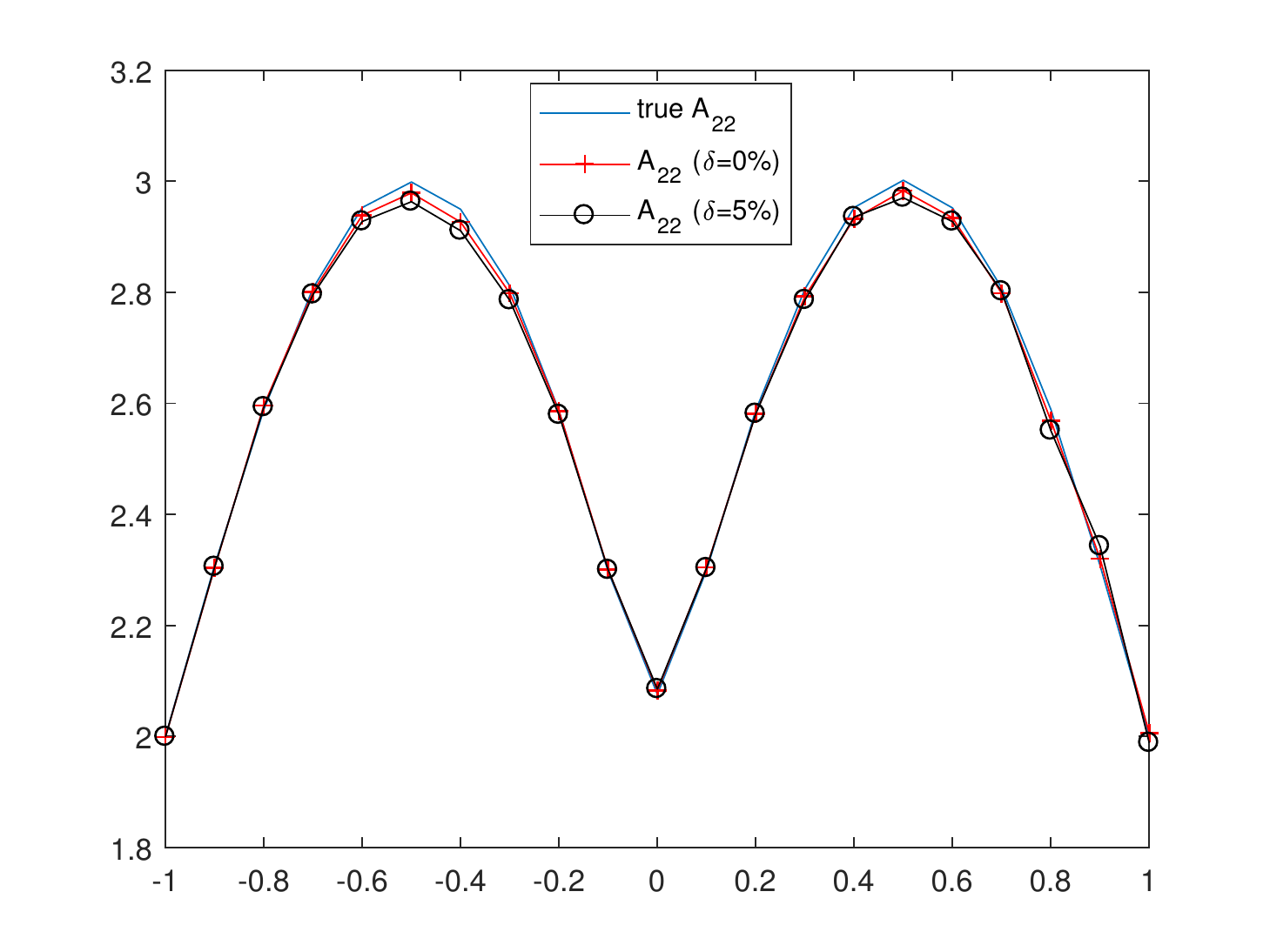}\\
\end{tabular}
\caption{The reconstructions for the conductivity tensor for Example 5: the top, middle and bottom rows refer to $A_{11}$, $A_{12}$ and $A_{22}$,
respectively. Column 1 is for exact conductivity, Column 2 for reconstruction with exact data, Column 3 for noisy data with $\delta=5\%$, and column 4 for the cross section along \{$x_2=-0.5$\}.}
\label{exam5}
\end{figure}

\begin{figure}[!htbp]
\centering
\setlength{\tabcolsep}{0pt}
\begin{tabular}{cccc}
 \includegraphics[width=\wsize,height=\hhsize,trim={1.5cm 1cm 1cm 0.5cm},clip]{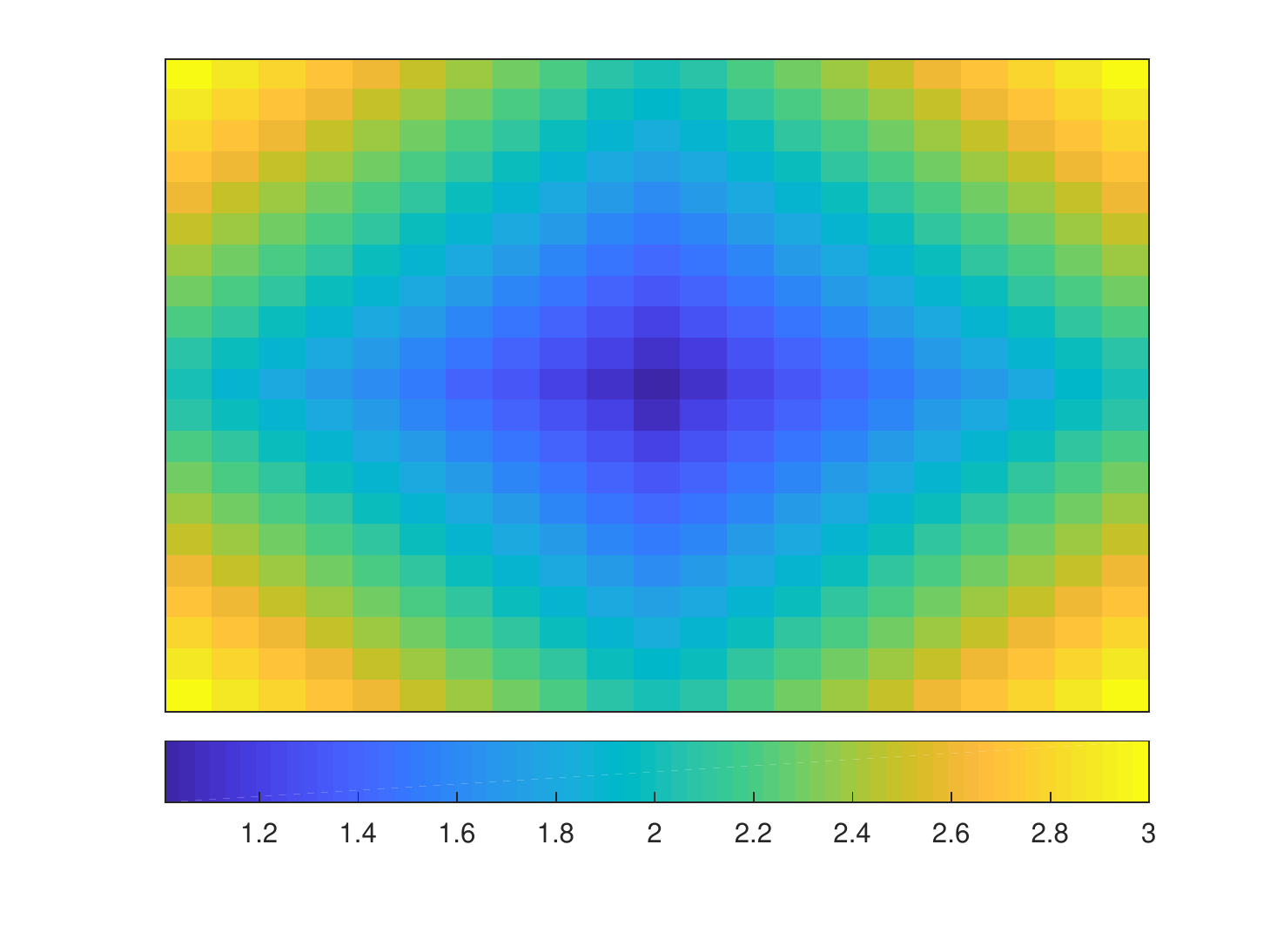}
&\includegraphics[width=\wsize,height=\hhsize,trim={1.5cm 1cm 1cm 0.5cm},clip]{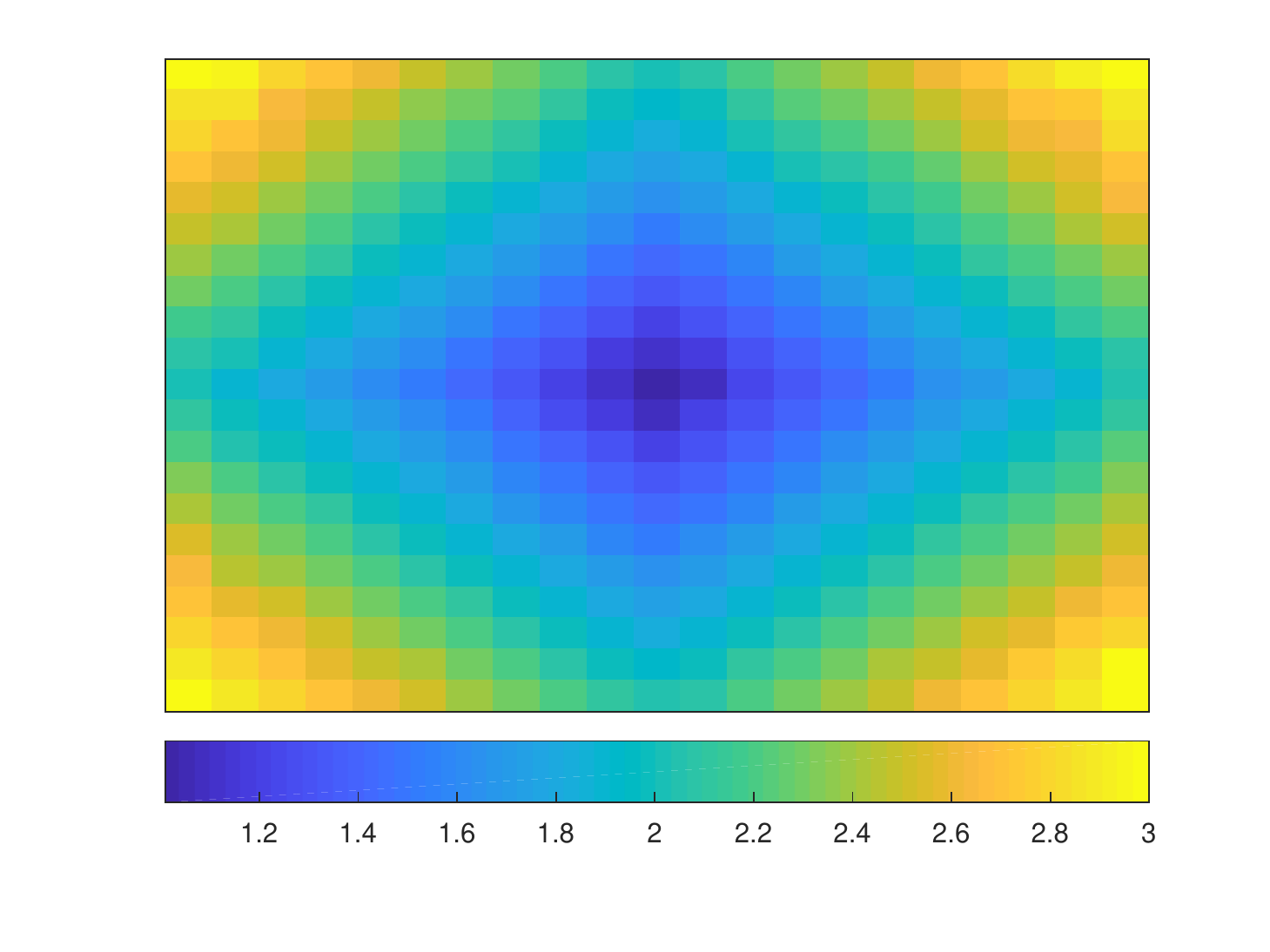}
&\includegraphics[width=\wsize,height=\hhsize,trim={1.5cm 1cm 1cm 0.5cm},clip]{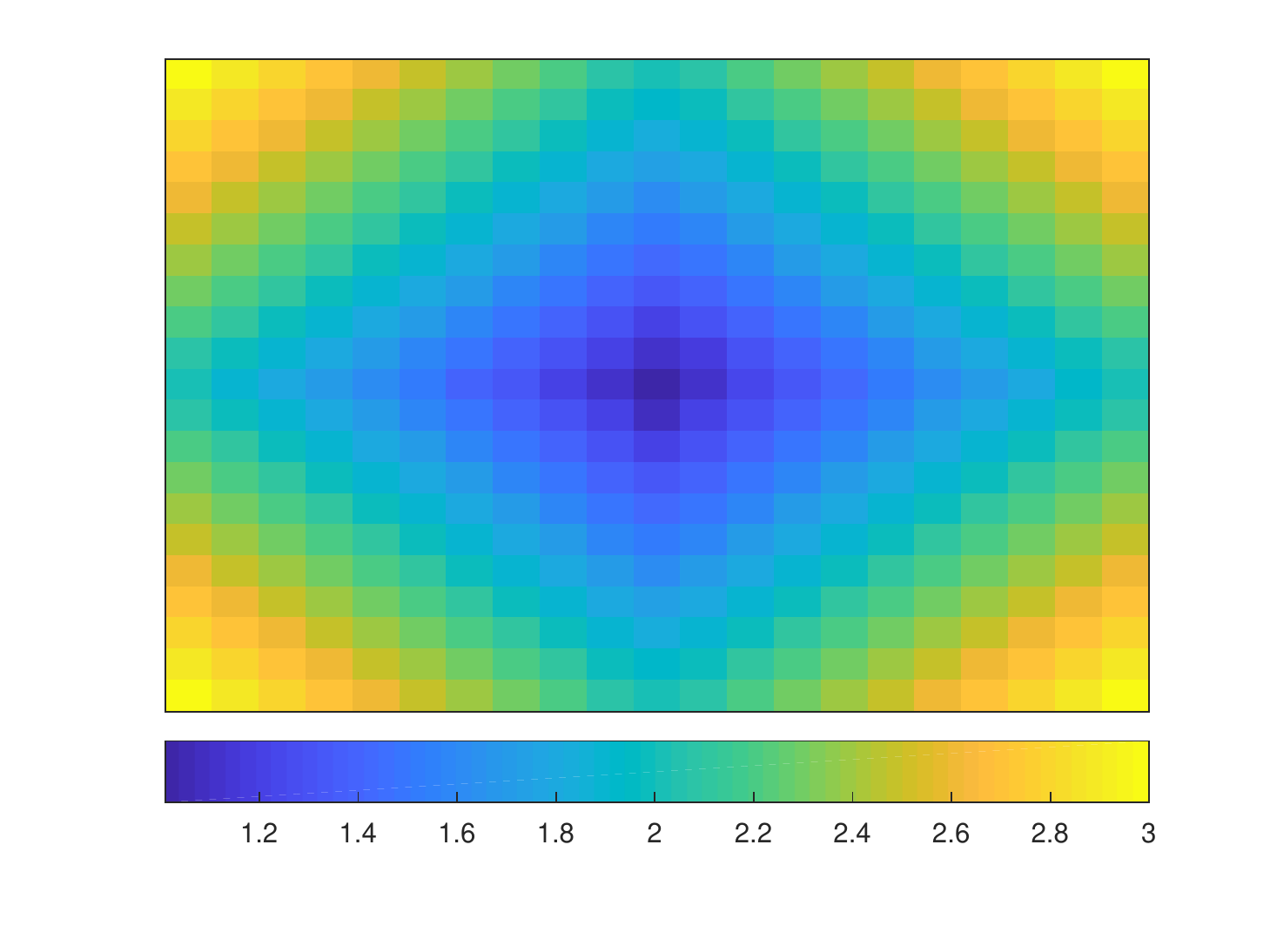}
&\includegraphics[width=\wsize,height=\hhsize,trim={1cm .5cm .5cm .5cm },clip]{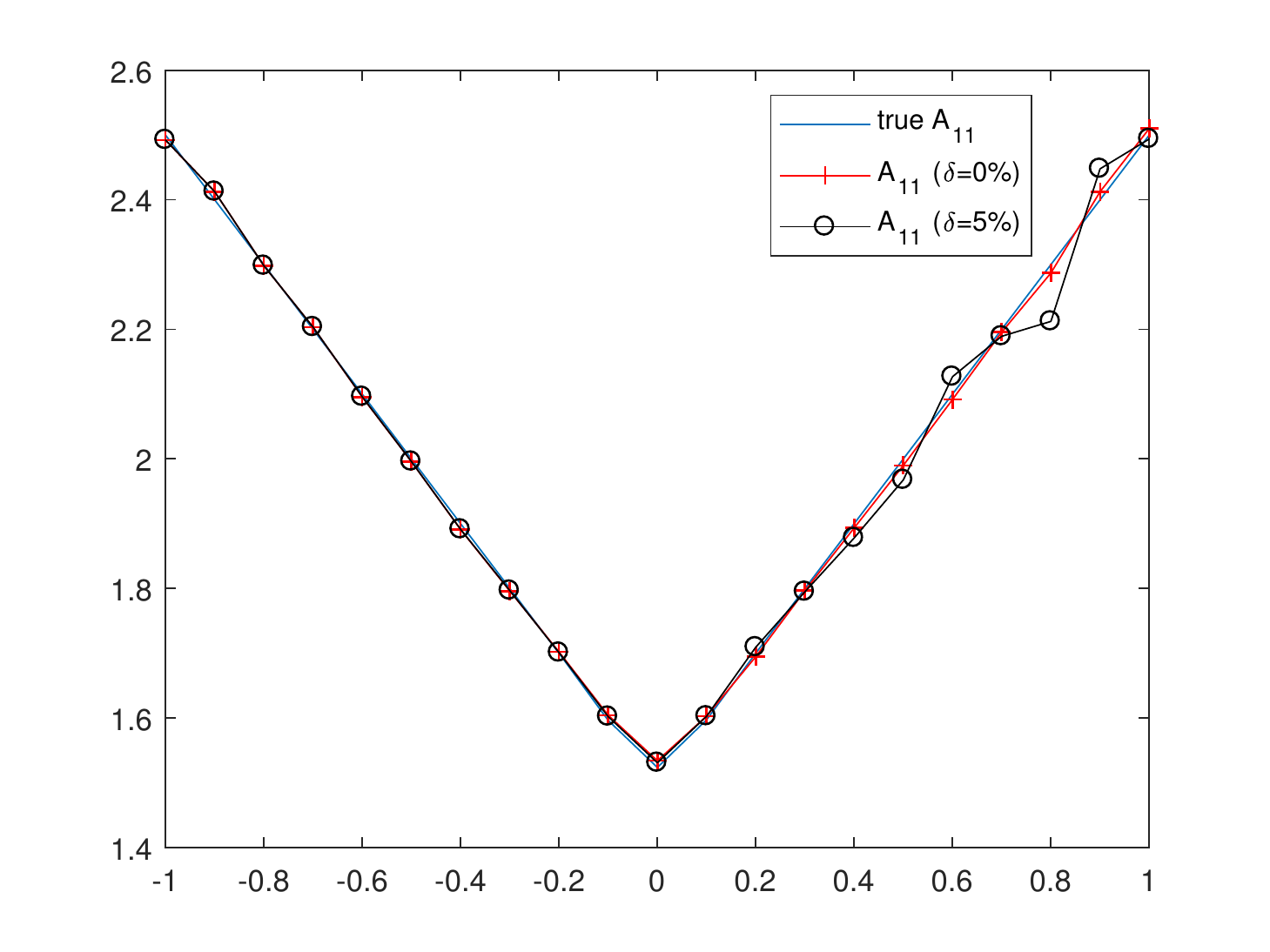}\\
 \includegraphics[width=\wsize,height=\hhsize,trim={1.5cm 1cm 1cm 0.5cm},clip]{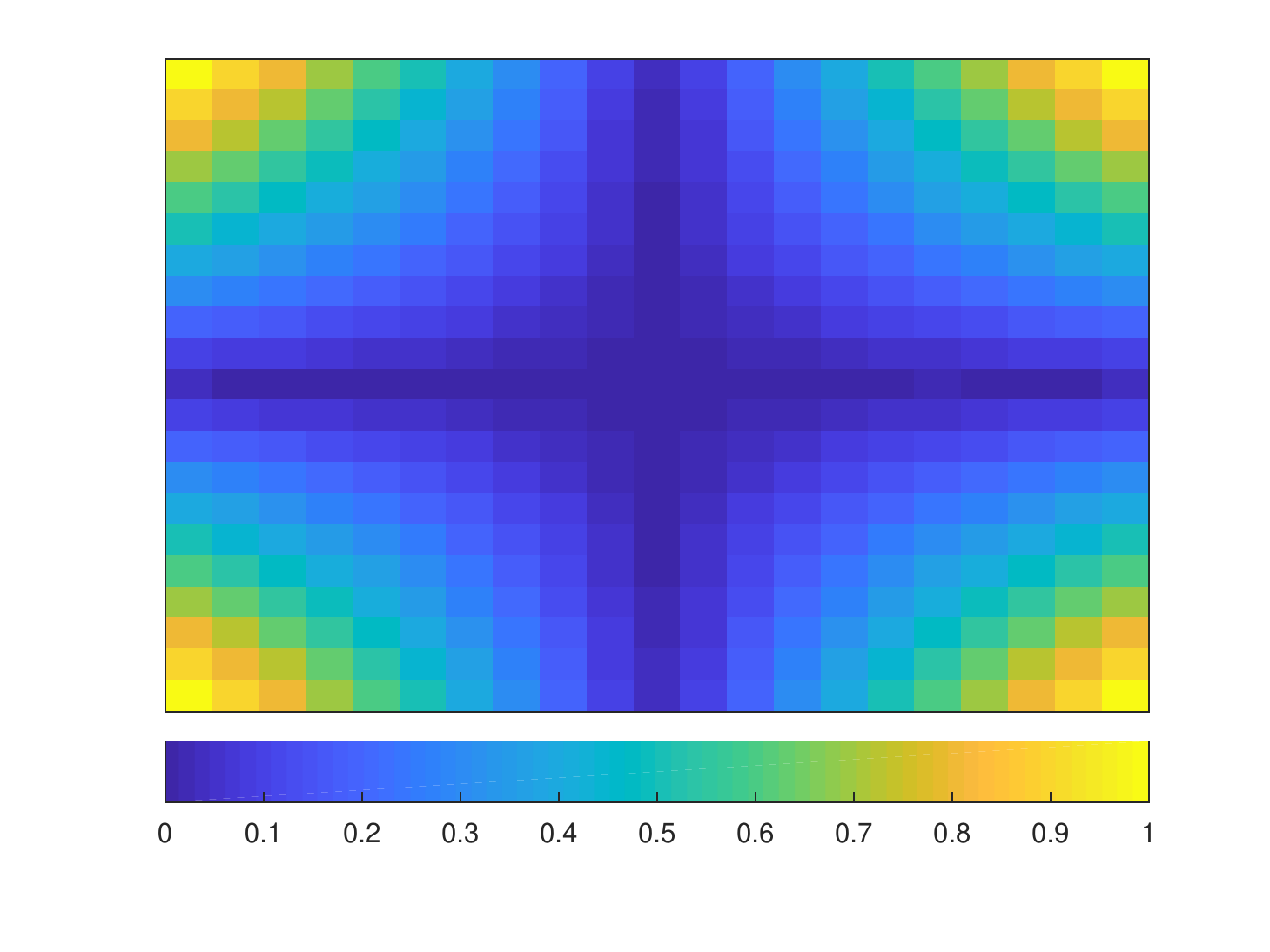}
&\includegraphics[width=\wsize,height=\hhsize,trim={1.5cm 1cm 1cm 0.5cm},clip]{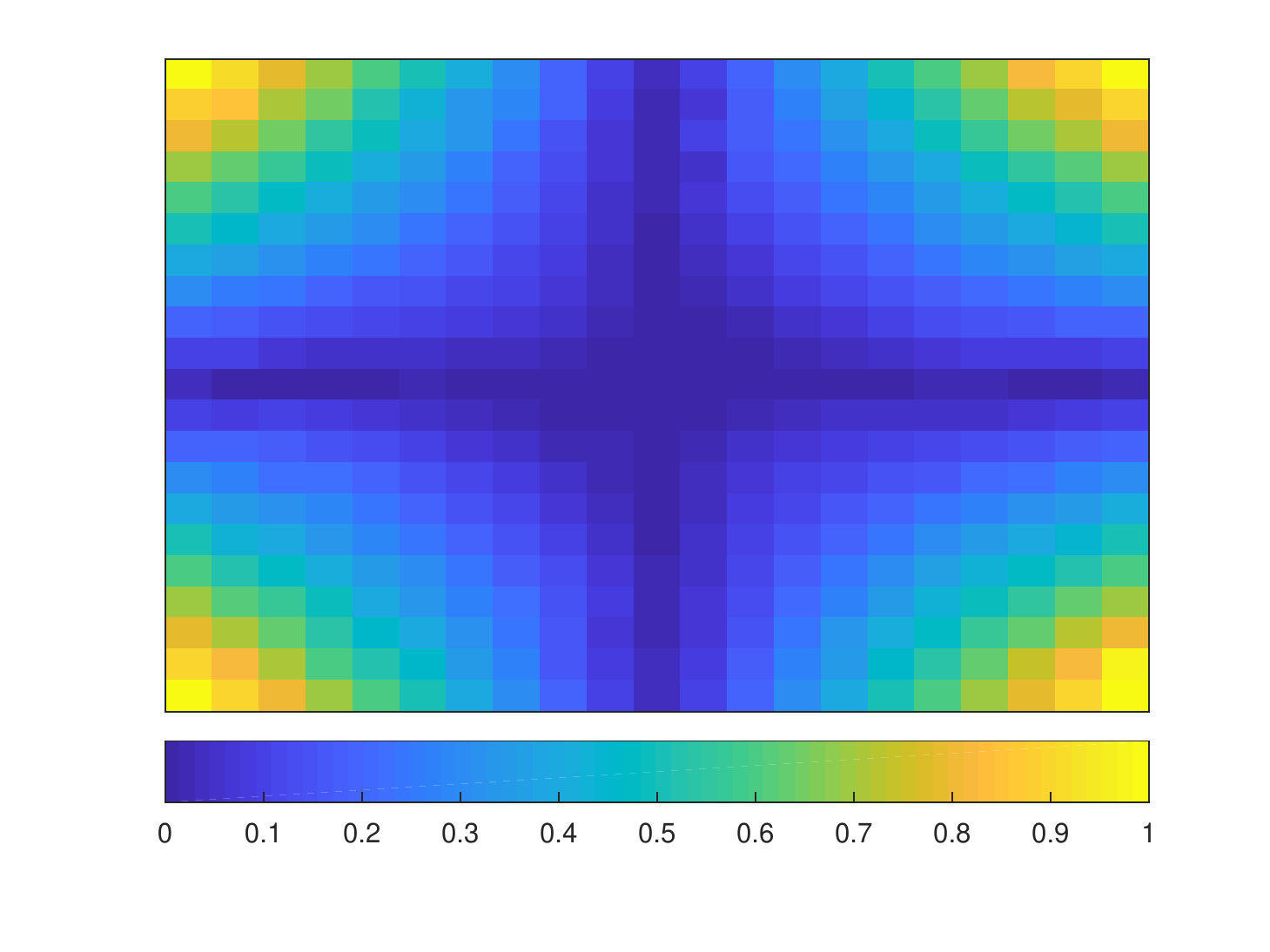}
&\includegraphics[width=\wsize,height=\hhsize,trim={1.5cm 1cm 1cm 0.5cm},clip]{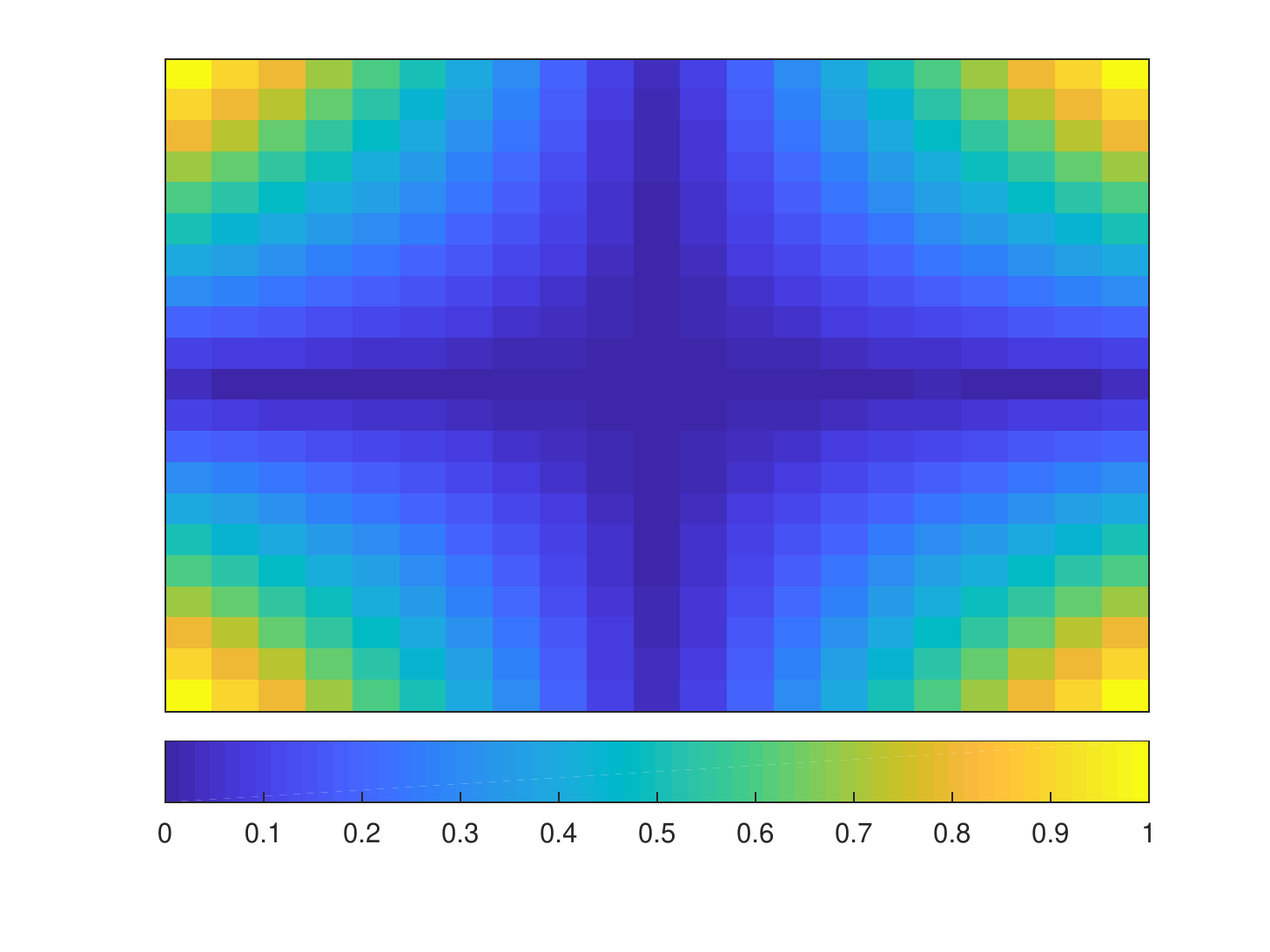}
&\includegraphics[width=\wsize,height=\hhsize,trim={1cm 0.5cm .5cm .5cm},clip]{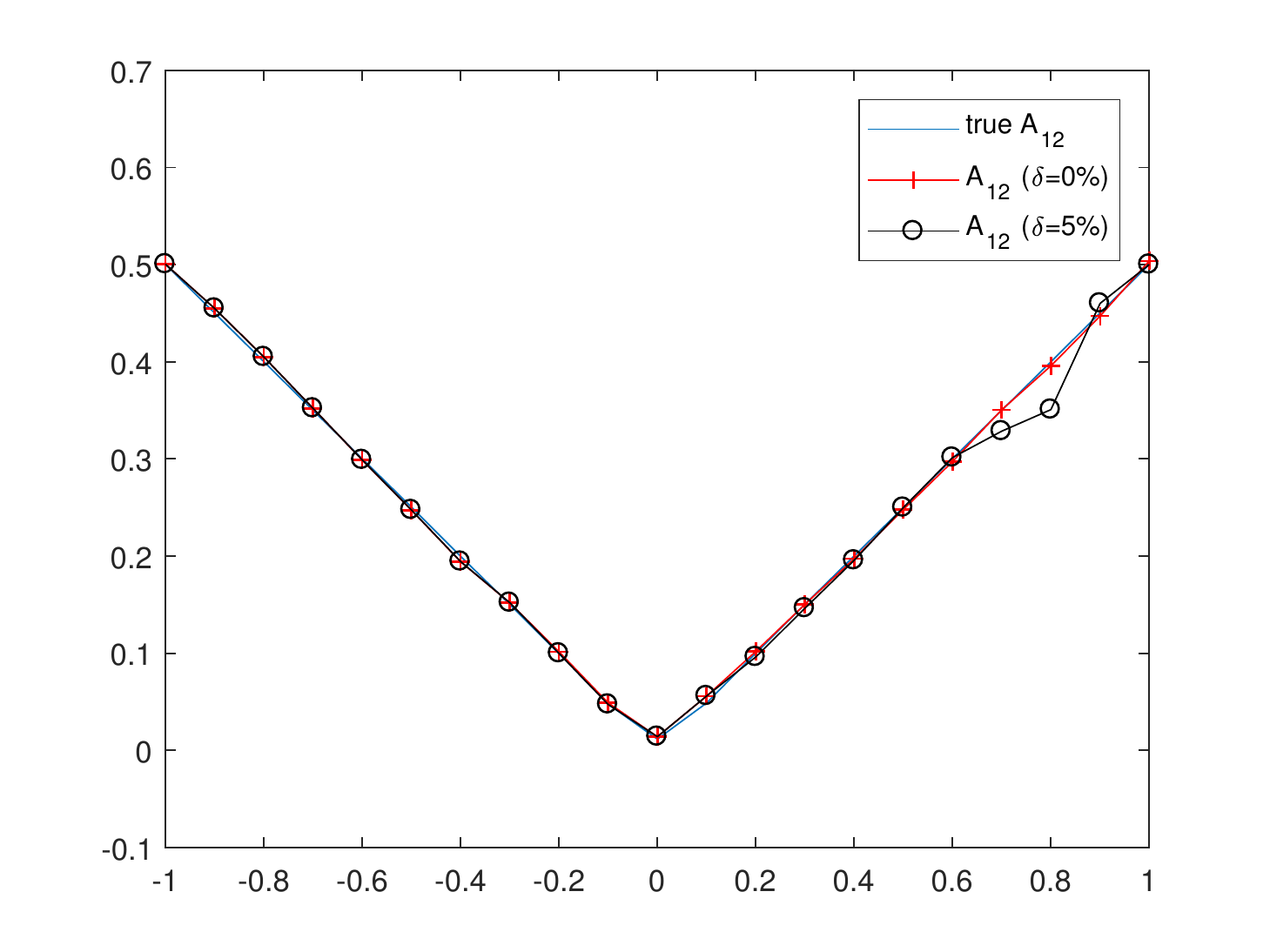}\\
 \includegraphics[width=\wsize,height=\hhsize,trim={1.5cm 1cm 1cm 0.5cm},clip]{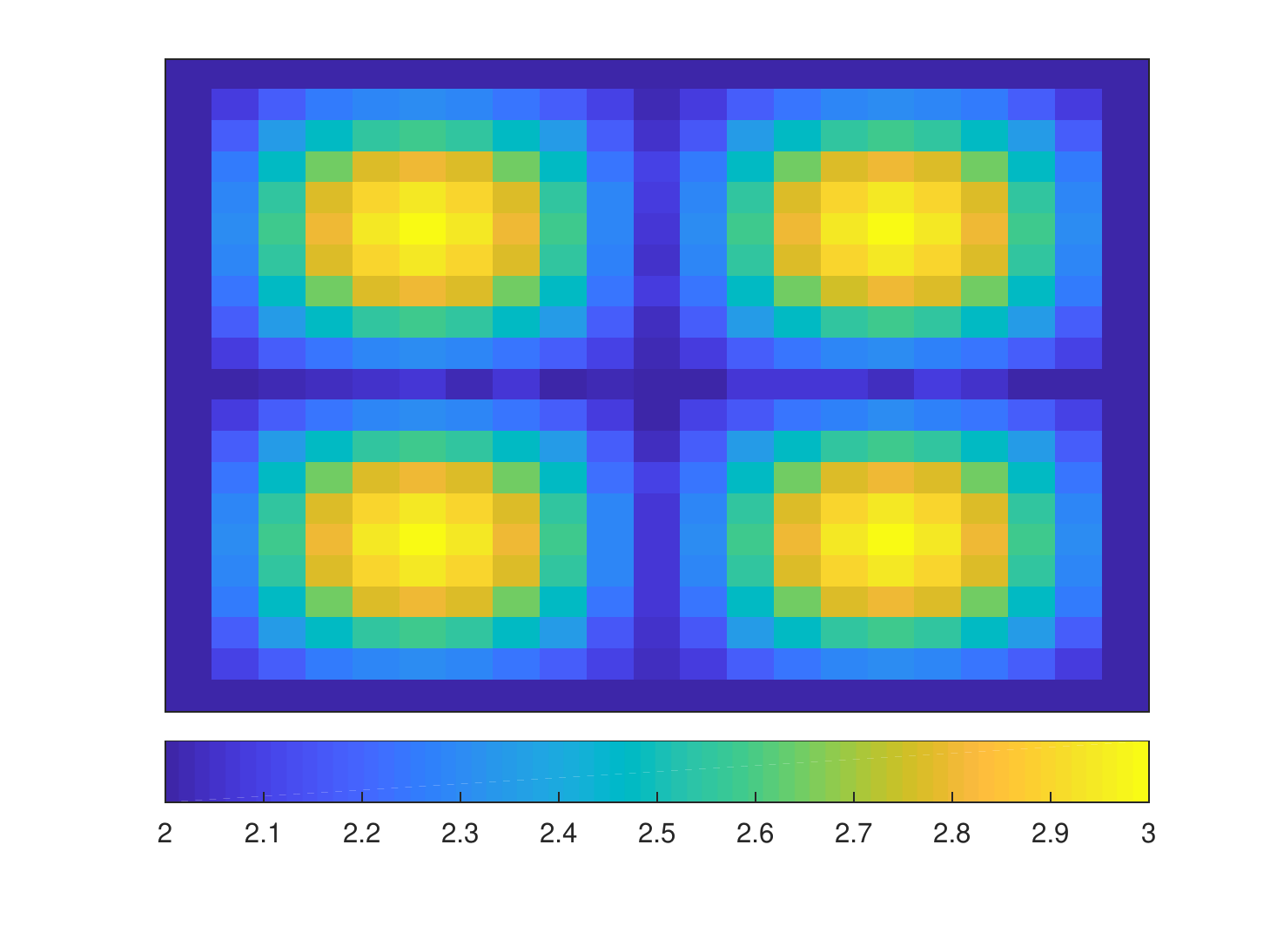}
&\includegraphics[width=\wsize,height=\hhsize,trim={1.5cm 1cm 1cm 0.5cm},clip]{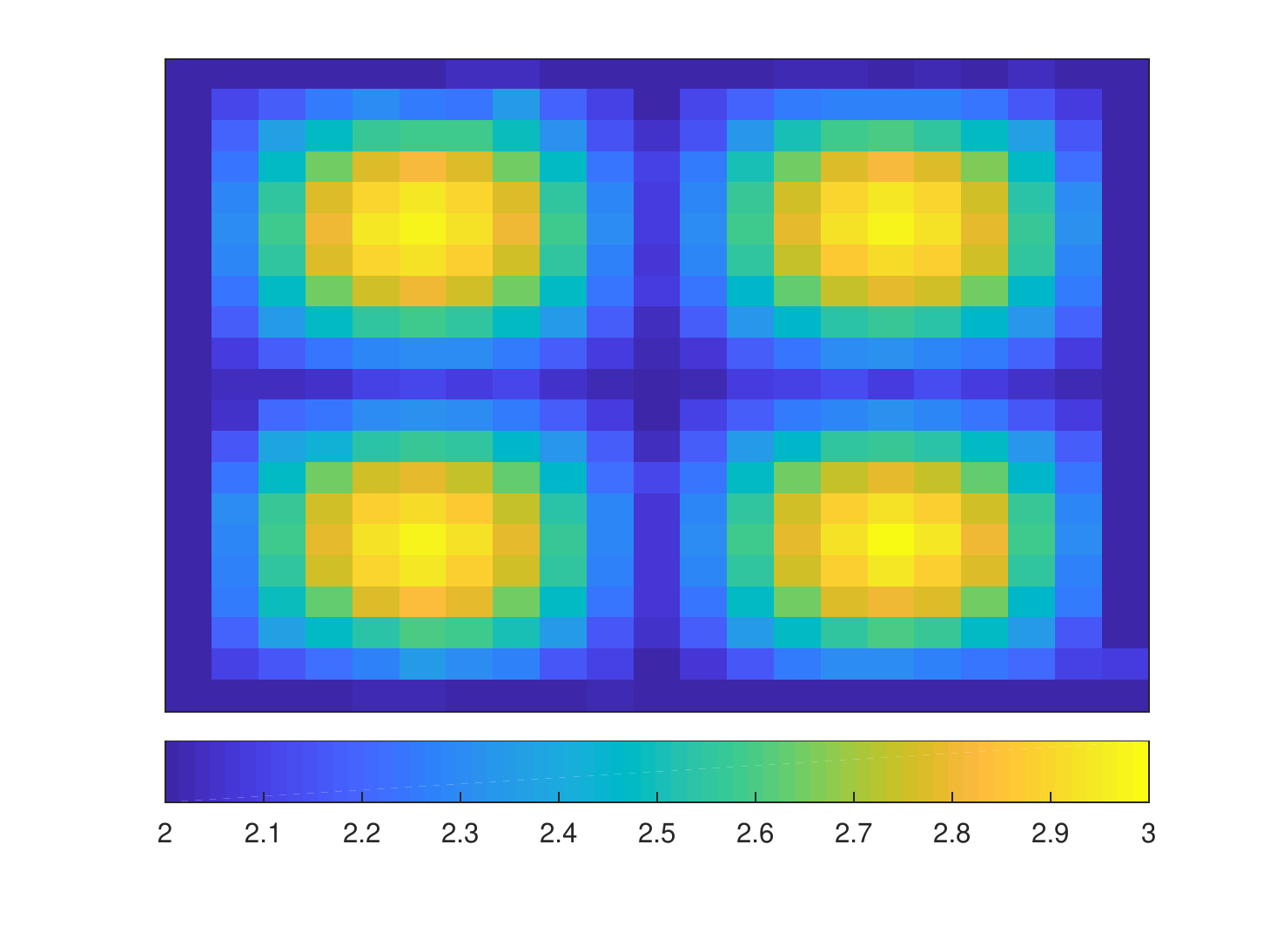}
&\includegraphics[width=\wsize,height=\hhsize,trim={1.5cm 1cm 1cm 0.5cm},clip]{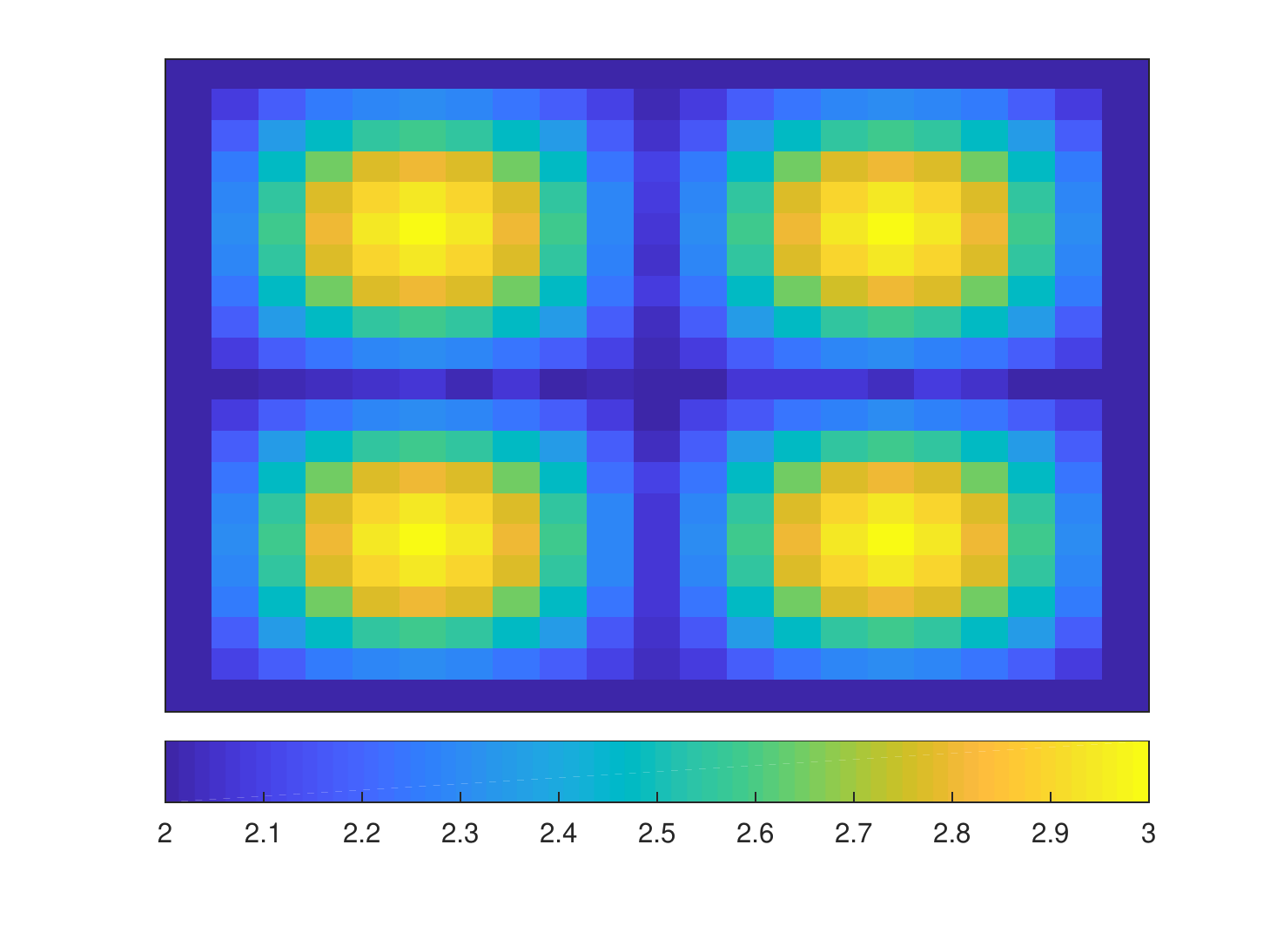}
&\includegraphics[width=\wsize,height=\hhsize,trim={1cm .5cm .5cm .5cm },clip]{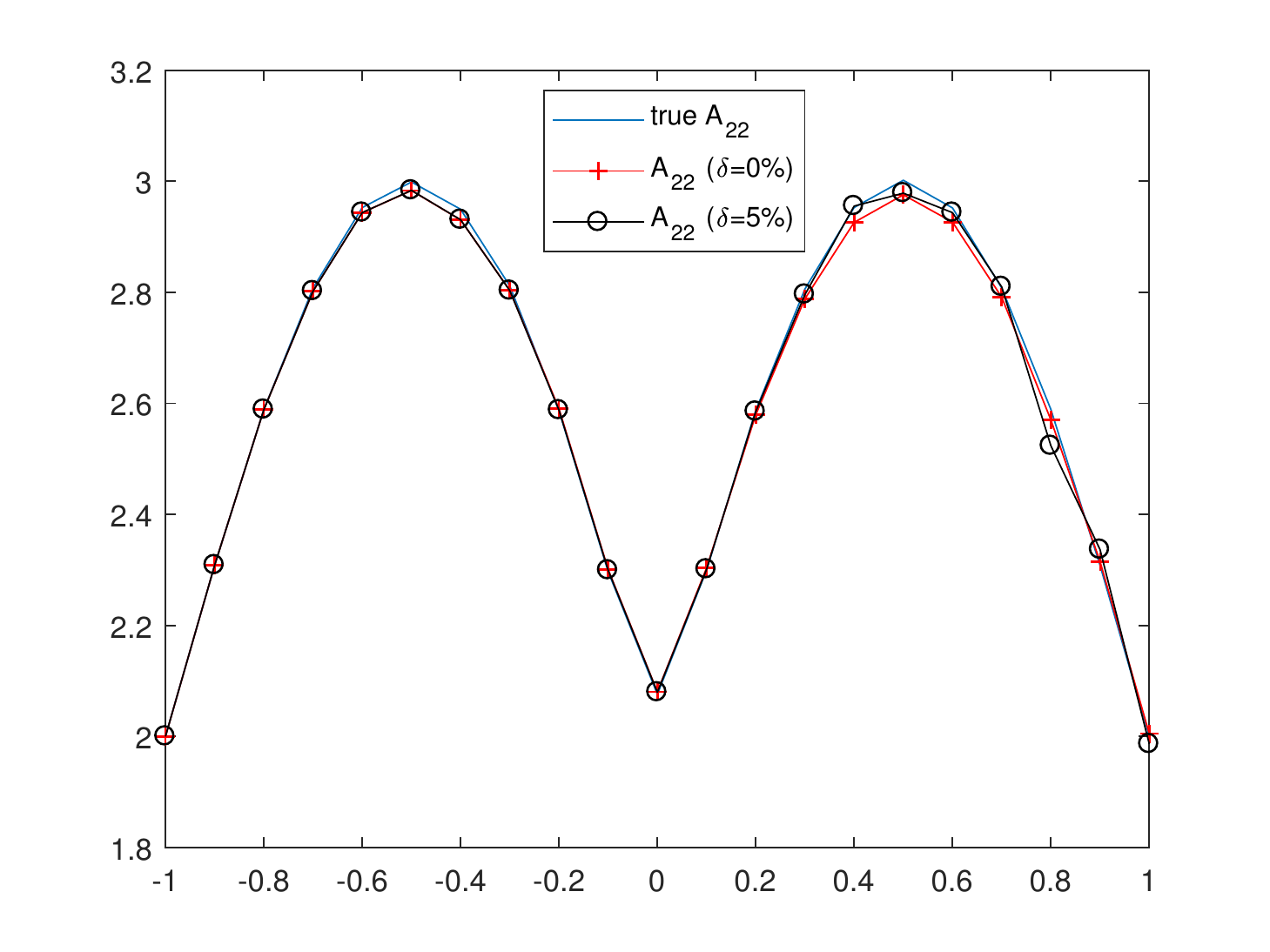}\\
\end{tabular}
\caption{The reconstructions for the conductivity tensor for Example 6: the top, middle and bottom rows refer to $A_{11}$, $A_{12}$ and $A_{22}$,
respectively. Column 1 is for exact conductivity, Column 2 for reconstruction with exact data, Column 3 for noisy data with $\delta=5\%$, and column 4 for the cross section along \{$x_2=-0.5$\}.}
\label{exam6}
\end{figure}

\begin{figure}[!htbp]
\centering
\setlength{\tabcolsep}{0pt}
\begin{tabular}{cccc}
 \includegraphics[width=\wsize,height=\hhsize,trim={1.5cm 1cm 1cm 0.5cm},clip]{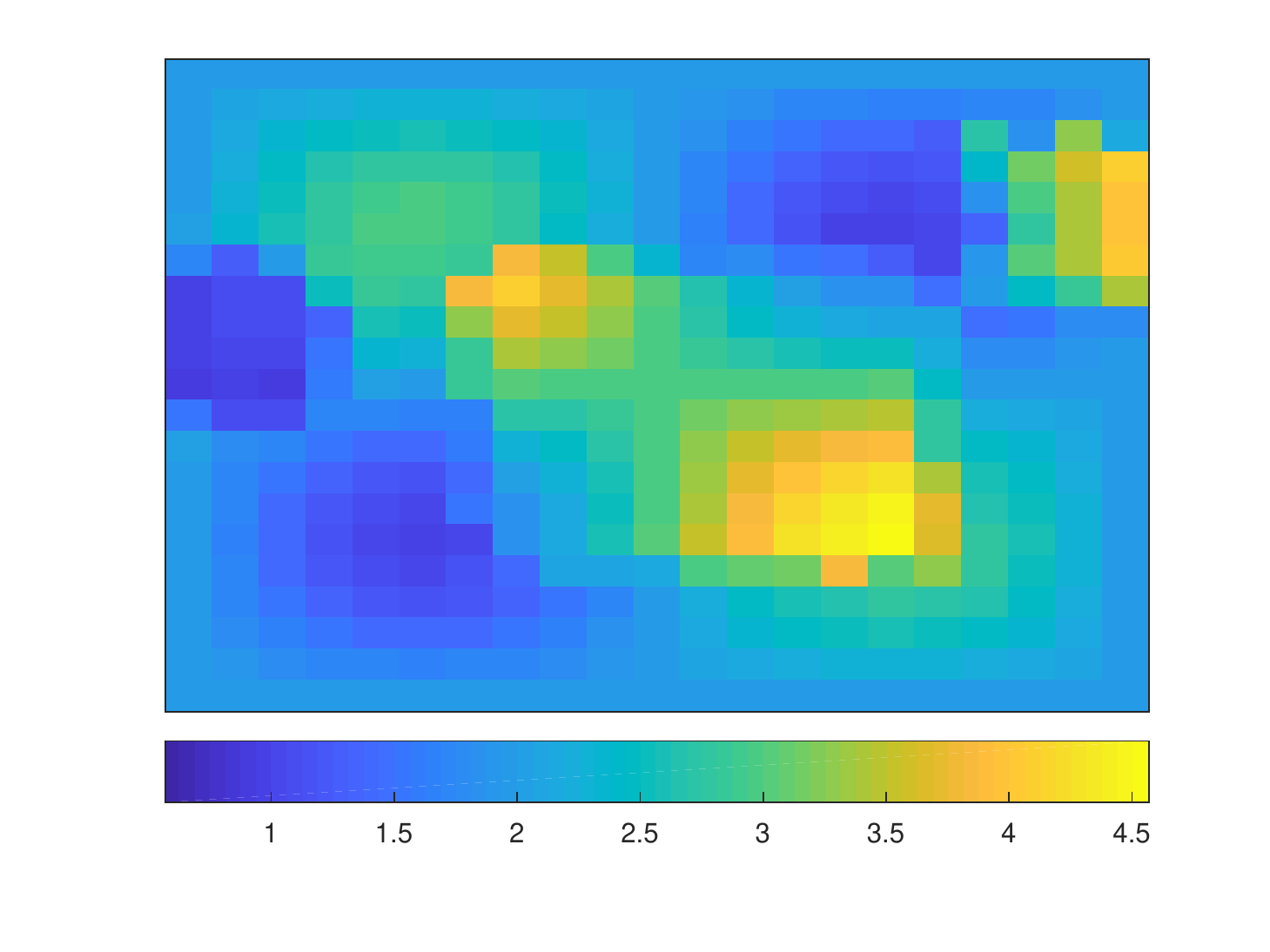}
&\includegraphics[width=\wsize,height=\hhsize,trim={1.5cm 1cm 1cm 0.5cm},clip]{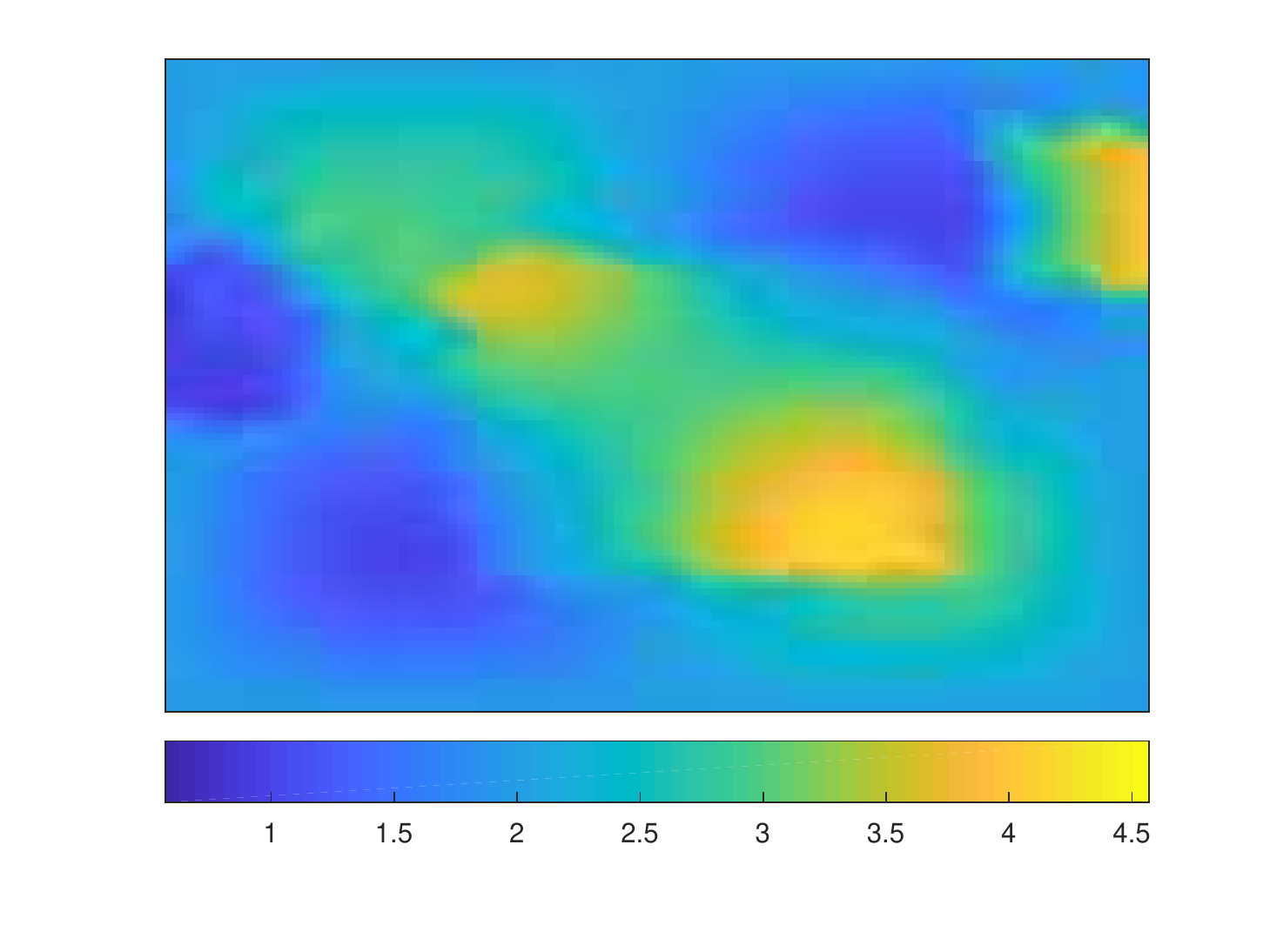}
&\includegraphics[width=\wsize,height=\hhsize,trim={1.5cm 1cm 1cm 0.5cm},clip]{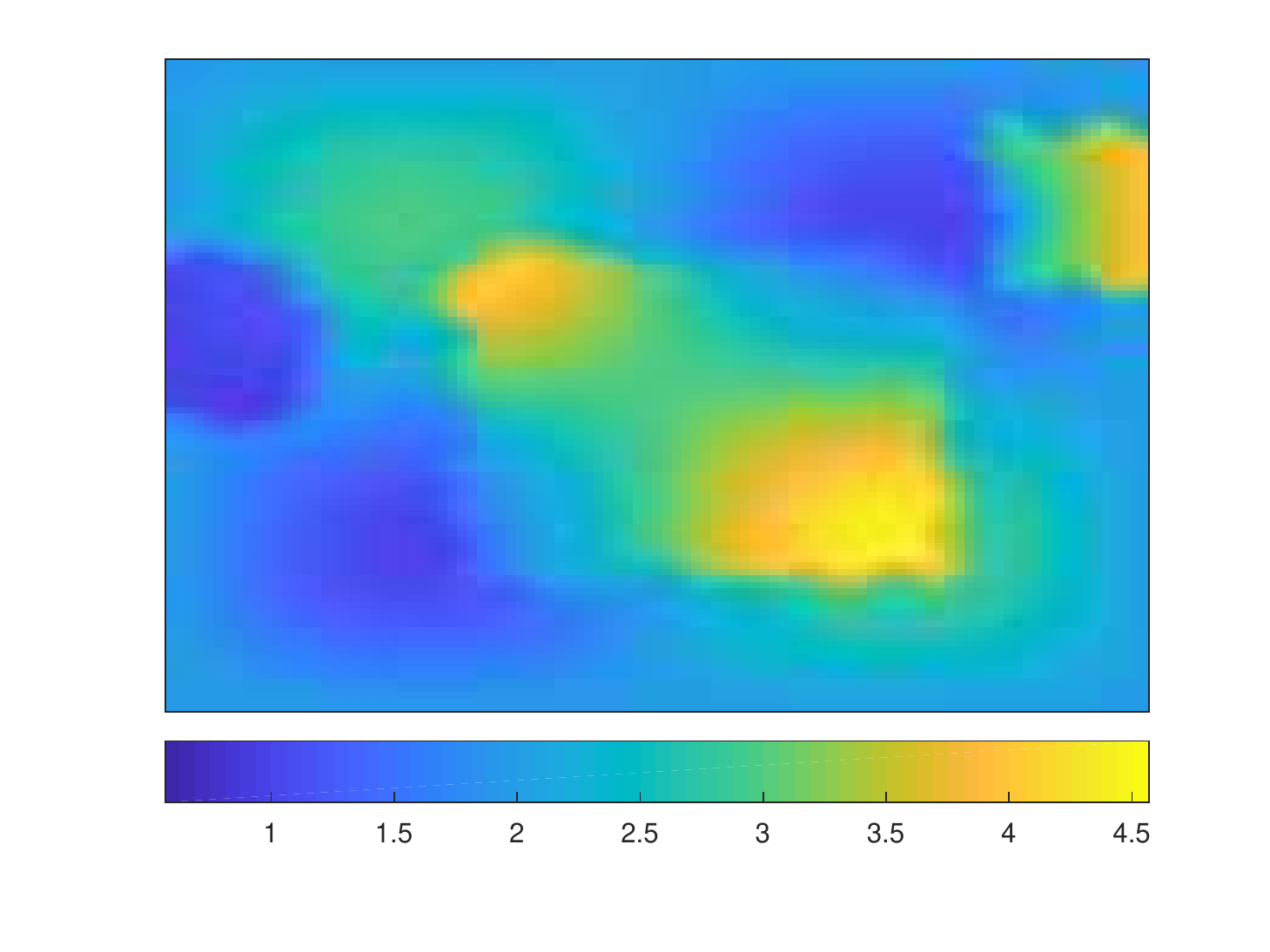}
&\includegraphics[width=\wsize,height=\hhsize,trim={1cm .5cm .5cm .5cm },clip]{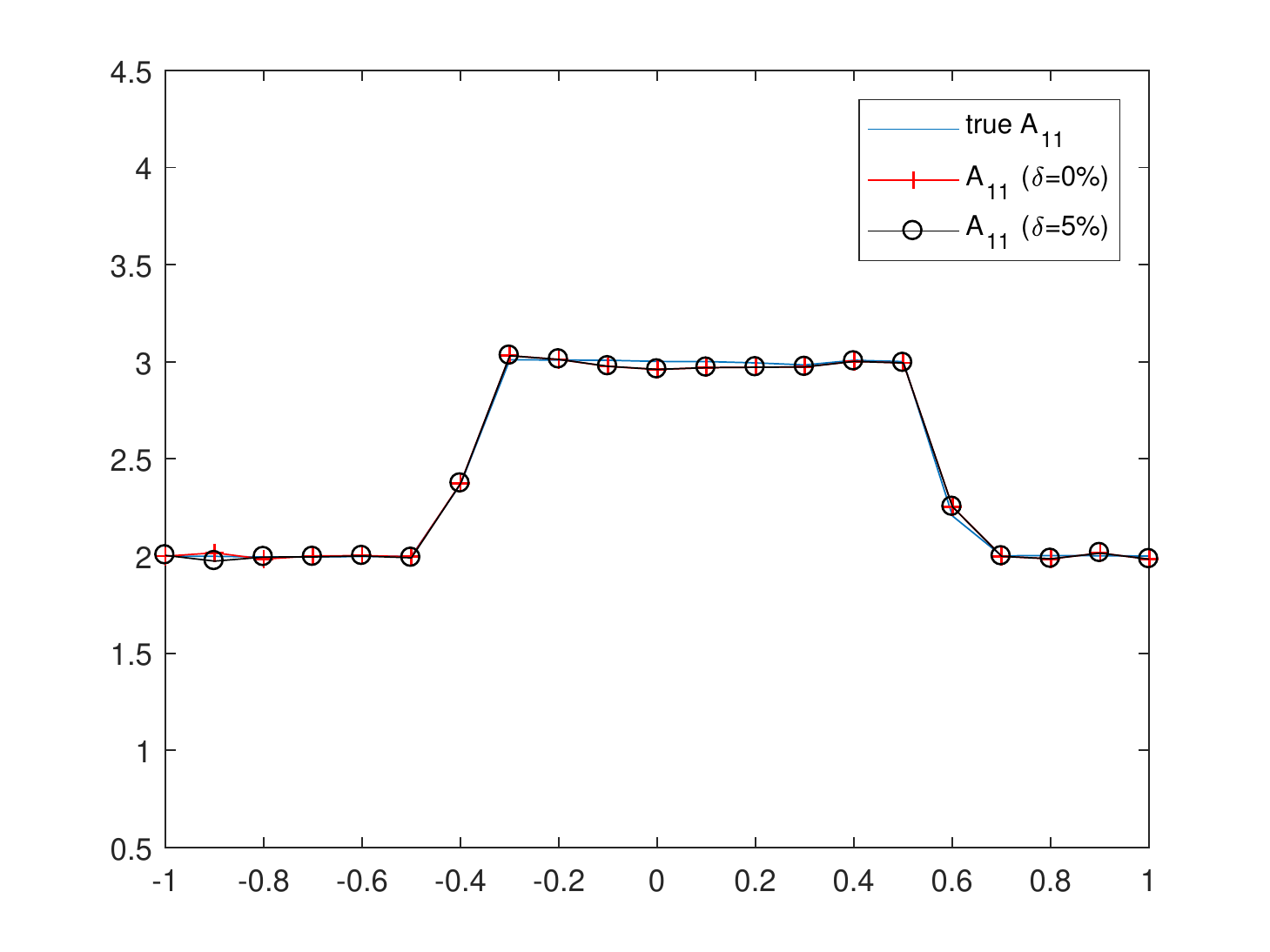}\\
 \includegraphics[width=\wsize,height=\hhsize,trim={1.5cm 1cm 1cm 0.5cm},clip]{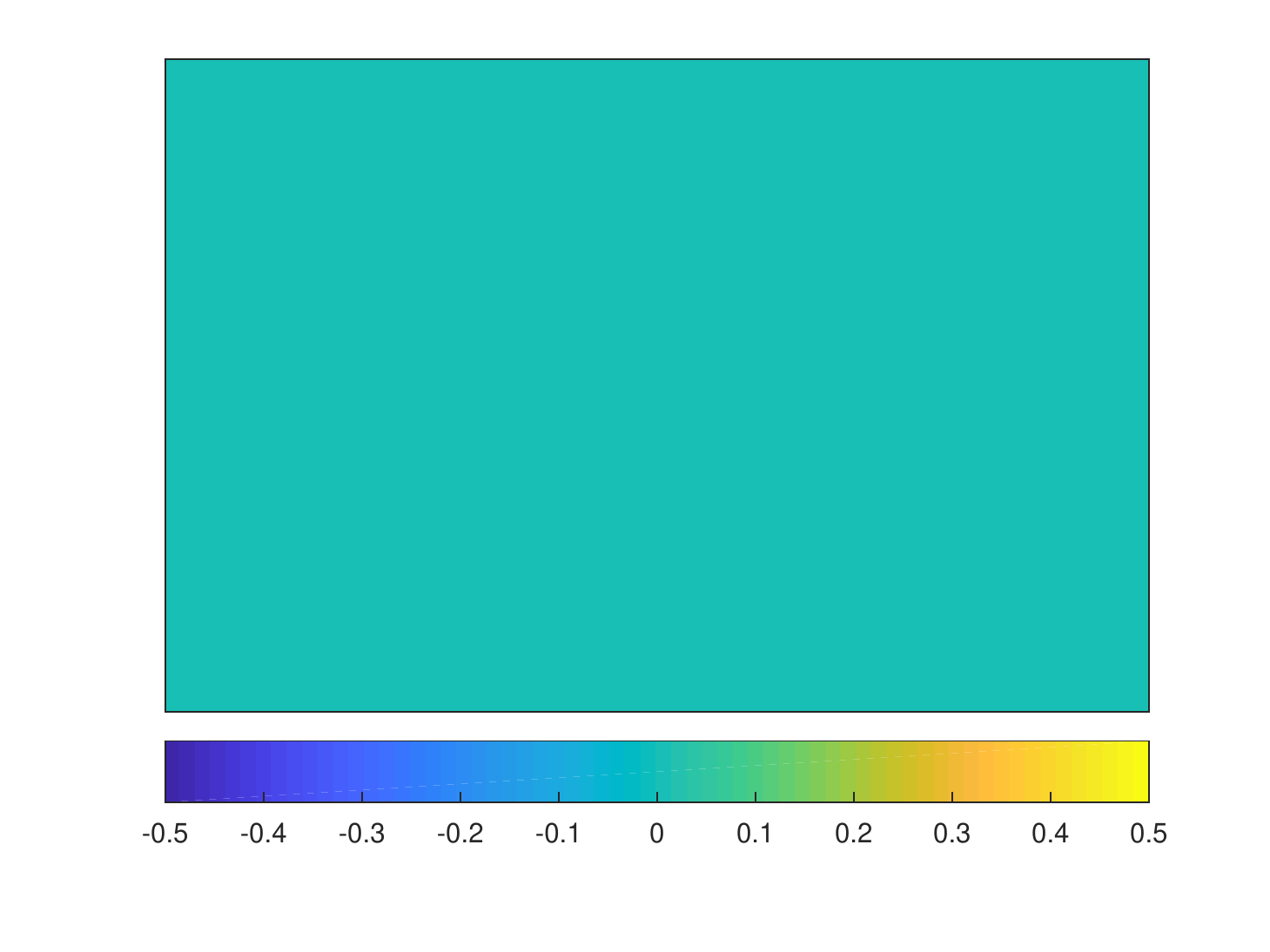}
&\includegraphics[width=\wsize,height=\hhsize,trim={1.5cm 1cm 1cm 0.5cm},clip]{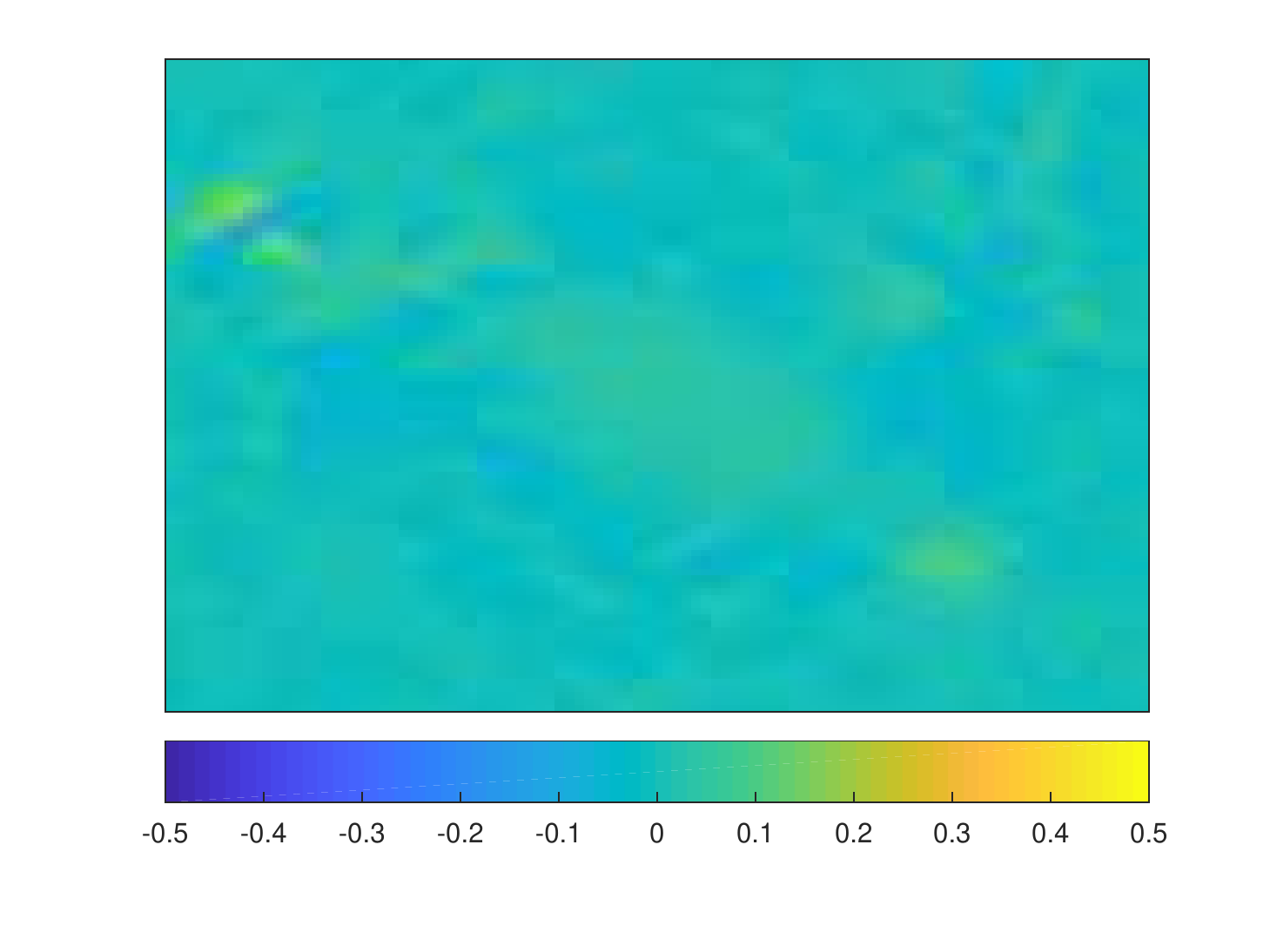}
&\includegraphics[width=\wsize,height=\hhsize,trim={1.5cm 1cm 1cm 0.5cm},clip]{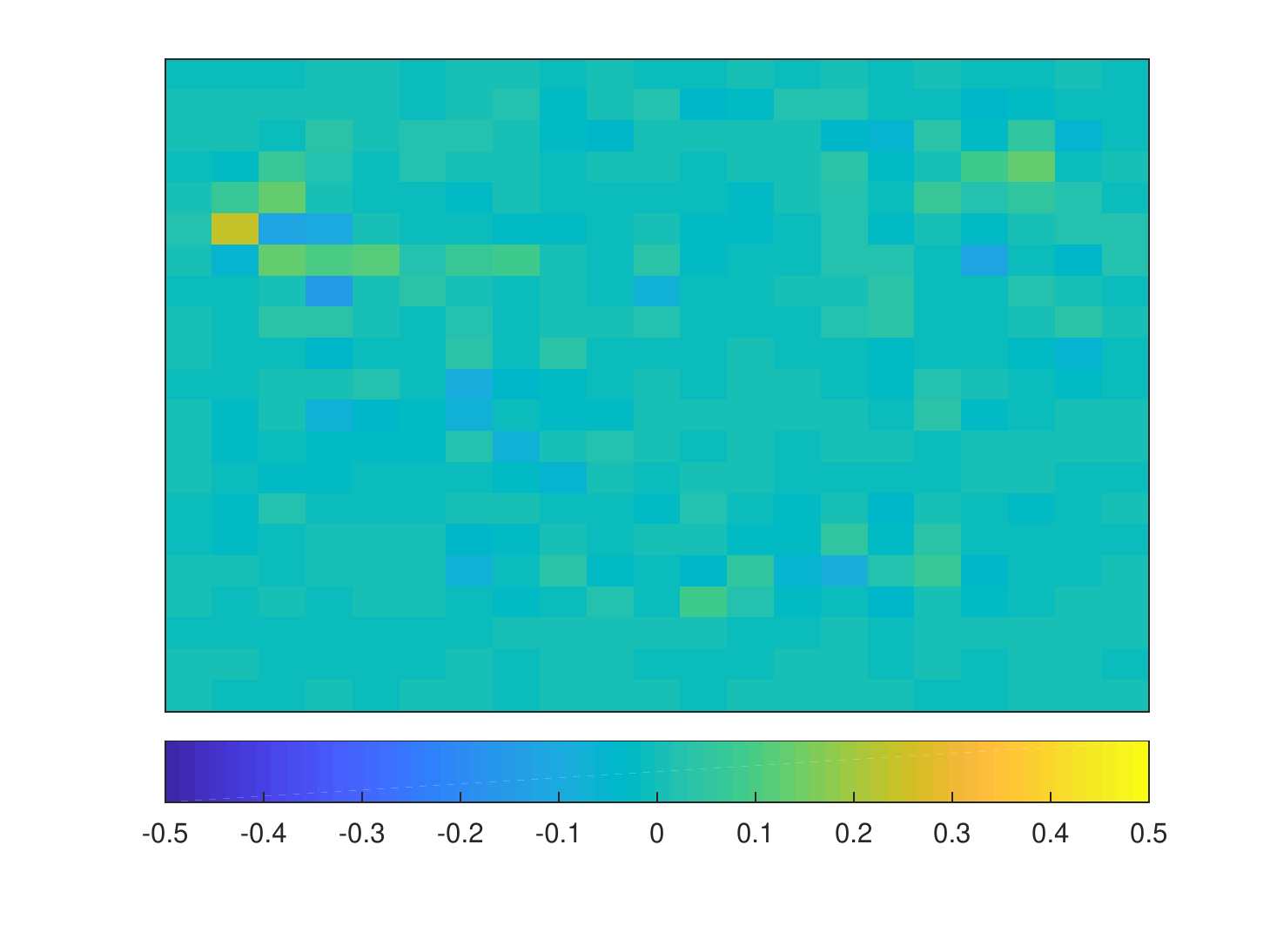}
&\includegraphics[width=\wsize,height=\hhsize,trim={1cm 0.5cm .5cm .5cm},clip]{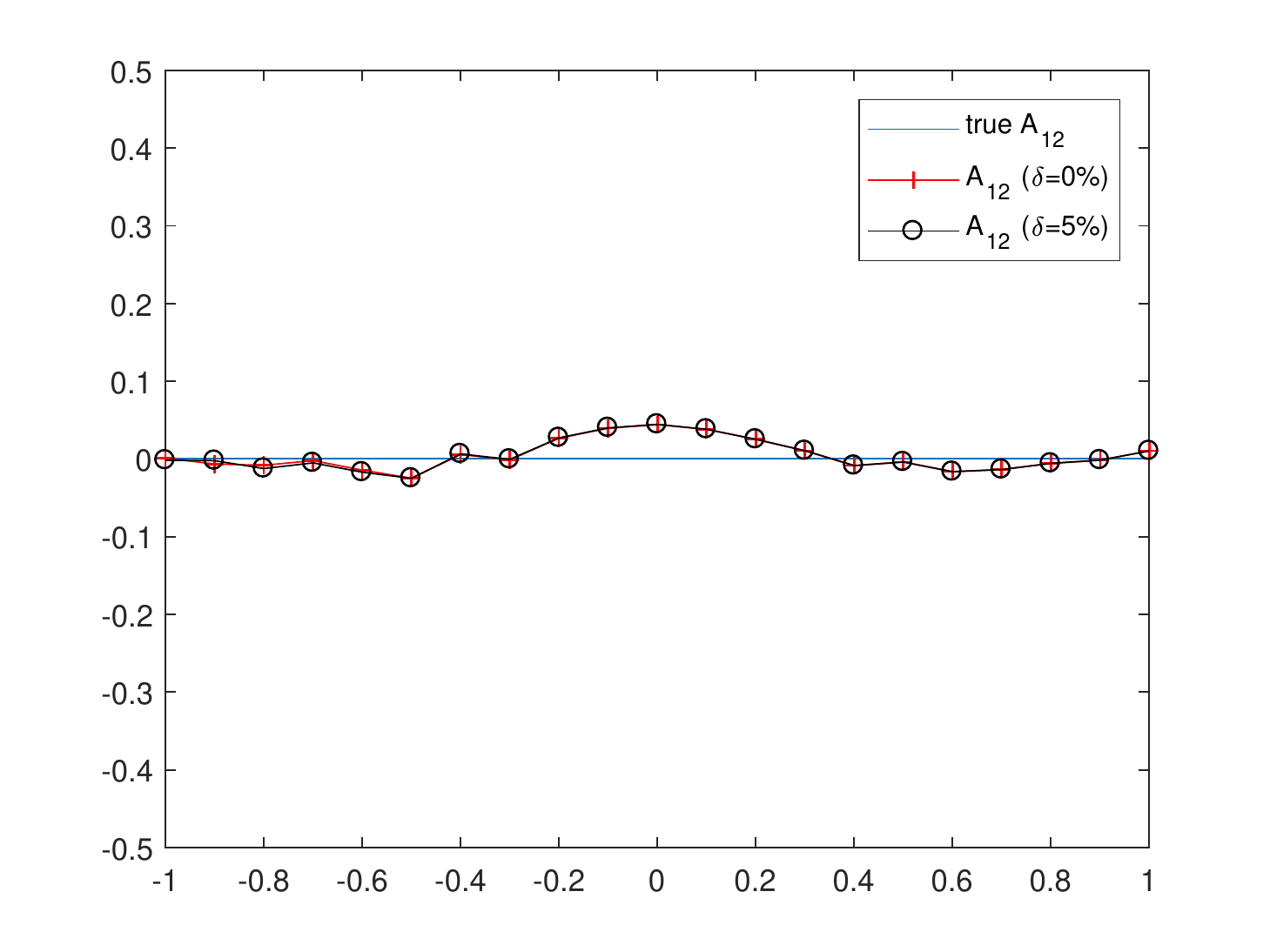}\\
 \includegraphics[width=\wsize,height=\hhsize,trim={1.5cm 1cm 1cm 0.5cm},clip]{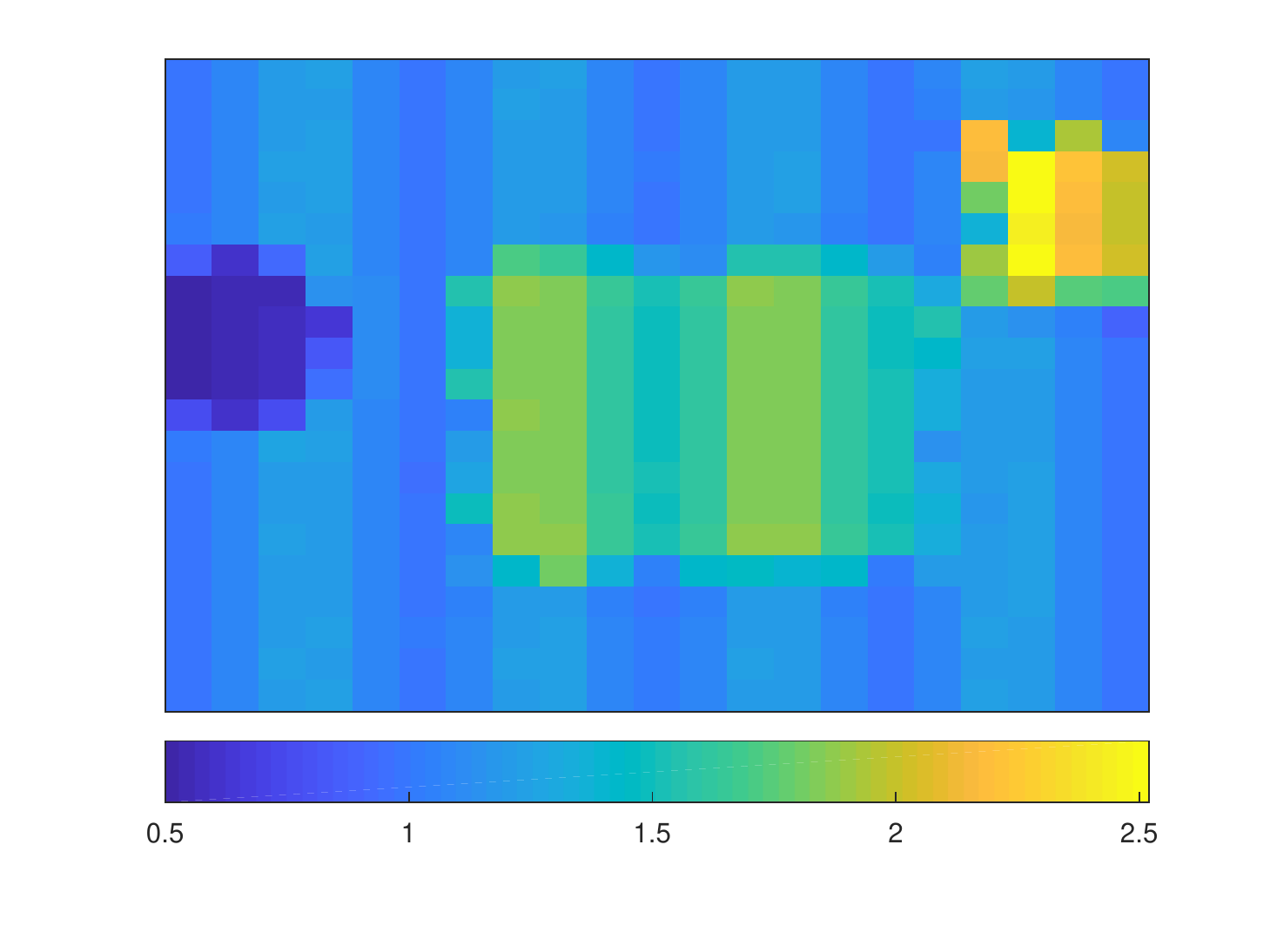}
&\includegraphics[width=\wsize,height=\hhsize,trim={1.5cm 1cm 1cm 0.5cm},clip]{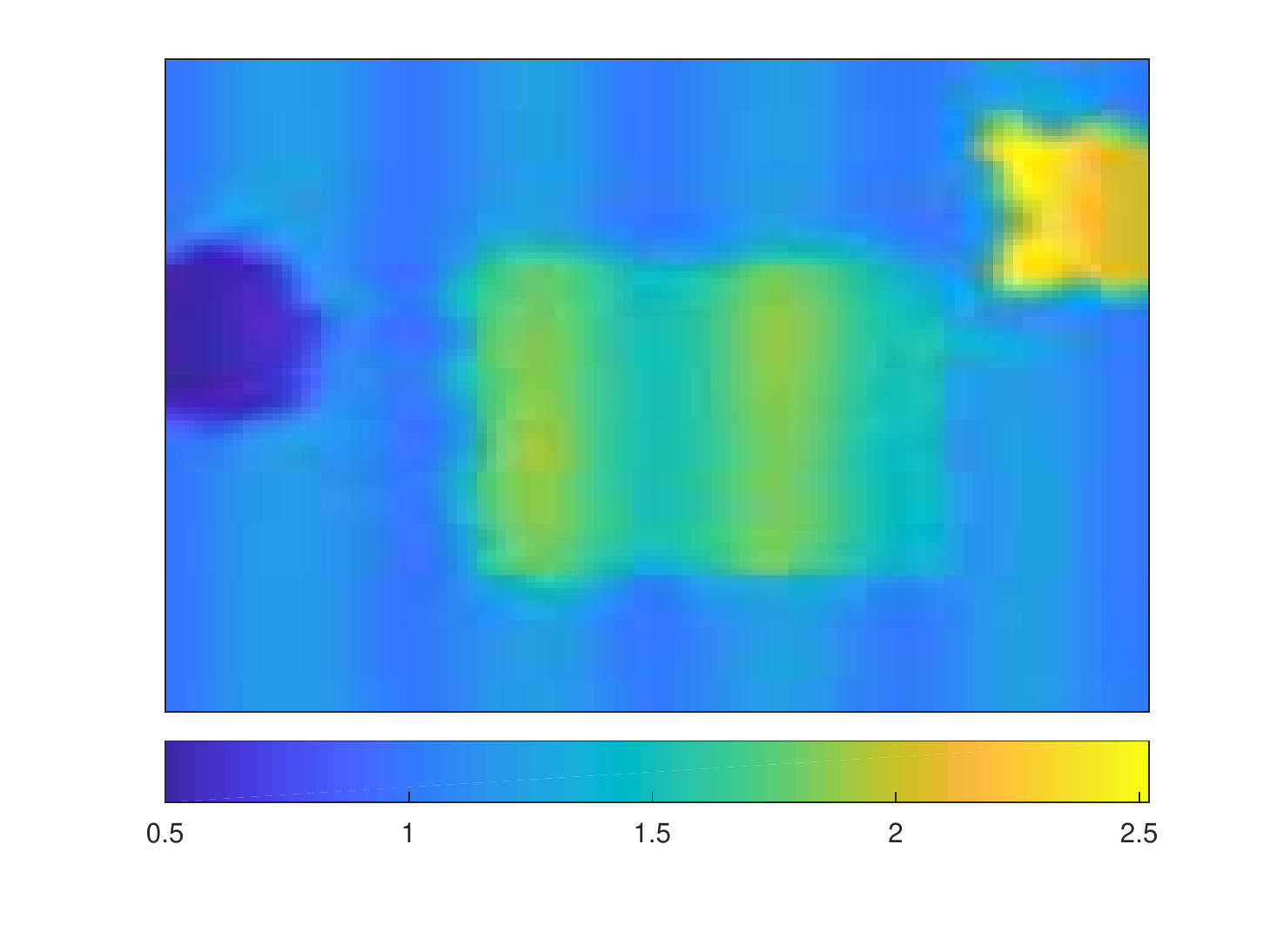}
&\includegraphics[width=\wsize,height=\hhsize,trim={1.5cm 1cm 1cm 0.5cm},clip]{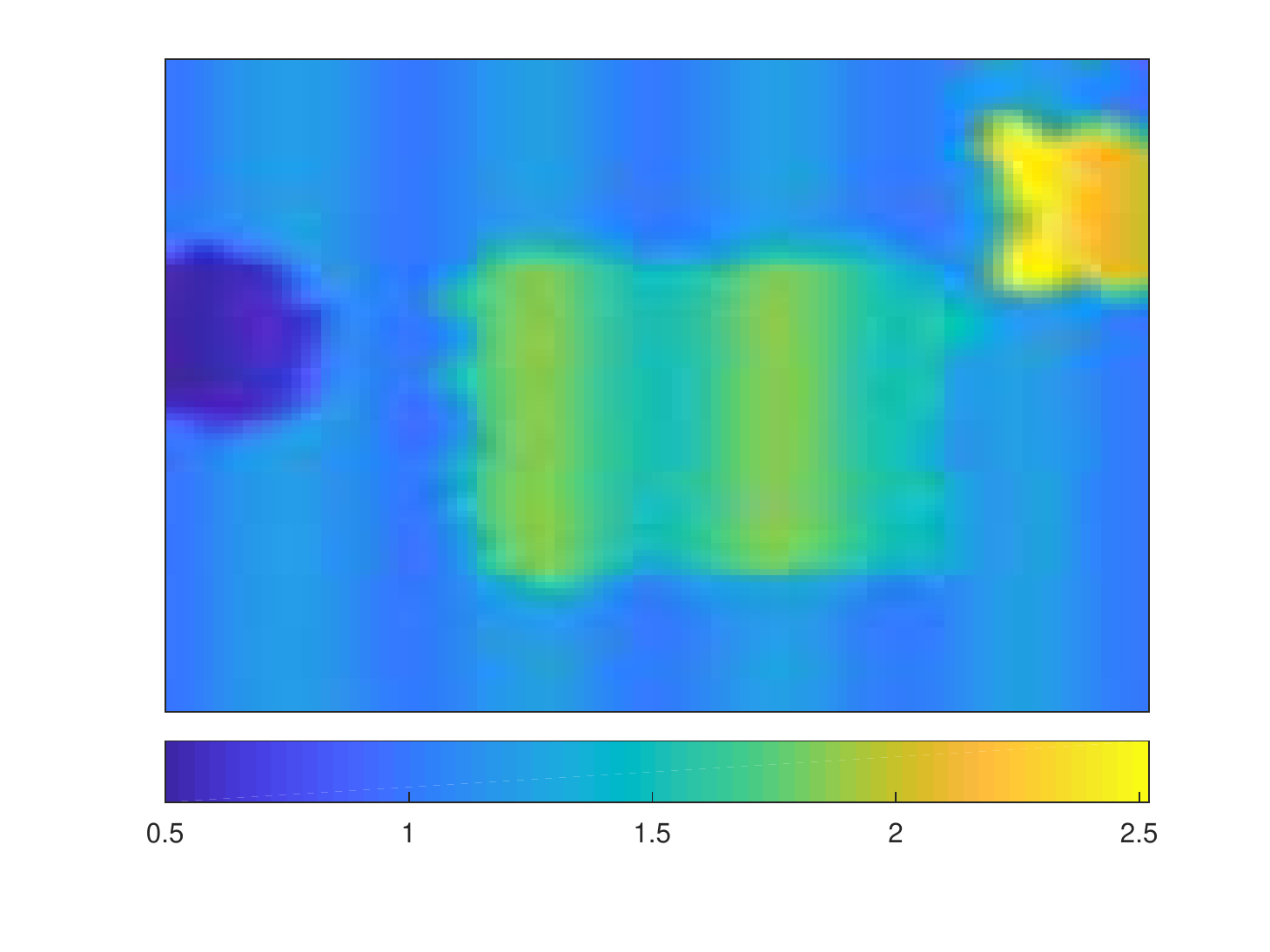}
&\includegraphics[width=\wsize,height=\hhsize,trim={1cm .5cm .5cm .5cm },clip]{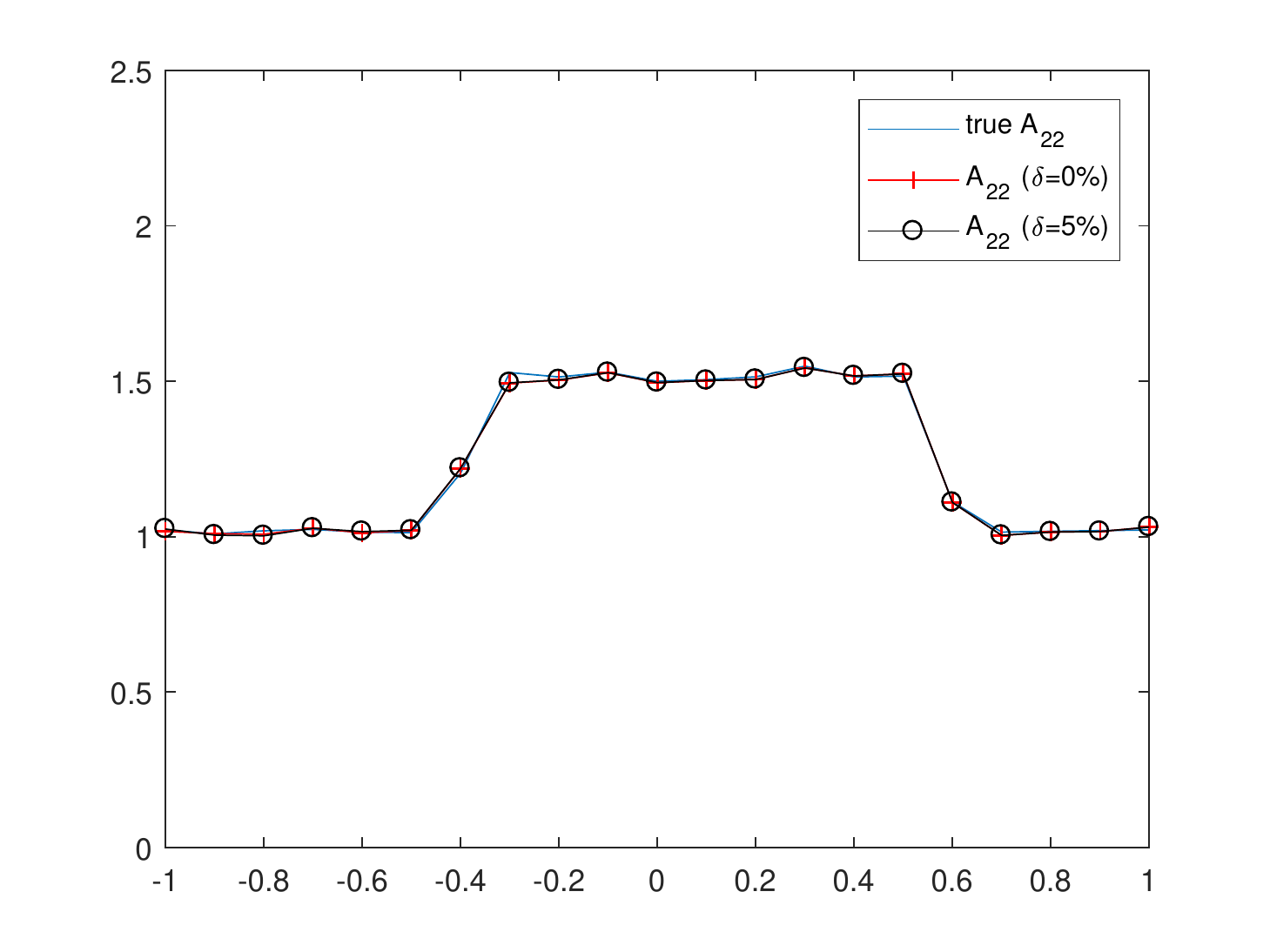}\\
\end{tabular}
\caption{The reconstructions for the conductivity tensor for Example 7: the top, middle and bottom rows refer to $A_{11}$, $A_{12}$ and $A_{22}$,
respectively. Column 1 is for exact conductivity, Column 2 for reconstruction with exact data, Column 3 for noisy data with $\delta=5\%$, and column 4 for the cross section along \{$x_1=0$\}.}
\label{exam7}
\end{figure}

\begin{figure}[ht!]
\centering
\subfigure[accuracy error $e$]{
\includegraphics[width=5.2cm,height=4.5cm,trim={0.5cm 0cm 1cm 0.5cm},clip]{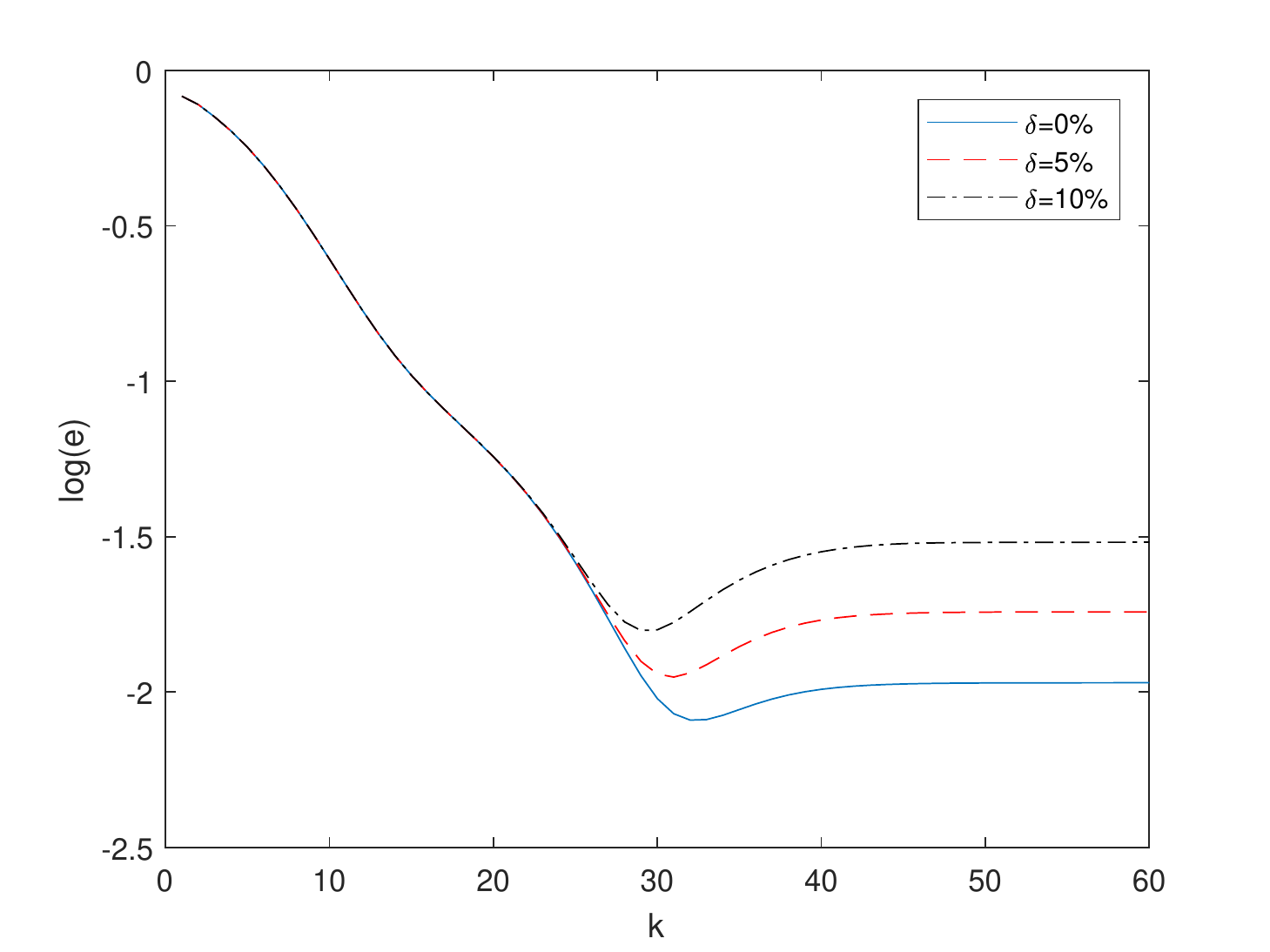}}
\hspace{0.5in}
\subfigure[residual $E$]{
\includegraphics[width=5.2cm,height=4.5cm,trim={0.5cm 0cm 1cm 0.5cm},clip]{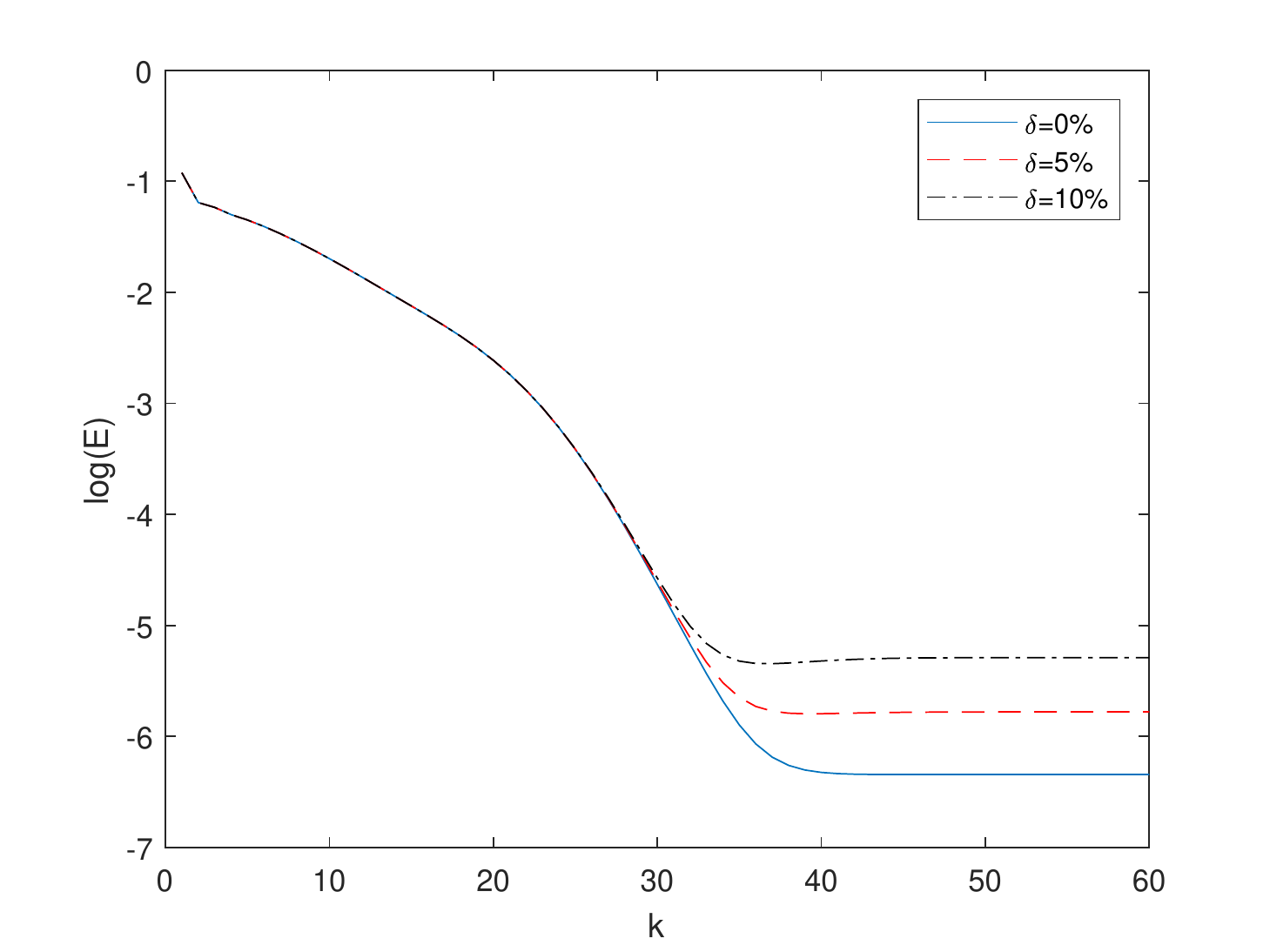}}
\hspace{0.5in}
\caption{The convergence of the aglorithm for Example 4 with noisy data.}
\label{err} 
\end{figure}

The accuracy of a reconstruction $\tilde A$ with respect to the exact one $A^\dag$
is measured by the relative error $e=\normmm{\tilde{A}-A^\dag}_2/\normmm{A^\dag}_2$.
The relative errors for Examples \ref{exam1}-\ref{exam7} with different noise levels are
given in Table \ref{table2}. It is clearly observed that for all examples, as
the noise level $\delta$ decreases, the reconstructions become more accurate.
For exact data, the reconstruction error $e$ is nonzero, since the mesh used
for the inversion step is not very refined, and there is an inevitable discretization error
(and also requires suitable regularization). One can also observe that the errors in the
reconstruction of the discontinuous conductivity tensor (Example 7) are much larger than
the continuous cases in Examples 1--6, since the sharp edges cannot be accurately
captured. This clearly indicates the numerical challenge
associated with recovering nonsmooth conductivity tensors.

To gain further insight into the reconstructions, we present in Figs. \ref{exam1}-\ref{exam7}
the exact conductivity tensor, the reconstructed tensors for exact and noisy data, and
a cross section. These plots show clearly that the reconstructions are fairly
accurate, even in the discontinuous case, agreeing well with the quantitative results in
Table \ref{table2}. Indeed, the kinks in the nonsmooth conductivity components can be
accurately resolved, showing clearly the accuracy of the proposed approach. A careful comparison
between the results shown in Figs. \ref{exam1}-\ref{exam2}
and Figs. \ref{exam3}-\ref{exam7} indicates that the reconstructions of an isotropic
conductivity tends to be more stable than the anisotropic case. This is consistent with
existing theoretical findings \cite[Theorem 2.4]{G.Bal2014}.
Nonetheless, Figs. \ref{exam5}-\ref{exam7} show that the approach can also stably and accurately
recover an anisotropic conductivity tensor of distinct features for both exact and noisy data (with
up to 10 percent of noise in the data).

Next we briefly examine the convergence behavior of the algorithm. In Fig. \ref{err}, we
present the convergence of the relative error and residual for Example 4 versus
iteration number (without stopping by the discrepancy principle). It is observed that as the noise
level $\delta$ decreases from $10\%$ to $0\%$, the accuracy of the reconstruction $\tilde{A}$
improves steadily, which partly agrees with Theorem \ref{stability}. Also from Fig. \ref{err},
we observe a robust convergence of the algorithm: it converges after about 30 iterations,
irrespective of the noise level $\delta$. Although not presented, these observations hold
also for other examples. These results confirm the efficiency of the projected
Newton algorithm for conductivity tensor recovery.

\begin{figure}[!htbp]
\centering \setlength{\tabcolsep}{0pt}
\begin{tabular}{cccc}
 \includegraphics[width=\wsize,height=\hhsize,trim={1.5cm 1cm 1cm 0.5cm},clip]{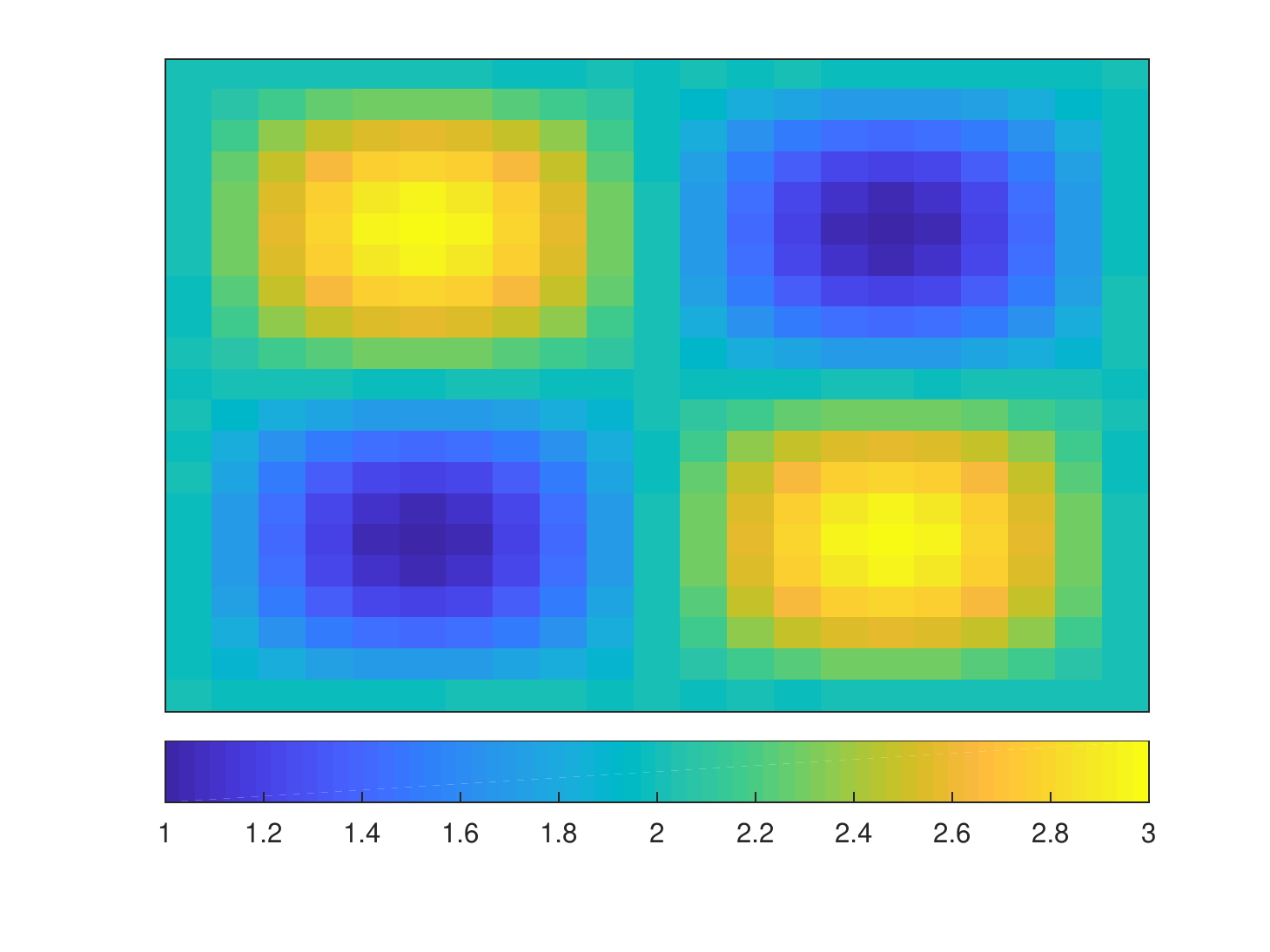}
&\includegraphics[width=\wsize,height=\hhsize,trim={1.5cm 1cm 1cm 0.5cm},clip]{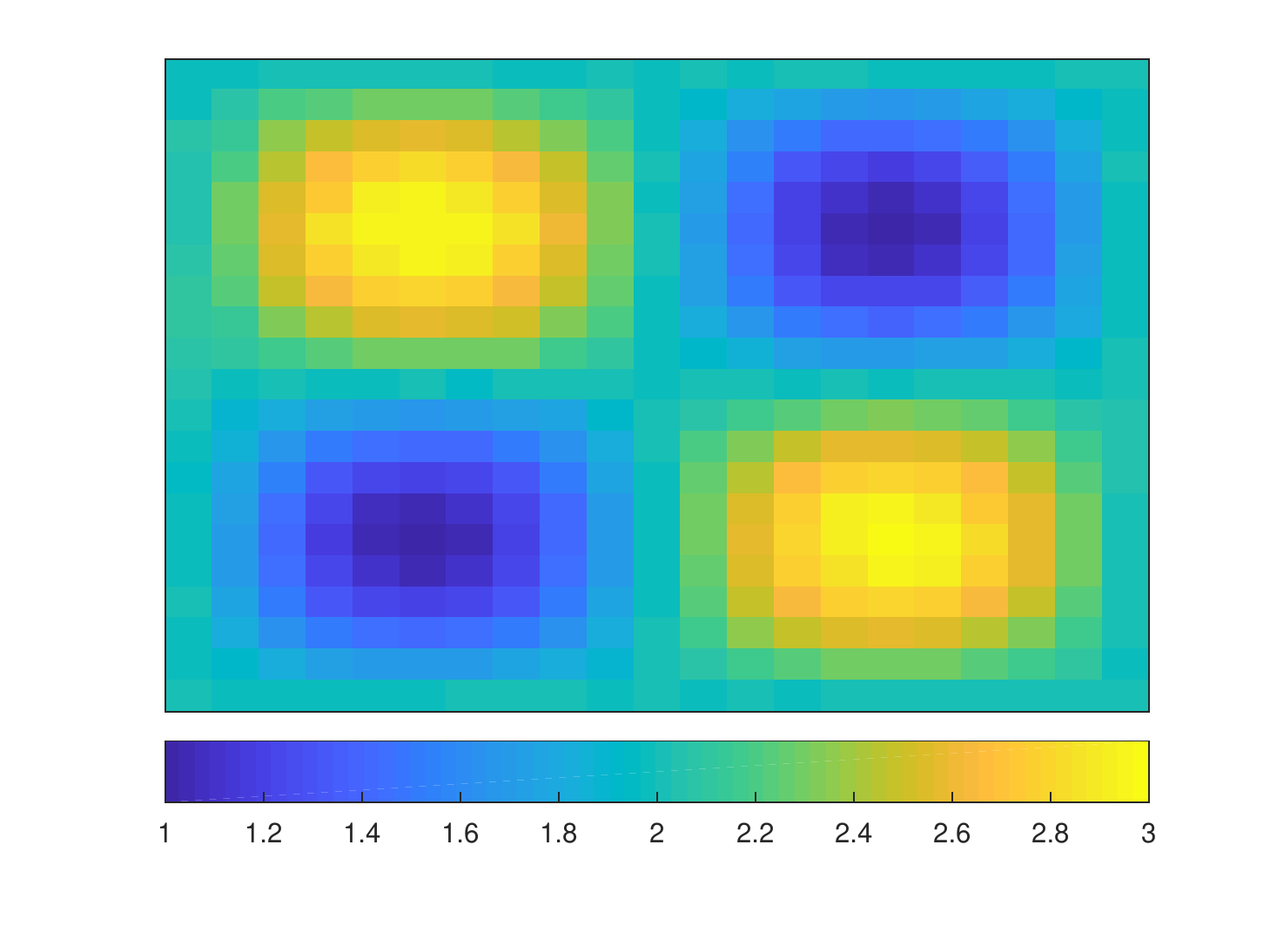}
&\includegraphics[width=\wsize,height=\hhsize,trim={1.5cm 1cm 1cm 0.5cm},clip]{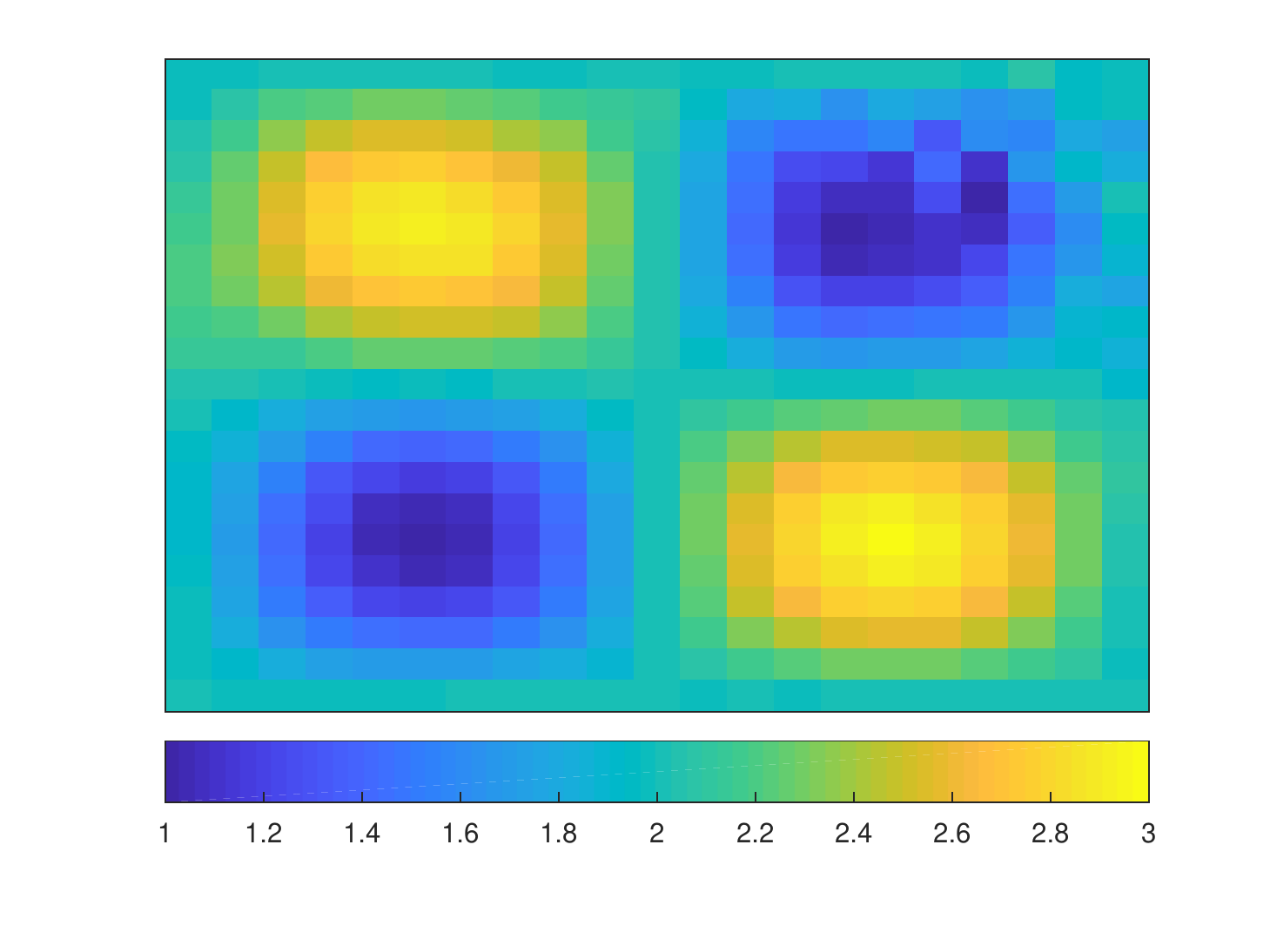}
&\includegraphics[width=\wsize,height=\hhsize,trim={1cm .5cm .5cm .5cm },clip]{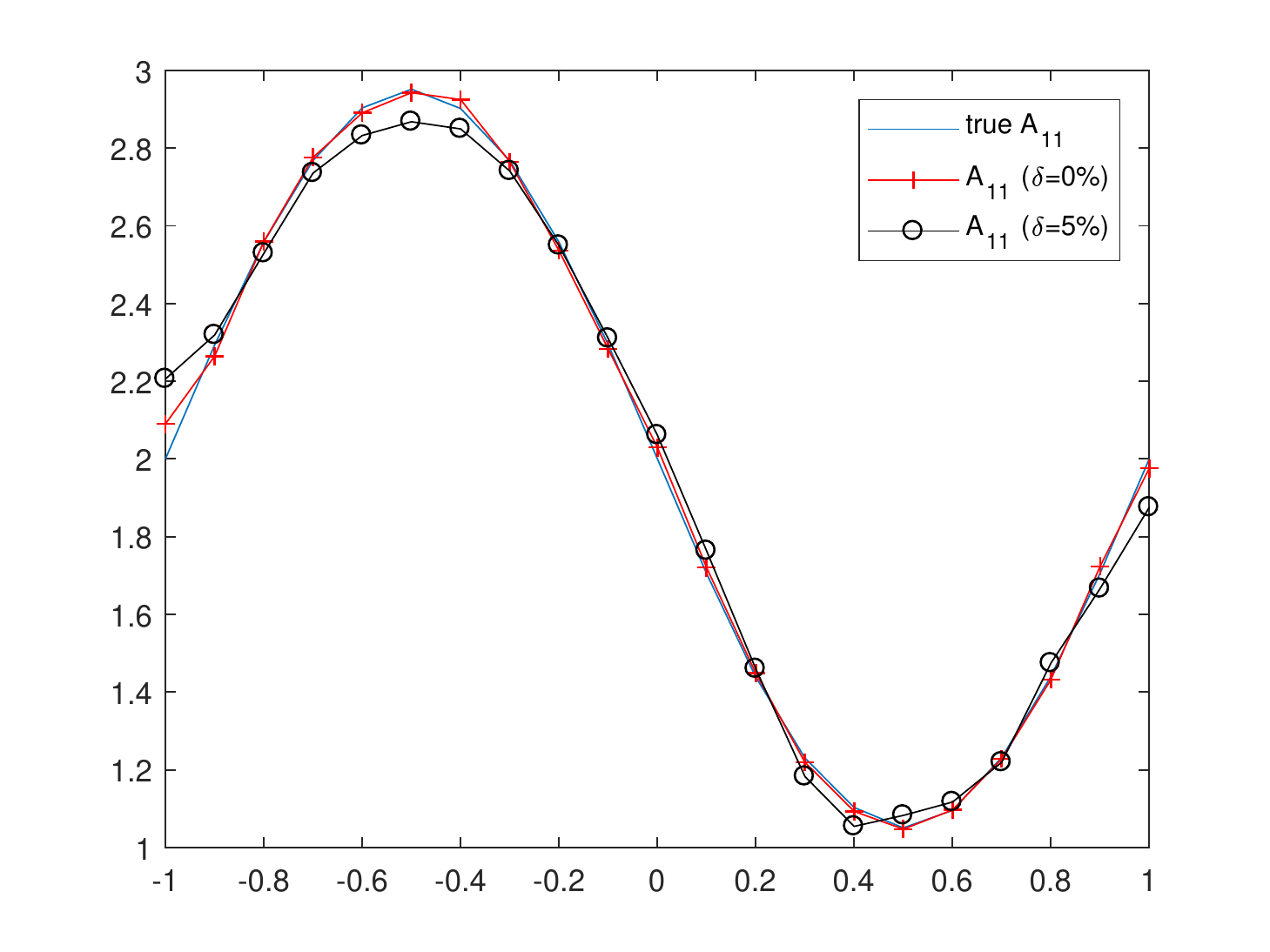}\\
 \includegraphics[width=\wsize,height=\hhsize,trim={1.5cm 1cm 1cm 0.5cm},clip]{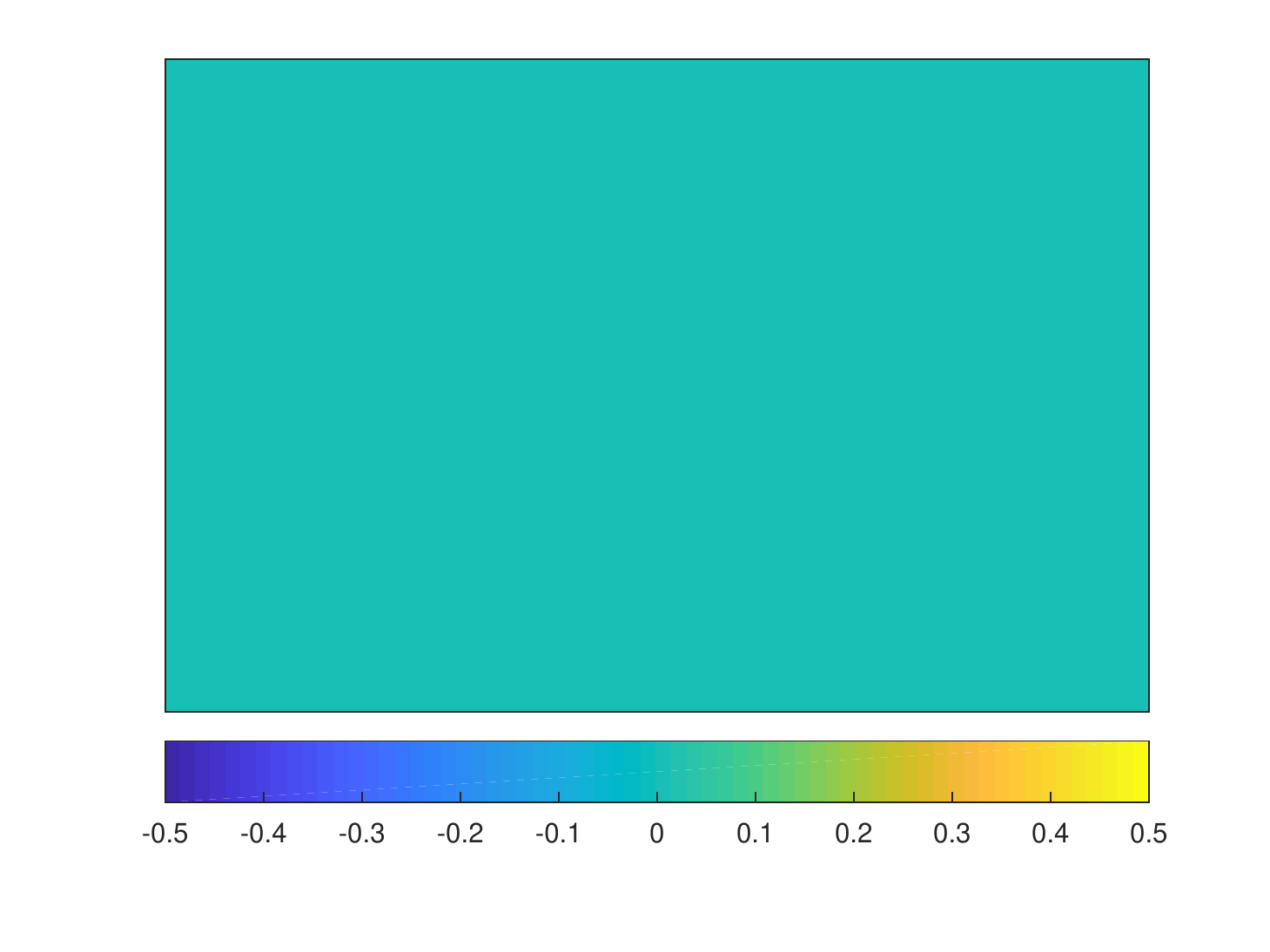}
&\includegraphics[width=\wsize,height=\hhsize,trim={1.5cm 1cm 1cm 0.5cm},clip]{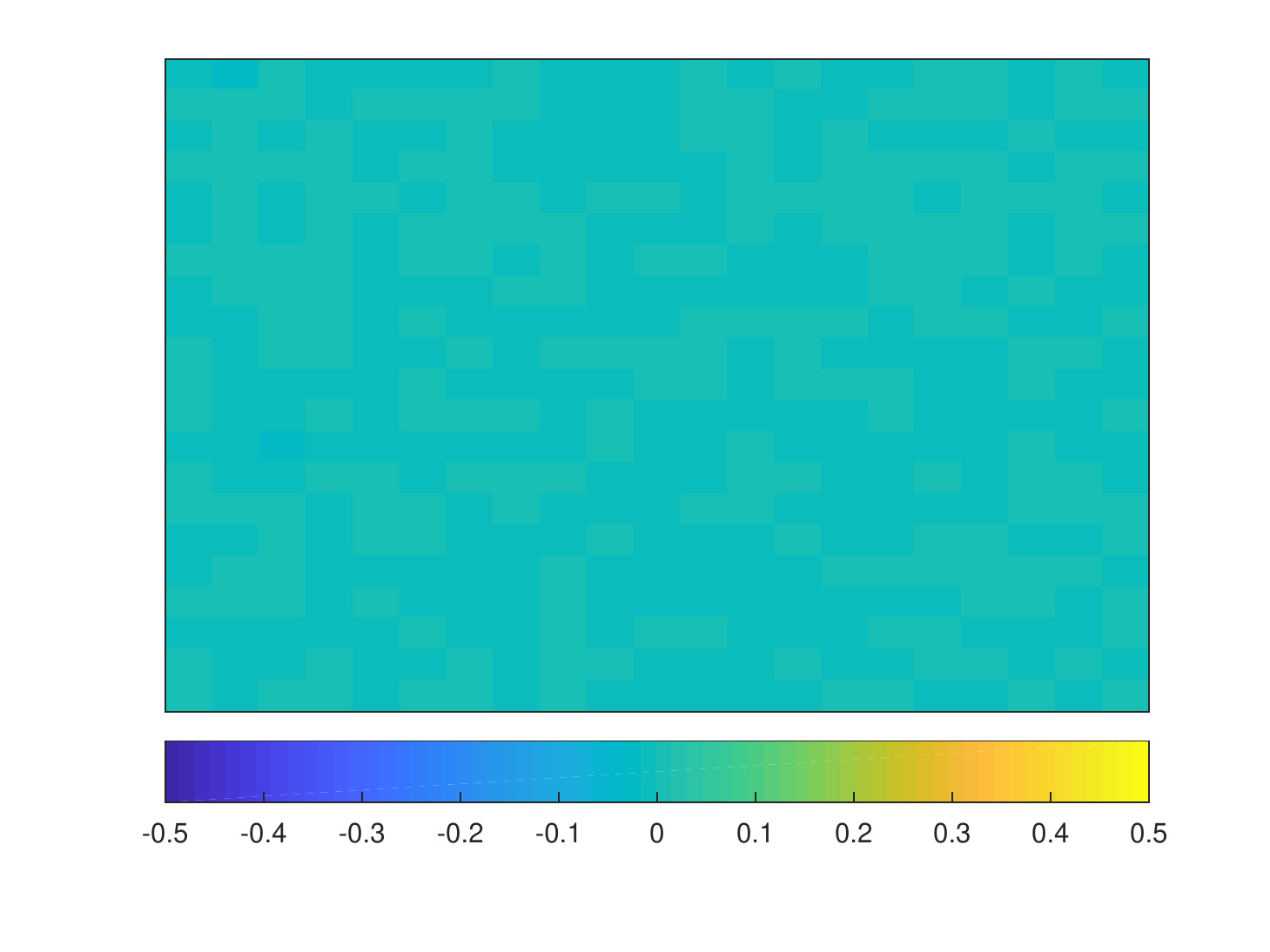}
&\includegraphics[width=\wsize,height=\hhsize,trim={1.5cm 1cm 1cm 0.5cm},clip]{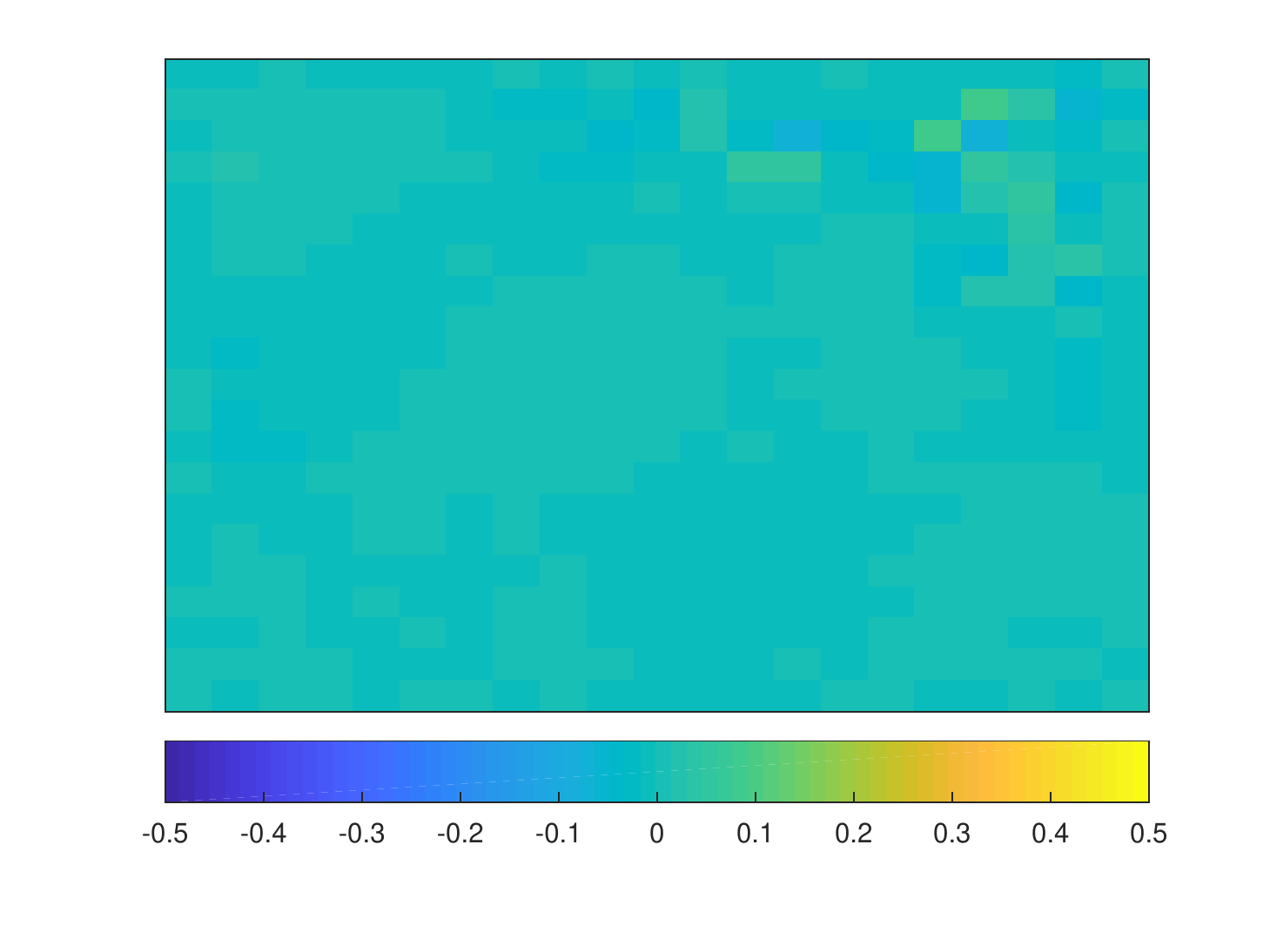}
&\includegraphics[width=\wsize,height=\hhsize,trim={1cm 0.5cm .5cm .5cm},clip]{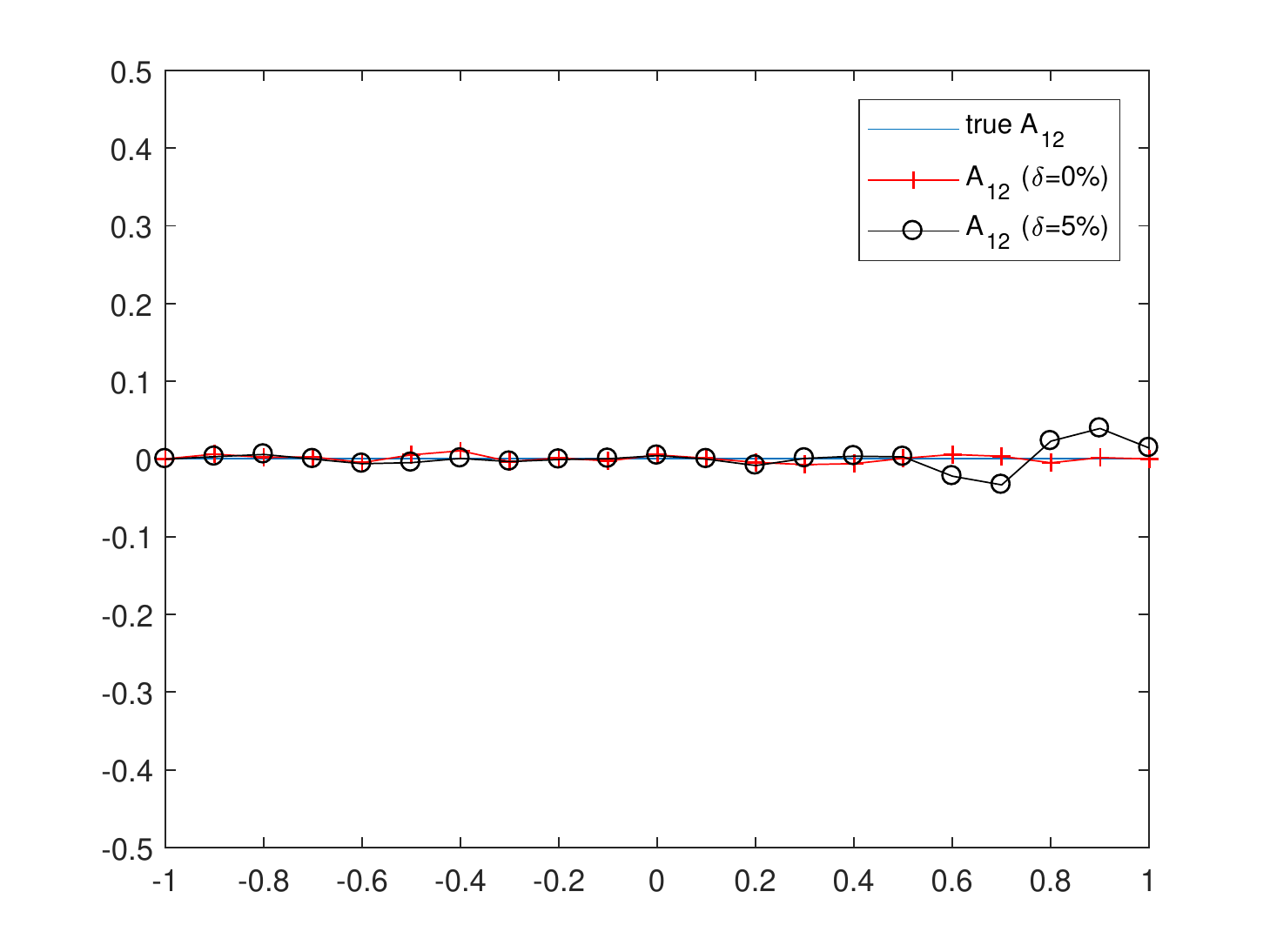}\\
 \includegraphics[width=\wsize,height=\hhsize,trim={1.5cm 1cm 1cm 0.5cm},clip]{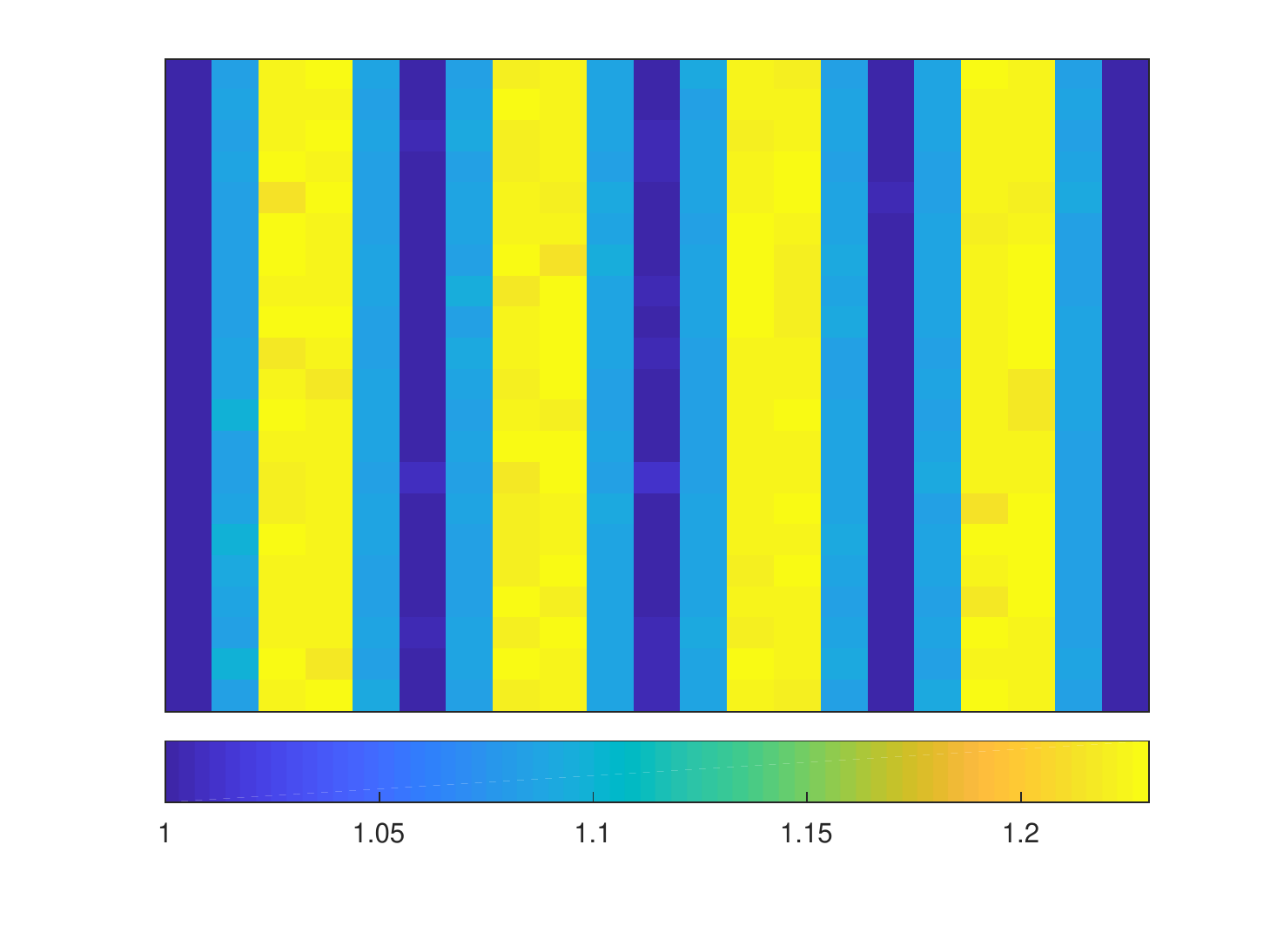}
&\includegraphics[width=\wsize,height=\hhsize,trim={1.5cm 1cm 1cm 0.5cm},clip]{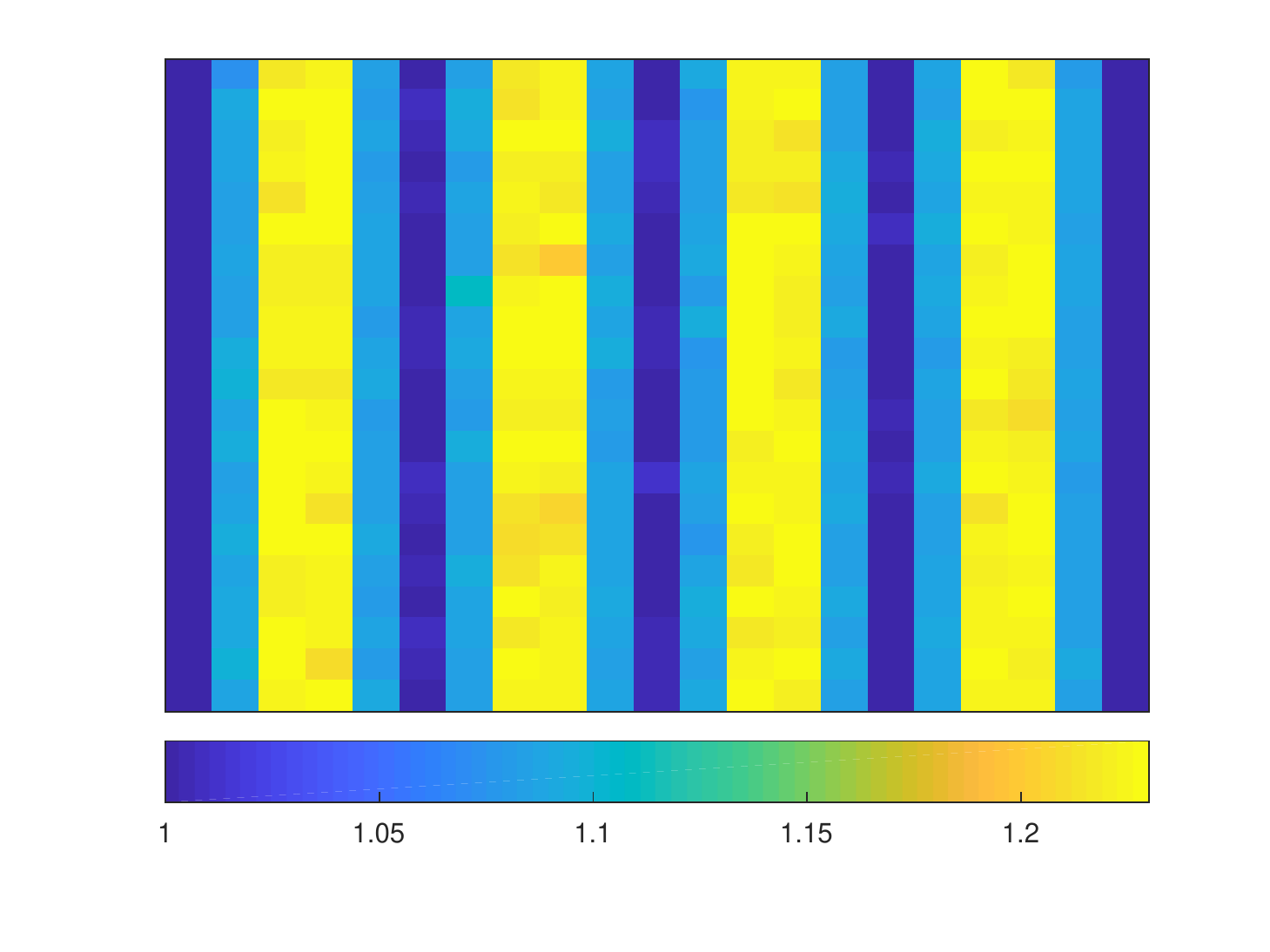}
&\includegraphics[width=\wsize,height=\hhsize,trim={1.5cm 1cm 1cm 0.5cm},clip]{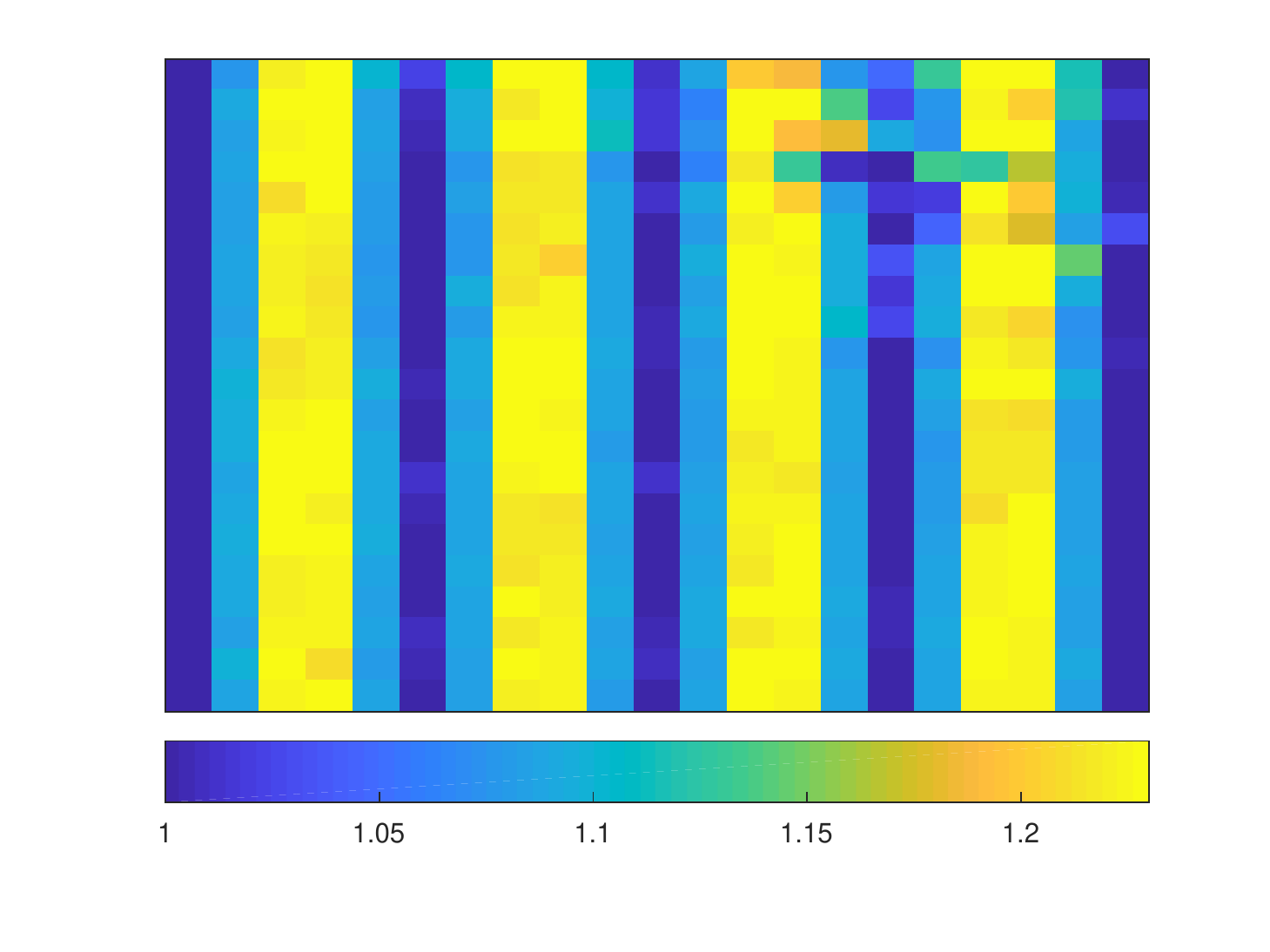}
&\includegraphics[width=\wsize,height=\hhsize,trim={1cm .5cm .5cm .5cm },clip]{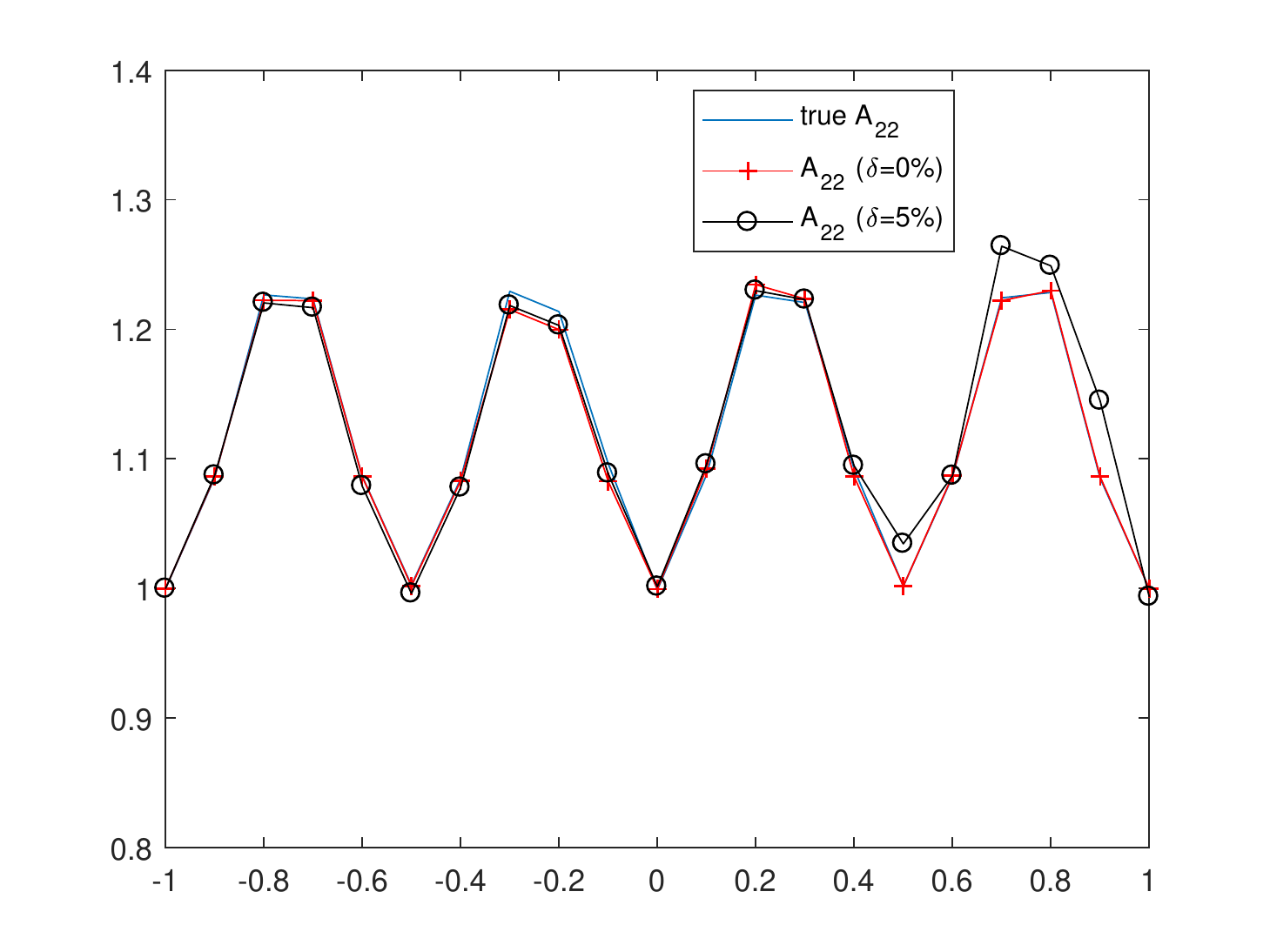}\\
\end{tabular}
\caption{The reconstructions for the conductivity tensor for Example 3 with three measurements: the top, middle and bottom rows refer to $A_{11}$, $A_{12}$ and $A_{22}$,
respectively. Column 1 is for exact conductivity, Column 2 for reconstruction with exact data, Column 3 for noisy data with $\delta=5\%$, and column 4 for the cross section along \{$x_2=-0.4$\}.}
\label{exam3data3}
\end{figure}

\begin{figure}[!htbp]
\centering     \setlength{\tabcolsep}{0pt}
\begin{tabular}{cccc}
 \includegraphics[width=\wsize,height=\hhsize,trim={1.5cm 1cm 1cm 0.5cm},clip]{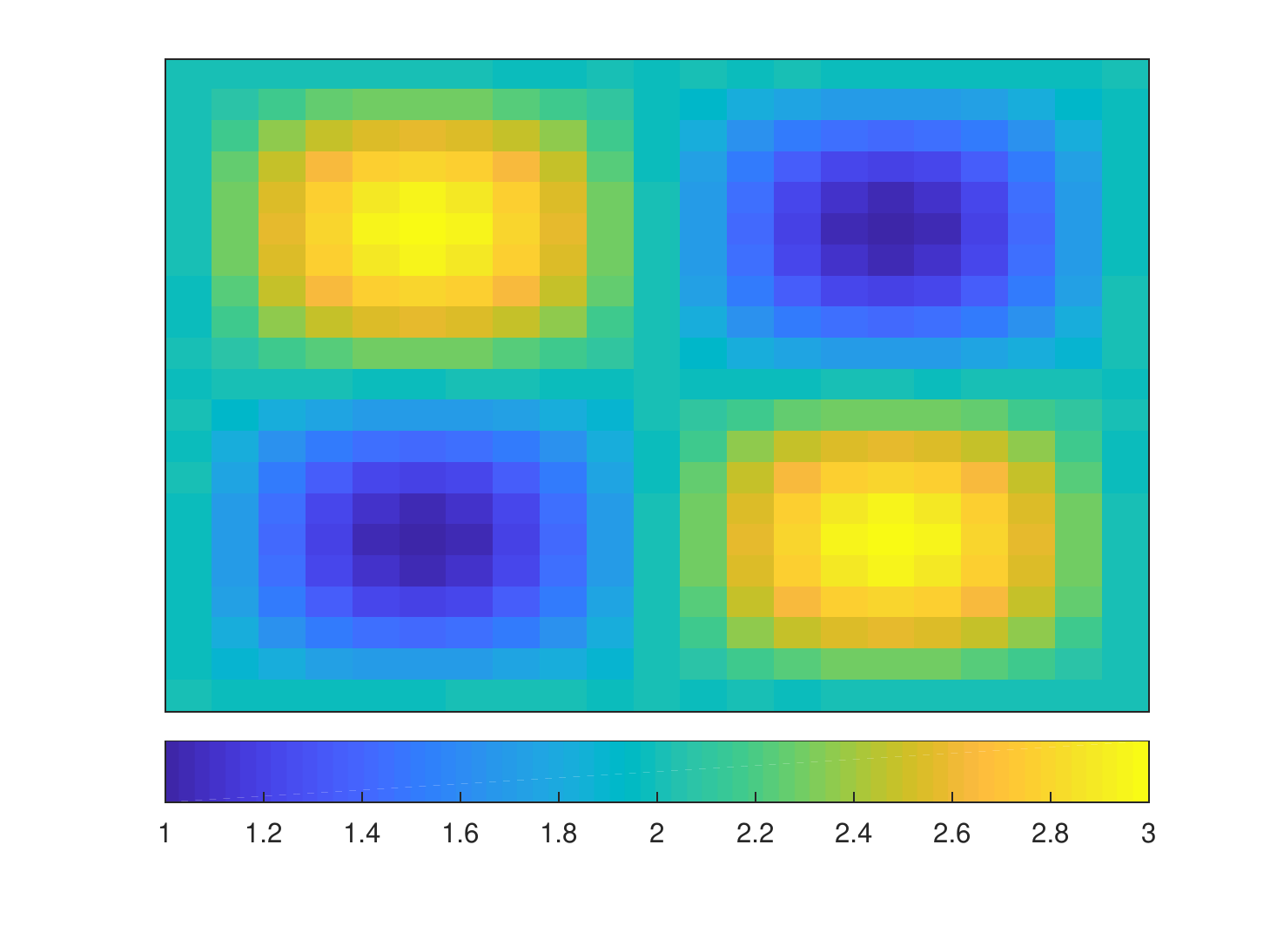}
&\includegraphics[width=\wsize,height=\hhsize,trim={1.5cm 1cm 1cm 0.5cm},clip]{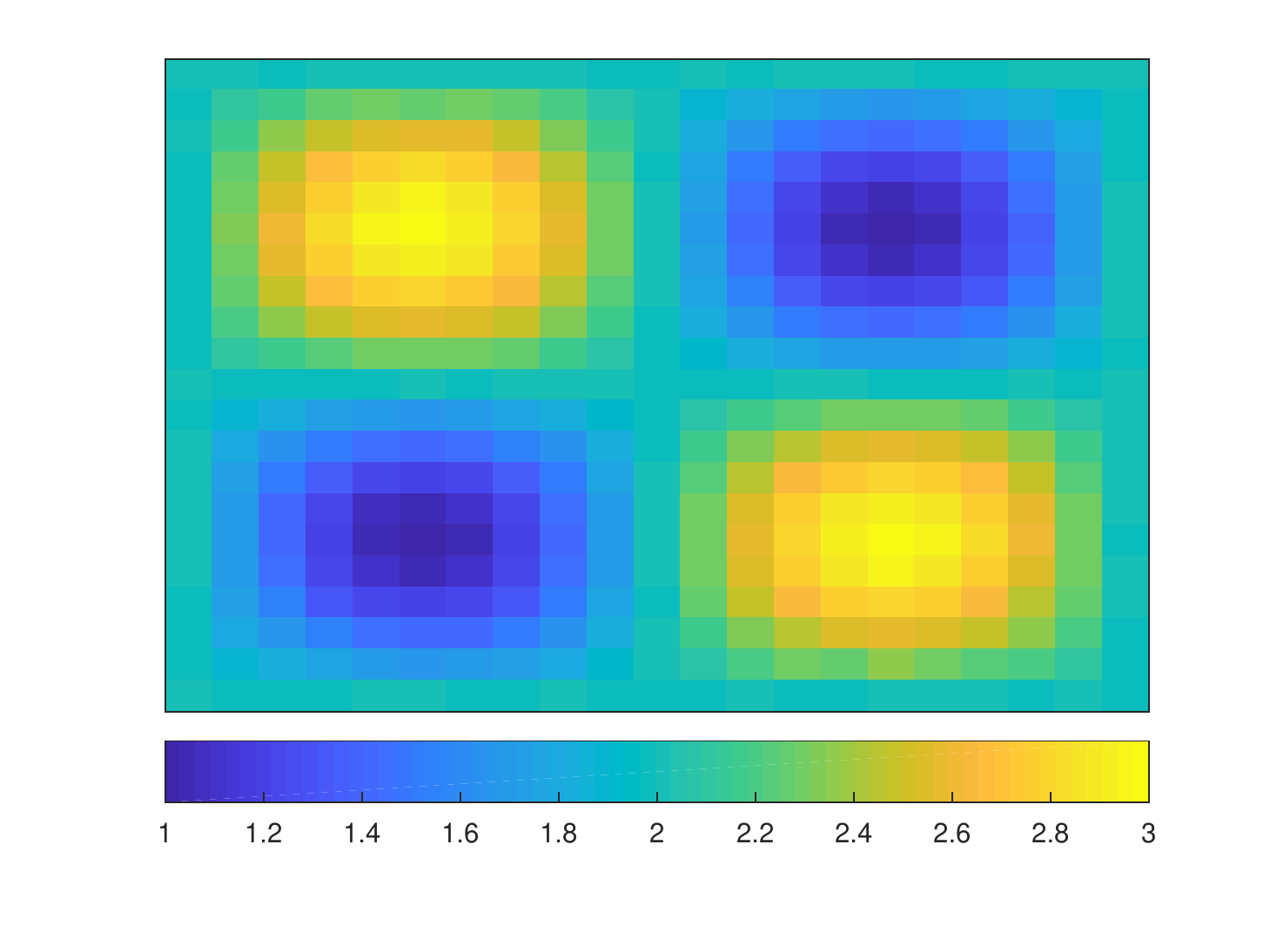}
&\includegraphics[width=\wsize,height=\hhsize,trim={1.5cm 1cm 1cm 0.5cm},clip]{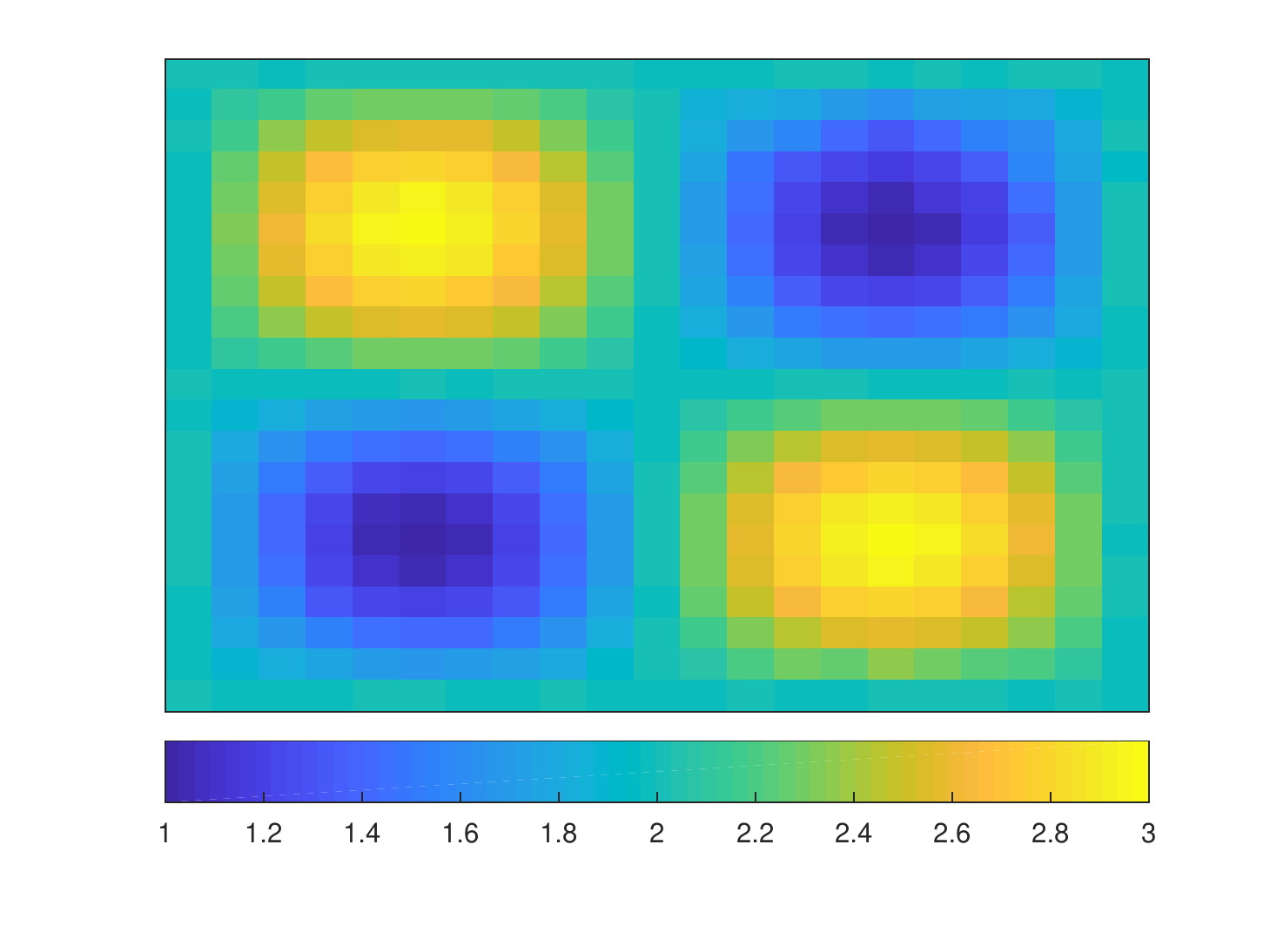}
&\includegraphics[width=\wsize,height=\hhsize,trim={1cm .5cm .5cm .5cm },clip]{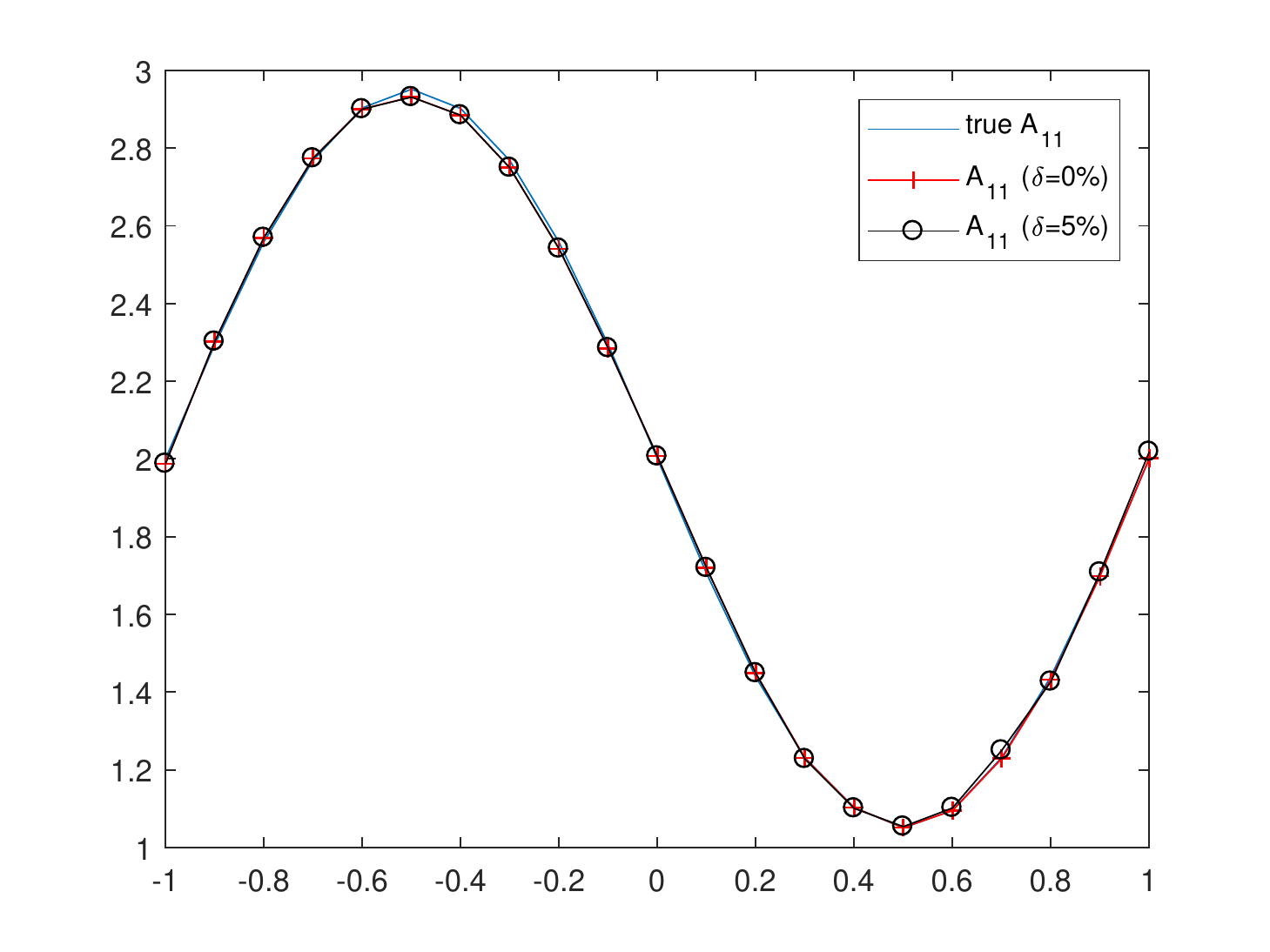}\\
 \includegraphics[width=\wsize,height=\hhsize,trim={1.5cm 1cm 1cm 0.5cm},clip]{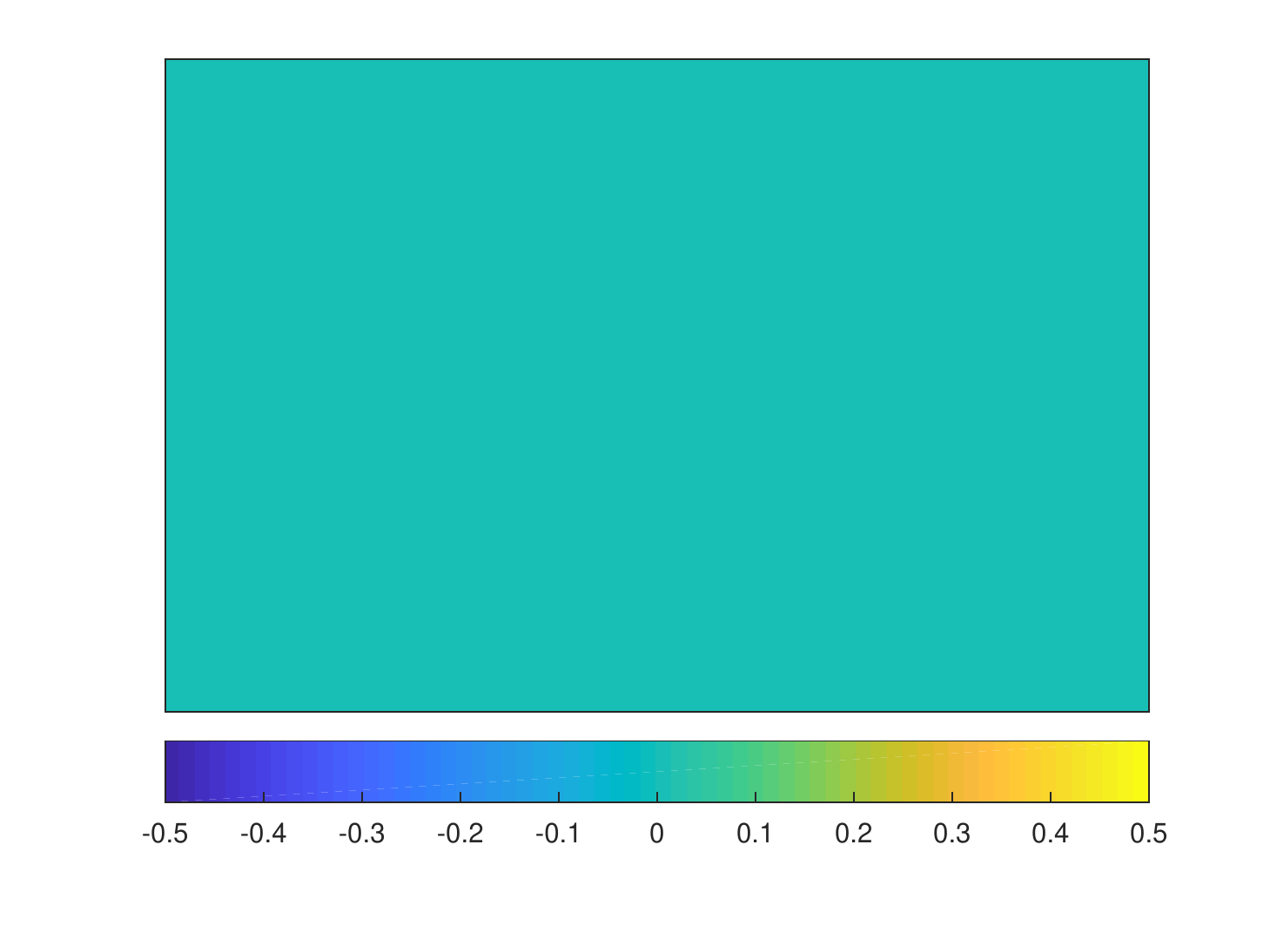}
&\includegraphics[width=\wsize,height=\hhsize,trim={1.5cm 1cm 1cm 0.5cm},clip]{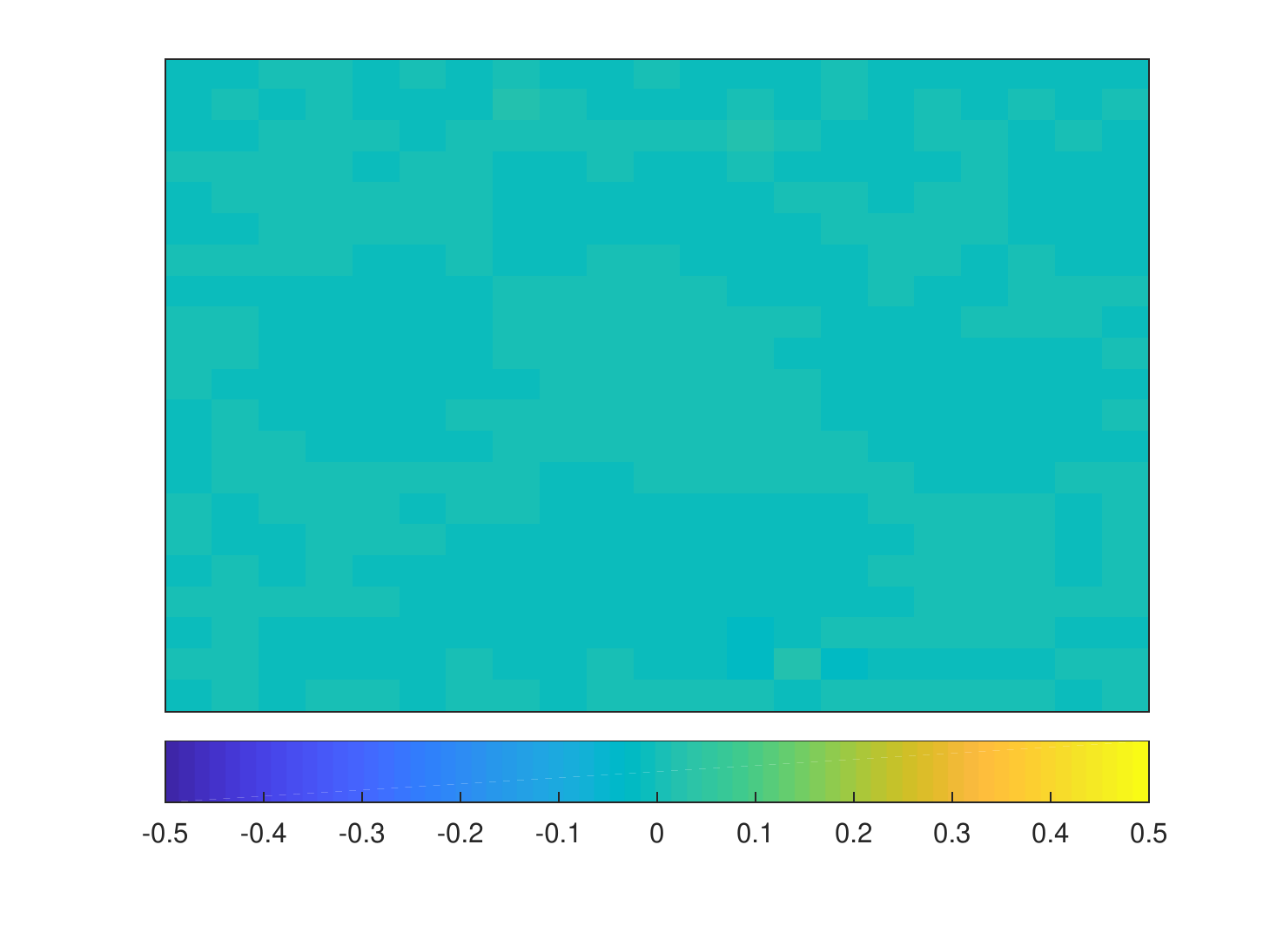}
&\includegraphics[width=\wsize,height=\hhsize,trim={1.5cm 1cm 1cm 0.5cm},clip]{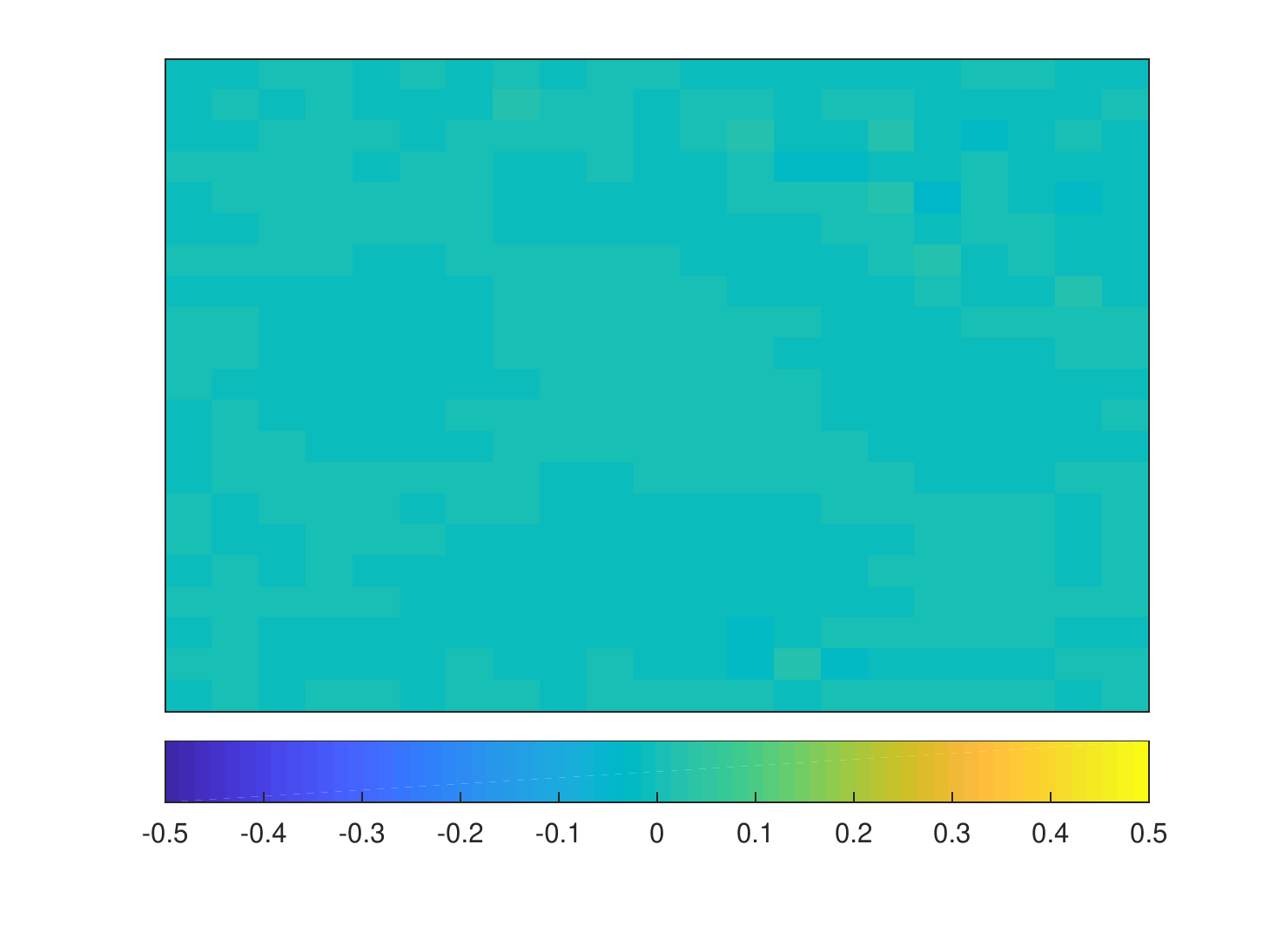}
&\includegraphics[width=\wsize,height=\hhsize,trim={1cm 0.5cm .5cm .5cm},clip]{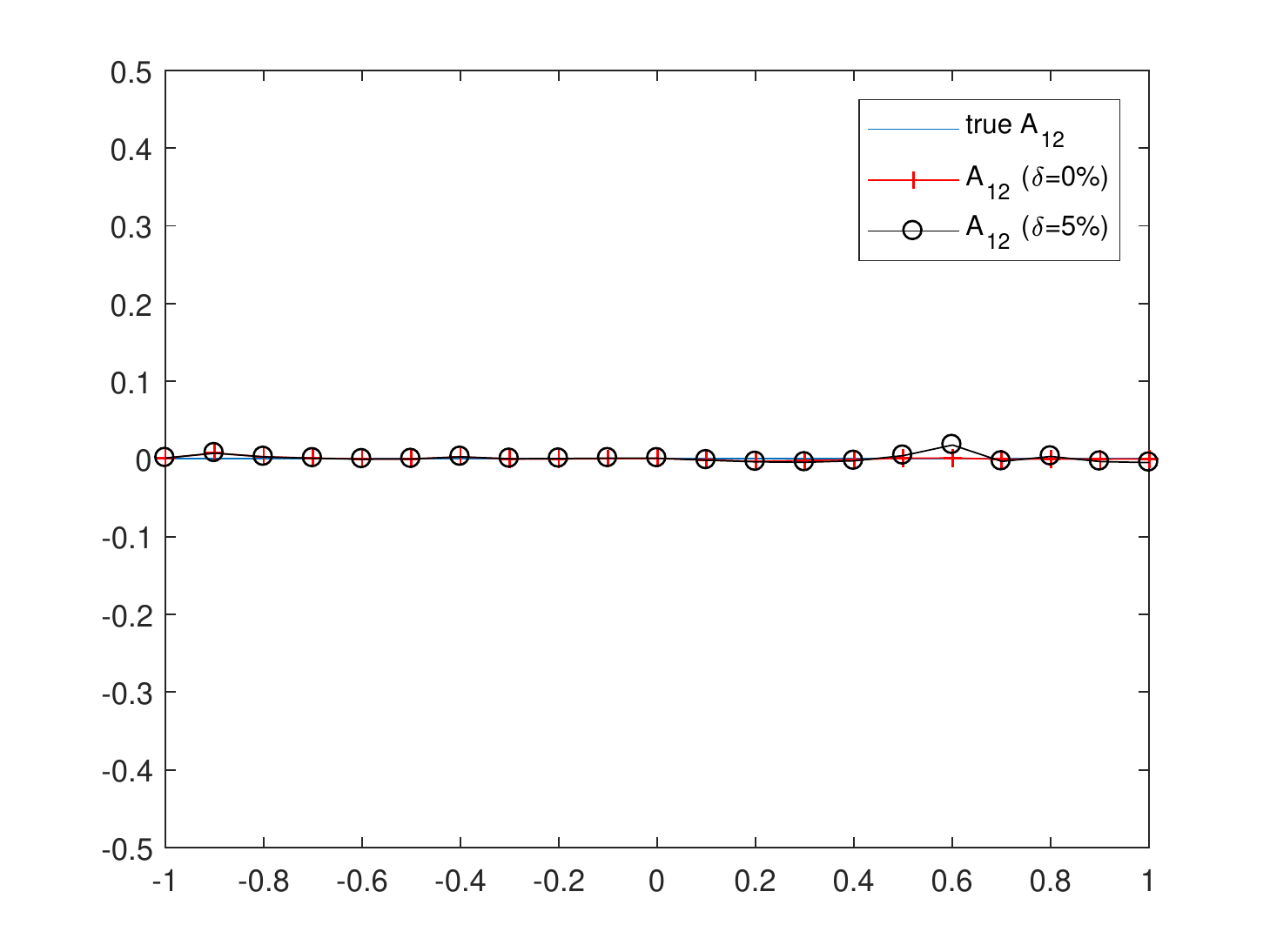}\\
 \includegraphics[width=\wsize,height=\hhsize,trim={1.5cm 1cm 1cm 0.5cm},clip]{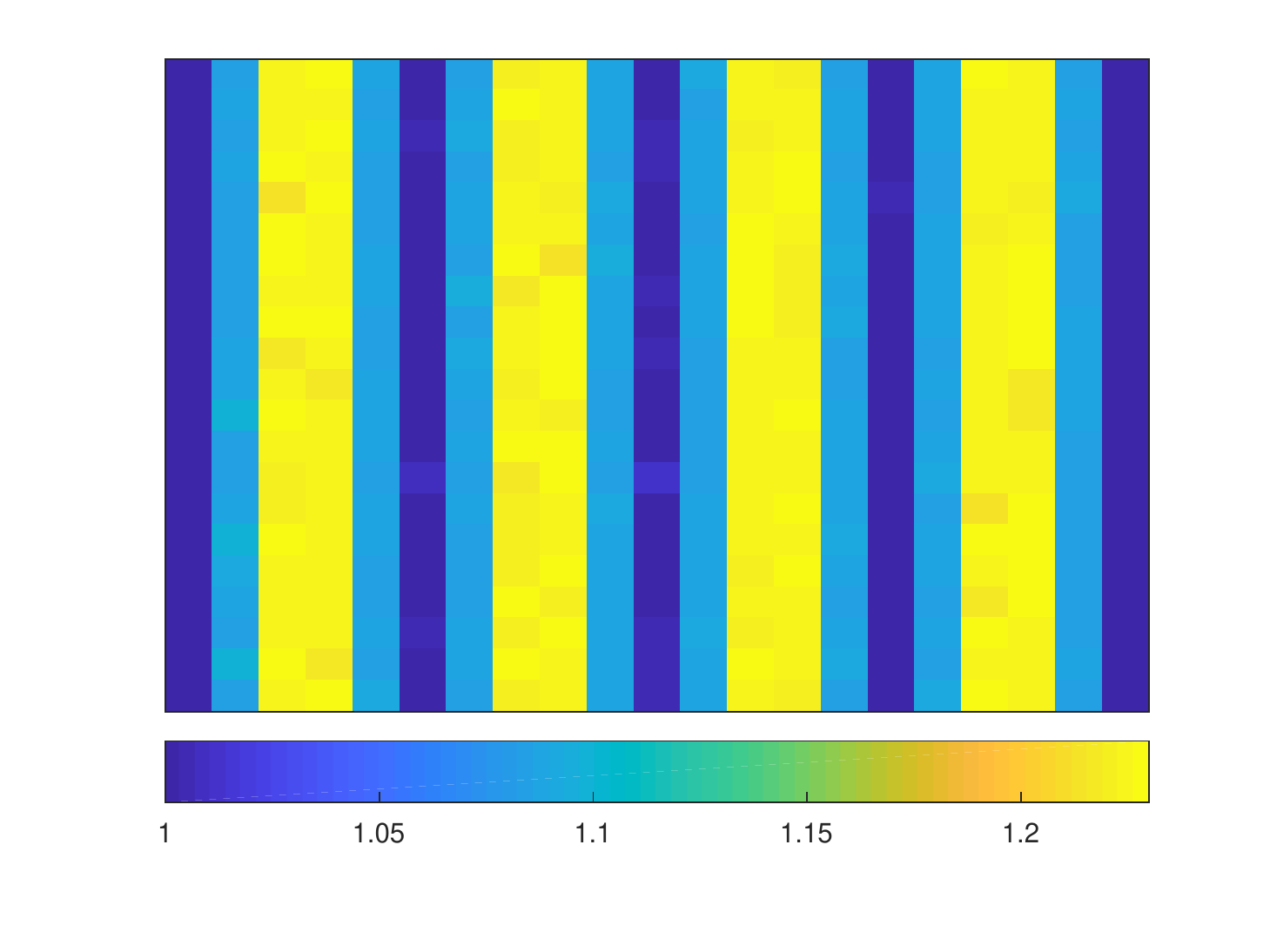}
&\includegraphics[width=\wsize,height=\hhsize,trim={1.5cm 1cm 1cm 0.5cm},clip]{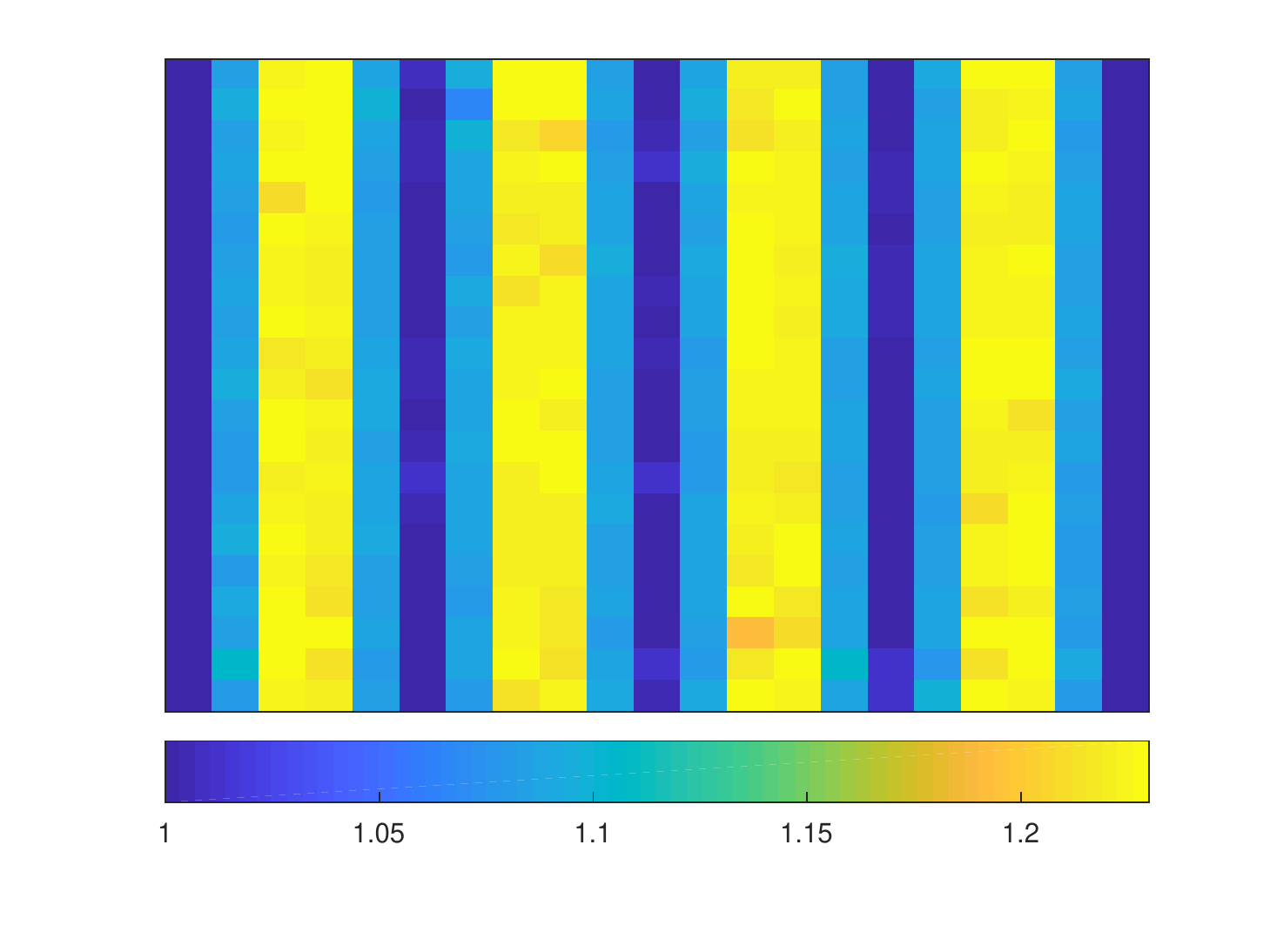}
&\includegraphics[width=\wsize,height=\hhsize,trim={1.5cm 1cm 1cm 0.5cm},clip]{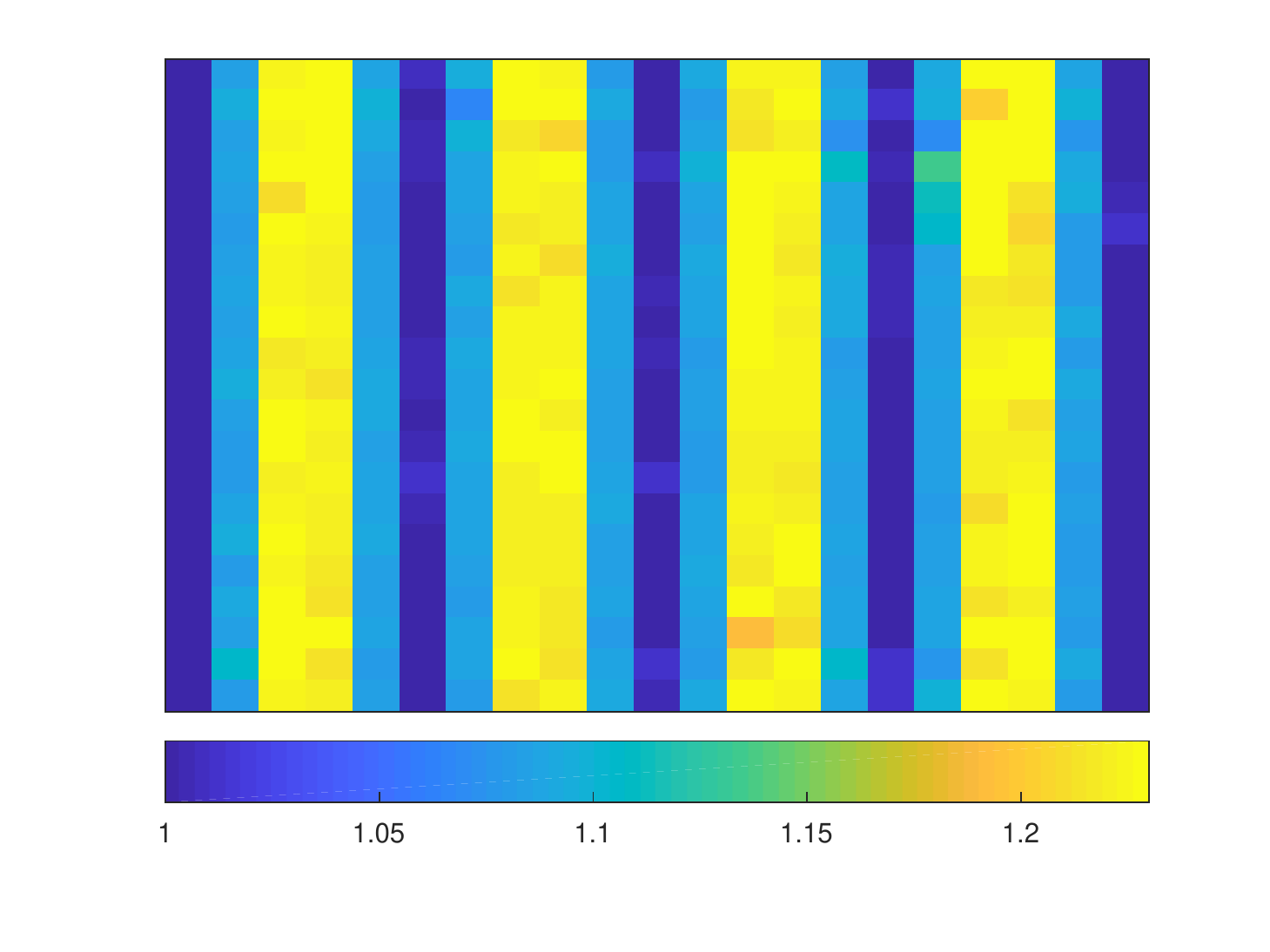}
&\includegraphics[width=\wsize,height=\hhsize,trim={1cm .5cm .5cm .5cm },clip]{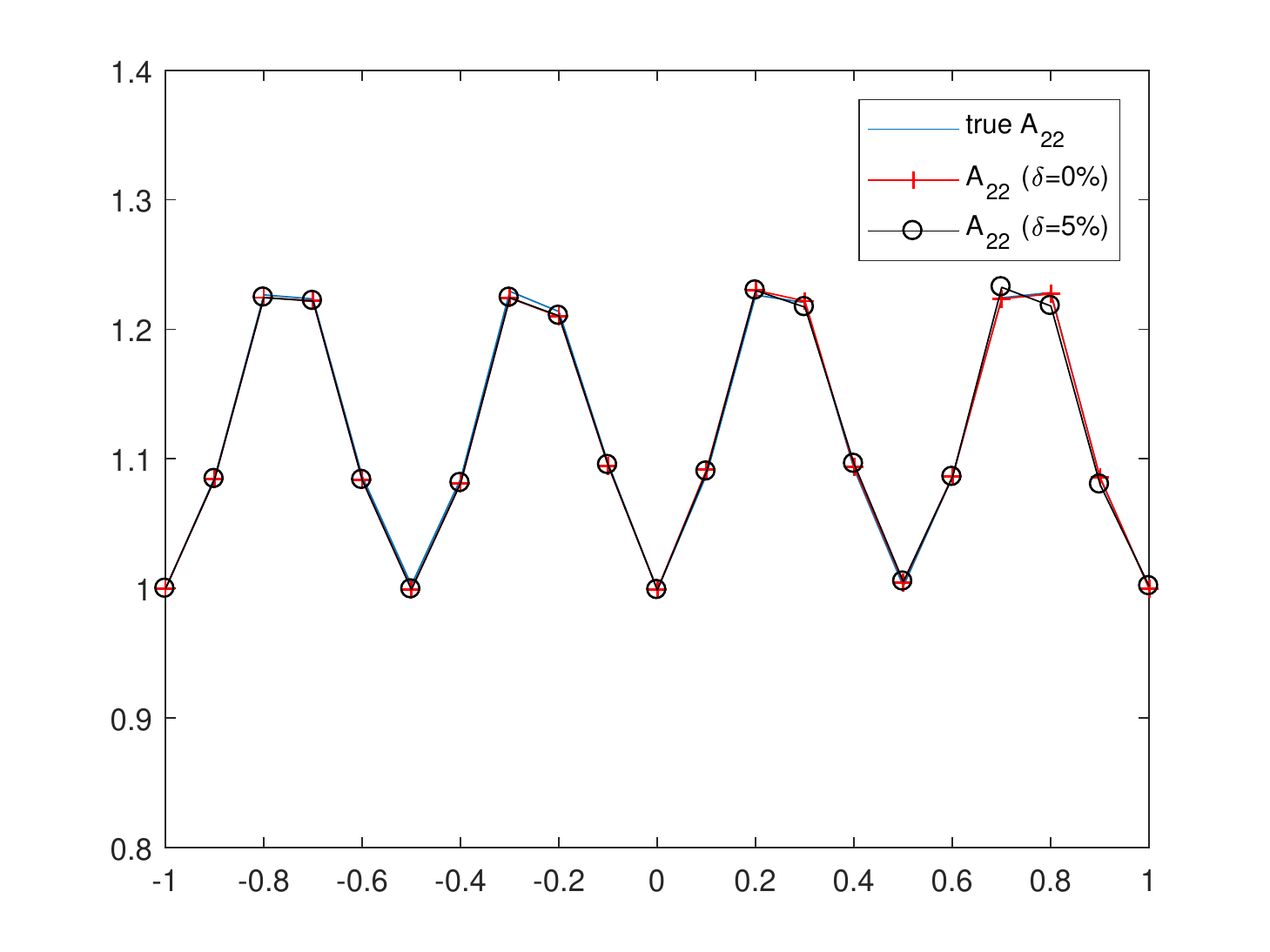}\\
\end{tabular}
\caption{The reconstructions for the conductivity tensor for Example 3 with nine measurements: the top, middle and bottom rows refer to $A_{11}$, $A_{12}$ and $A_{22}$,
respectively. Column 1 is for exact conductivity, Column 2 for reconstruction with exact data, Column 3 for noisy data with $\delta=5\%$, and column 4 for the cross section along \{$x_2=-0.4$\}.}
\label{exam3data9}
\end{figure}

Last, we examine the influence of the measurement number $L$ on the reconstruction quality.
In practice, the number $L$ of available measurements
may depend on the specific application. In Figs. \ref{exam1}-\ref{exam7}, we have always used
five measurements. In Figs. \ref{exam3data3} and \ref{exam3data9}, we repeat Example 3 with
three and nine measurements, respectively given by
\begin{itemize}
\item[(i)] three Dirichlet inputs
$(g^{1},g^2,g^3)(x)=(x_1+x_2,x_2+0.5x_1^2,x_1-0.1x_2^2)$
\item[(ii)] nine Dirichlet inputs
$(g^1,g^2,g^3,g^4,g^5,g^6,g^7,g^8,g^9)=(x_1+x_2,x_2+0.5x_1^2,x_1-0.1x_2^2,$
  $0.1(\cos(10x_2)- \cos(10x_1)),x_1x_2, 0.1\sin(10 x_1), x_1^2+x_2^2, x_1-0.5x_1^2, 0.1e^{x_1})$.
\end{itemize}
It is observed from Figs. \ref{exam3}, \ref{exam3data3} and \ref{exam3data9} that when the number of
measurements increases, the reconstruction becomes more accurate. Thus, a larger number of data
is beneficial for recovering conductivity tensors, as one might expect.

\section{Conclusions}
In this paper, we have studied both analytically and numerically the reconstruction of an
anisotropic conductivity tensor $A$ from the knowledge of multiple current densities of the form
$A\nabla u$, corresponding to different Dirichlet boundary excitations. We have established important analytical
properties, e.g., Lipschitz continuity and Frech\'{e}t differentiability, of the parameter-to-state
map  with respect to $L^p(\Omega)^{d,d}$-norms. We proposed and analyzed a regularized formulation of
Tikhonov type, which is suitable for both smooth and nonsmooth conductivity tensors. We proved
the existence of minimizers and consistency as the noise level tends to zero by means of
H-convergence, and derived the necessary optimality system and a Newton algorithm for
efficiently solving the coupled system. Also we established the convergence
of the finite element approximation by means of Hd-convergence. The extensive numerical
simulations clearly showcase the significant potentials of the approach.

\section*{Acknowledgments}
The authors would like to thank the two anonymous referees for several constructive
comments, which have led to an improved presnetation.

\appendix

\section{Proof of Theorem \ref{thm:diff-tikh}}
\begin{proof}
Since the penalty is quadratic, its gradient expression and differentiability follow directly. It suffices to compute
the gradient of the fidelity term. For any fix feasible direction $H$, let ${\rm I}_\ell(t) = \frac12\|{(A+tH)}\nabla u^\ell(A+tH)-\vec{h}^\ell\|^2_{L^2(\Omega)^d}-\frac12\|A\nabla u^\ell (A)-\vec{h}^\ell\|^2_{L^2(\Omega)^d}$.
Clearly,
\begin{align*}
{\rm I}_\ell(t)
&=\tfrac12\|(A+tH)\big(\nabla u^\ell(A+tH)-\nabla u^\ell(A)\big)+tH\nabla u^\ell\|_{L^2(\Omega)^d}^2\\
 &\quad +\big(A\nabla u^\ell(A)-\vec{h}^\ell, (A+tH)(\nabla u^\ell(A+tH)-\nabla u^\ell(A))+tH\nabla u^\ell\big)_{L^2(\Omega)^d}.
\end{align*}
Let $w^\ell = F_\ell^{\prime}(A)[H]$ be the linearized solution for $u^\ell$ at $A$ in
the direction of $H$, and $R^\ell$ the corresponding residual, i.e., $R^\ell=u^\ell(A+tH)-u^\ell(A)-tF_\ell^{\prime}(A)[H]$. Then
\begin{align*}
{\rm I}_\ell(t)&=\tfrac{t^2}{2}\|(A+tH)(\nabla w^{(i)} + t^{-1}\nabla R^\ell)+H\nabla u^\ell\|_{L^2(\Omega)^d}^2\\
     &\quad+\big(A\nabla u^\ell(A)-\vec{h}^\ell, (A+tH)(t\nabla w^\ell+ \nabla R^\ell) +tH\nabla u^\ell\big)_{L^2(\Omega)^d}.
\end{align*}
By the weak formulations for $p^\ell$ and $w^\ell$, we have
\begin{align*}
(A\nabla w^\ell,\nabla {p}^\ell)_{L^2(\Omega)^d}&=-(H\nabla u^\ell, \nabla {p}^\ell)_{L^2(\Omega)^d},\\
(A\nabla {p}^\ell,\nabla w^{(i)})_{L^2(\Omega)^d}&=-(A(A\nabla u^\ell-\vec{h}^\ell),\nabla w^\ell)_{L^2(\Omega)^d}.
\end{align*}
Subtracting these two identities gives
\begin{equation}\label{eq}
-(H\nabla u^\ell,\nabla{p}^\ell)_{L^2(\Omega)^d}=-(A(A\nabla u^\ell-\vec{h}^\ell), \nabla w^\ell)_{L^2(\Omega)^d}.
\end{equation}
This and the identity $\|R^\ell\|_{H^1(\Omega)}=O(t^2)$ from
the proof of Proposition \ref{prop:deriv} imply
\begin{align*}
\lim_{t\rightarrow 0^+}t^{-1}{\rm I}_\ell(t)&=
\lim_{t\rightarrow 0^+}\big(A\nabla u^\ell(A)-\vec{h}^\ell,(A+tH)(\nabla w^\ell +t^{-1}\nabla R^\ell)+H\nabla u^\ell\big)_{L^2(\Omega)^d}\\
&=
\lim_{t\rightarrow 0^+}\big(A\nabla u^\ell(A)-\vec{h}^\ell,tH(\nabla w^\ell+t^{-1}\nabla R^\ell)
+H\nabla u^\ell+At^{-1}\nabla R^\ell)_{L^2(\Omega)^d}\\
&\qquad+(H\nabla u^\ell,\nabla {p}^\ell)_{L^2(\Omega)^d}\\
&=\big(\nabla u^\ell\otimes \nabla {p}^\ell
+\nabla u^\ell\otimes (A\nabla u^\ell-\vec{h}^\ell),H\big)_{L^2(\Omega)^{d,d}},
\end{align*}
which directly gives the expression of $J_\gamma'(A)[H]$.
Next, we verify that it defines a bounded linear functional on $L^p(\Omega)^{d\times d}$. For each $\ell=1,
\ldots,L$, it follows from  \eqref{inp} that for any $p, q\geq1$ satisfying $\frac1p+\frac1q=\frac12$,
\begin{align*}
&\quad \big|(\nabla u^{\ell}\otimes \nabla \bar{p}^{\ell}+\nabla u^{\ell}\otimes (A\nabla u^{\ell}-\vec{h}^{\ell}), H)_{L^2(\Omega)^{d,d}}\big|\\
&=\big|(H\nabla u^{\ell},\nabla \bar{p}^{\ell})_{L^2(\Omega)^d}+(H\nabla u^{\ell},A\nabla u^{\ell}-\vec{h}^{\ell})_{L^2(\Omega)^d}\big|\\
&\leq \normmm{H}_p\big(\|\nabla u^{\ell}\|_{L^q(\Omega)^d}\|\nabla \bar{p}^{\ell}\|_{L^2(\Omega)^d}+
\|\nabla u^{\ell}\|_{L^q(\Omega)^d}
\|A\nabla u^{\ell}-\vec{h}^{\ell}\|_{L^2(\Omega)^d}\big)\\
&\leq \normmm{H}_p\big((\|f^{\ell}\|_{L^2(\Omega)}+\|g^{\ell}\|_{H^1(\Gamma)})\|\bar{p}^{\ell}\|_{H^1(\Omega)}
+(\|f^{\ell}\|_{L^2(\Omega)}+\|g^{\ell}\|_{H^1(\Gamma)})\|A\nabla u^{\ell}-\vec{h}^{\ell}\|_{L^2(\Omega)^d}\big).
\end{align*}
Hence, it is a bounded linear functional on $L^p(\Omega)^{d,d}$.
Last, we show that it is actually the Fr\'{e}chet derivative of $J_\gamma(A)$. Let
${\rm II}_\ell: ={\rm I}_\ell(1)-(\nabla u^{\ell}\otimes \nabla \bar{p}^{\ell}+
\nabla u^{\ell}\otimes (A\nabla u^{\ell}-\vec{h}^{\ell}), H)_{L^2(\Omega)^{d,d}}$. Indeed, we have
\begin{align*}
{\rm II}_\ell
&=\tfrac12\|(A+H)\big(\nabla u^{\ell}(A+H)-\nabla u^{\ell}(A)\big)+H\nabla u^{\ell}\|_{L^2(\Omega)^d}^2\\
&\quad +\big(A\nabla u^{\ell}(A)-\vec{h}^{\ell},(A+H)(\nabla u^{\ell} (A+H)- \nabla u^{\ell}(A))\big)_{L^2(\Omega)^d}-\big(H\nabla u^{\ell},\nabla\bar{p}^{\ell}\big)_{L^2(\Omega)^d}.
\end{align*}
This and \eqref{eq} lead to
\begin{align*}
|{\rm II}_\ell|&=\Big|\tfrac12\|(A+H)\big(\nabla u^{\ell}(A+H)-\nabla u^{\ell}(A)\big)+H\nabla u^{\ell}\|_{L^2(\Omega)^d}^2\\
  &\quad +(A\nabla u^{\ell}(A)-\vec{h}^{\ell},H(\nabla u^{\ell} (A+H)-\nabla u^{\ell}(A)))_{L^2(\Omega)^d}+(A\nabla u^{\ell}(A)-\vec{h}^{\ell}, A\nabla R)_{L^2(\Omega)^d}\Big|\\
&\leq \|A+H\|^2_{L^\infty(\Omega)^{d,d}} (\|\nabla u^{\ell}(A+H)-\nabla u^{\ell}(A)\|^2_{L^2(\Omega)^d}+
\|H\nabla u^{\ell}\|^2_{L^2(\Omega)^d})\\
&\quad +\|A\nabla u^{\ell}(A)-\vec{h}^{\ell}\|
_{L^2(\Omega)^d}\|H\|_{L^\infty(\Omega)^{d,d}}
\|\nabla u^{\ell}(A+H)-\nabla u^{\ell}(A)\|_{L^2(\Omega)^d}\\
&\quad +\|A\nabla u^{\ell}(A)-\vec{h}^{\ell}\|_{L^2(\Omega)^d}\|A\|_{L^\infty(\Omega)^{d,d}}\|\nabla R\|_{L^2(\Omega)^d}.
\end{align*}
This estimate, Lemma \ref{lem:LipCont} and Proposition \ref{prop:deriv} yield
the desired differentiability. This complete the proof of the theorem.
\end{proof}
\bibliographystyle{siam}
\bibliography{mybib}
\end{document}